\documentclass[12pt]{book}

\usepackage{epsfig, psfrag, color}
\usepackage[toc,page]{appendix}
\usepackage{multind}
\makeindex{Index}
\makeindex{Symbols}
\input epsf.sty
\epsfverbosetrue 

\usepackage{latexsym}
\usepackage{amsfonts}
\usepackage{amssymb}

\newcommand{\R}{\mathbb{R}}
\newcommand{\Q}{\mathbb{Q}}
\newcommand{\Z}{\mathbb{Z}}
\newcommand{\C}{\mathbb{C}}
\newcommand{\IP}{\mathbb{P}}
\newcommand{\IH}{\mathbb{H}}
\newcommand{\T}{\mathbb{T}}
\newcommand{\N}{\mathbb{N}}
\newcommand{\D}{\mathbb{D}}
\newcommand{\EE}{\mathbb{E}}
\newcommand{\IS}{\mathbb{S}}

\newcommand{\rs}{\mbox{$\widehat{\C}$}}

\def\SSS{{\mathcal S}}
\def\TTT{{\mathcal T}}
\def\AAA{{\mathcal A}}
\def\BBB{{\mathcal B}}
\def\CCC{{\mathcal C}}
\def\DDD{{\mathcal D}}

\def\FFF{{\mathcal F}}
\def\GGG{{\mathcal G}}
\def\HHH{{\mathcal H}}
\def\III{{\mathcal I}}

\def\MMM{{\mathcal M}}

\def\MMM{{\mathcal M}}

\def\YYY{{\mathcal Y}}

\def\PPP{{\mathcal P}}
\def\RRR{{\mathcal R}}
\def\SSS{{\mathcal S}}
\def\UUU{{\mathcal U}}
\def\VVV{{\mathcal V}}

\def\XXX{{\mathcal X}}

\newtheorem{thm}{Theorem}[section]
\newtheorem{defn}[thm]{Definition}

\newtheorem{prop}[thm]{Proposition}

\newtheorem{lemma}[thm]{Lemma}
\newtheorem{cor}[thm]{Corollary}

\newcommand{\qed}{\nopagebreak \begin{flushright}
        \rule{2mm}{2.5mm} \end{flushright}}

\newcommand{\implies}{\Rightarrow}
\newcommand{\bdry}{\partial}                     
\newcommand{\id}{\mbox{\rm id}}                  
\newcommand{\cl}{\overline}                      

\newcommand{\Aut}{\mbox{\rm Aut}}                        
\newcommand{\mod}{\mbox{\rm mod}}    
\newcommand{\intersect}{\cap}                    
\newcommand{\union}{\cup}                        
\newcommand{\la}{\langle}                        
\newcommand{\ra}{\rangle}                        
\newcommand{\diam}{\mbox{\rm diam}}					         
\newcommand{\dist}{\mbox{\rm dist}}
\newcommand{\interior}{\mbox{int}}     

\newcommand{\mtwo}[4]                            
{\mbox{$\left(\begin{array}{cc}                  
#1 & #2 \\
#3 & #4 
\end{array}
\right)$}}

\newcommand{\dettwo}[4]                          
{\mbox{$\left|\begin{array}{cc}                  
#1 & #2 \\
#3 & #4 
\end{array}
\right|$}}

\newcommand{\pf}{\noindent {\bf Proof: }}

\newcommand{\be}{\begin{enumerate}}
\newcommand{\eb}{\end{enumerate}}
\newcommand{\bi}{\begin{itemize}}
\newcommand{\ib}{\end{itemize}}
\newcommand{\bl}{\begin{list}}
\newcommand{\lb}{\end{list}}

\newcommand{\gap}{\vspace{5pt}}                 

\newcommand{\ibid}{\mbox{\it ibid}}

\newcommand{\PP}{\IP}
\newcommand{\HH}{\IH}
\newcommand{\qc}{quasiconformal }
\def\dis{\displaystyle}
\newcommand{\al}{\alpha} 
\newcommand{\cbar}{\rs} 
\newcommand{\roundness}{\mbox{\rm Round}}

\newcommand{\wtB}{\widetilde{B}}
\newcommand{\wtU}{{\widetilde{U}}}
\newcommand{\wtV}{\widetilde{V}}
\newcommand{\wtW}{\widetilde{W}}

\newcommand{\txi}{\tilde{\xi}}
\newcommand{\tv}{\tilde{v}}
\newcommand{\ta}{\tilde{a}}
\newcommand{\tx}{\tilde{x}}
\newcommand{\ty}{\tilde{y}}
\newcommand{\tu}{\tilde{u}}
\newcommand{\tw}{\tilde{w}}
\newcommand{\tr}{\tilde{r}}

\newcommand{\XX}{\mbox{$\mathfrak{X}$}}
\newcommand{\YY}{\mbox{$\mathfrak{Y}$}}

\newcommand{\te}{t_\epsilon}

\newcommand{\cwstar}{\mbox{\rm star}}

\begin{document}
\frontmatter
\pagestyle{empty}

\begin{center}
\vskip 2in

{ \LARGE Coarse expanding conformal dynamics}
\vskip 1in


\vfill

\begin{tabular}{ll} Peter Ha\"{\i}ssinsky & Kevin M. Pilgrim\\
LATP/CMI & Dept. of Mathematics\\ Universit\'e de Provence & Indiana University\\ 
39, rue Fr\'ed\'eric Joliot-Curie & Bloomington\\
13453 Marseille cedex 13 & IN 47405\\France & USA\\
{\it E-mail:} {\tt phaissin@cmi.univ-mrs.fr \hspace{0.5in}} & {\it E-mail:} {\tt pilgrim@indiana.edu}\end{tabular}

\end{center}

\newpage
\pagestyle{headings}

$\;$
\vspace{1in}

\begin{center}{\bf Abstract}\end{center}
Motivated by the dynamics of  rational maps, we introduce a class of topological dynamical systems satisfying 
certain topological regularity, expansion, irreducibility, and finiteness conditions.  
We call such maps ``topologically coarse expanding conformal" (top. cxc) dynamical systems.  
Given such a system $f: X \to X$ and a finite cover of $X$ by connected open sets, 
we construct a negatively curved infinite graph on which $f$ acts naturally by local isometries.  
The induced topological dynamical system on the boundary at infinity is naturally  conjugate to the dynamics of $f$.  This implies that $X$ inherits
metrics in which the dynamics of $f$ satisfies the Principle of the Conformal Elevator:  arbitrarily small balls may 
be blown up with bounded distortion to nearly round sets of definite size. This property is preserved under conjugation by a quasisymmetric map, 
and top. cxc dynamical systems on a metric space satisfying this property we call ``metrically cxc''.  
The ensuing results deepen the analogy between rational maps and Kleinian groups by extending it 
to analogies between metric cxc systems and hyperbolic groups.   
We give many examples and several applications. In particular, we provide a new interpretation of 
the characterization of rational functions among topological maps and 
of generalized Latt\`es examples among uniformly quasiregular maps. 
{\em Via} techniques in the spirit of those used to construct quasiconformal measures for hyperbolic groups,  
we also establish existence, uniqueness, naturality, and metric regularity properties 
for the measure of maximal entropy of such systems.

\vfill 



\newpage

$\;$
\vspace{2in}
\begin{center}
{\em This work is dedicated to Adrien Douady.}
\end{center}
\vspace*{4cm}

\noindent{``Que s'est-il pass\'e dans ta t\^ete ?\\
Tu as pris la poudre d'escampette\\
Sans explication est-ce b\^ete\\
Sans raison tu m'as plant\'e l\`a\\
Ah ah ah ah !''}

\gap

\gap

\noindent (Boby Lapointe, {\it La question ne se pose pas})

\newpage

\tableofcontents

\mainmatter
\pagestyle{headings}

\chapter{Introduction}

The classical  conformal dynamical systems include iterated rational maps and Kleinian groups acting on the Riemann sphere.  
The development of these two theories was propelled  
forward in the early 1980's by Sullivan's introduction of quasiconformal methods and of  a  ``dictionary'' 
between the two subjects \cite{sullivan:qcI}.   
Via complex analysis, many basic dynamical objects can be similarly defined and results similarly proven.  
There is a general deformation theory which specializes to both subjects and which yields deep finiteness results  \cite{ctm:ds:qciii}.  
Since then, the dictionary has grown to encompass a guiding heuristic whereby constructions, methods, 
and results in one subject suggest similar ones in the other.   
For example, in both subjects there are common themes in the combinatorial classification 
theories \cite{ctm:classification}, \cite{kmp:cds}, 
the fine geometric structure of the associated fractal objects \cite{ctm:kleinjohn}, \cite{ctm:hdimii}, 
\cite{stratmann:urbanski:parabolic}, \cite{stratmann:urbanski:diophantine}, 
and the analysis of certain geometrically infinite systems \cite{ctm:renorm:3mfds}.   
A Kleinian group uniformizes a hyperbolic three-manifold, and there is now a candidate three-dimensional 
object associated to a rational map \cite{lyubich:minsky:lamination}, \cite{kaimanovich:lyubich:harmonic}.   
Of course, essential and important differences between the two theories remain.  

Other examples of conformal dynamical systems include iteration of smooth maps of the  
interval to itself and discrete groups of M\"obius transformations acting properly discontinuously on higher-dimensional spheres.   
However, a theorem of Liouville \cite{rickman:qrmappings} asserts that any conformal map in 
dimensions $\geq 3$ is the restriction of a M\"obius transformation.  
Thus, there is no nonlinear classical theory of iterated conformal maps in higher dimensions.  

Two different generalizations of  conformal dynamical systems have been studied.  
One of these retains the Euclidean metric structure of the underlying space and keeps 
some regularity of the iterates or group elements,  but replaces their conformality with uniform quasiregularity.  
Roughly, this means that they are differentiable almost everywhere, and they distort the roundness of  
balls in the tangent space by a uniformly bounded amount.   
In dimension two, Sullivan showed that the Measurable Riemann Mapping Theorem and an averaging process imply 
that each such example is obtained from a rational map or a Kleinian group by a quasiconformal deformation \cite{DS3}.   
In higher dimensions, Tukia gave an example to show that  this fails \cite{tukia:not}.  
The systematic study of uniformly quasiconformal groups of homeomorphisms on $\R^n$ was begun by Gehring and Martin 
\cite{gehring:martin:qcgroupsi}.  They singled out a special class of such groups, the convergence groups, 
which are characterized by topological properties. The subsequent theory of such quasiconformal groups turns out to be quite rich.   
The study of iteration of uniformly quasiregular maps on manifolds is somewhat more recent; see  e.g. \cite{iwaniec:martin:uqr}.  
At present, examples of chaotic sets of such maps are either spheres or Cantor sets, 
and it is not yet clear how rich this subject will be in comparison with that of classical rational maps.  

A second route to generalizing classical conformal dynamical systems is to replace 
the underlying Euclidean space with some other metric space, and to replace the condition 
of conformality with respect to a Riemannian metric with one which makes sense for metrics given as distance functions.   
Technically, there  are many distinct such reformulations--some local, some global, some infinitesimal 
(quasim\"obius, quasisymmetric, quasiconformal).    
An important source of examples with ties to many other areas of mathematics is the following.   
A convex  cocompact Kleinian group acting on its limit set in the Riemann sphere generalizes to a negatively curved group  
(in the sense of Gromov) acting on its boundary at infinity.   
This boundary carries a natural topology and a natural quasisymmetry class of so-called {\em visual} metrics  
\cite{ghys:delaharpe:groupes}, \cite{bonk:schramm:embeddings}. 
With respect to such a  metric, the elements of the group act by uniformly quasim\"obius maps.   
Negatively curved groups acting on their boundaries thus provide a wealth of examples of generalized ``conformal'' dynamical systems.   

Tukia \cite{tukia:convergence_groups} generalized Gehring and Martin's notion of a convergence group 
from spheres to compact Hausdorff spaces, and Bowditch \cite{bowditch:characterization} 
then characterized negatively curved groups acting on their boundaries by purely topological conditions:
\index{Index}{Bowditch characterization}%
\begin{thm}[Characterization of boundary actions]
\label{thm:bowditch_characterization}
Let $\Gamma$ be a group acting on a perfect metrizable compactum $M$ by homeomorphisms. 
If the action on the space of triples is properly discontinuous and cocompact, 
then $\Gamma$ is hyperbolic, and there is a $\Gamma$-equivariant homeomorphism of $M$ onto $\bdry \Gamma$.
\end{thm}
Following Bowditch \cite{bowditch:convergence_groups}  and abusing terminology, we refer to such actions as uniform convergence groups. 
In addition to providing a topological characterization, the above theorem may be viewed as a uniformization-type result.  
Since the metric on the boundary is well-defined up to quasisymmetry, 
it follows that associated to any uniform convergence group action of $\Gamma$ on $M$, there is a preferred  
class of metrics on $M$ in which the dynamics is conformal in a suitable sense:
the action is uniformly quasim\"obius.   

Sullivan referred to convex cocompact Kleinian groups and their map analogs, hyperbolic rational maps, 
as expanding conformal dynamical systems.  Their characteristic feature is the following
principle which we may refer to as the {\em conformal elevator}:  
\begin{quote}
\label{quote:blowup}
{\em Arbitrarily small balls can be blown up via the dynamics to nearly round sets of 
definite size with uniformly bounded distortion, and vice-versa.}   
\end{quote}
This property is also enjoyed by negatively curved groups acting on their boundaries, and is 
the basis for many rigidity arguments in dynamics and geometry.   
Recalling the dictionary, we have then the following table:

\begin{center}
\begin{tabular}{|l|l|} \hline
{\bf Group actions} & {\bf Iterated maps} \\ \hline
Kleinian group & rational map \\ \hline
convex cocompact Kleinian group & hyperbolic rational map\\  \hline
uniform convergence group & ? \\ \hline
\end{tabular}
\end{center}

The principal goal of this work is to fill in the missing entry in the above table.   
To do this, we introduce {\em topological} and {\em metric coarse expanding conformal} (cxc) dynamical systems.   
We emphasize that topologically cxc systems may be locally non-injective, i.e. branched, on their chaotic sets.  
Metric cxc systems are topologically cxc by definition.  
Hyperbolic rational maps on their Julia sets and uniformly quasiregular maps on manifolds with good expanding properties are 
metric cxc.  
Thus, our notion includes both the classical and generalized Riemannian examples of expanding conformal 
dynamical systems mentioned above.   
As an analog of Bowditch's characterization, viewed as a uniformization result, we have the following result:
\begin{thm}[Characterization of metric cxc actions]
\label{thm:characterization_of_metric_cxc_actions}
Suppose $f: X \to X$ is a continuous map of a compact metrizable space to itself.  If $f$ is topologically cxc, 
then there exists a metric $d$ on $X$, unique up to quasisymmetry, such that with respect to this metric, 
$f$ is essentially metric cxc.
\end{thm}
(See Corollary \ref{thm:canonical_gauge}).  
In many cases (e.g. when $X$ is locally connected) we may drop the qualifier  ``essentially''  
from the conclusion of the above theorem.  In general, we cannot.  
It is unclear to us whether this is a shortcoming of our methods, or reflects some key difference between group 
actions and iterated maps;  see \S \ref{secn:properties_following_hyp}.
The naturality of the metric $d$ implies that quasisymmetry invariants of $(X, d)$ 
then become topological invariants of the dynamical system.    
Hence, tools from the theory of analysis on metric spaces may be employed.  
In particular, the conformal dimension (see \S \ref{secn:properties_following_hyp}) becomes 
a numerical topological invariant, distinct from the entropy.  
The  existence of the metric $d$ may be viewed as a generalization of the well-known 
fact that given a positively expansive map of a compact set to itself, there exists a canonical 
H\"older class of metrics in which the dynamics is uniformly expanding.  

Our class of metric cxc systems $f: (X,d)  \to (X,d)$ includes a large number of previously studied types of dynamical systems.    
A rational map is cxc on its Julia set with respect to the standard spherical metric if and only if it is a so-called 
semihyperbolic map (Theorem \ref{thm:newchar_semihyperbolic}).   
A metric cxc map on the standard two-sphere is quasisymmetrically conjugate to a semihyperbolic rational  
map with Julia set the  sphere (Theorem \ref{thm:nothing_new}).  
Using elementary Lie theory, we construct by hand the metric $d$ in the case when $X$ is a manifold and $f$ is an expanding map, 
and show that in this metric $f$ becomes locally a homothety (\S  \ref{secn:expanding_maps}).   
Theorems \ref{thm:cxc_implies_Lattes} and \ref{thm:Lattes_are_cxc} imply that 
uniformly quasiregular maps on Riemannian manifolds of dimension greater 
or equal to $3$ which are metric cxc are precisely the generalized Latt\`es examples of Mayer \cite{mayer:lattes}.  

Just as negatively curved groups provide a wealth of examples of non-classical ``conformal'' group actions, 
so our class of metric cxc maps provides a wealth of examples of non-classical ``conformal'' iterated maps as dynamical systems.   
The case of the two-sphere is of particular interest.  Postcritically finite branched coverings 
of the two-sphere to itself arising from rational maps were characterized combinatorially by Thurston \cite{DH1}.  
Among such branched coverings, those  which are expanding with respect to a suitable orbifold metric give examples 
of topologically cxc systems on the two-sphere.   Hence by our results, 
they are uniformized by a metric such that the dynamics becomes conformal.  
This metric, which is a distance function on the sphere, need not be quasisymmetrically equivalent to the standard one.   
A special class of such examples are produced from the finite subdivision rules on the sphere considered by Cannon, 
Floyd and Parry \cite{cfp:fsr}, \cite{cfkp:rationalmaps}; cf. \cite{meyer:origami}.  
These provide another source of examples of dynamics on the sphere which are conformal 
with respect to non-standard metrics.    
Conjecturally, given a negatively curved group with two-sphere boundary, 
the visual metric is always quasisymmetrically equivalent to the standard one, 
hence (by Sullivan's averaging argument and the Measurable Riemann mapping theorem) the action is isomorphic to that 
of a cocompact Kleinian group acting on the two-sphere.  
This is Bonk and Kleiner's reformulation of Cannon's Conjecture \cite{bonk:kleiner:qsparam}. 

In Theorem \ref{thm:characterization} below, we characterize in several ways when a topologically cxc map on the two-sphere, 
in its natural metric, is quasisymmetrically conjugate to a rational map.   This result was our original motivation. 
The natural metrics associated to a topologically cxc map $f: S^2 \to S^2$ are always linearly locally connected 
(Cor. \ref{cor:get_llc}).  If $f$ is not quasisymmetrically conjugate to a rational map, e.g. 
if $f$ is postcritically finite and has a Thurston obstruction, 
then Bonk and Kleiner's characterization of the quasisymmetry class of the standard two-sphere 
\cite{bonk:kleiner:qsparam} allows us to conclude indirectly that these natural metrics are never Ahlfors 2-regular.   
Recent results of Bonk and Meyer \cite{bonk:meyer:subdivisions}, \cite{bonk:icm:qcgeom} and the authors \cite{haissinsky:pilgrim:cxcIII} suggest that in general, 
Thurston obstructions manifest themselves directly as metric obstructions to Ahlfors 2-regularity in a specific and natural way.  
Differences with the group theory emerge:  we give an example of a metric cxc map on a $Q$-regular two-sphere of 
Ahlfors regular conformal dimension $Q>2$ which is nonetheless not $Q$-Loewner.   
In contrast, for hyperbolic groups, Bonk and Kleiner (\cite{bonk:kleiner:conf_dim}, Theorem 1.3)  
have shown that if the Ahlfors regular conformal dimension is attained,  then the metric is Loewner. 

As mentioned above, the dictionary is rather loose in places.  
From the point of  view of combinatorics and finiteness principles, a postcritically finite subhyperbolic rational map $f$ is a 
reasonable analog of a cocompact Kleinian group $G$.  By Mostow rigidity, $G$  is determined up to M\"obius conjugacy 
by the homotopy type of the associated quotient three-manifold.  
This is turn is determined by the isomorphism type of $G$.  
Since $G$ as a group is finitely presented, a finite amount of combinatorial data determines the geometry of Kleinian group $G$.  
For the analogous rational maps, Thurston \cite{DH1} showed that they are determined up to M\"obius conjugacy by their homotopy 
type, suitably defined.  Recently, Nekrashevych \cite{nekrashevych:book:selfsimilar} introduced tools from the theory 
of automaton groups that show that these homotopy types are again determined by a finite amount of group-theoretical data. 
In a forthcoming work \cite{kmp:ph:cxcii}, we introduce a special class of metric cxc systems that enjoy similar 
finiteness principles.  From the point of view of analytic properties, however, our results suggest that another candidate 
for the analog of a convex cocompact Kleinian group is a so-called semi-hyperbolic rational map, 
which is somewhat more general (\S \ref{secn:semihyperbolic}) and which allows non-recurrent branch points 
with infinite orbits in the chaotic set.   

Our construction of a natural metric associated to a topologically cxc system $f$ proceeds via identifying the chaotic set $X$ 
of the system as the boundary at infinity of a locally finite, negatively curved graph $\Gamma$ with a preferred basepoint.   
By metrizing $\Gamma$ suitably and using the Floyd completion to obtain the metric on the boundary, 
the dynamics becomes quite regular.  The map $f$ behaves somewhat like a homothety:  
there exists a constant $\lambda>1$ such that if $f$ is injective on a ball $B$, then on the smaller ball 
$\frac{1}{4}B$ it multiplies distances by $\lambda$.  
In particular, $f$ is Lipschitz, and (Theorem \ref{thm:construction} and \ref{thm:BPI}) $X$ becomes a  BPI-space
in the sense of David and Semmes \cite{david:semmes:dreams}.  
By imitating the Patterson-Sullivan construction of conformal measures \cite{patterson:measures} as generalized by 
Coornaert \cite{coornaert:patterson-sullivan}, we construct a natural measure $\mu_f$ on 
the boundary with a perhaps remarkable coincidence of properties.  
The measure $\mu_f$ is quasiconformal with constant Jacobian, is the unique measure of maximal entropy $\log \deg(f)$, 
describes the distribution both of backwards orbits and of periodic points, and satisfies Manning's formula 
relating Hausdorff dimension, entropy, and Lyapunov exponents (\S\S \ref{secn:measure_theory} and \ref{secn:properties_following_hyp}). 
Thus, all variation in the distortion of $f$ is ironed out to produce a metric in which the map is in some sense a piecewise homothety, 
much like a piecewise linear map of the interval to itself with constant absolute value of slope.  
In this regard, our results may be viewed as an analog 
of the Milnor-Thurston theorem asserting that a unimodal map with positive topological entropy is semiconjugate 
to a tent map whose slope is the exponential of the entropy \cite{milnor:thurston}.  
Our estimates generalize those of Misiurewicz-Przytycki 
\cite{misiurewicz:przytycki:smooth} and Gromov \cite{gromov:entropy}.  

By way of contrast, Zdunik \cite{zdunik:singular} shows that among rational maps, only the usual family of 
exceptions (critically finite maps with parabolic orbifold) has the property that the measure of maximal entropy 
is equivalent to the Hausdorff measure in the dimension of the Julia set.   
Our construction, however, yields a metric with this coincidence for any rational map which is suitably expanding.

It turns out (Theorem \ref{thm:newchar_semihyperbolic}) that $f$ is semihyperbolic if and only if $\Gamma$ is 
quasi-isometric to the convex hull of the Julia set of $f$ in hyperbolic three-space. 
Lyubich and Minsky \cite{lyubich:minsky:lamination} give a similar three-dimensional characterization of 
this family of maps using hyperbolic three-manifold laminations.  
Analogously, convex cocompact Kleinian groups are characterised by the property 
that their Cayley graphs are quasi-isometric to the convex hull of their limit sets in $\IH^3$.

In summary, we suggest the following enlargement of the above dictionary:

\renewcommand{\arraystretch}{2}
\begin{center}
\begin{tabular}{|l|l|} \hline
{\bf Group actions} & {\bf Iterated maps} \\ \hline
Kleinian group & rational map \\ \hline
convex cocompact Kleinian group & (semi) hyperbolic rational map\\  \hline
uniform convergence group & topologically cxc map \\ \hline
uniform quasim\"obius convergence group & metric cxc map\\ \hline
Cayley graph $\Gamma$ & graph $\Gamma$ \\ \hline 
visual metric & visual metric \\ \hline 
quasiconformal measure $\mu$ & canonical measure $\mu_f$ \\ \hline 
\parbox{3in}{Cannon Conjecture\\ on groups with sphere boundary}  & \parbox{2in}{Thurston's Theorem
\\ characterizing rational maps}  \\ \hline 
\parbox{3in}{Cannon's, Bonk-Kleiner's \\ Characterization Theorems\\ of cocompact Kleinian groups}  & \parbox{2in}{Characterization Theorem for cxc maps \\ on the standard $\IS^2$}\\ \hline 
\end{tabular}
\end{center}
\gap

Our basic method is the following.  Since we are dealing with noninvertible mappings whose chaotic sets are possibly 
disconnected, we imagine the repellor $X$ embedded in a larger, nice  space $\XX_0$ and we suppose that 
$f: \XX_1 \to \XX_0$ where $\cl{\XX}_1 \subset \XX_0$.  We require some regularity on $f$:  
it should be a finite branched  covering.  Our analysis proceeds as follows:
\begin{quote}
\label{quote:pov}
{\em We suppose that the repellor $X$ is covered by a finite collection $\UUU_0$ of open, connected subsets.   
We pull back this covering by iterates of $f$ to obtain a sequence $\UUU_0, \UUU_1, \UUU_2, \ldots$ of coverings of $X$.  
We then examine the combinatorics and geometry of this sequence.}
\end{quote}
The collection of coverings $\{\UUU_n\}$ may be viewed as a discretization of Pansu's quasiconformal structures \cite{pansu:confdim}.  
This motivates our use of the adjective ``coarse'' to describe our metric dynamical systems.  
\gap

\noindent{\bf Contents.}  
In Chapter 2, we begin with the topological foundations needed to define topologically cxc mappings.  
We give the definitions of topologically and metric cxc mappings, prove metric and dynamical regularity properties 
of the repellor, and prove that topological conjugacies between metric cxc systems are quasisymmetric.

In Chapter 3,  we construct the graph $\Gamma$ associated to topologically cxc maps (and to more general maps as well) 
and discuss its geometry and the relation of its boundary with the repellor.  We construct the natural measure and study 
its relation to equidistribution, entropy, and Hausdorff dimension.  
The chapter closes with those properties enjoyed specifically by metric cxc mappings.  

Chapter 4 is devoted to a discussion of examples. It contains a proof of
the topological characterization
of semihyperbolic rational maps among cxc mappings on the two-sphere
(Theorem \ref{thm:characterization}).   We also discuss maps with
recurrent branching and we very briefly point out some formal similarities
between our constructions and analogous constructions in $p$-adic
dynamics.

We conclude with an appendix in which we briefly recall those facts concerning hyperbolic groups 
and convergence groups which served as motivation for this work. 
\gap

\noindent{\bf Notation.}  The cardinality of a set $A$ is denoted by $\#A$ and its closure by $\cl{A}$.  Given a metric space, if $B$ denotes a ball of radius $r$ and center $x$, the notation $\lambda B$ is used for the ball of center $x$ and radius $\lambda r$.  The diameter of a set $A$
is the supremum of the distance between two points of $A$
\index{Index}{diameter} and is denoted by \index{Symbols}{$\diam\, A$}$\diam\, A$.
The Euclidean $n$-sphere, regarded as a metric space, is denoted by $\IS^n$;
\index{Symbols}{$\IS^n$}%
 we use the notation $S^n$
 \index{Symbols}{$S^n$}%
  for the underlying topological space.  
Generally, we will write, for two positive functions, $a\lesssim b$ or $a\gtrsim b$ if there is a universal constant $u>0$ such that $a\le ub$. The notation $a\asymp b$ will
mean $a\lesssim b$ and $a\gtrsim b$.  As usual, we use the symbols $\forall$ and $\exists$ for the quantifiers ``for every'' and ``there exists'' when convenient.  We denote by $\N$ the set of natural numbers $\{0, 1, 2, \ldots\}$, and by $\R_+$ the non-negative real numbers.

\gap
\noindent{\bf Acknowledgements.}  We are grateful for the many opportunities given to present 
our results while this manuscript was in preparation.  We benefited from 
many helpful discussions and encouragment.   In particular we thank M. Bonk, J.-Y. Briend,  
G. Havard, M. Lyubich, V. Mayer, M. Misiurewicz, and the anonymous referee.

The first author 
was partially supported
by the  project ANR ``Cannon'' (ANR-06-BLAN-0366);
he thanks Indiana University for its hospitality.
The second author was supported by U.S. National Science Foundation, Division of Mathematical Sciences grant \#0400852; 
he thanks the Universit\'e de Provence and the LATP for its hospitality where part of the research took place.  
Both authors are also grateful to the IHP which hosted them during the trimester on Dynamical systems Sept.-Nov. 2003.

 \chapter{Coarse expanding conformal dynamics}

The following setup is quite common in the dynamics of noninvertible maps. One is
given a nice, many-to-one map 
\[ f: \XX_1 \to \XX_0\] where $\XX_0$ and $\XX_1$ are nice spaces and $\XX_1 \subset \XX_0$.  One studies  the typically complicated set $X$ of nonescaping points, i.e. points $x \in \XX_1$ for which $f^n(x) \in \XX_1$ for all $n \geq 0$.   We are particularly interested in  maps  for which the restriction of $f$ to $X$ need not be locally injective.   
For those readers unused to noninvertible dynamics, we suggest assuming  
that $\XX_0=\XX_1=X$ upon a first reading.   

A basic method for analyzing such systems is to consider the behavior of small open connected sets of $\XX_0$ under backward, instead of forward, iteration.   For this reason, it is important to have some control on restrictions of iterates of the form $f^k: \wtU \to U$, where $U$ is a small open connected subset of $\XX_0$, and $\wtU$ is a connected component of $f^{-k}(U)$.   Hence it is reasonable to assume that $\XX_0, \XX_1$ are at least locally connected.  The nonescaping set $X$ itself, however, may be disconnected and non-locally connected.   To rule out topological pathology in taking preimages, we impose some tameness restrictions on $f$ by assuming that $f: \XX_1 \to \XX_0$ is a so-called {\em branched covering} between suitable topological spaces.  
When $\XX_0$ is a metric space it is tempting to ask for control over inverses images of metric balls instead of connected open sets.  However, this can be awkward since balls in $\XX_0$ might not be connected.    

We focus  on those topological dynamical systems with good expanding properties.  
However, a map $f: X \to X$ which is not locally injective is never positively expansive, and neither is the induced map on the natural extension.   Thus, notions of expansiveness in this category need to be defined with some care.

\section{Finite branched coverings} 
\label{secn:fbc}

There are have been many different definitions of ramified coverings and branched coverings, most of which coincide in the
context of manifolds (cf. e.g. \cite{fox,edmonds:fbc, dinh:sibony:dynamique}. We define here the notion of finite branched coverings
which suits our purpose: it generalises the topological properties of rational maps of the Riemann sphere, and behaves well
for their dynamical study (e.g. pull-backs of Radon measures are well-defined).

\gap

Suppose $X, Y$ are locally compact Hausdorff spaces, and let $f: X \to Y$ be a finite-to-one continuous map.  
The {\em degree}
\index{Index}{degree}\index{Symbols}{\ensuremath{\deg(f)}}%
of $f$ is  
\[ \deg(f) = \sup\{\# f^{-1}(y): y \in Y\}.\] For $x \in X$, the {\em local degree} 
\index{Index}{degree, local}\index{Symbols}{\ensuremath{\deg(f;x)}}%
of $f$ at $x$ is
\[ \deg(f; x) = \inf_U \sup\{\#f^{-1}(\{z\}) \intersect U : z \in f(U)\}\] where $U$
ranges over all neighborhoods of $x$.   

\begin{defn}[finite branched covering]
The map $f$ is a {\em finite branched covering}
\index{Index}{finite branched covering (fbc)}%
 (abbreviated fbc) provided $\deg(f)<\infty$ and 
\bi 
\item[(i)] $$
\sum_{x \in f^{-1}(y)} \deg(f;x) = \deg f$$ holds for each $y \in Y$;
\item[(ii)] for every $x_0 \in X$, there are compact neighborhoods
$U$ and $V$ of $x_0$ and $f(x_0)$ respectively such that
\[ \sum_{x \in U, f(x)=y} \deg(f; x) = \deg(f; x_0)\]
for all $y\in V$.\ib  
\end{defn}
We note the following two consequences of (ii): the restriction $f:f^{-1}(V)\cap U\to V$
is proper and onto and $f^{-1}(\{f(x_0)\})\cap U=\{x_0\}$.

\gap

The composition of  fbc's is an fbc, and the degrees of fbc's multiply under compositions.  In particular, local degrees of fbc's multiply under compositions. 

Given an fbc $f: X \to Y$, a point $y \in Y$ is a {\em principal value} if $\# f^{-1}(y) = \deg(f)$.   Condition (ii) implies that if $x_n \to x_0$, then $\deg(f; x_n) \leq \deg(f; x_0)$.  It follows that  
the {\em branch set} $B_f = \{x \in X : \deg(f; x)>1\}$ \index{Symbols}{$B_f$} is closed.    The set \index{Symbols}{$V_f$} of {\em branch values} is defined as $V_f = f(B_f)$. 
Thus $Y-V_f$ is the set of principal values.

\begin{lemma}Let $X,Y$ be Hausdorff locally compact topological spaces. An fbc $f:X\to Y$ of
degree  $d$ is open, closed, onto and proper:  the inverse image of a compact subset is compact
and the image of an open set is open. Furthermore, $B_f$ and $V_f$ are nowhere dense.
\end{lemma}

Since the spaces involved are not assumed to be
metrizable, we are led to use filters instead of sequences
in the proof \cite{bour:top_gen}. Recall that a {\it filter base} (or filter basis) of a set $S$ is a collection  
$B$ of subsets of $S$ with the following properties:
   1. the intersection of any two sets of $B$ contains a set of $B$;
   2. the subset $B$ is non-empty and the empty set is not in $B$.

\gap

\pf The map is onto by definition.

\gap

 {\noindent\bf Claim.} {\it For any $x\in X$, let $U(x)$ and $V(x)$ be the neighborhoods of $x$ 
and $f(x)$ given
by (ii).  If $\FFF$ denotes the set of neighborhoods of $f(x)$ contained in $V(x)$, then
$f^{-1}(\FFF)\cap U(x)$ is a filter base converging to $x$. }

\gap

\noindent{\sc Proof of Claim.} Fix $x$ and $y=f(x)$.
Let $\FFF$ be the set of neighborhoods of $y$ contained in $V(x)$.
Since $x$ is  accumulated by $f^{-1}(\FFF)\cap U(x)$,
it follows that if $f^{-1}(\FFF)\cap U(x)$ is not convergent to $x$, then there is another accumulation
point $x'$ of $f^{-1}(\FFF)\intersect U(x)$ in $U(x)$, since $U(x)$ is compact. 
By continuity of $f$, this
implies that $f(x')=y$, so that $x'=x$ since $f^{-1}(\{y\})\cap U(x)=\{x\}$.
Thus, the family of sets $f^{-1}(\FFF)\cap U(x)$ is a filter base converging to $x$. 

This ends the proof of the claim.

\gap

For any $y\in Y$, the set $W(y)=\cap_{f(x)=y} V(x)$ is a compact neighborhood of $y$
since $y$ has finitely many preimages. Let $(N(x))_{x\in f^{-1}(\{y\})}$ be 
compact neighborhoods of $x\in f^{-1}(\{y\})$ which are pairwise disjoint.
It follows from the claim that there is a compact neighborhood $V(y)\subset  W(y)$ 
of $y$  such
that $f^{-1}(V(y))\cap U(x)\subset N(x)$ for all $x\in f^{-1}(\{y\})$.
Therefore, (ii) holds for each pair $(N(x),V(y))$.

\gap

Let $\Omega\subset X$ be an open set, and let us consider $x\in\Omega$ and $y=f(x)$.
We choose a compact neighborhood $N'(x)\subset N(x)\cap\Omega$.
It follows from the claim that a neighborhood
$V'(y)\subset V(y)$ exists  such that $f^{-1}(V'(y))\cap N(x) \subset N'(x)$.  
So, for any $y'\in V'(y)$, by (ii) 
$$\sum_{x'\in f^{-1}(\{y'\})\cap N(x)} \deg(f;x')=\deg(f;x)\ge 1.$$ Hence, $y' = f(x')$ for some $x' \in N'(x) \intersect \Omega$.  Thus $V'(y) \subset f(\Omega)$.  
This establishes that $f$ is open.

\gap

Let us fix $y\in Y$ and let us consider $y'\in V(y)$. Then 
$$\renewcommand{\arraystretch}{2}
\begin{array}{ll} d=\dis\sum_{f(x)=y} \deg(f;x) & =\dis\sum_{f(x)=y}\left( \dis\sum_{x'\in f^{-1}(\{y'\})\cap N(x)} \deg(f;x') \right)\\
& =\dis\sum_{x'\in f^{-1}(\{y'\})\cap (\cup_{f(x)=y} N(x))} \deg(f;x')\,.\end{array}$$
This implies that $f^{-1}(\{y'\})\subset \cup_{f(x)=y} N(x)$. Using the relative compactness and the continuity of
$f$, it follows that the filter base $f^{-1}(\FFF)$ is finer than the set of neighborhoods of $f^{-1}(\{y\})$, 
where $\FFF$ is any filter base converging to $y$.

Let $K\subset Y$ be a compact set and set $L=f^{-1}(K)$. Let $\FFF$ be a filter base in $L$. 
Since $f(L)$ is compact,
there is some accumulation point $y$ in $K$ of $f(\FFF)$. We claim that at least one preimage 
of $y$ is accumulated
by $\FFF$. If it was not the case, then, for any $x\in f^{-1}(\{y\})$, 
there would be some $F_x\in\FFF$ with $x\notin \overline{F_x}$. 
We could
therefore find a compact neighborhood $N'(x)$ of  $x$ such that 
$N'(x)\cap \overline{F_x}=\emptyset$. The claim implies the existence of
 some neighborhood $V'(y)\subset V(y)$ such that 
$f^{-1}(V'(y))\subset (\cup_{f(x)=y}N'(x))$.

Since the fibers are finite, $\cap_{f(x)=y} F_x$ contains a set $F_y\in\FFF$ 
and 
$f^{-1}(V'(y))\cap F_y=\emptyset$.
Hence $V'(y)\cap f(F_y)=\emptyset$, which contradicts that $y$ was an accumulation point
of $f(\FFF)$.
Therefore $f$ is proper. 
\gap

Let us prove that $f$ is closed. Let $Z\subset X$ be a closed set, and let $\FFF$ be any filter base in $f(Z)$ tending to some $y\in Y$.
Fix a compact neighborhood $V$ of $y$ such that $f^{-1}(V)$ is compact,
and consider $\FFF'=\{F\cap V,\ F\in\FFF\}$ to be the trace of $\FFF$ in $V$: this remains a filter base in $f(Z)$ converging to $y$.

Note that, according to what we proved above, $f^{-1}(\FFF')$
is a filter base finer than the set of neighborhoods
of $f^{-1}(\{y\})$. But the trace $\FFF_1=\{F\cap Z, F\in\FFF'\}$ 
of $f^{-1}(\FFF')$ in $Z$ remains
a filter base as well (since it is nonempty),
and  $Z\cap f^{-1}(V)$ is compact: this implies
that $\FFF_1$ accumulates a point $x \in Z\cap f^{-1}(y)$, so $y\in f(Z)$.
Hence $f$ is closed. 

\gap

 The set $V_f$ cannot have interior since $f$ has bounded multiplicity.  Indeed, if $V_f$ had interior,  we could construct a decreasing sequence of open sets $W(y_n)\subset V(y_n)\cap V(y_{n-1})\subset V_f$, so
we would have $p(y_{n+1})\ge p(y_n)+1\ge n+1$, where $p:Y\to \N\setminus\{0\}$ denotes the map that counts the number of preimages of points in $Y$.

Therefore, $B_f$ cannot have interior either since $f$ is an open mapping.\qed

\gap

Many arguments are done using pull-backs of sets and restricting to connected components. It is therefore necessary to work with fbc's defined on sets $X$ and $Y$ enjoying more properties.  The lemma below summarizes results proved in  \cite[\S 2]{edmonds:fbc}.

\begin{lemma}
\label{lemma:edmonds}
Suppose $X$ and $Y$ are locally connected, connected, Hausdorff spaces and $f: X \to Y$ is a finite-to-one, closed, open, surjective, continuous
finite branched covering  map.  
\be
\item If $V \subset Y$ is open and connected,
and $U \subset X$ is a connected component of $f^{-1}(V)$, then $f|U: U \to V$ 
is also a finite branched covering.

\item If $y \in Y$, and $f^{-1}(\{y\}) = \{x_1, x_2, \ldots, x_k\}$, then there exist arbitrarily small connected open neighborhoods $V$ of $y$ such that 
\[ f^{-1}(V)=U_1 \sqcup U_2 \sqcup \ldots \sqcup U_k\]
is a disjoint union of connected open neighborhoods $U_i$ of $x_i$ such that $f|U_i: U_i \to V$ is an fbc of degree $\deg(f; x_i)$, $i=1, 2, \ldots, k$.  

\item If $f(x)=y$, $\{V_n\}$ is sequence of nested open connected sets with $\intersect_n V_n = \{y\}$, and if $\wtV_n$ is the component of $f^{-1}(V_n)$ containing $x$, then $\intersect_n \wtV_n = \{x\}$. 

\eb
\end{lemma}

\section{Topological cxc systems}  
\label{secn:def_top_cxc}   In this section, we state the topological axioms underlying the definition of a cxc system.

Let $\XX_0, \XX_1$\index{Symbols}{$\XX_0$, $\XX_1$} be Hausdorff, locally compact, locally connected topological spaces, each with finitely many connected components.  We further assume that  $\XX_1$ is an open subset of $\XX_0$ and that $\cl{\XX_1}$ is compact in $\XX_0$.    
Note that this latter condition implies that if $\XX_0 = \XX_1$, then $\XX_0$ is compact.  

Let $f: \XX_1 \to \XX_0$ be a  finite branched covering map of degree $d \geq 2$, and for $n \geq 0$ put 
\[ \XX_{n+1} = f^{-1}(\XX_n).\]
Lemma \ref{lemma:edmonds} (1) implies that $f: \XX_{n+1} \to \XX_n$ is again an fbc of degree $d$.  Since $f$ is proper, $\cl{\XX_{n+1}}$ is compact in $\XX_n$, hence in $\XX_0$.  

The {\em nonescaping set}\index{Index}{nonescaping set}, or {\em repellor}\index{Index}{repellor}, of $f: \XX_1 \to \XX_0$ is
\[ X = \{ x \in \XX_1 | f^n(x) \in \XX_1 \;\forall n > 0\} = \bigcap_n \cl{\XX_n}.\]
\index{Symbols}{$X$}
See Figure 2.1. 
\begin{figure}
\label{fig:Repellor}
\psfragscanon
\psfrag{1}{$\mathfrak{X}_0$}
\psfrag{2}{$\mathfrak{X}_1$}
\begin{center}
\includegraphics[width=3in]{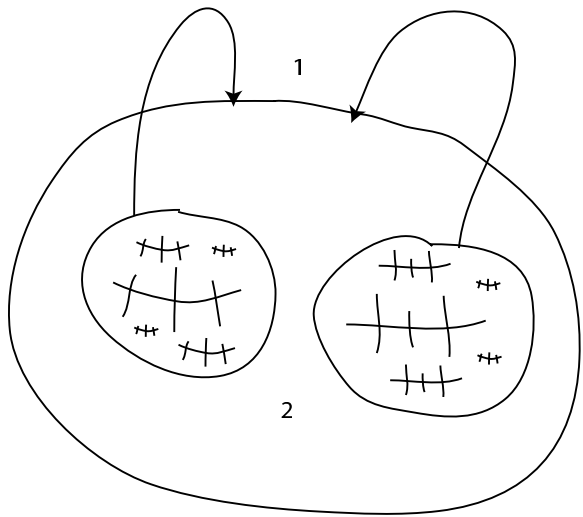}
\end{center}
\caption{{\sf While $\mathfrak{X}_0$ and $\mathfrak{X}_1$ have finitely many components, the repellor $X$ may have uncountably many components.}}
\end{figure}

We make the technical assumption that the restriction $f|X: X \to X$ is also an fbc of degree equal to $d$.  
This implies that  $\#X \geq 2$.   Also,  $X$ is {\em totally invariant}:  $f^{-1}(X) = X = f(X)$.  The definition of the nonescaping set and the compactness of $\cl{\XX_1}$ implies that given any open set $Y$ containing $X$, $\XX_n \subset Y$ for all $n$ sufficiently large.

The following is the essential ingredient in this work.  Let $\UUU_0$ be a finite cover of $X$ by open, connected subsets of $\XX_1$ whose intersection with $X$ is nonempty.  A  {\em preimage}\index{Index}{preimage}
 of a connected set $A$ is defined as a connected component of $f^{-1}(A)$. Inductively, set, for $n\ge 0$,
\[ \UUU_{n+1}=f^{-1}(\UUU_n)=\{\wtU : \exists U \in \UUU_n \ \mbox{with}\ \wtU \ \mbox{a preimage of}\ U\}.\]
\index{Symbols}{$\UUU_n$}
That is, the elements of $\UUU_n$ are the connected components of $f^{-n}(U)$,  where $U$ ranges over $\UUU_0$.

We denote by \index{Symbols}{$\mathbf{U}$}$\mathbf{U} = \union_{n \geq 0} \UUU_n$the collection  of all such open sets thus obtained.  

We say $f: (\XX_1,X) \to (\XX_0, X)$ is {\em topologically coarse expanding conformal with repellor $X$}
\index{Index}{topologically cxc}%
 provided there exists  a finite covering $\UUU_0$ as above, such that the following axioms hold.  

\be  

\item {\bf [Expansion]}  
\index{Index}{[Expansion] Axiom}%
The mesh of the coverings $\UUU_n$ tends to zero as $n \to \infty$.  That is, for any finite open cover $\YYY$ of $X$ by open sets of $\XX_0$, there exists $N$ such that for all $n \geq N$ and all $U \in \UUU_n$, there exists $Y \in \YYY$ with $U \subset Y$.  

\item {\bf [Irreducibility]}  
\index{Index}{[Irreducibility] Axiom}%
The map $f: \XX_1 \to \XX_0$ is {\em locally eventually onto near $X$}:  for any $x \in X$ and any neighborhood $W$ of $x$ in $\XX_0$, there is  some $n$ with $f^n(W) \supset X$

\item {\bf [Degree]} 
\index{Index}{[Degree] axiom}%
The set of degrees of maps of the form $f^k|\wtU: \wtU \to U$, where $U \in \UUU_n$, $ \wtU \in \UUU_{n+k}$, and $n$ and $k$ are arbitrary, has a finite maximum, denoted $p$.  
\eb

Axiom [Expansion] is equivalent to saying that, when $\XX_0$ is a metric space, the diameters of the elements of $\UUU_n$ tend to zero as $n \to \infty$.  Axiom [Irreducibility] implies that $f: X \to X$ is topologically exact; we  give a useful, alternative characterization below.    

These axioms are reminiscent of the following properties of a group $G$ acting on a compact Hausdorff space $X$; see Appendix \ref{appendix:hyperbolic_groups} and \cite{bowditch:convergence_groups}.  
Axiom [Irreducibility] is analogous to $G$ acting minimally on $X$. Axiom  [Expansion] is analogous to $G$ acting properly discontinuously on triples, i.e. that $G$ is a {\em convergence group}.  Axiom [Degree] is analogous to $G$ acting cocompactly on triples; we will see later that this condition implies good regularity properties of metrics and measures 
associated to cxc systems.  

Together, a topologically cxc system we view as the analog, for iterated maps, of a {\em uniform convergence group}.   

The elements of $\UUU_0$ will be referred to as {\em level zero good open sets}\index{Index}{good open sets}%
.  While as subsets of $\XX_0$ they are assumed connected, their intersections with the repellor $X$ need not be.   Also, the elements of $\mathbf{U}$, while connected, might nonetheless be quite complicated topologically--in particular they need not be contractible.  

If $\XX_0 = \XX_1 = X$, then the elements of $\mathbf{U}$ are connected subsets of $X$. \gap

\section{Examples of topological cxc maps}

\subsection{Rational maps}  
\label{subsecn:rational_maps}
Let $\rs$
\index{Symbols}{$\rs$}%
 denote the Riemann sphere, and let $f: \rs \to \rs$ be a rational function 
 \index{Index}{rational map}%
 of degree $d \geq 2$
for which the critical points either converge under iteration to attracting cycles, or land on a repelling periodic cycle 
(such a function is called {\em subhyperbolic}\index{Index}{rational map, subhyperbolic}%
).  For such maps, every point on the sphere belongs either to the {\em Fatou set} and converges to an attracting cycle, or belongs to the {\em Julia set}\index{Index}{Julia set}%
  $J(f)$.  One may find a small closed neighborhood $V_0$ of the attracting periodic cycles such that $f(V_0) \subset \interior(V_0)$.   Set $\XX_0 = \rs - V_0$ and $\XX_1 = f^{-1}(\XX_0)$.    Then $f: \XX_1 \to \XX_0$ is an fbc of degree $d$, the repellor $X=J(f)$, and $f|X: X \to X$ is an fbc of degree $d$.  

Let $\UUU_0$ be a finite cover of $J(f)$ by open spherical balls contained in $\XX_1$, chosen so small that each ball contains at most one forward iterated image of a critical point.  The absence of periodic critical points in $J(f)$ easily implies that the local degrees of iterates of $f$ are uniformly bounded at such points, and so Axiom [Degree] holds.  Since $J(f)$ can be characterized as the locus of points on which the iterates fail to be locally a normal family, Montel's theorem implies that Axiom [Irreducibility] holds.   Finally, $f$ is uniformly expanding near $X$ with respect to a suitable orbifold metric, and Axiom [Expansion] holds; see \cite{shishikura:tan:mane}, Thm. 1.1(b).    

{\em Latt\`es maps} are a special class of subhyperbolic rational maps defined as follows.  Fix a lattice $\Lambda \subset \C$.  The quotient $\TTT = \C/\Lambda$ is a complex torus, and the quotient of this torus by the involution $z \mapsto -z$ is the Riemann sphere, $\rs$.  Let $\pi: \C \to \rs$ be the composition of the two projections.  A rational map $f: \rs \to \rs$ is called a Latt\`es map
\index{Index}{Latt\`es map}%
 if the degree of $f$ is at least two and there is an affine map $L: \C \to \C$ such that $f \circ \pi = \pi \circ L$; see \cite{milnor:lattes}.  Away from the finite set of critical points, Latt\`es  maps are uniformly expanding on the whole Riemann sphere with respect to the length structure induced by pushing forward the Euclidean metric on the complex plane.   The Julia set of such a map is the whole sphere.  
\gap

\subsection{Smooth expanding partial self-covers} Let $\XX_0$ be a connected complete Riemannian manifold, $\XX_1 \subset \XX_0$ an open submanifold with finitely many components which is compactly contained in $\XX_0$.  Let $f: \XX_1 \to \XX_0$ be a $C^1$ covering map which is {\em expanding}\index{Index}{expanding}%
, i.e. there are constants $c>0, \lambda>1$ such that whenever $f^n(x)$ is defined, $||Df^n_x(v)|| > c\lambda^n||v||$.  
If $X$ denotes the set of nonescaping points, then $f:\XX_1 \to \XX_0$ is topologically cxc with repellor $X$--we may take $\UUU_0$ to be a finite collection of small balls centered at points of $X$.  

One may argue as follows.  Since $X$ is compact, there is a uniform lower bound $r$ on the injectivity radius of $\XX_1$ at points $x \in X$.  Thus, for each $x \in X$, the ball $B(x,r)$ is homeomorphic to an open Euclidean ball; in particular, it is contractible.  Let $\UUU_0$ be a finite open cover of $X$ by such balls.  Since $f: \XX_1 \to \XX_0$ is a covering map, all iterated preimages $\wtU$ of elements $U \in \UUU_0$ map homeomorphically onto their images, so the [Degree] Axiom holds with $p=1$.   Since $f$ is expanding, the diameters of the elements of $\UUU_n$ tend to zero exponentially in $n$, so the [Expansion] Axiom holds.  The restriction $f|X: X \to X$ is clearly an f.b.c.   To verify the [Irreducibility] Axiom, we use an alternative characterization given as Proposition \ref{prop:elementary}(2) below.   Suppose $x \in X$ and $x_0\in \XX_0$.  Since $X \union \{x_0\}$ is compact, there exists $L>0$ such that for all $n$, there exists a path $\gamma_n$ of length a
 t most $L$ joining $f^n(x)$ and $x_0$.  
Let $\widetilde{\gamma}_n$ denote the lift of $\gamma_n$ based at $x$.  The other endpoint $\tilde{x}_n$ of $\widetilde{\gamma}_n$ lies in $f^{-n}(x_0)$. By expansion, the length of $\widetilde{\gamma}_n$ tends to zero.  Hence $\tilde{x}_n \to x$ and so $x$ belongs to the set $A(x_0)$ of accumulation points
of $\cup_{n\ge 0}f^{-n}(x_0)$.

Following Nekrashevych \cite{nekrashevych:book:selfsimilar}, we will refer to the topologically cxc system $f: \XX_1 \to \XX_0$ a {\em smooth expanding partial self-covering}.  A common special case is when $\XX_1, \XX_0$ are connected and the homomorphism $\iota_*: \pi_1(\XX_1) \rightarrow \pi_1(\XX_0)$ induced by the inclusion map $\iota: \XX_1 \hookrightarrow \XX_0$ induces a surjection on the fundamental groups.  In this case, the preimages of $\XX_0$ under $f^{-n}$ are all connected, and the repellor $X$ itself is connected.

One can generalize the above example so as to allow branching, by e.g. working in the category of orbifolds\index{Index}{orbifold}%
; see \cite{nekrashevych:book:selfsimilar}.   

\section{Elementary properties}

\noindent{\bf Conjugacy.}  Suppose $f: \XX_1 \to \XX_0$ and $g: \YY_1 \to \YY_0$ are f.b.c.'s with repellors $X$, $Y$ as in the definition of topologically cxc.  A homeomorphism $h: \XX_0 \to \YY_0$ is  called a {\em conjugacy}
\index{Index}{conjugacy}%
  if it makes the diagram
\[ 
\begin{array}{ccc}
(\XX_1, X) & \stackrel{h}{\longrightarrow} & (\YY_1, Y)
\\
f \downarrow & \; & \downarrow g \\
(\XX_0, X) & \stackrel{h}{\longrightarrow} &(\YY_0, Y) \\
\end{array}
\]
commute.  (Strictly speaking, we should require only that $h$ is defined near $X$; however, we 
will not need this more general point of view here.)  

It is clear that the property of being topologically cxc is closed under conjugation:  if $\UUU_0$ is a set of good open sets at level zero for $f$, then $\VVV_0 = \{V=h(U) | U \in \UUU_0\}$ is a set of good open sets at level zero for $g$.  

\gap

Suppose $\XX_1, \XX_0$ are Hausdorff, locally compact, locally connected topological spaces, each with finitely many connected components, $\XX_1 \subset \XX_0$ is open,  and $\cl{\XX_1} \subset \XX_0$. 

The proofs of the following assertions are straightforward consequences of the definitions.  
\begin{prop} 
\label{prop:elementary}
Suppose $f: \XX_1 \to \XX_0$ is an fbc of degree $d\ge 2$ with nonescaping set $X$ and let $\UUU_0$ be a finite open cover of $X$.
\be
\item The condition that $f|_X: X \to X$ is an fbc of degree $d$ implies that the set $V_f \intersect X$ is nowhere dense in $X$.  

\item Axiom [Expansion] implies 
\be
\item $\mathbf{U}$ is a basis for the subspace topology on $X$.  In particular, if $U \intersect X$ is connected for all $U \in \mathbf{U}$, then $X$ is locally connected.

\item For  distinct $x, y \in X$, there is an $N$ 
such that for all $n>N$, and all $U \in \UUU_n$, $\{x, y\} \not\subset 
U$.

\item There exists $N_0$ such that for all $U_1', U_2' 
\in \UUU_{N_0}$, $U_1' \intersect U_2' \neq \emptyset \implies \exists 
U \in \UUU_0$ with $U_1' \union U_2' \subset U$.

\item Periodic points are dense in $X$.  
\eb

\item Axiom [Irreducibility] 
\be
\item holds if and only if for each $x_0 \in \XX_0$, the set $A(x_0)$ of limit points of $\union_{n \geq 0} f^{-n}(x_0)$ equals the nonescaping set $X$.  

\item  implies that either $X=\XX_0=\XX_1$, or $X$ is nowhere dense in $\XX_0$.

\item together with  $f|_X: X \to X$ is an fbc of degree $d$, implies that $X$ is perfect, 
\index{Index}{perfect}%
i.e. contains no isolated points.\eb

\item The class of topologically cxc systems is closed under taking Cartesian products.  
\eb
\end{prop}

In the remainder of this section, we assume $f: \XX_1 \to \XX_0$ is topologically cxc with repellor $X$ and level zero good open sets $\UUU_0$.  

To set up the next statement, given $U \in \UUU_n$ mapping to $U \in \UUU_0$ under $f^n$, denote by \index{Symbols}{$d(U)$}$d(U)=\deg(f^n|U)$ if $n \geq 1$ and $d(U)=1$ if $n=0$.   

\begin{prop}[Repellors are fractal]
\label{prop:repellors_are_fractal}
For every $x \in X$, every neighborhood $W$ of $x$, every $n_0 \in \N$, and every $U \in \UUU_{n_0}$, there exists a preimage $\wtU \subset f^{-k}(U)$ with $\cl{\wtU}\subset W$ and $\deg(f^k:\wtU \to U) \leq \frac{p}{d(U)}$ where $p$ is the maximal degree obtained in the [Degree] Axiom. 
\end{prop}

\pf  Let $\YYY$ be an open cover of $X$ with the property that (i) $W \in \YYY$ and (ii) there exists a neighborhood $W' \subset W$ of $x$ such that for all $Y \neq W$ in $\YYY$, $Y \intersect W' = \emptyset$.  Axiom [Expansion] then implies that there exists $n_1 \in \N$ such that for all $n \geq n_1$, any element of $\UUU_{n_1+n_0}$ containing $x$ is contained in $W$.  Axiom [Irreducibility] implies that there exists $n_2$ such that $f^n(W')=X \supset U$ for all $n \geq n_2$.  Hence for  $k = \max\{n_1, n_2\}$,  there is a preimage $\wtU$ of $U$ under $f^{-k}$ contained in $W$. The assertion regarding degrees follows immediately from the multiplicativity of degrees under compositions.  
\qed
\gap

\noindent{\bf Post-branch set.}  
The {\em post-branch set}
\index{Index}{post-branch set}\index{Symbols}{$P_f$}%
 is defined by 
\[ P_f = X \intersect \cl{\bigcup_{n>0} V_{f^n}}.\]

\begin{prop}
\be
\item A point $x \in X$ belongs to $X-P_f$ if there exists $U \in \mathbf{U}$ such that all preimages of $U$ under iterates of $f$ map by degree one onto $U$. 
\item The post-branch set is a possibly empty, closed, forward-invariant, nowhere dense subset of $X$.
\eb
\end{prop}

Without further finiteness hypotheses on the local topology of $\XX_0$, we do not know if the converse of (1) holds, i.e. if every point in the complement of the post-branch set has a neighborhood over which all preimages under all iterates map by degree one, as is the case for e.g. rational maps.
\gap   

\pf {\bf 1.}  If such a $U$ exists, then $U \intersect V_{f^n}=\emptyset$ for all $n$ and so $x \not\in P_f$.

{\bf 2.}  All but the last assertion are clear.  To show $P_f$ is nowhere dense, let $x \in X$ and let $W$ be any neighborhood of $x$ in $\XX_0$.  Let $\wtU \subset W$ be the element of $\mathbf{U}$ given by Proposition \ref{prop:repellors_are_fractal} applied with a $U$ chosen so that $d(U)=p$.  Then all further preimages of $\wtU$ map by degree one and so $\wtU \intersect V_{f^n}=\emptyset$ for all $n$.  Hence $\wtU \intersect P_f = \emptyset$.  Finally, since $\wtU \intersect X \neq \emptyset$ we conclude that $W \intersect (X-P_f) \neq \emptyset$ and so $P_f$ is nowhere dense in $X$.\qed

\section{Metric cxc systems} In this section, we state the definition of metric cxc systems; we will henceforth drop the adjective, metric.  
\gap

\noindent{\bf Roundness.} Let $Z$ be a metric space and let $A$ be a bounded, proper  subset of $Z$ with nonempty interior.  Given $a \in \interior(A)$, let 
\[ L(A,a)=\sup\{|a-b|: b \in A\}\]
and
\[ l(A,a)=\sup\{ r : r \leq L(A,a) \; \mbox{ and } \; B(a,r) \subset A\}\]
denote, respectively, the {\em outradius}\index{Index}{outradius}
and {\em inradius}\index{Index}{inradius}
of $A$ about $a$.  
While the outradius is intrinsic, the inradius depends on how $A$ sits in $Z$.  The condition $r \leq L(A,a)$ is necessary to guarantee that the outradius is at least the inradius. The\index{Index}{roundness}\index{Symbols}{$\roundness(A,a)$}
{\em roundness of $A$ about $a$} is defined as 
\[ \roundness(A,a) = L(A,a)/l(A,a) \in [1, \infty).\]
One says $A$ is {\em $K$-almost-round} if $\roundness(A,a) \leq K$ for some $a \in A$, and this implies that for some $s>0$, 
\[ B(a,s) \subset A \subset  \cl{B(a,Ks)}.\]
\gap

\noindent{\bf Metric cxc systems.}  
\index{Index}{metrically cxc}%
A key feature of many conformal dynamical systems is the fact that small balls can be blown up using the dynamics to sets of definite size which are uniformly $K$-almost-round and such that ratios of diameters are distorted by controlled amounts.  Below, we formulate abstract versions of these properties which make sense in arbitrary metric spaces.  

Suppose we are given a topological cxc system $f:\XX_1 \to \XX_0$ with level zero good neighborhoods $\UUU_0$,  and that $\XX_0$ is now endowed with a metric  
compatible with its topology.  The resulting metric dynamical system equipped with the covering $\UUU_0$  is called {\em
coarse expanding conformal}, \index{Index}{cxc}%
abbreviated cxc, provided there exist
\bi
\item continuous, increasing embeddings $\rho_{\pm}:[1,\infty) \to [1,\infty)$, the {\em forward and backward
roundness distortion functions}\index{Index}{roundness distortion function}\index{Symbols}{$\rho_{\pm}$}%
, and

\item increasing homeomorphisms $\delta_{\pm}:[0,1] \to [0,1]$, the {\em forward and backward relative
diameter distortion functions}
\index{Index}{diameter distortion function}\index{Symbols}{$\delta_{\pm}$}%
\ib
satisfying the following axioms:
\be
\setcounter{enumi}{3}

\item {\bf [Roundness distortion]} For all $n, k \in \N$ and for all
\[ U \in \UUU_n, \;\;\wtU \in \UUU_{n+k}, \;\; \tilde{y} \in \wtU, \;\;  
y
\in U\]
if
\[ f^{\circ k}(\wtU) = U, \;\;f^{\circ k}(\tilde{y}) = y \]
then the {\em backward roundness bound}
\begin{equation}
\label{eqn:backward_roundness_bound}
  \roundness(\wtU, \tilde{y}) <
\rho_-(\roundness(U,y))
\end{equation}
and the {\em forward roundness bound}
\begin{equation}
\label{eqn:forward_roundness_bound}
\roundness(U, y) <
\rho_+(\roundness(\wtU,
\tilde{y})).
\end{equation}
hold.  

In other words: for a given element of $\mathbf{U}$,
iterates of $f$ both forward and backward distorts its roundness by an
amount independent of the iterate.

\item {\bf [Diameter distortion]} For all $n_0, n_1, k \in \N$ and for all
\[ U \in \UUU_{n_0}, \;\;U' \in \UUU_{n_1}, \;\;\wtU \in \UUU_{n_0+k},
\;\;\wtU'
\in
\UUU_{n_1+k}, \;\; \wtU' \subset \wtU, \;\; U' \subset U\]
if
\[ f^k(\wtU) = U, \;\;f^k(\wtU') = U'\]
then
\[ \frac{\diam\wtU'}{\diam\wtU} < \delta_-\left(\frac{\diam U'}{\diam
U}\right)\]
and
\[ \frac{\diam U'}{\diam U} < \delta_+\left(\frac{\diam \wtU'}{\diam
\wtU}\right)\,.\]

In other words:  given two nested elements of $\mathbf{U}$, iterates of  
$f$
both forward and backward distort their relative sizes by an amount
independent of the iterate.

As a consequence, one has then also the {\em backward upper and lower relative diameter
bounds}:
\begin{equation}
\label{eqn:brdb}
\delta_+^{-1}\left(\frac{\diam U'}{\diam U}\right) <
\frac{\diam\wtU'}{\diam \wtU} < \delta_-\left(\frac{\diam U'}{\diam
U}\right)
\end{equation}
and the {\em forward upper and lower relative diameter bounds}:
\begin{equation}
\label{eqn:frdb}
\delta_-^{-1}\left(\frac{\diam\wtU'}{\diam\wtU}\right) < \frac{\diam
U'}{\diam U} < \delta_+\left(\frac{\diam\wtU'}{\diam\wtU}\right).
\end{equation}
\eb

Axiom [Expansion] implies that the maximum diameters of the elements of $\UUU_n$ tend to zero uniformly in $n$.   Since $\UUU_0$ is assumed finite, each covering $\UUU_n$ is finite, so for each $n$ there is a  minimum diameter of an element of $\UUU_n$.  Since $X$ is perfect and, by assumption, each $U \in \mathbf{U}$ contains a point of $X$, each $U$ contains many points of $X$ and so has positive diameter.   Hence there exist decreasing positive sequences $c_n, d_n \to 0$ such that the {\em diameter bounds} hold:
\begin{equation}
\label{eqn:diameter_bounds}
0 < c_n \leq  \inf_{U \in \UUU_n} \diam U \leq \sup_{U \in \UUU_n} \diam U \leq  d_n.
\end{equation}
\index{Symbols}{$c_n$}%
\index{Symbols}{$d_n$}%

\section{Metric regularity of cxc systems}
\label{secn:metric_regularity}

Suppose now $\XX_0, \XX_1$ are metric spaces.  Let $f: \XX_1 \to \XX_0$ be an f.b.c. as in the previous section with repellor $X$ and level zero good neighborhoods $\UUU_0$.   Throughout this subsection, we assume that the topological axiom [Expansion] and the metric axioms [Roundness distortion] and [Relative diameter distortion] are satisfied.   We  assume neither axiom [Irreducibility] nor [Degree].

In this section, we derive metric regularity properties of the elements of the coverings $\UUU_n$ and the repellor $X$.
\gap 

\noindent{\bf A word regarding notation and the strategy of the proofs.}  In this and the following section, 
$U$ will always denote an element of $\mathbf{U}=\union_n \UUU_n$.  
Generally, $\widetilde{(\cdot )}$ denotes an inverse image of $(\cdot)$ under some iterate of $f$.  
Often, but not always, $U'$ denotes an element of $\mathbf{U}$ which is contained in $U$.   \
Many of the statements of the propositions below make reference to an element $U$ of $\mathbf{U}$.  
The typical proof consists of renaming $U$ as $\wtU$, 
mapping $\wtU$ forward via some iterate to an element $U$ of definite size, making estimates, 
and then transporting these estimates to $\wtU$ via the distortion functions.  
Our estimates, and the implicit constants will always be independent from the iterates
of the map, unless otherwise stated.

\gap

We first resolve a technicality.

\begin{prop}
\label{prop:technicality}
Let $D_0$ denote the minimum diameter of a connected component of $\XX_0$.  Then for any ball $B(a,r)$ in $\XX_0$ where $r\leq D_0/2$, we have $\diam B(a,r) \geq r$.  
\end{prop}

\pf  Fix $\epsilon>0$, and let $C$ denote the component of $\XX_0$ containing $a$.  Pick $p,q\in C$ with $|p-q|>D-\epsilon$.  Then at least one of $|a-p|, |a-q|$ is at least
$ (D-\epsilon)/2$, say $|a-p|$.  Since $C$ is connected, the function $y \mapsto |a-y|$ has an image which contains $[0,(D-\epsilon)/2]$.  Thus for any $s\leq (D-\epsilon)/2$, there exists $y\in C$ with $|a-y|=s$.  Letting $\epsilon \to 0$ proves that $B(a,r)$ has diameter at least $r$.  
\qed

{\em When dealing with balls below, we shall always assume that $r<D_0/2$.}   

\gap

{\noindent\bf Lebesgue number.} 
Let $\UUU$ be a finite covering of a compact metric space $X$ by open sets. The {\em Lebesgue number}
\index{Index}{Lebesgue number}%
 $\delta$ of the covering is the supremum
over all radii $r$ such that, for any point $x\in X$, there is some element $U\in\UUU$ which contains $B(x,r)$.   
Since the covering is finite, $\delta$ is positive. 

\begin{prop}[Uniform roundness]
\label{prop:uniform_roundness}
There exists $K>1$ such that 
\be
\item For all $\tx \in X$ and $n \in \N$, there exists $U \in
\UUU_n$ such
that $U$ is $K$-almost round 
with respect to $\tx$, i.e. $(\exists r>0)$
\[ B(\tx, r) \subset U \subset \cl{B(\tx, Kr)}.\]
\item  For all $n \in \N$ and for all $U \in \UUU_n$, there exists $
\tx
\in X$ such that $\roundness(U,\tx) < K$.
\eb
\end{prop}

\pf
\noindent {\bf (1).}  Denote the set we are looking for by $\wtU$ instead of $U$.  Let $\delta$ be the Lebesgue number of the covering $\UUU_0$ 
and 
$\Delta = \sup_{U \in \UUU_0} \diam U$. Then given any
$x \in X$, there exists $U \in \UUU_0$ such that
\[ B(x,\delta) \subset U \subset B(x,\Delta)
\implies \roundness(U,x) < K_1:=\frac{\Delta}{\delta}.\]
Now let $\tilde{x} \in X$ and $n \in \N$ be arbitrary.  
Set $x=f^n(\tilde{x})$
and let $U \in \UUU_0$ be the element constructed as in the previous paragraph.
Let $\wtU \in \UUU_n$ be the component of $f^{-n}(U)$ containing $\tx$.  
By the backward
roundness bound (\ref{eqn:backward_roundness_bound}),
\[ \roundness(\wtU, \tilde{x}) < \rho_-(K_1).\]

\noindent {\bf (2).}  Denote the given element of $\UUU_n$ by $\wtU$ instead of $U$.  For each $U \in \UUU_0$, choose $x_U \in U$ arbitrarily.  
Let
\[ K_2 = \max_{U \in \UUU_0} \roundness(U, x_U).\]  
Given $\wtU \in \UUU_n$ arbitrary, let $U=f^n(\wtU)$, and let
$\tilde{x}_U \in f^{-n}(x_U) \intersect \wtU$.  
By the backward roundness
bound,
\[ \roundness(\wtU, \tilde{x}_U) < \rho_-(K_2).\]
Thus the conclusions of the lemma are satisfied with 
\[ K=\max\left\{\rho_-(K_1),
\rho_-(K_2)\right\}.\]
\qed

In the lemma below, let $K$ denote the constant in Proposition \ref{prop:uniform_roundness} and $c_n$ the constants giving the lower bound on the diameters of the elements of $\UUU_n$, cf. (\ref{eqn:diameter_bounds}).  

\begin{prop}[Lebesgue numbers]
\label{prop:lebesgue_numbers}
 For all $n\in \N$, all $x \in X$, and all $0< r < \frac{c_n}{2K})$, there exists $U \in \UUU_n$ and $s>r$ such that 
\[ B(x,r) \subset B(x,s) \subset U \subset B(x,Ks).\]
In particular, the Lebesgue number of the covering $\UUU_n$ is at least $\delta_n = \frac{c_n}{2K}$.
\end{prop}

\pf Given $n$ and $x$, by Proposition \ref{prop:uniform_roundness} there is $s>0$ and $U \in \UUU_n$ with 
\[ B(x,s) \subset U \subset B(x,Ks).\]
Thus $c_n < \diam U < 2Ks$ so that $ \frac{c_n}{2K} < s$, whence $ r<s$.

\qed

The next statement says that two elements of covers which intersect over $X$ have roughly the same 
diameter as soon as their levels are close.

\begin{prop}[Local comparability]
\label{prop:local_comparability}
There exists a constant $C>1$ such that
for all $x \in X$, all $n \in \N$, all $U \in \UUU_n$, and all $U'
\in \UUU_{n+1})$ we have:  
if $U \intersect U'\intersect X \neq \emptyset$ then
\[ \frac{1}{C} < \frac{\diam U'}{\diam U} < C.\]
That is, two elements of $\mathbf{U}$ at consecutive levels which
intersect at a point of $X$ are nearly the same size.
\end{prop}

\pf By Axiom [Expansion] there exists $n_0 \in \N$ such that
$2(d_{n_0}+d_{n_0+1})$ is less 
than the Lebesgue number of the covering $\UUU_0$.  
Thus there exist $r>0$ and $n_0$ such that whenever 
$U \in
\UUU_{n_0}$ and $U' \in \UUU_{n_0+1}$ contain a common point $x \in X$,
there exists $V \in \UUU_0$ depending on the pair $U, U'$  such that
\[ U \union U' \subset B(x, r)  \subset V.\]
By renaming as usual, let $\wtU \in \UUU_{n_0+n},
\wtU' \in \UUU_{n_0+n+1}$ denote respectively the sets $U, U'$ as in the statement of the lemma,  
and suppose $\tilde{x} \in \wtU \intersect \wtU' \intersect X$.  
Set
$U = f^{\circ n}(\wtU)$, $U' = f^{\circ n}(\wtU')$, $x=f^{\circ n}(\tilde{x})$  
and
let
\[ S = \sup_{U \in \UUU_{n_0}, U' \in \UUU_{n_0+1}} \max \left\{
\frac{\diam U}{\diam V}, \frac{\diam U'}{\diam V}\right\}.\]
Note that $S$ depends only on the integer $n_0$.  
If $\wtV$ denotes the preimage of $V$ under $f^{-n}$ containing $\wtU \union \wtU'$, then the backwards relative diameter bounds (\ref{eqn:brdb}) imply 
\[ \delta_+^{-1}(S)  \frac{\diam \wtU'}{\diam \wtV} < \delta_-(S)\]
and
\[  \delta_+^{-1}(S) < \frac{\diam \wtU}{\diam \wtV} < \delta_-(S).\]
Dividing yields,
\[ \frac{\delta_+^{-1}(S)}{\delta_-(S)} < \frac{\diam\wtU'}{\diam\wtU} <
\frac{\delta_-(S)}{\delta^{-1}_+(S)}.\] 
Since $S$ and $n_0$ are  
independent of
$n$, the conclusion follows with
\[ C = \max
\left\{
\left(
\frac{\delta_+^{-1}(S)}{\delta_-(S)}
\right)
^{\pm 1}
,\;\;
\frac{
\sup\{\diam U | \; U \in \union_0^{n_0}\UUU_n\}
}
{
\inf\{\diam U | \; U \in \union_0^{n_0}\UUU_n\}
}
\right\}.
\]
\qed

The following lemma  shows that cxc systems are truly expanding in a natural metric sense, and that the
$\delta_-$ function depends essentially only on the relative levels of the sets involved.

\begin{prop}[Contraction implies exponential contraction]
\label{prop:exponential_contraction}
Constants  $C'>0$ and $\theta \in (0,1)$ exist such that, for any $n,k\ge 0$, any $U'\in \UUU_{n+k}$ and any $U\in\UUU_n$,
if $U'\cap U\cap X\ne \emptyset$, then $$\frac{\diam U'}{\diam U}\le C' \theta^{k}\,.$$
In particular, in the upper diameter bounds (\ref{eqn:diameter_bounds}), one may assume $d_n = C'd_0\theta^{n}$.
\end{prop}

\pf The diameters
of the elements of $\UUU_0$ are bounded from below by the constant $c_0$.
Since the diameters of the elements of $\UUU_n$ tend uniformly to zero (by the [Expansion] Axiom), 
and the backwards relative diameter distortion function $\delta_-$  is a homeomorphism, there exists
$N_0 \in \N$ with the following property:
\[ (\forall U' \in \UUU_{N_0})(\exists U \in \UUU_0) \mbox{ such that } 
U' \subset U\;
\mbox{ and }\delta_-\left(\frac{\diam U'}{\diam U}\right) < 
\frac{1}{2}.\]
Now let $k \in \N$ be arbitrary and let $\wtU' \in \UUU_{k+N_0}$.
Let $U' = f^k(\wtU')$, 
let $U \supset U'$ be as above, 
and let $\wtU$ be the component of $f^{-N_0}(U)$ containing $\wtU'$.  
Then by the
backwards relative diameter bounds (\ref{eqn:brdb}), 
\[ \diam \wtU' < \frac{1}{2} \diam \wtU.\]
Thus, for any $k\in\N$, for any $U'\in\UUU_{N_0+k}$, there exists $U\in\UUU_k$ such that $U'\subset U$ and
$\diam U'\le (1/2)\diam U$.

\gap

Let us set $\theta=2^{-1/N_0}$ and $C'=2C^{N_0-1}$ where $C$ is given by Proposition \ref{prop:local_comparability}.

\gap

Let $n,k\ge 0$, and let us fix $U'\in\UUU_{n+k}$ and $U\in\UUU_n$ such that $U\cap U'\cap X\ne \emptyset$.
There are  integers $a\ge 0$, $b\in\{0,\ldots, N_0-1\}$ such that
$k= a\cdot N_0+b$. 

Define inductively $U_j\in \UUU_{jN_0+b}$, $j=0,\ldots, a$, such that $U_a=U'$, $U_{j+1}\subset U_j$ and $\diam U_{j+1}\le (1/2)\diam U_j$.
It follows that $$\diam U'\le \left(\frac{1}{2}\right)^a\diam U_0\le C^b\left(\frac{1}{2}\right)^a\diam U$$
by Proposition \ref{prop:local_comparability}. But $$2^{-a}=\theta^k 2^{b/N_0}\le 2\theta^k$$
so the proposition follows.\qed

The lemma below shows that in $\XX_0$, a possibly disconnected ball $B(x,r)$ with $x \in X$ can be both enlarged and shrunk to obtain a pair of elements $U, U'$ of $\mathbf{U}$ whose levels are equal up to an additive constant and whose diameters in $\XX_0$ are equal  to the diameter of $U$ up to a multiplicative constant.

\begin{prop}[Balls are like connected sets (BLC)]
\label{prop:BLC}
There exist constants $L>1$ and $n_0 \in \N$ such that for all $x \in X$ and $r<\delta_0$, there exist levels $m$ and $n$ and sets 
$U \in \UUU_n$, $U' \in \UUU_{m}$ such that $|m-n|\le n_0$ and
\[ B(x, r/L) \subset U' \subset B(x,r) \subset U \subset B(x,Lr).\]
\end{prop}

\pf We will first find $U$ and $L$ so that $B(x,r) \subset U \subset B(x,Lr)$, where $L=4KC$, $K$ is the roundness constant from Proposition \ref{prop:uniform_roundness}, and $C$ is the constant from Proposition \ref{prop:local_comparability}.  

Let $\delta_0$ denote the Lebesgue number of $\UUU_0$.  Given $x$ and $r<\delta_0$, the number 
\[ n = \sup\{i : \exists U, \exists i \mbox{ with } B(x, r) \subset U \in
\UUU_i \mbox{ and
}\roundness(U,x) < K\}\] 
exists.  (The set is nonempty by Proposition \ref{prop:lebesgue_numbers} and finite by the [Expansion] Axiom.)
Suppose $U \in \UUU_n$ and $B(x,r) \subset U$.  
We must bound $\diam U$ from above.

By Proposition [\ref{prop:uniform_roundness}, Uniform roundness], there exists $V \in \UUU_{n+1}$ for which
$\roundness(V, x) < K$.
Thus, $B(x, s) \subset V \subset B(x, Ks)$ for some $s$.  Since $n$ is
maximal, $s<r$, and so
$\diam V < 2Ks < 2Kr$.  By  Proposition 
[\ref{prop:local_comparability}, Local comparability], 
$\diam U <  C\diam V < C2Ks<2KCr$ and so $U
\subset  B(x, 4KCr)$ as required.  Thus, we have found $U$.  

The same argument applied to $B(x,r/L)$ produces $U'$ such that $B(x,r/L)\subset U'\subset B(x,r)$. 

Assume $U' \in \UUU_m$ and $U \in \UUU_n$.  If $m=n+k$ where $k \geq 0$, Proposition \ref{prop:exponential_contraction} implies $k \leq -\log(2L^2C')/\log \theta$.   If  $n=m+k$ where $k\geq 0$, then another application of the proposition (with the roles of $U$ and $U'$ reversed) yields $k \leq -\log(2C'/\log\theta)$.  The factors of two arise since the diameter of a ball of radius $r$ is bounded below by $r$, not $2r$ (Proposition \ref{prop:technicality}).  
\qed

Recall that a metric space is {\em uniformly perfect}
\index{Index}{uniformly perfect}  if there is a positive constant $\lambda  <1$ such that $B\setminus (\lambda B)$ is non-empty for every ball $B$ of radius at most the diameter of the space.

\begin{prop}
\label{prop:X_is_up}
We have $\diam (U \intersect X) \asymp \diam U$ for all $U \in \mathbf{U}$.  As a consequence, 
the repellor $X$ is uniformly perfect.
\end{prop}

\pf Recall that $X$ is perfect, i.e. contains no isolated points.  

We first claim that there exists some constant $c\in (0,1)$ such that 
\begin{equation}
\label{eqn:diamUintersectX}
\forall U \in \mathbf{U}, \ c \cdot \diam U \leq  \diam(U \intersect X).
\end{equation}

There exists $n_0$ large enough such that for each $U \in \UUU_0$, we 
can choose points $a, b \in U \intersect X$ and neighborhoods $U'_a, 
U'_b \in \UUU_{n_0}$ of $a$ and $b$, respectively, which are disjoint 
and contained in $U$.  We now assume such a choice has been fixed.
Given $\wtU \in \UUU_k$, let $U = f^k(\wtU)$ and let $a, b, U'_a, U'_b$ 
be as in the previous paragraph.   Choose arbitrarily $\tilde{a} \in 
\wtU \intersect f^{-k}(a)$ and let $\wtU'_{\tilde{a}} \in\UUU_{n_0+k}$ 
be the unique component of $f^{-k}(U'_a)$ containing $\tilde{a}$.  
Similarly, define $\tilde{b}$ and $\wtU'_{\tilde{b}}$.   Then 
$\wtU'_{\tilde{a}}$ and $\wtU'_{\tilde{b}}$ are disjoint and are 
contained in $\wtU$.  Each contains an element of $X$, since $X$ is 
totally invariant.  Thus, $\diam (\wtU\intersect X)$ is at least as large 
as the radius $r$ of the largest ball centered at $\tilde{a}$ and 
contained in $\wtU'_{\tilde{a}}$.  By the definition of roundness
\[ r > \frac{1}{2} \diam \wtU'_{\tilde{a}} \cdot 
\roundness(\wtU'_{\tilde{a}}, \tilde{a})^{-1}.\]
The backward relative diameter distortion bounds (\ref{eqn:brdb}) imply
\[  \diam \wtU'_{\tilde{a}} > \diam\wtU \cdot \delta_+^{-1}\left( 
\frac{\diam U'_a}{\diam U}\right).\]
The backward roundness distortion bound 
(\ref{eqn:backward_roundness_bound}) implies
\[ \roundness(\wtU'_{\tilde{a}}, \tilde{a}) < \rho_-(\roundness(U_a, 
a)).\]
Since $\UUU_0$ is finite, $r/\diam\wtU$ is therefore bounded from below by a positive
constant $c$ independent of $k$.

We now prove the proposition.   Let $B=B(x,r)$ be any ball centered at a point $x \in X$.  By Proposition [\ref{prop:BLC}, BLC] there exists $U \in \mathbf{U}$ with $U \subset B$ and $B(x, r/L) \subset U \subset B(x,r)$.  By (\ref{eqn:diamUintersectX}) there exists $x_1, x_2 \in X$ with $|x_1-x_2| > c\cdot \diam U > cr/L$.  At least one $x_i$ must lie outside of $B(x,cr/2L)$, so $X$ is $\lambda$-uniformly perfect where $\lambda = \frac{c}{2L}$.  
\qed

\begin{defn}[Linear local connectivity]
\label{defn:llc}
 Let $\lambda \geq 1$.  
A metric space $Z$ is $\lambda$-{\em linearly locally connected}
\index{Index}{linear locally connected}%
 if the following two conditions hold:
\be
\item if $B(a,r)$ is a ball in $Z$ and $x,y \in B(a,r)$, then there exists a continuum $E \subset B(a, \lambda r)$ containing $x$ and $y$; 
\item if $B(a,r)$ is a ball in $Z$ and $x,y \in Z-B(a,r)$, then there exists a continuum $E \subset Z-B(a, r/\lambda)$ containing $x$ and $y$.  
\eb 
\end{defn}

Propositions \ref{prop:BLC} and \ref{prop:X_is_up} imply immediately that (i) if $U \intersect X$ is connected for all $U \in \mathbf{U}$, then condition (1) above holds, and (ii) if $X \setminus (U \intersect X)$ is connected for all $U \in \mathbf{U}$, then condition (2) holds.  We obtain immediately 

\begin{cor}
\label{cor:get_llc}
If, for all $U \in \mathbf{U}$, the sets $U \intersect X$ and $X \setminus (U \intersect X)$ are connected, then $X$ is linearly locally connected.  
\end{cor}

Unlike the preceding results in this section, the following lemma uses the [Degree] Axiom.  
Recall that a metric space is {\em doubling}
\index{Index}{doubling}%
 if there is a positive integer $C_d$ such that any set of finite diameter
can be covered by $C_d$ sets of at most half its diameter (cf. \cite{heinonen:analysis}, \S\,10.13).

\begin{prop}[cxc implies doubling]
\label{prop:cxc_implies_doubling}
If axiom [Degree] is satisfied, then $X$ is a doubling metric space.
\end{prop}

\pf It follows from Proposition \ref{prop:exponential_contraction} that an integer $k_0$ exists such that, for any $n\ge 0$, 
any $U\in \UUU_n$, and any $U'\in\UUU_{n+k_0}$, we have $\diam U'\le (1/4L)\diam U$ as soon as $U'\cap U\cap X\ne \emptyset$. 

>From the finiteness of $\UUU_0$, it follows  that any $U\in\UUU_0$ can be covered by $N$ sets of level $k_0$.

Let $E\subset X$, and $x\in E$. If its diameter is larger than the Lebesgue number of $\UUU_0$, then it can be covered by a uniform number
of sets of half its diameter. Otherwise, one can find a level $n$ and a set $\wtU\in \UUU_n$ such that 
$$E\subset B(x,\diam E)\subset \wtU\subset B(x, L\diam E)$$ by Proposition \ref{prop:BLC}.

Let us cover $f^n(\wtU)$ by  $N$ sets $U_1',\ldots, U_N'$ of level $k_0$. Axiom [Degree] implies that $\wtU$, so $E$ as well,
is covered by at most $pN$ sets $\wtU'_j$ of level $n+k_0$. Thus, $$\diam \wtU_j'\le \frac{1}{4L}\diam \wtU\le \frac{2L}{4L}\diam E,$$
and so we may take $C_d\le pN$.\qed

From Assouad's theorem (see \cite{heinonen:analysis}, Thm. 12.1) we obtain 

\begin{cor}
If axiom [Degree] is satisfied, then $X$ is quasisymmetrically embeddable in some Euclidean space $\R^n$.
In particular, $X$ has finite topological dimension.  
\end{cor}

The definition of a quasisymmetric map is given below in \S 2.8.

\section{Dynamical regularity}
\label{secn:dynamical_regularity}

Suppose again that $\XX_0, \XX_1$ are metric spaces.  Let $f: \XX_1 \to \XX_0$ be an f.b.c. as in the previous section with repellor $X$ and level zero good neighborhoods $\UUU_0$.   Throughout this subsection, we again assume that the topological axiom [Expansion] and the metric axioms [Roundness distortion] and [Relative diameter distortion] are satisfied.   We  assume neither axiom [Irreducibility] nor [Degree]. 

Recall that a subset $A$ of a metric space $Z$ is {\em $c$-porous} 
\index{Index}{porous}%
when
every ball of radius $r<\diam Z$ contains a ball of radius $cr$ which does not meet $A$.  A subset is {\em porous} if it is $c$-porous for some $c>0$.

\begin{prop}  If axiom [Degree] is satisfied, then the post-branch set $P_f = \union_{n>0}f^n(B_f)$ is porous,  and the sets $B_{f^n} \intersect X$, $n=1, 2, 3, \ldots$ are porous with porosity constants independent of $n$.  
\end{prop}

\pf 
Axiom [Degree] implies there exists $n_0$ and  $U_{n_0} \in \UUU_{n_0}$ so that the degree 
$\deg(f^{n_0}| U_{n_0})$   is maximal. 
Then all iterated preimages $\wtU_{n_0}$ of $U_{n_0}$ map by degree one onto $U_{n_0}$.  
So $U_{n_0}$ and any iterated preimage $\wtU_{n_0}$ lie in the complement of the post-branch set.  
By the first assertion of Proposition \ref{prop:repellors_are_fractal}, 
for every element $U$ of $\UUU_0$, there is a $k(U) \in \N$ and a preimage $U'$ of $U_{n_0}$ under $f^{-k(U)}$ which is contained in $U$.   Let 
\[ c_0 = \min_{U \in \UUU_0} \frac{\diam U'}{\diam U}.\]

Let $B(x,r)$ be a small ball in $\XX_0$ centered at a point $x \in P_f$.  By Proposition [\ref{prop:BLC},  BLC] there exists some $n$ and $\wtU \in \UUU_n$ such that $B(x, r/L) \subset \wtU \subset B(x,r)$.  Let $U=f^n(\wtU) \in \UUU_0$.  Then by the previous paragraph, $U \supset U'$ where $U' \subset X-P_f$.  If $\wtU'$ is any preimage of $U'$ under $f^n$ which is contained in $\wtU$, then the forward invariance of $P_f$ implies $\wtU' \subset X-P_f$.  By the backward lower relative diameter distortion bounds (\ref{eqn:brdb}), 
\[ \diam \wtU' > \delta^{-1}_+(c_0) \diam \wtU > \delta^{-1}_+(c_0)r/L = c_1 r.\]
Since good open sets are uniformly $K$-almost round (Proposition \ref{prop:uniform_roundness}), 
$\wtU' \supset B(y, c_1r/K)$ for some $y \in X$  and so $P_f$ is $c$-porous  where $c=c_1/K$.  

We merely sketch the second assertion.   Suppose $B(x,r)$ is a small ball centered at a point $x \in B_{f^k}\intersect X$.  Then for some $n$, $r \asymp \diam \wtU$ where $\wtU \in \UUU_{n+k}$.  Let $U=f^k(\wtU)$.  Since $f^k(x) \in P_f$ and $P_f$ is porous, there is some $U' \subset U$ with $\diam U' \asymp \diam U$ and $U' \subset X-P_f$.  If $\wtU'$ is any preimage of $U'$ under $f^k$ which is contained in $\wtU'$, then $\wtU' \subset X-B_{f^k} \intersect X$, and the backwards relative diameter distortion bounds again imply $\diam \wtU' \asymp \diam \wtU \asymp r$.  Since $\wtU'$ is $K$-almost round this implies that $X-B_{f^k}\intersect X$ is uniformly porous as a subset of $X$.  
\qed

The next lemma shows that the roundness distortion control of iterates of $f$, which was assumed to hold only for the sets in $\mathbf U$, in fact extends to any iterate of
$f$ and any ball of small enough radius. 

\begin{prop}[cxc is uniformly weakly quasiregular]
\label{prop:uwqr}
There is a constant $H<\infty$ and a sequence of radii $\{r_n\}_{n=1}^\infty $ decreasing to $0$ such that, for any
iterate $n$, for any $x\in X$, and any $r\in(0,r_n)$, 
$$\roundness(f^n(B(x,r)),f^n(x))\le H.$$
\end{prop}

\pf 
Let $r_n = \frac{c_n}{2L}$, let $n \in \N$ be arbitrary, and fix $r<r_n$ and $\tx \in X$.  By Proposition \ref{prop:BLC}, there exist $\wtU \in \UUU_m$ and $\wtU' \in \UUU_{m+n_0}$ such that 
\[ B(\tx, r/L) \subset \wtU' \subset B(\tx,r) \subset \wtU \subset  B(\tx, Lr).\]
Thus $\diam \wtU \leq 2Lr < c_n$ and so $m>n$ since the sequence $(c_k)$ is decreasing.  Set as usual $U=f^n(\wtU)$, $U'=f^n(\wtU')$, and $x=f^n(\tx)$.  
Now,  
\[ \roundness(\wtU', \tx), \roundness(\wtU, \tx) < L\]
and
\[ \frac{1}{2L^2} < \frac{\diam \wtU'}{\diam \wtU} \leq  1.\]
By the forward roundness (\ref{eqn:forward_roundness_bound}) and relative diameter (\ref{eqn:frdb}) bounds, 
\[ \roundness(U',x), \roundness(U,x) < \rho_+(L)\]
and
\[ \delta_+\left(\frac{1}{2L^2}\right) < \frac{\diam  U'}{\diam U} \leq  1.\]
Moreover, 
\[ U' \subset f^n(B(\tx, r)) \subset U.\]
It follows easily that $\roundness(f^n(B(\tx, r)),x)$ is bounded by a constant independent of $x$, $n$, and $r$.

\qed

\section{Conjugacies between cxc systems}  

Given an increasing homeomorphism $\eta:\R_+\to\R_+$, a homeomorphism $\varphi:X\to Y$ between metric
spaces is said to be {\em $\eta$-quasisymmetric}
\index{Index}{quasisymmetric map}%
 if, for any distinct triples $x_1$, $x_2$ and $x_3$ in $X$,
$$\frac{|\varphi(x_1)-\varphi(x_2)|}{|\varphi(x_1)-\varphi(x_3)|}\le \eta\left(\frac{|x_1-x_2|}{|x_1-x_3|}\right)$$
holds.
A homeomorphism $h$ between metric spaces is {\em weakly quasisymmetric}
\index{Index}{quasisymmetric map, weakly}%
 if it distorts the roundness of balls by a uniform factor, i.e. 
\[ \roundness(h(B(x,r)), h(x)) \leq H\]
for all $x \in X$ and $r \leq \diam X$.  An $\eta$-quasisymmetric map is $\eta(1)$-weakly quasisymmetric.

We start with a result which will enable us to promote weak quasisymmetry to
the usual strong quasisymmetry.

\begin{thm} 
\label{thm:weak_qs}
Let $X,Y$ be two uniformly perfect,  doubling, compact metric 
spaces. Let $h:X\to Y$ be a homeomorphism. If both $h$ and $h^{-1}$
are weakly quasisymmetric, then $h$ and $h^{-1}$ are both 
quantitatively quasisymmetric.\end{thm}

In the above theorem, the term ``quantitative'' means that the function $\eta$ occuring in the definition can be taken to depend only on the constants in the definition of weak quasisymmetry and on the metric regularity constants. 

In the proof, we adapt the argument of Theorem 10.19 of \cite{heinonen:analysis}.

\medskip

\pf The assumptions imply respectively
\bi
\item[(a)] there is a constant $\lambda>1$ such that, for any ball $B$ in $X$ or
$Y$ with non-empty complement, $B\setminus (1/\lambda)B\ne\emptyset$;
\item[(b)] there are constants $C,\beta>0$ such that any set of diameter $d$ in $X$
or $Y$ can be covered by at most $C\epsilon^{-\beta}$ sets of diameter at most
$\epsilon d$;
\item[(c)] there is a constant $H$ such that
$$\left\{\begin{array}{ll} 
\hbox{if } a,b,x\,\in X \hbox{ and } |a-x|\le |b-x| &\hbox{then } |h(a)-h(x)|\le H 
|h(b)-h(x)|\\
\hbox{if } c,d,y\,\in Y \hbox{ and } |c-y|\le |d-y| &\hbox{then } |h^{-1}(c)-h^{-1}(y)|\le H
 |h^{-1}(d)-h^{-1}(y)|.\end{array}\right.
$$
\ib

Choose $\te\in (0,1)$ small enough so that $\te \lambda\le 1/3$.
Let $a,b,x\,\in X$ and set $$t=\frac{|a-x|}{|b-x|} \quad \hbox{and} \quad 
t'=\frac{|h(a)-h(x)|}{|h(b)-h(x)|}\,.$$

Let us assume that $t<\te$.  Property (a) implies there are points $b_0,\ldots,b_s$ 
such that 
$b_j\in B(x,\te^j|b-x|)\setminus B(x,(\te^j/\lambda)|b-x|)$, where $s$ is the 
least integer such that $\te^s<t$.

It follows that if $i<j$ then 
$$\frac{|b_i-b_j|}{|b-x|}\ge\frac{|b_i-x|}{|b-x|}-\frac{|x-b_j|}{|b-x|}$$
so that $$\frac{|b_i-b_j|}{|b-x|}\ge (\te^i/\lambda) -\te^j\ge 
(\te^i/\lambda)(1-\lambda \te)>0$$
and these points are all pairwise disjoint.

Furthermore, from the definition of $s$, we have $$\frac{\log (1/t)}{\log 
(1/\te)}\le s\,.$$

Let $0\le i < j\le s-1$; then $|a-b_j|\le 2|x-b_j|$ and $$|b_i-b_j|\ge 
(\te^{j-1}/\lambda)(1-\lambda\te)|b-x|\ge 2|x-b_j|\,.$$
Hence $|a-b_j|\le |b_i-b_j|$ and it follows from property (c) that $$|h(a)-h(b_j)|\le 
H|h(b_i)-h(b_j)|.$$
Similarly,  $|x-b_j|\le |b_i-b_j|$ implies that $$|h(x)-h(b_j)|\le 
H|h(b_i)-h(b_j)|\,.$$
Therefore $$|h(a)-h(x)|\le 2H|h(b_i)-h(b_j)|.$$

It follows that the balls $B(h(b_j),(1/5H)|h(a)-h(x)|)$ are pairwise disjoint. 
Indeed, if $y\in B(h(b_j),(1/5H)|h(a)-h(x)|)$,
then $$|y-h(b_i)|\ge |h(b_i)-h(b_j)| - |y-h(b_j)|\ge (3/5H)|h(a)-h(x)|$$ so that 
$y\notin B(h(b_i),(1/5H)|h(a)-h(x)|)\,.$
Furthermore they are contained in $B(h(x), 2H|h(x)-h(b)|)$, so the doubling 
property (b) implies
$$s\le C \left(\frac{t'}{5H}\right)^{-\beta}$$ from which we deduce that $t'$ is 
bounded by a function of $t$ which decreases
to $0$ with $t$.

Therefore,  there is a homeomorphism $\eta:[0,1]\to [0,\eta(1)]$ such that 
$\eta(1)\ge 1$ and if $|a-x|\le |b-x|$ then
$$|h(a)-h(x)|\le \eta\left(\frac{|a-x|}{|b-x|}\right) |h(b)-h(x)|\,.$$
Similarly, if $|c-y|\le |d-y|$ then
$$|h^{-1}(c)-h^{-1}(y)|\le \eta\left(\frac{|c-y|}{|d-y|}\right) 
|h^{-1}(d)-h^{-1}(y)|\,.$$

Let us assume now that $t \ge 1/\eta^{-1}(1)$. It follows that $$|h(b)-h(x)|\le 
\eta(1/t)|h(a)-h(x)|\le |h(a)-h(x)|,$$
whence $$|b-x|\le \eta\left(\frac{|h(b)-h(x)|}{|h(a)-h(x)|}\right) |a-x|\,.$$
It follows that $$t'\le 1/\eta^{-1}(1/t).$$
This establishes that $f$ is quasisymmetric, and $f^{-1}$ as well.
\qed

The main result of this section is 

\begin{thm}[Invariance of cxc]
\label{thm:invariance_of_cxc}
Suppose $f: (\XX_1, X) \to (\XX_0, X)$ and $g:(\YY_1, Y) \to (\YY_0, Y)$ are two topological cxc systems which are conjugate via a homeomorphism $h: \XX_0 \to \YY_0$, where $\XX_0$ and $\YY_0$ are metric spaces. 
\be
\item If $f$ is metrically cxc and $h$ is quasisymmetric, then $g$ is metrically cxc, quantitatively.
\item If $f, g$ are both metrically cxc, then $h|_X: X \to Y$ is quasisymmetric, quantitatively. 
\eb
\end{thm}

In the proof below, we use subscripts to indicate the dependence of the metric regularity constants on the system, 
e.g. $\delta_{\pm, f}, \delta_{\pm, g}$, etc.  

\gap

\pf {\bf (1).}  Suppose first that  
$h$ is $\eta$-quasisymmetric.  
Then
\bi
\item {\bf [Roundness quasi-invariant]} 
The map $h$ sends $K$-almost-round sets  
with respect to $x$ to
$\eta(K)$-almost-round sets with respect to $h(x)$; 
\item {\bf [Relative distance distortion]} for all $A,B \subset X$ with $A \subset B$, 
\[ \frac{1}{2\eta\left(\frac{\diam B}{\diam A}\right)} \leq \frac{\diam  
h(A)}{\diam h(B)} \leq
\eta\left(2\frac{\diam A}{\diam B}\right)\]
(see \cite{heinonen:analysis}, Prop. 10.8).
\ib
The topological axioms ([Expansion], [Irreducibility], [Degree]) are invariant under topological conjugacies.  
Axiom [4, Roundness distortion] follows immediately from property 
[Roundness
quasi-invariant] above.  Thus, it suffices to check   
Axiom [5, Diameter distortion].
Let us use small  letters to denote sets and drop   ``$\diam $'' for ease of readability.  
Let $\cl{\eta}(t)= 1/(\eta^{-1})(1/t)$, and notice that $h^{-1}$ 
is $\cl{\eta}$-quasisymmetric.

We have
\[
\renewcommand{\arraystretch}{2}
\begin{array}{cccr}
\frac{\tilde{v}'}{\tilde{v}} & < & \eta(2\frac{\tilde{u}'}{\tilde{u}}) & \mbox{ 
[Rel. dist. distortion]}
\\
\frac{\tilde{u}'}{\tilde{u}} & < & \delta_{-,f}( \frac{u'}{u}) & \mbox{ def.
$\delta_-$} \\
\frac{\tilde{v'}}{\tilde{v}} & < & \eta(2\delta_{-,f}( \frac{u'}{u})) &
\mbox{ $\eta$ increasing} \\
\frac{u'}{u} & < & \cl{\eta}(2\frac{v'}{v}) &\mbox{ [Rel. dist. distortion]}
\end{array}
\]
Thus, $$\frac{\tilde{v}'}{\tilde{v}} 
<\eta\left(2\delta_{-,f}\left(\cl{\eta}\left(2\frac{v'}{v}\right)\right)\right).$$
Now define 
$$\delta_{-,g}(t) =\eta(2\delta_{-,f}(\cl{\eta}(2t))).$$
This is a composition of homeomorphisms, hence a homeomorphism, and so it satisfies the requirements.  
Finding $\delta_{+,g}$ is
accomplished similarly\,:
$$\renewcommand{\arraystretch}{2}\begin{array}{ll}
\frac{v'}{v} & < \eta(2\frac{u'}{u})\\
 & <  \eta(2\delta_{+,f}(\frac{\tilde{u}'}{\tilde{u}}) )\\
& <  \eta(2\delta_{+,f}( \cl{\eta}(2\frac{\tilde{v}'}{\tilde{v}} ))).

\end{array}
$$

\noindent{\bf (2).} Now suppose $g$ is metrically cxc.   By Propositions \ref{prop:cxc_implies_doubling} and \ref{prop:X_is_up}, $X$ and $Y$ are doubling and uniformly perfect.   Therefore, 
it suffices to
show $h$ and $h^{-1}$  are weakly quasisymmetric (cf. Theorem \ref{thm:weak_qs}). Since the setting is
symmetric with respect to $f$ and $g$, it is enough to prove that $h$ is  weakly quasisymmetric.
To show this, it suffices to show (since $h$ and its inverse are uniformly
continuous) that if
$B=B(\tilde{x}, r)$ is a   sufficiently
small ball, then its image $h(B)$  is almost round with respect to  
$\tilde{y}=h(\tilde{x})$, with
roundness constant independent of $B$.
Our proof below follows the usual method (see \cite{sullivan:ihes1982}): 
given a small ball $B$, we  
use the dynamics and the
distortion axioms to blow it up to a ball of definite size and bounded  
roundness.
By compactness, moving over to $Y$ via $h$ distorts roundness by a  
bounded amount.  We then pull
back by the dynamics and apply the distortion axioms again.

Our argument is slightly tricky, since we must trap balls, which are  
possibly disconnected, inside
connected sets in order to apply the pullback step and make sense of  
the ``lift'' of a ball.
We will accomplish this as follows.  Let $\mathbf{U} = \{\UUU_n\}_{n=0}^\infty , \mathbf{V} = \{\VVV_n\}_{n=0}^\infty $ be  
the sequences of good open sets for
$f$ and $g$, respectively.  We are aiming for the following diagram:
{\small 
\begin{equation}
\label{eqn:aimfor}
\begin{array}{ccc}
\widetilde{U}' \subset B \subset \widetilde{U} & \stackrel{h}{\to} &  
\widetilde{V}' \subset
h(\widetilde{U}') \subset h(B) \subset h(\widetilde{U}) \subset  
\widetilde{V} \\
f^n \downarrow & \; & \downarrow g^n \\
U' \subset f^n(B) \subset U & \stackrel{h}{\to} & V' \subset h(U')  
\subset h(f^n(B))=g^n(h(B)) \subset h(U)
\subset V
\end{array}
\end{equation}
}
Below, we indicate by subscripts the dependence on the map of the metric regularity constants $K, C, L, c_n, d_n$, etc.  defined in the previous two sections .  

The diameters of elements of   $\VVV_0$ are bounded from
below.  
Since $\XX_1$ is relatively compact, $h|_{\XX_1}: \XX_1 \to \YY_1$ is uniformly continuous.  Hence there  
exists $\delta_0 > 0$ such
that
\begin{equation}
\label{eqn:delta0}
\diam E < \delta_0 \implies \diam(h(E)) < \epsilon_0 =  \mbox{Lebesgue \# of $\VVV_0$}.
\end{equation}

\noindent{\bf Finding $\wtU, \wtU'$.}  The [Expansion] Axiom implies that there exists $N_0$ such that $d_{N_0, f} <  \delta_0$.   Let $B=B(\tx, r)$ where $r < c_{N_0, f}/(2L_f)$.  By Proposition \ref{prop:BLC}, there exists $n_{0,f}$ and $m \in \N$, $\wtU \in \UUU_m$, and $\wtU' \in \UUU_{m+n_{0,f}}$ such that 
\[ B(\tx, r/L_f) \subset \wtU' \subset B \subset \wtU \subset B(\tx, L_f r).\]
Thus $\diam \wtU \leq  2L_f r \leq  c_{N_0,f}$ and so $m=N_0+n$ where $n \geq 0$.  
\gap

\noindent{\bf Finding $U', U$.}  Let as usual $U=f^n(\wtU)$, $U'=f^n(\wtU')$, $x=f^n(\tx)$.  Then $U \in \UUU_{N_0}$ and $U' \in \UUU_{N_0+n_{0,f}}$.  
\gap

\noindent{\bf Finding $V$.}  Let $y=h(x)$.  Since $U \in \UUU_{N_0}$ and $d_{N_0,f} < \delta_0$, the bound (\ref{eqn:delta0}) implies $\diam(h(U)) < \epsilon_0$ and so there exists $V \in \VVV_0$ with $h(U) \subset V$.
\gap

\noindent{\bf Finding $V'$.}  The forward roundness bound 
(\ref{eqn:forward_roundness_bound}) implies that 
\[ \roundness(U',x) < \rho_{+,f}(L_f).\]
Hence 
\[ U' \supset B(x,s'), \;\;\; \mbox{ where } s' = \frac{c_{N_0+n_{0,f}}}{2\rho_{+,f}(L_f)}.\]
Since $X$ is compact, $h(B(x,s')) \supset B(y,t')$ where 
\[ t' = \inf\{|h(x)-h(a)| : x \in X, |a-x| = s'\}.\]
The [Expansion] Axiom implies that there exists $k_0$ such that $d_{k_0,g} < t'/2$.  Proposition \ref{prop:uniform_roundness} implies that there exists $V' \in \VVV_{k_0}$ such that $\roundness(V', y) < K_g$.   Then 
\[  V' \subset h(U') \subset h(f^n(B)) \subset h(U) \subset V\]
where 
\[ \roundness(V', y), \roundness(V,y) \leq \min\left\{\frac{d_{0,g}}{t'}, K_g\right\} =: R\]
and
\[ \frac{\diam V'}{\diam V} < \frac{c_{k_0,g}}{d_{0,g}} =: D.\]
\gap

\noindent{\bf Finding $\wtV, \wtV'$.}  Let $\wtV, \wtV'$ denote the preimages of $V$ and $V'$, respectively, containing $\ty=h(\tx)$.  We have now achieved the situation summarized in (\ref{eqn:aimfor}).  
\gap

\noindent{\bf Conclusion.}  The backwards roundness bound 
(\ref{eqn:backward_roundness_bound}) and backwards relative diameter distortion bounds (\ref{eqn:brdb}) imply
\[ \roundness(\wtV, \ty), \roundness(\wtV',\ty) < \widetilde{R} = \rho_{-,g}(R)\]
and 
\[ \frac{\diam\wtV'}{\diam\wtV} > \widetilde{D} = \delta_{+,g}(D).\]
Hence $\roundness(h(B), h(\tx)) < 2\widetilde{R}^2/\widetilde{D}$ and the proof is complete.
\qed

\chapter{Geometrization}

In this chapter, we assume  we are given the data of an fbc $f:\XX_1 \to \XX_0$ with repellor $X$ as in the beginning of \S\,\ref{secn:def_top_cxc}   and a finite cover $\UUU$ which together satisfy Axiom [Expansion].  Given a suitably small parameter $\varepsilon>0$, we will associate to this data a metric $d_\varepsilon$ on the repellor $X$ such that in this metric, $f$ acts very much like a piecewise linear map of an interval with constant absolute value of slope:  it sends balls of radius $r$ onto balls of radius $e^\varepsilon r$.  In so doing, we promote our topological dynamical system to a non-classically conformal one.  We will also show that the quasisymmetry class of this metric is natural, in that it does not depend on the choice of the open cover $\UUU$, so long as Axiom [Expansion] is satisfied.   This means our topological dynamical system has a canonically associated conformal gauge.  

The metric $d_\varepsilon$ arises as a visual metric on the boundary at infinity $\bdry \Gamma$ of a certain Gromov hyperbolic space $\Gamma$ associated to the data.   The map $f$ induces a $1$-Lipschitz map $F: \Gamma \to \Gamma$ and so extends to a map $F:  \cl{\Gamma} \to \cl{\Gamma}$ of the compactification.   We show that on the boundary $\bdry \Gamma$, the map $F$ is conjugate via a natural homeomorphism $\phi_f$  to our original dynamical system $f: X \to X$.  
This approach follows not only Thurston's philosophy that Topology implies a natural Geometry, but also Gromov's point of view that coarse notions capture enough information to determine Geometry.

We then exploit the existence of the extension $F: \cl{\Gamma} \to \cl{\Gamma}$.  
Assuming in addition Axiom [Irreducibility], we construct, using the Patterson-Sullivan method of Poincar\'e series, a natural measure $\mu_f$ which is invariant, quasiconformal, ergodic, and mixing, and which governs the distribution of preimages of points and of periodic points.  

When the original map $f: \XX_1 \to \XX_0$ is metrically cxc with respect to a metric $d$, we show that the conjugacy $\phi_f$ constructed above is quasisymmetric, and we prove that the measure $\mu_f$ is the unique measure of maximal entropy $\log \deg(f)$, and
 is Ahlfors regular of exponent 
$\frac{1}{\varepsilon} \log \deg(f)$.
\gap

This chapter is organized as follows. In the first section, we review the basic geometric theory of 
unbounded metric spaces, emphasizing hyperbolicity and compactifications. Section 2 is devoted to the construction of the 
space $\Gamma$, and we establish its first properties. In section 3, the hyperbolicity of $\Gamma$ is proved, and its
naturality is established.   In Section 4 we study measure-theoretic properties 
and Hausdorff dimension.   
In Section 5 we assume that the original system is topologically or metrically cxc, and we refine the results of the preceding sections. 

\section{Compactifications of quasi-starlike spaces}
\label{scn:compactification}

A  metric space  $(X,d)$ is said to be {\em proper}
\index{Index}{proper metric space}%
 if, for all $x\in X$, the function $y\mapsto d(x,y)$ is proper, meaning
that closed balls of finite radius are compact. A {\em geodesic curve}
\index{Index}{geodesic curve}%
 is a continuous function $\gamma:I\to X$ such that
$d(\gamma(t),\gamma(t'))=|t-t'|$ for all $t,t'\in I$ and  where $I$ is an interval. We will often not distinguish between the
function $\gamma$ and its image in $X$. The space $X$ is said to be {\em geodesic}
\index{Index}{geodesic metric space}%
 if any pair of points can be joined by a geodesic.

\gap

\noindent{\bf Rectifiable curves and integration.} 
We refer to \cite[Chapter 1]{vaisala:lectures_qc} and \cite[Chapter 7]{heinonen:analysis} for what follows.
Let $I=[a,b]\subset \R$ be
an interval, $a \le b$. A {\em rectifiable curve}
\index{Index}{rectifiable curve}%
$\gamma: I\to X$ is a continuous map of bounded variation i.e.,
\begin{equation}\label{eqn:rect}
\sup_{a = s_0 < s_1 \ldots <  s_n=b} \sum_{0\le j<n} d(\gamma(s_{j+1}),\gamma(s_j))<\infty\,,\end{equation} 
where the supremum is taken over all subdivisions of $[a,b]$ with $s_0=a$ and $s_n=b$. The supremum in (\ref{eqn:rect}) is the
{\em length}\index{Index}{length of a curve} $\ell(\gamma)$ of $\gamma$. When $\gamma$ is rectifiable, it can be parametrized by
arclength\index{Index}{arclength parametrization}, 
i.e. there is an increasing and continuous function $s:[a,b]\to [0,\ell(\gamma)] $  and a curve
$\gamma_s:[0,\ell(\gamma)]\to X$ such that $\gamma=\gamma_s\circ s$ and, for all
$0\le c\le d\le \ell(\gamma)$, $\ell(\gamma_s|_{[c,d]})= |d-c|$.
We then say that $\gamma_s$ is the parametrization of $\gamma$
by {\em arclength}. Note that geodesic curves are already
parametrized by arclength by definition.

Let $\gamma$ be a rectifiable curve, and $\rho$ a non-negative Borel
function defined on $\gamma$. We set
$$\int_{\gamma} \rho=\int_{\gamma} \rho(x)ds(x)=\int_0^{\ell(\gamma)} \rho(\gamma_s(t))dt$$ where
$\gamma_s$ is the parametrization by arclength.

\gap
\noindent{\bf Quasi-starlike spaces.}
\index{Index}{quasi-starlike}%
Fix a base point $o\in X$. A {\em ray
\index{Index}{ray}%
 based at $o$} is a geodesic curve $\gamma:\R_+\to X$ such that $\gamma(0)=o$.
 Let $\RRR$\index{Symbols}{$\RRR$} be the set of geodesic curves starting at $o$, and 
let $\RRR_\infty$\index{Symbols}{$\RRR_\infty$} be the set of rays based at $o$.
The space $(X,o)$ is  {\em $K$-quasi-starlike (about $o$)} if, for any $x\in X$, there is a ray $\gamma\in\RRR_\infty$ such that
$d(x,\gamma)\le K$.

\gap

In this section, we assume that $(X,d)$ is a geodesic, proper, $K$-quasi-starlike space about a point $o$. 
For convenience, we write  $d(x,y)=|x-y|$ and $|x|=|x-o|$.

\gap

{\noindent\bf Hyperbolic spaces.} The {\em Gromov product}\index{Index}{Gromov product} of two points $x,y\in X$ is defined by\index{Symbols}{$(x\mid y)$} $(x|y)=(1/2)(|x|+|y|-|x-y|)$.
The metric space $X$ is {\em Gromov hyperbolic}
\index{Index}{hyperbolic space}%
\index{Index}{Gromov hyperbolic space}%
 if there is some constant $\delta\ge 0$ such that $$(x|z)\ge \min\{(x|y),(y|z)\}-\delta$$
for any points $x,y,z\in X$.  (By Proposition 1.2 of \cite{coornaert:delzant:papadopoulos}, this definition agrees with the more common one 
in which the above inequality is required to hold for all $x, y, z$ and $o$ instead of just at a single basepoint $o$.)  
Let us note that in such a space, $(x|y)$ equals  $d(o,[x,y])$ up to a universal additive constant, where $[x,y]$ is any geodesic segment
joining $x$ to $y$. We refer to \cite{coornaert:delzant:papadopoulos} and to \cite{ghys:delaharpe:groupes} for more information on Gromov hyperbolic spaces.
\gap

{\noindent\bf Compactification.} 
\index{Index}{compactification}%
Here, we do not assume $X$ to be hyperbolic.
We propose to compactify $X$ using the method of W.\,Floyd \cite{floyd:completion}. Let $\varepsilon >0$, and, for $x\in X$, 
define\index{Symbols}{$\rho_{\varepsilon}$}
$\rho_\varepsilon(x)=\exp (-\varepsilon|x|)$.

\gap

If $\gamma$ is rectifiable, we set\index{Symbols}{$\ell_{\varepsilon}(\gamma)$} $$\ell_{\varepsilon}(\gamma)=\int_{\gamma}\rho_{\epsilon}\,.$$

For $x,y\in X$, define\index{Symbols}{$d_{\varepsilon}$, $\mid\cdot\mid_{\varepsilon}$ } 
$$d_\varepsilon(x,y)=|x-y|_\varepsilon=\inf_\gamma \int_\gamma \rho_\varepsilon=\inf_\gamma \ell_{\varepsilon}(\gamma)$$ where the infimum is taken over all rectifiable
curves which join $x$ to $y$. Thus, $|x-y|_\varepsilon\le |x-y|$.

\gap

The space $(X,|\cdot|_\varepsilon)$ is not complete since if $\gamma\in\RRR_\infty$  and if $t'>t$ then
$$|\gamma(t)-\gamma(t')|_\varepsilon\le \int_{t}^{t'}e^{-\varepsilon s}ds\le  e^{-\varepsilon t}/\varepsilon\,.$$
Therefore $\{\gamma(n)\}_n$ is a non-convergent $d_\varepsilon$-Cauchy sequence.

\gap

{\noindent\bf Definition.} Let $\overline{X_\varepsilon}$ be the completion of $(X,d_\varepsilon)$, 
and set  $\partial X_\varepsilon= \partial_\varepsilon X= \overline{X_\varepsilon}\setminus X$.\index{Symbols}{$\overline{X_\varepsilon}$, $\partial X_\varepsilon$}  
Thus, $\cl{X}_\varepsilon$ is also a length space.

\gap

\noindent{\bf Visual metrics.} If $X$ is Gromov hyperbolic, then, for $\varepsilon$ small enough, $\partial X_\varepsilon$ coincides with the Gromov boundary of 
$X$ and $d_\varepsilon$ is a so-called {\em visual distance}
\index{Index}{visual metric}%
 (cf. \cite[Chap.\,11]{coornaert:delzant:papadopoulos} and \cite[Chap.\,4]{bonk:heinonen:koskela}). That is, we may extend the definition of the Gromov product to the boundary
$\partial X_\varepsilon$ and then, $|x-y|_\varepsilon\asymp e^{-\varepsilon (x|y)}$ holds if 
$0<\varepsilon\le \varepsilon_0(\delta)$ for some constant $\varepsilon_0(\delta)>0$ which depends only on $\delta$. 
When dealing with a hyperbolic space $X$, we will write indifferently $$\partial_\varepsilon X,\
\partial_\infty X,\ \partial X$$ to denote its boundary. In any case and unless specified, metrics
on the boundary will always be visual metrics $d_{\varepsilon}$ as above for some fixed
parameter $\varepsilon>0$ small enough.

\gap

{\noindent\bf Topology on $\RRR$.} If $\gamma\in\RRR$ then $\gamma$ is geodesic for $d_\varepsilon$. 
Indeed, let $\gamma$ be a curve starting from $o$ that is parametrized by arclength. 
It follows that $|\gamma(t)|\le t$ for all $t\in[0,\ell(\gamma)]$ with equality for all $t$ if $\gamma \in\RRR$.
Therefore, 
\begin{equation}\label{eqn:lowebound_rays}
\ell_\varepsilon(\gamma)=\int_0^{\ell(\gamma)}e^{-\varepsilon |\gamma(t)|}dt\ge \int_0^{\ell(\gamma)}e^{-\varepsilon t}dt\ge  \frac{1}{\varepsilon}\left(1-e^{-\varepsilon\ell(\gamma)}\right)\,. \end{equation}
We have equality when $\gamma$ is geodesic for $d_0$.

\gap

For $\gamma\in\RRR$, the limit in $\cl{X}_\varepsilon$ of $\gamma(t)$ at $\ell(\gamma)$ exists since any sequence $(\gamma(t_n))$ with $t_n\to\ell(\gamma)$ is a $d_\varepsilon$-Cauchy 
sequence.
Let us define $$\pi(\gamma)=\lim_{t\to\ell(\gamma)}\gamma(t).$$

\gap
The {\em Hausdorff distance} between two closed nonempty subsets $A,B$ of a compact metric space $(Z,d)$ is given by 
\[ d_H(A,B)=\max\{ \sup_{a \in A}\ \inf_{b \in B}d(a,b), \sup_{b \in B}\ \inf_{a \in A}d(a,b)\}.\]  
\index{Index}{Hausdorff distance}%
The Hausdorff distance turns the set of nonempty closed subsets of a compact metric space into a compact metric space.  

The closure of each element $\gamma\in\RRR$ is  compact in $\overline{X_\varepsilon}$, so 
the Hausdorff distance with respect to the metric $d_\varepsilon$ between the closures of rays defines a metric, and hence a topology, on the set of rays $\RRR$. 

\begin{lemma}\label{toprayon} The set $\RRR$ is compact and the map $\pi:\RRR\to\overline{X_\varepsilon}$ is continuous and surjective.
Furthermore, $\RRR_\infty$ is closed in $\RRR$ and $\pi|_{\RRR_\infty} : \RRR_\infty \to \bdry X_\varepsilon$ is also surjective.\end{lemma}

\pf Let $(\gamma_n)$ be a sequence in $\RRR$.  
Suppose first that $\lim\inf\ell(\gamma_n)=L<\infty$, where $\ell$ denotes the $d_0$-length of a geodesic in $X$.  
Regard each $\gamma_n$ as a function from an interval into $X$.  Since $X$ is proper and $\gamma_n(0)=o$, the basepoint,  for all $n$, the Arzela-Ascoli theorem implies that after passing to a subsequence we may assume $(\gamma_n)_n$ converges uniformly on $[0,L]$ 
to a continuous map $\gamma: [0,L] \to X$.   
Since $|a-b|=d_0(\gamma_n(a),\gamma_n(b))\to d_0(\gamma(a),\gamma(b))$ for all $a, b \in [0,L]$, the curve $\gamma$ is a $d_0$-geodesic 
and so $\gamma \in \RRR$.  It follows easily that $\gamma_n \to \gamma$ in the Hausdorff topology on $\RRR$ 
with respect to the metric $d_\varepsilon$.  

Suppose now that $\lim\inf\ell(\gamma_n)=\infty$.  
Again, the Arzela-Ascoli theorem and a diagonalization argument shows that we may assume after passing to a subsequence that 
$\gamma_n \to \gamma$ uniformly on compact subsets of $[0,\infty)$, where $\gamma: [0,\infty) \to X$ is a $d_0$-geodesic, i.e. 
an element of $\RRR_{\infty}$ (by this, we mean that on any compact subset of $[0,\infty)$, $\gamma_n$ 
is defined on this subset for all $n$ sufficiently large, and for such $n$ the convergence is uniform on this subset).  
We now prove $\gamma_n \to \gamma$ in the Hausdorff topology with respect to the metric $d_\varepsilon$.  
It suffices to show that $\gamma_n \to \gamma$ uniformly when regarded as maps from $[0,\infty)$ to $(\cl{X}_\varepsilon, d_\varepsilon)$.  Observe that for any $L \in [0,\infty]$, any geodesic ray $\alpha: [0, L] \to X$, and any $t_0 \le L$, the length of the tail  satisfies
\begin{equation}
\label{eqn:tail}
 \ell_\varepsilon(\alpha|_{[t_0,L]}) = \int_{t_0}^L e^{-s\varepsilon}ds \leq \frac{1}{\varepsilon}e^{-\varepsilon t_0}.
\end{equation}
Fix now $\eta>0$ and choose $t_0$ so large that $\frac{1}{\varepsilon}e^{-\varepsilon t_0}<\eta/3$.  
Next, using uniform convergence, choose $n_0$ so large that $\sup_{0 \leq t \leq t_0}|\gamma_n(t)-\gamma(t)|_{\varepsilon} < \eta/3$ 
for all $n \geq n_0$.  
Then for all $t>t_0$ and all $n>n_0$ large enough so that $\gamma_n(t)$ is defined, 
by the triangle inequality and (\ref{eqn:tail}) we have 
\[ |\gamma_n(t)-\gamma(t)|_\varepsilon \leq |\gamma_n(t)-\gamma_n(t_0)|_{\varepsilon} + |\gamma_n(t_0)-\gamma(t_0)|_{\varepsilon} + |\gamma(t_0)-\gamma(t)|_{\varepsilon} < \eta.\]
We have thus established that $\RRR$ is compact and that $\RRR_\infty$ is closed in the Hausdorff topology with respect to $d_\varepsilon$.  

The continuity of the map $\pi$ follows from entirely analogous arguments.  It remains only to show that $\pi: \RRR_\infty \to \bdry X_\varepsilon$ is surjective.  If $x\in\partial X_\varepsilon$, then
there is a sequence $(x_n)$ in $X_\varepsilon$ which converges to $x$. Let $\gamma_n$ be geodesic segments joining $o$ to $x_n$, so that $\pi(\gamma_n)=x_n$.  
By compactness of $\RRR$, there exists $\gamma \in \RRR$ such that after passing to a subsequence, $\gamma_n \to \gamma$.  Clearly $\gamma \in \RRR_\infty$.  Since $\pi$ is continuous, $\pi(\gamma)=x$.  
\qed

\begin{lemma}\label{identif} The following hold $$X_\varepsilon=B_\varepsilon(o,1/\varepsilon), \partial X_\varepsilon=\{x,\ |x|_\varepsilon=(1/\varepsilon)\} \quad   \hbox{ and } \quad 
B(o,R)= B_\varepsilon(o,(1/\varepsilon)(1-e^{-\varepsilon R}))\,.$$\end{lemma}

\pf
Let $x\in X$ and $\gamma\in\RRR$ be a geodesic curve joining $o$ to $x$.
According to (\ref{eqn:lowebound_rays}), $\gamma$ is also
geodesic for $d_{\varepsilon}$.
Therefore $$|x|_\varepsilon=\frac{1}{\varepsilon}(1-e^{-\varepsilon |x|})<\frac{1}{\varepsilon}\,.$$
This implies that $B(o,R)= B_\varepsilon(o,(1/\varepsilon)(1-e^{-\varepsilon R}))$ and $X_\varepsilon\subset B(o,1/\varepsilon)$.


\gap

Let $x\in \partial X_\varepsilon$. There is a sequence $(x_n)$ of $X$ such that $x_n$ converges to $x$. Since $X$ is a proper
space, it follows that $|x_n|\to\infty$. Furthermore,
$$|x_n|_\varepsilon\ge (1/\varepsilon)(1-e^{-\varepsilon|x_n|})$$ so that $|x|_\varepsilon=1/\varepsilon$. This establishes the lemma.\qed

\bigskip

\begin{figure}
\label{fig:Shadow}
\psfragscanon
\psfrag{1}{$o$}
\psfrag{2}{$x$}
\psfrag{3}{$R$}
\psfrag{4}{$\mho(x,R)$}
\psfrag{5}{$\mho_\infty(x,R)$}
\psfrag{6}{$\bdry_\varepsilon X$}
\begin{center}
\includegraphics[width=3in]{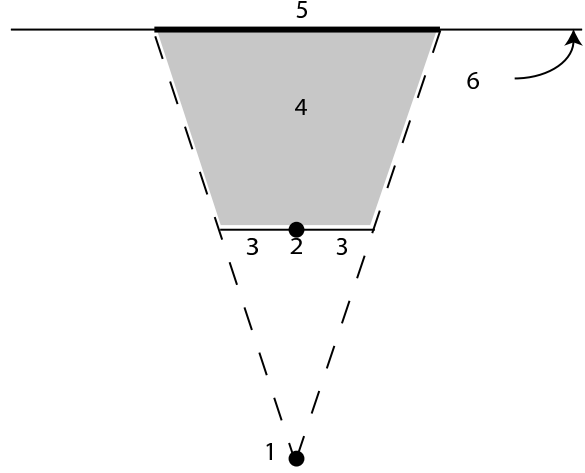}
\end{center}
\caption{{\sf The shadow $\mho(x,R)$ cast by the ball of radius $R$ about $x$.}}
\end{figure}

{\noindent\bf Shadows.} Let $x\in X$, $R>0$. The {\em shadow}
\index{Index}{shadow}%
 $\mho(x,R)$\index{Symbols}{$\mho(x,R)$, $\mho(x)$, $\mho_{\infty}(x,R)$} of $B(x,R)$ is the set of points $y$ in $\overline{X}_\varepsilon$
for which there is a $d_0$-geodesic curve joining $o$ to $y$ which intersects $\cl{B(x,R)}$. 
See Figure \ref{fig:Shadow}.
Let 
$\mho_\infty(x,R)=\mho(x,R)\cap\partial X_\varepsilon$ be its trace on $\partial X_\varepsilon$.
When $R=1$ we employ the notation $\mho(x)$ for $\mho(x,1)$.

\begin{lemma}\label{omdiam} For any $x, R$, there is a constant $C_R>0$ such that $$\diam_\varepsilon \mho(x,R)\le C_R e^{-\varepsilon |x|}\,.$$\end{lemma}

\pf Let $y\in\mho(x,R)$. There is a geodesic segment $[o,y]$ and a point $p\in B(x,R)\cap [o,y]$.
Therefore, $$|y-x|_\varepsilon\le |y-p|_\varepsilon +|x-p|_\varepsilon.$$
Since $|x|-R\le |p|\le |x|+R$, it follows that $|x-p|_\varepsilon\le Re^{\varepsilon R} e^{-\varepsilon|x|}$.
Let $[p,y]$ denote the subsegment of the segment $[o,y]$ joining $p$ to $y$.  We have 
$$|p-y|_\varepsilon \leq \ell_\varepsilon([p,y]) \leq \int_{|p|}^{|y|}e^{-\varepsilon t} dt \leq \frac{1}{\varepsilon}e^{-\varepsilon |p|} \leq \frac{1}{\varepsilon}e^{-\varepsilon(|x|-R)} = \frac{e^{R\varepsilon}}{\varepsilon}e^{-\varepsilon |x|}.$$
This establishes the estimate.\qed

{\noindent\bf Remark.} Shadows are almost round subsets of the boundary.  More precisely:  if $X$ is Gromov hyperbolic and $K$-quasi-starlike, then, for a fixed $R$  
which is chosen large enough, there is 
a constant $C=C(\varepsilon,R,K)$ such that, for any $x\in X$, there is a boundary point $a\in\partial X$ such that 
$$B(a,(1/C)e^{-\varepsilon|x|})\subset \mho_\infty(x,R)\subset B(a,Ce^{-\varepsilon|x|})\,.$$
A proof of this fact can be found in \cite[Prop. 7.4]{coornaert:patterson-sullivan}. Furthermore, the family $\{\interior(\mho(x,R)),\,x\in X,\,R>0\}$ defines a basis of neighborhoods
in $\overline{X_\varepsilon}$ for points at infinity.

\gap

{\noindent\bf Distance to the boundary.} If $x\in X_\varepsilon$, we let $\delta_\varepsilon(x)=\mbox{\rm dist}_\varepsilon(x,\partial X_\varepsilon)$.

\begin{lemma}\label{distbo} If $X$ is  $K$-quasi-starlike, then for all $x \in X$, 
$$ \frac{e^{-\varepsilon|x|}}{\varepsilon}\le \delta_\varepsilon(x)\le C_{K,\varepsilon}\frac{e^{-\varepsilon|x|}}{\varepsilon}\,.$$\end{lemma}

\pf Let $x\in X$. We start with a first coarse estimate\,:
$$\delta_\varepsilon(x)\ge \int_{|x|}^\infty e^{-\varepsilon t}dt =\frac{1}{\varepsilon} e^{-\varepsilon |x|}\,.$$
If there is a ray $\gamma\in\RRR_\infty$ such that $x\in\gamma$, then 
$$\delta_\varepsilon(x)=\int_{|x|}^{\infty} e^{-\varepsilon t}dt= \frac{e^{-\varepsilon|x|}}{\varepsilon}=\frac{\rho_\varepsilon(x)}{\varepsilon}\,.$$
In general, since $X$ is $K$-quasi-starlike, there is a ray $\gamma\in\RRR_\infty$ and a point $p\in B(x,K)\cap \gamma$.
Therefore, $|x|-K\le |p|\le |x|+K$ and $|x-p|_\varepsilon\le C_K e^{-\varepsilon |x|}|x-p|$. Then
$$\delta_\varepsilon(x)\le |x-p|_\varepsilon+ \delta_\varepsilon(p)\le  C_K e^{-\varepsilon |x|}+\frac{e^{-\varepsilon(|x|-K)}}{\varepsilon}\le C_{K,\varepsilon}e^{-\varepsilon |x|}\,.$$\qed

{\noindent\bf Quasi-isometries versus quasisymmetries.} A {\em quasi-isometry}
\index{Index}{quasi-isometry}%
 $f:X\to Y$ between two metric spaces is
a map for which there are constants $\lambda\ge 1$ and $c>0$ such that 
\be
\item {[bi-Lipschitz in the large]} for any $x,x'\in X$, $$\frac{1}{\lambda}|x-x'|-c\le |f(x)-f(x')|\le \lambda |x-x'|+c,$$
\item {[nearly surjective]} for any $y\in Y$, there is some $x\in X$ such that $|f(x)-y|\le c$.
\eb
We note that if $f:X\to Y$ is a quasi-isometry, then there exists a quasi-isometry $g:Y\to X$ such that
$|g\circ f(x)-x|\le C$ for some constant $C<\infty$.

It is well-known that if $\Phi:X\to Y$ is a quasi-isometry between two hyperbolic spaces, then it extends as
a quasisymmetric homeomorphism $\varphi:\partial X\to \partial Y$, if we endow the boundaries with visual metrics; see \cite[Prop.\,4.5]{paulin:determined}
in the general setting of hyperbolic metric spaces.
For the converse, we have

\begin{thm}[M.\,Bonk \& O.\,Schramm]\label{thm:bsch_funct} Let $X$, $Y$ be two quasi-starlike hyperbolic spaces. For any quasisymmetric homeomorphism
 $\varphi:\partial X\to\partial Y$, there is a quasi-isometric map $\Phi:X\to Y$ which extends $\varphi$.\end{thm}

For a proof, see \cite{paulin:determined} or Theorem 7.4 and Theorem 8.2 in \cite{bonk:schramm:embeddings}.

\section{Spaces associated to finite branched coverings}
\label{secn:spaces_associated}

Suppose $f: \XX_1 \to \XX_0$ is a finite branched covering with repellor $X$, and all the conditions on $\XX_0, \XX_1, f$, and $X$  
stated at the beginning of \S\,\ref{secn:def_top_cxc} are satisfied.  
We assume furthermore that we are given a finite open covering $\UUU=\UUU_0$ of $X$ by connected subsets of $\XX_0$.

\gap

Under these assumptions, we prove 

\begin{thm}\label{thm:construction} The pair $(f,\UUU)$ defines a proper, geodesic, unbounded,  quasi-starlike, metric space $\Gamma$
together with a continuous map $F:\Gamma\to\Gamma$ with the following property.  
For any $\varepsilon >0$, let $(\cl{\Gamma}_\varepsilon, d_\varepsilon)$ denote the metric space giving the 
compactification of $\Gamma$ as constructed in \S\,\ref{scn:compactification}.  Then  
\be
\item There exists a continuous map 
$$\phi_f:X\to\partial_{\varepsilon}\Gamma$$ such that $\phi_f\circ f=F\circ\phi_f$.
\item The map $F$ extends as a Lipschitz map $F:\cl{\Gamma}_\varepsilon \to \cl{\Gamma}_\varepsilon$ sending the boundary to the boundary.
\item Balls are sent to balls:  $F(B(\xi,re^{-\varepsilon}))=B(F(\xi),r)$ holds  for any $\xi\in\overline{\Gamma}\setminus\{o\}$ 
and any $r\in (0,|F(\xi))|_\varepsilon)$.
\item If $(f,\UUU)$ satisfies the [Expansion] axiom, then there is some $\varepsilon_0>0$ such that, 
for any $\varepsilon\in (0,\varepsilon_0)$, the map
$\phi_f:X\to\partial_{\varepsilon}\Gamma$ is a homeomorphism.
\eb
\end{thm}

In the above theorem, $\Gamma_\varepsilon$ and $\bdry_\varepsilon \Gamma$ are defined as in the previous section, 
and $B(\zeta, r)$ denotes the ball of radius $r$ about $\zeta$ in $\cl{\Gamma_\varepsilon}$.    

Note the similarity of this statement with the case of hyperbolic groups; cf. Appendix B.

In the next section, we investigate more closely the geometry of $\cl{\Gamma}_\varepsilon$.  

\gap

{\noindent\bf Definition of $\Gamma$.} 
\index{Symbols}{$\Gamma(f,\UUU)$}%
\index{Symbols}{$\Gamma$}%
From the data consisting of the map $f:\XX_1\to\XX_0$ and the cover $\UUU$, we will define an infinite graph $\Gamma$  equipped with a distinguished basepoint.   
Our construction is quite similar to that employed by Elek \cite[\S\,3]{elek:lp_confdim} and by Bourdon and Pajot \cite[\S\,2.1]{bourdon:pajot:besov}.  
However, in our setting, $\Gamma$ is defined using topological instead of metric data, and it will be used later on to construct metrics associated with topological dynamical systems.

The set $V(\Gamma)$ of vertices is the union of the elements of $\mathbf{U}=\union_{n \geq 0} \UUU_n$, together with a base vertex $o=X$.  
It will be convenient to reindex the levels as follows.  For $n \in \N$ set \index{Symbols}{$S(n)$}
\[ S(n)=\left\{
\begin{array}{lr}
\UUU_{n-1} & \mbox{ if }  n \geq 1,\\
\{o\} &  \mbox{ if } n=0.
\end{array}
\right. 
\]
For $n \in \N$ and a vertex $W \in S(n)$, we set $|W|=n$.\index{Symbols}{\ensuremath{\mid W\mid}}  Thus, $V(\Gamma)=\union_{n\geq 0}S(n)$.  

\begin{figure}
\label{fig:PicOfSigma01}
\begin{center}
\includegraphics[width=3in, angle=90]{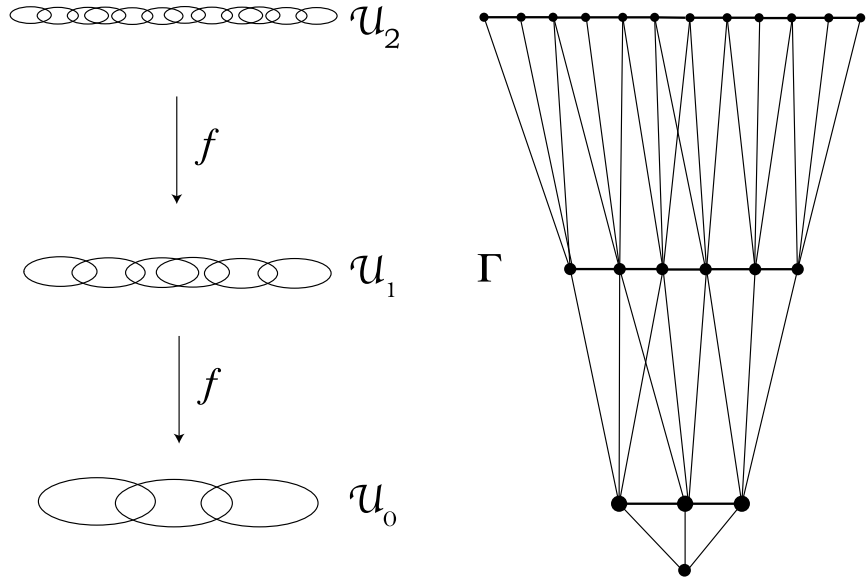}
\end{center}
\caption{{\sf Definition of the graph $\Gamma$.}}
\end{figure}

Two vertices $W_1, W_2$ are joined by an edge if 
\[ \left|   |W_1| - |W_2| \right| \leq 1 \;\;\mbox{ and }\;\; W_1\intersect W_2 \intersect X \neq \emptyset.\]
See Figures \ref{fig:PicOfSigma01} and 
3.3.

\begin{figure}
\label{fig:PathInGamma}
\psfragscanon
\psfrag{1}{in $\XX_1$}
\psfrag{2}{in $\Gamma$}
\begin{center}
\includegraphics[width=3in]{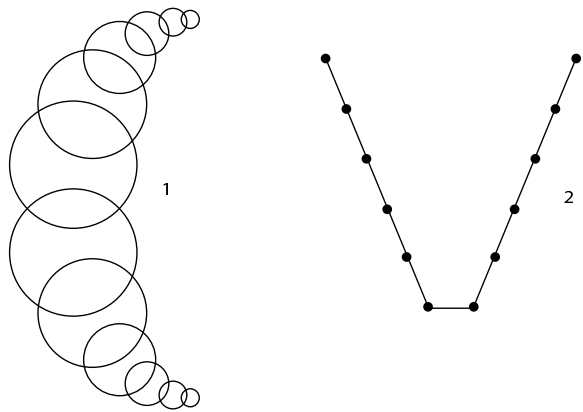}
\end{center}
\caption{{\sf The vertices of a parameterized edge-path in $\Gamma$, at right, correspond to a sequence of elements of $\mathbf{U}$ in which consecutive elements intersects in points of the repellor, $X$.}}
\end{figure}
This definition forbids loops from a vertex to itself and multiple edges between vertices, so $\Gamma$ is indeed a graph as claimed.   The graph $\Gamma$ is turned into a geodesic metric space in the usual way by decreeing that each closed edge is isometric to the Euclidean unit interval $[0,1]$.  Since each $S(n)$ is finite, the valence at each vertex is bounded (though not necessarily uniformly so) and so $\Gamma$ is proper.    Since as subsets of $\XX_0$, any vertex $W \in S(n)$ intersects a set $W' \in S(n-1)$, any vertex $W$ can be joined to the basepoint $o$ by a geodesic ray in $\Gamma$.  Hence $\Gamma$ is connected.  It is also $1/2$-quasi-starlike, since any point of an edge joining vertices at the same level lies within distance at most $1/2$ from a geodesic ray emanating from $o$.  By construction, $S(n)$ is the 
sphere of radius $n$ about the origin $o$.  
\gap

\gap

{\noindent\bf Action of the finite branched covering.} 
If $n \geq 2$ and $W \in S(n)=\UUU_{n-1}$, then as subsets of $\XX_0$,
the set $f(W) \in S(n-1)=\UUU_{n-2}$, so $f$ induces a map 
\index{Symbols}{$F$}$F: \union_{n\geq 2} S(n) \to \union_{n\geq 2} S(n-1)$.  Define $F(W)=o$ for all $W \in S(1) \union S(0)$; 
thus $F$ is defined on the vertex set $V(\Gamma)$.  To extend $F$ over edges, observe that if $n \geq 1$ and if (as subsets of $\XX_0$) the sets
$\widetilde{W}, \widetilde{W}'$ are distinct inverse images of $W$, then $\widetilde{W}, \widetilde{W}'$, being distinct components of the inverse image, cannot intersect.
Thus, if $W_1, W_2$ are joined by an edge, the definition of edges given above then implies that either 
\be
\item $F(W_1) \neq F(W_2)$ and $F(W_1), F(W_2)$ are joined by an edge, or
\item $|W_1|, |W_2| \leq 1$ and $F(W_1)=F(W_2)$.
\eb 
Letting $E$ be the union of edges joining pairs of elements of $S(1)$, properties (1) and (2) above imply 
that $F$ extends naturally to a continuous map $F: \Gamma \to \Gamma$
which collapses $\cl{B(o, 1)}  \union S(0) \union E \to \{o\}$, and which otherwise sends all edges homeomorphically onto their images.

\gap

\noindent{\bf Properties of $F$.}
\begin{itemize}
\item[--] $F$ is $1$-Lipschitz.  
\item[--] $F$ decreases levels by one: $|F(\xi)|=|\xi|-1$ for all $|\xi|\ge 1$.
\item[--] $F$ sends rays to rays:  $F:(\RRR,\RRR_\infty)\to(\RRR,\RRR_\infty)$ 
\item[--] $F$ has the path lifting property for paths which avoid the base point $o$:   
any path $\gamma$ in $\Gamma\setminus\{o\}$ can be lifted by $F^{-1}$.   

Once a basepoint has been chosen, the only ambiguity in defining the lift arises from  vertices corresponding to a component on which  $f$ is non-injective.   In the sequel of the paper, we will use this property without mentioning it explicitly.

Lifts preserve lengths: if $\gamma'$ is a lift of a curve $\gamma$, then $\ell(\gamma)=\ell(\gamma')$.

\item[--] $F$ maps shadows onto shadows:  for any $|\xi|\ge 2$, $F(\mho(\xi))=\mho(F(\xi))$.

To see this, note that since $F$ maps rays to rays, it follows that $F(\mho(\xi))\subset\mho(F(\xi))$. For the converse, let $\zeta\in \mho(F(\xi))$ and let us consider a geodesic
curve $\gamma$ joining $F(\xi)$ to $\zeta$. The function $t\mapsto |\gamma(t)|$ is strictly monotone. Since $F$ has the lifting property, there is a 
strictly monotonic geodesic curve $\gamma'$ starting from $\xi$ such that $F(\gamma')= \gamma$. This curve can be extended geodesically to the base point $o$.
It follows that $\mho(F(\xi))\subset F(\mho(\xi))$, which proves the claim.
\end{itemize}

Let $e$ be an edge in $\Gamma$ at distance at least $1$ from the origin $o$. Then, since $F|_e$ is injective, 
$$\ell_{\varepsilon}(F(e))=\int_{F(e)}e^{-\varepsilon |x|}ds(x)=\int_{e}e^{-\varepsilon |F(x)|}ds(x)=
 \int_e e^{-\varepsilon (|x|-1)}ds(x)= e^\varepsilon \ell_{\varepsilon}(e)\,.$$
Let now  $\xi,\zeta\in \Gamma\setminus\overline{B(o,1)}$ and let $\gamma$ be a geodesic segment joining these points.  
Since two edges can be mapped to a common one, and since $\gamma$ may contain the origin,  
$$d_\varepsilon(F(\xi),F(\zeta))\le\int_{F(\gamma)}e^{-\varepsilon |\gamma(t)|}dt\le
 \int_\gamma e^{-\varepsilon (|\gamma(t)|-1)}dt\le e^\varepsilon d_\varepsilon(\xi,\zeta)\,.$$
\gap

Therefore $F$ is uniformly continuous and so extends to an $e^\varepsilon$-Lipschitz map $F:\cl{\Gamma}_\varepsilon\to \cl{\Gamma}_\varepsilon$.

\begin{prop}\label{boulesF} For any $\xi\in \overline{ \Gamma}_\varepsilon$, and 
$r<|F(\xi)|_\varepsilon$, $F(B_\varepsilon(\xi,re^{-\varepsilon}))=B_\varepsilon(F(\xi),r)$.\end{prop}

Hence $F$ is an open mapping.

\gap

\pf  We already know that $F(B_\varepsilon(\xi,re^{-\varepsilon}))\subset B_\varepsilon(F(\xi),r)$.  
Let us
consider $\zeta'\in B_\varepsilon(F(\xi),r)$ and $\gamma'$
a $d_\varepsilon$-geodesic curve joining $F(\xi)$ to $\zeta'$. Since $r<|F(\xi)|_\varepsilon$, it follows that $\gamma'$ avoids $o$.
We let $\gamma$ be a  
lift of
$\gamma'$  which joins $\xi$ to a point $\zeta\in  
\overline{\Gamma_\varepsilon}$.
It follows that
$$|\xi-\zeta|_\varepsilon\le \ell_\varepsilon(\gamma)=\int_\gamma \rho_\varepsilon(\xi)ds(\xi)=e^{-\varepsilon}\int_\gamma \rho_\varepsilon(F(\xi))ds(\xi)=e^{-\varepsilon}\int_{\gamma'}
\rho_\varepsilon(\xi)ds(\xi)$$
so $|\xi-\zeta|_\varepsilon\le e^{-\varepsilon} |F(\xi)-\zeta'|_\varepsilon\le  e^{-\varepsilon} r$ and $\zeta\in  
B_\varepsilon(\xi,
e^{-\varepsilon} r)$.\qed

\gap

The following proposition says that if $F^n$ is injective on a ball, then it is a similarity on the ball of one-fourth the size.

\begin{prop}
\label{prop:tentlike}
Suppose $B=B(\xi,r) \subset \cl{\Gamma}_\varepsilon$ and $F^n|_B: B \to B(F(\xi), e^{n\varepsilon} r)$ is a homeomorphism.  Then for all $\zeta_1, \zeta_2 \in B(\xi, r/4)$, 
\[ |F^n(\zeta_1) - F^n(\zeta_2)|_\varepsilon = e^{n\varepsilon}|\zeta_1-\zeta_2|_\varepsilon.\]
\end{prop}

\pf We first claim that the above equality holds when $\zeta_1=\xi$ and $\zeta=\zeta_2$ is an arbitrary point in $B$.  The upper bound is clear.  To show the lower bound, 
notice that $F^{-n}:B(F^n(\xi),re^{\varepsilon n})\to B(\xi,r)$ is well defined, and let $\gamma\subset B(F^n(\xi),re^{\varepsilon n})$ be a  curve joining $F^n(\xi)$ to $F^n(\zeta)$.
It follows that $F^{-n}(\gamma)$ is a curve joining $\xi$ to $\zeta$ inside $B$, so the proof of Proposition \ref{boulesF} shows that 
\[\ell_{\varepsilon}(F^{-n}(\gamma))=e^{-\varepsilon n}\ell_{\varepsilon}(\gamma)\,.\] 
Since $F^n|_B$ is a homeomorphism, the claim follows.

The proposition follows immediately by applying the claim to the ball centred at $\zeta_1$ of radius $|\zeta_1-\zeta_2|_\varepsilon$, 
which by hypothesis is contained in $B$ and hence maps homeomorphically onto its image under $F^n$.  
\qed

{\noindent\bf Comparison of $X$ and $\partial \Gamma$.} For any $x\in X$ and $n \in \N$, let $W_n\in S(n)$ contain $x$.   The sequence $(W_n)$ defines a ray $\gamma_x$ in $\RRR_\infty$ such that $\gamma_x(n)=W_n$. There is a natural map $\phi_f: X \to \bdry_\varepsilon \Gamma$ defined by 
$\phi_f(x)=\pi (\gamma_x)$.
In other words:  the sequence $(W_n)$ is a Cauchy sequence in $\overline{\Gamma}_\varepsilon$, and we let $\phi_f(x)$ be its limit.  This map is well defined\,:
if $(W_n')$ is another sequence contained in a ray $\gamma_x'$, then $d(W_n,W_n')\le 1$ since $x \in W_n \intersect W'_n \intersect X$, so  
$\pi (\gamma_x)=\pi (\gamma_x')$. 
Furthermore, $F\circ \phi_f=\phi_f\circ f$
on $X$.

\begin{prop}\label{compvp} The map $\phi_f:X\to \partial_\varepsilon \Gamma$ is continuous and onto.
\end{prop}

\pf  To prove surjectivity, suppose $\xi \in \bdry\Gamma_\varepsilon$.  By Lemma \ref{toprayon}, there exists a ray $\gamma \in \RRR_\infty$ such that $\pi(\gamma) = \xi$.  For $k \in \N$ let $W_k = \gamma(k)$, so that $W_k \in S(k)$.  Then $W_k \in \cl{\Gamma}_\varepsilon$ and $\xi = \lim W_k$.  But each $W_k$ is also a subset of $\XX_0$ whose intersection with the repellor $X$ contains some point $w_k$.  Since $X$ is compact, there exists a limit point $x$ of $(w_k)$.  

We claim $\phi_f(x) = \xi$.  By definition $\phi_f(x) = \lim V_n$, where $V_n$ is an arbitrary element of $S(n)$ which as a subset of $X$ contains $x$ and where the limit is in $\cl{\Gamma}_\varepsilon$.  Then 
for each $n \in \N$, since $V_n$ is open and $w_k \to x$, there exists $k(n) \in \N$ such that $W_k \intersect V_n \intersect X \neq \emptyset$ for all $k \geq k(n)$.  By the definition of shadows, $W_k \subset \mho(V_n)$.  By Lemma \ref{omdiam},  $|W_{k(n)} - V_n|_\varepsilon \to 0$ as $k \to \infty$.   Hence $\xi = \lim W_k = \lim W_{k(n)} = \lim V_n = \phi_f(x)$ as required.   

To prove continuity, suppose $x_k \to x \in X$.  For all $n \in \N$ choose $W_n \in S(n)$ containing $x$, 
so that $\xi = \phi_f(x)=\lim W_n \in \mho_\infty(W_n)$.  Then for all $n \in \N$ there exists $k(n)$ such that 
$x_k \in W_n$ for all $k \geq k(n)$.  By the definition of $\phi_f$, $\phi_f(x_k) \in \mho_\infty(W_n)$.  
By Lemma \ref{omdiam}, $|\phi_f(x_k) - \xi|_\varepsilon \leq Ce^{-\varepsilon n} \to 0$ as $n \to \infty$ and
 so $\phi_f(x_k) \to \xi = \phi_f(x)$.  \qed

We now turn to the proof of Theorem  \ref{thm:construction}.  
We first prove the existence of a preliminary metric in which the diameters of the sets $\phi_f(U)$, $U \in \UUU_n$, tend to zero exponentially fast in $n$.  

\begin{thm}\label{thm:exp_decay}
Suppose Axiom [Expansion] holds.  Then there exists a metric on the 
repellor $X$ and constants $C>1, \theta < 1$ such that for all $n \geq 
0$,
\[ \sup_{U \in \UUU_n} \diam U \leq C\theta^n.\]
\end{thm}

The proof is standard and mimics the proof of a preferred H\"older structure given a uniform structure; see \cite[Chap.\,II]{bour:top_gen}.  

\gap

\pf Let $N_0$ be given by Proposition \ref{prop:elementary}, 2(c) and 
put $g=f^{N_0}$, $\VVV_0 = \union_{j=0}^{N_0-1}\UUU_j$, $\YY_1 = 
f^{-N_0}\XX_0$, $\YY_0=\XX_0$, and $\VVV_n = g^{-n}\VVV_0$.  Then $g$ 
is a finite branched covering, the repellor of $g$ is $X$ 
(by total invariance), and the mesh of $\VVV_n$ tends to zero.  The 
conclusion of the above proposition (applied $U_1' = U_2' = V')$ and the definition of 
$g$ implies
\begin{equation}
\label{eqn:nested}
\forall V' \in \VVV_n,\ \exists V \in \VVV_{n-1}\; \mbox{ with } \; V' 
\subset V.
\end{equation}   Then (\ref{eqn:nested}) and conclusion 2(b) of Proposition \ref{prop:elementary} imply immediately 
that for any distinct $x, y \in X$, the quantity
\[ [x|y] = \max\{ n : \forall 1 \leq i \leq n, \exists V_i \in \VVV_i 
\; \mbox{ with } \{x,y\} \subset V_i\}\]
is finite.   For $x=y$ set $[x,y]=\infty$.  The statement (\ref{eqn:nested}), Proposition \ref{prop:elementary} 2(c) and
the definition of $g$ imply  
for any triple $x, y, z \in X$,
\[ [x|z] \geq \min\{[x|y], [y|z]\} -1.\]
Fix $\epsilon>0$ small, and define
\[ \varrho_\epsilon(x,y) = \exp(-\epsilon [x|y]).\]
Then  $\varrho_\epsilon(x,y) = 0$ if and only 
if $x=y$, and indeed $\varrho_\epsilon$ satisfies all properties of a 
metric save the triangle inequality.  Instead, we have
\[ \varrho_\epsilon(x,z) \leq e^\epsilon\max\{\varrho_\epsilon(x,y),
\varrho_\epsilon(y,z)\}.\]

There is a standard way to extract a metric bi-Lipschitz equivalent to $\varrho_\varepsilon$; 
see \S\,\ref{secn:expanding_maps} for an outline.
If $\epsilon<\frac{1}{2}\log 2$ then Proposition 7.3.10 of 
\cite{ghys:delaharpe:groupes} implies that there is a metric 
$d_\epsilon$ on $X$ satisfying
\[ (1-2\sqrt{2})\varrho_\epsilon(x,y) \leq d_\epsilon(x,y) \leq \varrho_\epsilon(x,y).\]
Letting $\diam_\epsilon$ denote the diameter with respect to $d_\epsilon$,
it is clear from the definitions that $V \in \VVV_n$ implies 
$\diam_\epsilon V  \leq \exp(-n\epsilon)$.  It is then easy to check 
that taking $\theta = \exp(-\epsilon/N_0)$ and
\[ C = \max\{ \diam_\epsilon U : U \in \union_{i=0}^{N_0-1}\UUU_i\}\]
will do.
\qed

\gap

We may now prove Theorem \ref{thm:construction}.

\gap

\pf (Theorem \ref{thm:construction}) The graph $\Gamma$, the map $\phi_f$ and $F$ satisfy the three first points of the 
theorem by construction and  according to  Proposition \ref{compvp} and \ref{boulesF}.

We assume from now on that [Expansion] holds.
It follows from Theorem \ref{thm:exp_decay} that there is a metric $d_X$ and constants $C>0$ and $\theta\in  
(0,1)$ such
that, for any $W\in S(n)$, $\diam\,W\le C\theta^n$, where the diameter is with respect to the metric $d_X$.
Let $x,y\in X$.
Let $\gamma: (-\infty, \infty) \to \Gamma$ be a curve joining $\phi_f(x)$ to $\phi_f(y)$ such that for all $n \in \Z$, $\gamma(n)$ is a vertex of $\Gamma$ corresponding to an open set $W_n$, and such that $\gamma|_{[n, n+1]}$ traverses a closed edge of $\Gamma$ exactly once.   Then 
\[ \ell_\varepsilon(\gamma) = \sum_{n\in\Z}  \ell_\varepsilon(\gamma|_{[n, n+1]}) \asymp \sum_{n\in\Z} e^{-\varepsilon |\gamma(n)|} = \sum_{n\in\Z}  e^{-\varepsilon|W_n|}= \sum_{n\in\Z}  
\frac{ e^{-\varepsilon
|W_n|}}{\diam W_n}\diam W_n\,.\] 

If $\varepsilon>0$ is small enough, then 
$e^{-\varepsilon |W_n|}\ge \theta^{|W_n|}\ge  (1/C)\diam W_n$.
Furthermore, there are points
$z_n$ such that $z_n\in W_n\cap W_{n+1}$.  For all $k\in\N$, we let $\gamma_k$ be the subcurve 
of $\gamma$ joining $W_{-k}$ to $W_k$. Then
$$\begin{array}{ll}\ell_\varepsilon(\gamma_k) & \ge (1/C)\dis\sum_{|n|\le k} \diam W_n \\&\\
& \ge (1/C)\dis\sum_{|n|\le k} d_X(z_n,z_{n+1})  \\&\\
& \ge (1/C)d_X(z_{-k}, z_k)\,.\end{array}$$
Since $\gamma_k$ is a subset of $\gamma$, $\{ z_{-k}, z_k\}$ tends to $\{x,y\}$; this implies that $\ell_\varepsilon(\gamma)\ge (1/C)d_X(x,y)$,
where $C$ is independent of $\gamma$. Therefore, $d_\varepsilon(\phi_f(x),\phi_f(y))\ge (1/C)d_X(x,y)$.

\qed

\section{Geometry of $\Gamma$} 

Let $f: \XX_1 \to \XX_0$ and $\UUU$ satisfy the conditions listed at the  beginning of \S\,\ref{secn:spaces_associated}.  The main result of this section is

\begin{thm}\label{homeoimphyp} If $(f,\UUU)$ satisfies [Expansion], then $\Gamma$ is Gromov hyperbolic.  
If $(f, \VVV)$ also satisfies [Expansion], then $\Gamma(f, \UUU)$ is quasi-isometric to $\Gamma(f,\VVV)$.  
If $g=f^n: \XX_{n-1} \to \XX_0$, then $\Gamma(g,\UUU)$ is quasi-isometric to $\Gamma(f,\UUU)$.
\end{thm}

Hence, as long as the expansion axiom is satisfied, the quasi-isometry class of $\Gamma$ is an invariant of the conjugacy class of $f: \XX_1 \to \XX_0$.  

The proof of Thm \ref{homeoimphyp} will follow from both Proposition  \ref{phomeoimphyp} and Proposition \ref{homeoimpqi}.

\gap

The naturality results given above are analogous to those enjoyed by Cayley graphs of finitely generated groups.

\subsection{Metric estimates}

We start by gathering information on the geometry of balls, and how they interact with the coverings.

Our main estimates are the following, which assert that the elements $\phi_f(W)$ enjoy geometric properties with respect to the metric $d_\varepsilon$ similar to those enjoyed by the sets $U$ with respect to a metric for a cxc map; compare Propositions \ref{prop:uniform_roundness} (Uniform roundness) and \ref{prop:BLC} (Balls are like connected sets).  
\gap

\noindent{\bf Notation.}  For an element $W\in S(n)$, regarded as a subset of $\XX_1$, we denote by $\phi_f(W)$ the set $\phi_f(W\intersect X)$.   

\begin{prop}\label{ba3} Fix $(f,\XXX_1,\XXX_0,X,\UUU)$, and let us consider the graph $\Gamma$.
We assume that $\phi_f:X\to\partial_{\varepsilon}\Gamma$ is a homeomorphism for all $\varepsilon>0$ small enough.
\be
\item There is some constant $C>1$ such that, for all $W\in V\setminus\{o\}$, there is a point $\xi\in\phi_f(W)$ so that
$$B_{\varepsilon}(\xi,(1/C)e^{-\varepsilon|W|}) \subset   
\phi_f(W)\subset\mho_\infty(W)\subset B_\varepsilon(\xi,C e^{-\varepsilon|W|}).$$

\item There is a radius $r_1>0$ such that, for any $n\ge 1$ and for any $\xi\in\partial_\varepsilon \Gamma$, there
is some $W\in S(n)$ so that $B_{\varepsilon}(\xi,r_1 e^{-\varepsilon n})\subset \phi_f(W)$.

\item A maximal radius $r_0>0$ exists such that, for any $r\in (0,r_0)$ and any  $\xi\in\partial_\varepsilon \Gamma$,
there exist $W$ and $W'$ in $\mathbf{U}$  such that $|W-W'|=O(1)$, $$\phi_f(W')\subset B_\varepsilon(\xi,r)\subset \phi_f(W),$$
and $$\max\{\roundness(\phi_f(W),\xi),\roundness(\phi_f(W'),\xi)\}=O(1).$$
\eb\end{prop}

Since the set $\XX_0$ is not assumed to be endowed with a metric, we shall use uniform structures \cite[Chap.\,II]{bour:top_gen}.

Since $\XX_1$ has compact closure in $\XX_0$, there is a unique uniform structure on $\cl{\XX_1}$ compatible with its topology.  
We consider the uniform structure on $\XX_1$ induced by the one on $\cl{\XX_1}$. Let
us recall that an {\em entourage}
\index{Index}{entourage}%
 $\Omega$ is a neighborhood of the diagonal of $\XX_1\times\XX_1$.
If $x\in\XX_1$, then $\Omega(x)=\{y\in\XX_1,\ (x,y)\in\Omega\}$.

\begin{prop}\label{ba2} Given an entourage $\Omega$, there is some constant $r=r(\Omega)>0$ such that, 
whenever $U\in S(1)$, $u\in U\cap X$ 
and $\Omega(u)\subset U$, then, for any $n\ge 1$, any $\wtU\in S(n)$ such that $f^{n-1}(\wtU)=U$, and any 
preimage $\tu\in\wtU\cap f^{-(n-1)}(\{u\})$, the ball $B_{\varepsilon} (\phi_f(\tu), re^{-\varepsilon n})$ is contained in 
$\phi_f(\wtU)$.\end{prop}

Let us first prove some lemmata.

\begin{lemma}\label{1ba2} Let $\gamma:\R\to \Gamma$ be a curve such that $\gamma(\Z)\subset V\setminus\{o\}$ and which connects
two points $u$ and $v$ from the boundary. Let $r\in (0,1/\varepsilon)$.  If $\ell_{\varepsilon}(\gamma)<r$ then 
$$\overline{\cup_{n\in\Z}\phi_f(\gamma(n))}\subset B_\varepsilon(u,r).$$\end{lemma}

See 
Figure 3.4.

\begin{figure}
\label{fig:Lemma3_3_4}
\psfragscanon
\psfrag{1}{$u$}
\psfrag{2}{$v$}
\psfrag{3}{$r$}
\psfrag{4}{$\phi_f$}
\psfrag{5}{$X$}
\psfrag{6}{$\bdry_\varepsilon\Gamma$}
\psfrag{7}{$\gamma$}
\begin{center}
\includegraphics[width=6in]{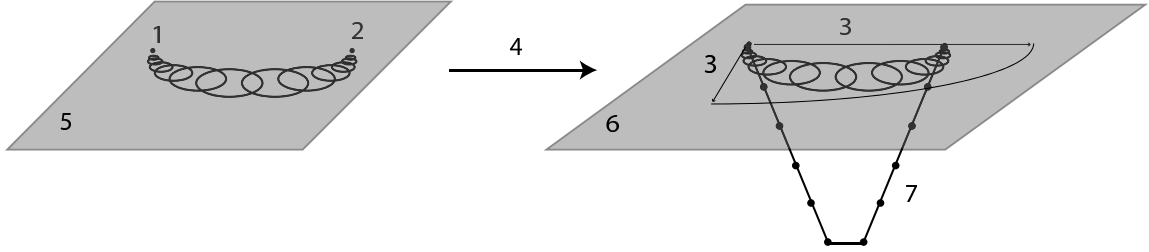}
\end{center}
\caption{{\sf  For simplicity, the figure is drawn as if $X=\XX_1$.  The points $u, v \in X$ are shown at left.  At right, the edge-path $\gamma$ in $\Gamma$ has closure in $\cl{\Gamma}_\varepsilon$ joining $\phi_f(u)$ to $\phi_f(v)$ and has length $< r$.  The conclusion of the lemma says that the image in $\bdry_\varepsilon\Gamma$ under $\phi_f$ of the open sets in $\mathbf{U}$ comprising the vertices of $\gamma$ is contained in the ball of radius $r$ about $\phi_f(u)$.}}
\end{figure}

\pf For any fixed $n$ and any $z\in \phi_f(\gamma(n))$, 
$$|z-u|_\varepsilon\le |z-\gamma(n)|_\varepsilon + |\gamma(n)-u|_\varepsilon\,.$$
Since $z\in\phi_f(\gamma(n))$, there is a geodesic ray $[\gamma(n),z)$ contained in some ray
in $\RRR_\infty$, so that 
$$|z-\gamma(n)|_\varepsilon= 
\mbox{\rm dist}_\varepsilon(\gamma(n),\partial_\varepsilon \Gamma)
\le |\gamma(n)-v|_\varepsilon\le \ell_\varepsilon(\gamma|_{[n,\infty[})$$ and 
$|\gamma(n)-u|_\varepsilon\le \ell_\varepsilon(\gamma|_{]-\infty,n]})$ so that 
$|z-u|_\varepsilon\le \ell_\varepsilon(\gamma)< r$. Therefore
$z\in \overline{B_\varepsilon(u,\ell_\varepsilon(\gamma) )}$ and 
$\overline{\phi_f(\gamma(n))}\subset\overline{B_\varepsilon(u,\ell_\varepsilon(\gamma) )}$ for all $n\in\Z$.
Hence $$\overline{\cup_{n\in\Z}\phi_f(\gamma(n))}\subset B_\varepsilon(u,r).$$\qed

{\noindent\bf Definition.} Given $u\in\partial_\varepsilon \Gamma$ and $r\in (0,1/\varepsilon)$, we let $V(u,r)$ be the set
of all vertices of $\Gamma$ contained in curves of $d_\varepsilon$-length less than $r$ joining $u$ to another boundary point.

\gap

It follows from the lemma above that $\phi_f(U)\subset  B_\varepsilon(u,r)$ for any $U\in V(u,r)$.

\begin{lemma}\label{3ba2} Let  $\Omega$ be an entourage of $\XX_1$.
There is a radius $r>0$ which depends only on $\Omega$ such that, for any $u\in\partial_{\varepsilon}\Gamma$, we have
 $W\subset\Omega(\phi_f^{-1}(u))$ 
whenever $W\in V(u,r)$.
\end{lemma}

\pf 
The uniform continuity of $\phi_f^{-1}$ provides us with a radius $r$ (independent from $u$) such that, 
for any $u\in\partial_{\varepsilon} \Gamma$,  $\phi_f^{-1}(B_\varepsilon(u,r))\subset \Omega(\phi_f^{-1}(u))$.
It follows from Lemma \ref{1ba2} that if $W\in V(u,r)$, then $\phi_f(W)\subset B_\varepsilon(u,r)$ so that $W\subset\Omega(\phi_f^{-1}(u))$.
\qed

\gap

We are now ready for the proofs of the Propositions.

\gap

\pf (Proposition \ref{ba2}). Let $\Omega$ be an entourage of $\XX_1$, $U\in S(1)$, $u\in U\cap X$, satisfy 
$\Omega(u)\subset U$. Let us choose another 
entourage $\Omega_0$ such that
$\overline{\Omega_0}\subset \Omega$. 

Choose $n\ge 1$, $\wtU\in S(n)$ and $\tu\in f^{-(n-1)}(\{u\})\cap\wtU$ such that $f^{n-1}(\wtU)=U$. 
Let us also consider the constant
$r>0$ given by Lemma \ref{3ba2} applied to $\Omega_0$.

\gap

Let $\tv\in \phi_f^{-1}(B_\varepsilon(\phi_f(\tu),re^{-\varepsilon (n-1)}))$ and $\gamma$ be a curve joining $\phi_f(\tu)$ to 
$\phi_f(\tv)$
of $d_\varepsilon$-length less than $re^{-\varepsilon (n-1)}$. Set $$K=\overline{\cup_{n\in\Z}\gamma(n)}\subset \XX_0.$$
Then $K$ is a continuum by definition which joins $\tu$ to $\tv$. Therefore, $f^{n-1}(K)$ joins
$u$ to $f^{n-1}(\tv)=v$, and $F^{n-1}(\phi_f(K))\subset B_\varepsilon(\phi_f(u),r)$. 
By Lemma \ref{3ba2},
$f^{n-1}(\gamma(k))$ is in $\Omega_0(\phi_f^{-1}(u))$ for any $k\in\Z$, so that $f^{n-1}(K)\subset U$.
It follows that $K\subset \wtU$ since $f^{n-1}:\wtU\to U$ is proper and $K$ is connected.\qed

\gap

\pf (Proposition \ref{ba3}) 
Let $\Omega$ be an entourage such that, for any $x\in X$, there is some $U\in S(1)$ such that $\Omega(x)\subset U$.
\be
\item 
Let $n$ be the level of $W$ and pick some $x'\in (X\cap f^{n-1}(W))$.

Let $x\in f^{-(n-1)}(\{x'\})\cap W$; it follows from Proposition \ref{ba2} that
 $\phi_f(W)$
will contain the ball $B_{\varepsilon}(\xi, r e^{-\varepsilon n})$ where $r=r(\Omega)$ and $\xi=\phi_f(x)$.

\gap

Furthermore, Lemma \ref{omdiam} implies that 
$\diam_{\varepsilon}\phi_f(W)\asymp \diam_{\varepsilon}\mho(W)\asymp e^{-\varepsilon|W|}$.
It follows that there is some constant $C>1$ such that, for all $W\in V$, there is a point $\xi\in\phi_f(W)$ so that
$$B_{\varepsilon}(\xi,(1/C)e^{-\varepsilon|W|}) \subset   
\phi_f(W)\subset\mho_\infty(W)\subset B_\varepsilon(\xi,C e^{-\varepsilon|W|}).$$

\item Similarly, Proposition \ref{ba2} implies that, 
for any $n\ge 1$, there is some $W\in S(n)$ such that $\phi_f(W)$
will contain the ball $B_{\varepsilon}(\xi,r_1 e^{-\varepsilon n})$ where $r_1=r(\Omega)$ is given by the proposition.

\item  Fix $r\in (0,\delta)$ and $\xi\in\partial \Gamma$, where $\delta$ is the Lebesgue number of $S(1)$ in 
$\partial_\varepsilon\Gamma$. 

It follows from point 1. above that, for any $n$
and any $W\in S(n)$, $\diam_\varepsilon\phi_f(W)\asymp e^{-\varepsilon n}$. 

Moreover, from point 2., 
there is some
$W$ such that $\roundness(\phi_f(W),\xi)=O(1)$ and $B_{\varepsilon}(\xi,r_1 e^{-\varepsilon n})\subset\phi_f(W)$.
Let $m\ge n$ be so that the diameter of any element of $S(m)$ is at most $r$. It follows from the diameter control above
that $m$ may be chosen so that $|m-n|=O(1)$.  Point 2. provides us with an element $W'\in S(m)$ so that
$\roundness(\phi_f(W'),\xi)=O(1)$ and $\phi_f(W')\subset  B_\varepsilon(\xi,r)$.\qed
\eb

As a consequence of Proposition \ref{ba3} and its proof, we obtain the following.
For $n \in \N$, let $\VVV_n = \{\phi_f(U): U \in \UUU_n\}$.  Thus for each $n$, the collection $\VVV_n$ is a covering of $\partial_\varepsilon\Gamma$ by open sets which, in general, need not be connected.  

\begin{prop}
\label{prop:F_is_nearly_cxc} 
The map $F: \bdry_\varepsilon\Gamma \to \bdry_\varepsilon\Gamma$ and the sequence of coverings $\VVV_n, n=0, 1, 2,\ldots$ together satisfy Axioms 
[Expansion], [Roundness distortion] and [Relative diameter distortion].
\end{prop}

\noindent{\bf Remark:}  If in addition axiom [Degree] is satisfied, it would be tempting to assert that 
$F: \bdry_\varepsilon\Gamma \to \bdry_\varepsilon\Gamma$ is also metrically cxc.  
However, even though $\overline{\Gamma_{\varepsilon}}$ 
is locally connected (since shadows define connected neighborhoods of points at infinity), the boundary $\bdry_\varepsilon\Gamma$ need not be (locally) connected.  
Our definition of metrically cxc is not purely intrinsic to the dynamics on the repellor $X$ since we require that the covering $\UUU_0$ consists of connected sets which are contained in an a priori larger space $\XX_1$.  Unfortunately, in general we do not know how to modify the definition of $\Gamma$ so that $F: \cl{\Gamma}_\varepsilon \to \cl{\Gamma}_\varepsilon$ 
becomes a finite branched covering map on an open connected neighborhood of $\bdry_\varepsilon\Gamma$.  
If this were possible, it seems likely that one could then establish a variant of Proposition \ref{prop:F_is_nearly_cxc} 
in which the conclusion asserted that the model dynamics was indeed metrically cxc.
\gap

\pf The [Expansion] axiom follows from Lemma \ref{omdiam}.
The forward and backward relative diameter distortion bounds follow immediately from Proposition \ref{ba3}.  Since $F$ maps round balls in the metric $d_\varepsilon$ to round balls, the forward roundness function $\rho_+$ may be taken to be the identity.   We claim that we may take the backward roundness distortion function to be linear.

First, suppose $F^n:(\wtV,\txi) \to (V,\xi)$ where $V=\phi_f(W)$ and $W\in S(k)$. 
 Suppose $B(\xi,r) \subset V \subset B(\xi, Kr)$.  Then $Kr\asymp e^{-\varepsilon k}$.
Proposition \ref{ba2} shows that $B(\txi, c e^{-\varepsilon n}r) \subset \wtV$ for some uniform constant $c>0$. 
 By Proposition \ref{ba3}, $\diam_\varepsilon(\wtV) \asymp e^{-(n+k)\varepsilon}$.  Hence 
\[ \roundness(\wtV, \txi) \lesssim\frac{e^{-(n+k)\varepsilon}}{e^{-n\varepsilon}r} \asymp K\asymp \roundness(V, \xi).\]
\qed

\begin{prop}
\label{prop:degree_fails}
Suppose axiom [Expansion] holds.  Let $Y$ denote the set of points $y 
$ in $X$ such that there exists an element $U'$ of $\mathbf{U}$  
containing $y$ such that all iterated preimages $\wtU'$ of $U'$ map  
by degree one onto $U'$.

If Axiom [Degree] fails, and if $Y \intersect X$ is dense in $X$,  
then $\bdry_\varepsilon\Gamma$ fails to be doubling.
\end{prop}

\noindent{\bf Remarks:}
\be
\item We have always $f^{-1}(Y) \subset Y$.  If Axiom  
[Irreducibility] holds and $Y$ is nonempty then $Y$ is dense in $X$,  
so the above proposition implies that $\bdry_\varepsilon\Gamma$ fails  
to be doubling.

\item It is reasonable to surmise that $Y =X-P_f$ -- this is the case  
e.g. for rational maps.  However, we have neither a proof nor counterexamples.
\eb

\pf
Suppose Axiom [Degree] fails.  It follows easily that then there  
exists some $U\in\UUU_0$ such that for all $p \in \N$, there exists  
$n \in \N$ and a preimage $\wtU \in \UUU_n$ of $U$ such that $f^n:  
\wtU \to U$ has degree at least $p$.  The assumption and axiom  
[Expansion] imply that there exist $U' \subset U$, $U' \in \UUU_N$,   
independent of $p$ and of $\wtU$, such that $\wtU$ contains at least $p 
$ disjoint preimages $\wtU'$ of $\wtU$.

By Proposition \ref{ba3}, $\phi_f(\wtU)$ and $\phi_f(\wtU')$ are  
uniformly almost round, $\diam_\varepsilon(\phi_f(\wtU)) \asymp \exp(- 
\varepsilon n)$, and $\diam_\varepsilon(\phi_f(\wtU)) \asymp \exp(- 
\varepsilon (n+N))$.  So at least $p$ balls of radius $C'\cdot  \exp(- 
\varepsilon(n+N))$ are needed to cover a ball of radius $C \exp(- 
\varepsilon n)$, where $C', C$ are independent of $n$.  Therefore $ 
\bdry_\varepsilon\Gamma$ fails to be doubling.

\qed

We close this section with the following consequence of Proposition \ref{ba3} which will be useful in our  characterization rational maps; cf. Definition \ref{defn:llc} and Corollary \ref{cor:get_llc}.

\begin{cor}\label{cor:llc}  If for each $W\in \cup S(n)$, the sets $\phi(W \intersect X)$ 
and $X\setminus \phi(W \intersect X)$ are connected, then $\partial\Gamma$ is 
linearly locally connected.\end{cor}

\pf Let us fix $B_{\varepsilon}(\xi,r)$. Proposition \ref{ba3}, point 3, implies the existence of vertices $W,W'$
such that $||W|-|W'||\le C_1$ for some universal constant $C_1$ and such that $$\phi(W)\subset B_{\varepsilon}(\xi,r)\subset
\phi(W').$$
Therefore, $\diam_{\varepsilon} \phi(W)\asymp\diam_{\varepsilon} B_{\varepsilon}(\xi,r)\asymp \diam_{\varepsilon}\phi( W')$.

If $\zeta,\zeta'\in B_{\varepsilon}(\xi,r)$, then they are connected by $\phi(W')$ which is connected by assumption.
Similarly, if $\zeta,\zeta'\notin B_{\varepsilon}(\xi,r)$ then they are joined within $X\setminus\phi(W)$.
\qed

\subsection{Hyperbolicity}
We are now ready to prove the first part of  Theorem  \ref{homeoimphyp}.  

\begin{prop}\label{phomeoimphyp} If $\phi_f:X\to \partial_{\varepsilon}\Gamma$ is a homeomorphism, then $\Gamma$ is hyperbolic.\end{prop}

The proof is an adaptation of Proposition 2.1 in \cite{bourdon:pajot:besov} and its main step 
 is given by the following lemma.

\begin{lemma}\label{lsehypimphyp} For any $W,W'\in V$, 
$$\diam_{\varepsilon}(\phi_f(W)\cup\phi_f(W'))\asymp e^{-\varepsilon (W|W')}\,.$$\end{lemma}

\pf We assume that $|W'|\ge |W|$.
We let $n\in\N\cup\{\infty\}$ be the smallest integer such that 
$$ \mbox{dist}_\varepsilon(\phi_f(W),\phi_f(W'))\ge r_1 e^{-\varepsilon n}$$
where $r_1$ is the constant given by Proposition \ref{ba3}.
Let $(\xi,\xi')\in \overline{\phi_f(W)}\times \overline{\phi_f(W')}$ satisfy 
$ \mbox{dist}_\varepsilon(\phi_f(W),\phi_f(W'))=|\xi-\xi'|_\varepsilon$.
Let $m=\min\{|W|,n\}-1$; there is some $C\in S(m)$ such that $B_{\varepsilon}(\xi,r_1e^{-\varepsilon m})\subset \phi_f(C)$, 
so that
$W,W'\in\mho(C)$ and $$\diam_{\varepsilon}(\phi_f(W)\cup\phi_f(W'))\ge 
\max\{\mbox{dist}_\varepsilon(\phi_f(W),\phi_f(W')), \diam_\varepsilon\phi_f(W)\}.$$
The maximum of $\mbox{dist}_\varepsilon(\phi_f(W),\phi_f(W'))$ and of $\diam_\varepsilon\phi_f(W)$ is at least
of order $e^{-\varepsilon m}$.
Hence
$$\diam_{\varepsilon}(\phi_f(W)\cup\phi_f(W'))\gtrsim  e^{-\varepsilon |C|}\,.$$
Since $(W|W')\ge |C|$, it follows that $$\diam_{\varepsilon}(\phi_f(W)\cup\phi_f(W'))\gtrsim e^{-\varepsilon (W|W')}\,.$$

\gap

For the other inequality, we let $\{W_j\}_{0\le j\le |W-W'|}$ be a geodesic chain which joins $W$ to $W'$.
For convenience, set $m=|W|$, $m'=|W'|$ and $D=|W-W'|$.
Then $$\begin{array}{ll}\diam_{\varepsilon}(\phi_f(W)\cup\phi_f(W'))& \le \dis\sum_{0\le j\le D}\diam_{\varepsilon}(\phi_f(W_j))\\ &\\
&\le \dis\sum_{0\le j\le k}\diam_{\varepsilon}(\phi_f(W_j)) +\dis\sum_{k+1\le j\le D}\diam_{\varepsilon}(\phi_f(W_j)) \\ &\\
&\lesssim \dis\sum_{0\le j\le k} e^{-\varepsilon(m-j)} +\dis\sum_{0\le j\le D-(k+1)}e^{-\varepsilon(m'-j)}  \\ &\\
&\lesssim e^{-\varepsilon(m-k-1)} +e^{-\varepsilon(m'-D+k))} \end{array}$$
Choosing $k=(1/2)(D +m-m')$, one gets
$$\diam_{\varepsilon}(\phi_f(W)\cup\phi_f(W'))\lesssim e^{-\varepsilon(1/2)(m'+m-D)} \lesssim e^{-\varepsilon(W|W')}\,.$$
The lemma is established.\qed

\gap

\pf (Prop.\,\ref{phomeoimphyp}) It follows from Lemma \ref{lsehypimphyp} that if $W_1,W_2,W_3$ are three
vertices, then $$\begin{array}{ll}e^{-\varepsilon(W_1|W_3)} & \lesssim\diam_{\varepsilon}(\phi_f(W_1)\cup\phi_f(W_3))\\
& \lesssim\diam_{\varepsilon}(\phi_f(W_1)\cup\phi_f(W_2))+\diam_{\varepsilon}(\phi_f(W_2)\cup\phi_f(W_3))\\
& \lesssim e^{-\varepsilon(W_1|W_2)}+e^{-\varepsilon(W_1|W_3)}\\
& \lesssim \max\{e^{-\varepsilon(W_1|W_2)},e^{-\varepsilon(W_2|W_3)}\}\end{array}$$
so that there is a constant $c$ such that
$$(W_1|W_3)\ge \min\{(W_1|W_2),(W_2|W_3)\} -c\,.$$
This proves the hyperbolicity of $X$.\qed

The hyperbolicity of $\Gamma$ implies that the homeomorphism and quasisymmetry type of $\partial\Gamma$ does not depend on the chosen parameter $\varepsilon>0$ provided that it  is small enough. It also implies that, for such a parameter $\varepsilon$, $$|\xi-\zeta|_\varepsilon\asymp e^{-\varepsilon (\xi|\zeta)}$$ for points on the boundary.

\gap 

We turn now to the second part of Theorem \ref{homeoimphyp} --- the  uniqueness of the quasi-isometry type of $\Gamma=\Gamma(f, \UUU)$.   

\begin{prop}\label{homeoimpqi} Assume that $f:(\XX_1,X)\to (\XX_0,X)$ is a finite branched covering of degree $d$. Let
$S_j(1)$, $j=1,2$, be finite coverings. We denote by $\Gamma_j$, $F_j$, $\varepsilon_j$ and
$\phi_j:X\to\partial \Gamma_j$  the graph, dynamics, weight and projection map
associated to $S_j(1)$. If both coverings satisfy {\bf\rm [Expansion]}
and if $\phi_1$ are $\phi_2$ are both homeomorphisms, then $\Gamma_1$ is quasi-isometric to $\Gamma_2$.
\end{prop}

By Theorem \ref{thm:bsch_funct}, it is enough to show that $\bdry \Gamma_{j}$, $j=1,2$ are quasisymmetrically equivalent.
We will actually prove the stronger statement that the boundaries are {\it snowflake equivalent} i.e., for any $x$, $y$ and $z$ in $X$,
$$\left(\frac{|x-z|_1}{|x-y|_1}\right)^{\varepsilon_2}\asymp \left(\frac{|x-z|_2}{|x-y|_2}\right)^{\varepsilon_1}\,,$$ where  $|\cdot|_j$ denotes the metric on the repellor $X$ obtained by pulling back the metric $d_{\varepsilon_j}$ on $\bdry_{\varepsilon_j}\Gamma$ via the homeomorphism $\phi_j$.
Without further combinatorial finiteness or uniformity properties, it seems difficult to work directly with the graphs $\Gamma_j$ to show that they are quasi-isometric.    

\gap

We start with some lemmata relating the combinatorics and geometry of the
 two graphs $\Gamma_j$, $j=1,2$ constructed from different choices of coverings.  To avoid cumbersome duplication, and to keep the statements symmetric, we denote by $j \mapsto j^*$ the involution of $\{1,2\}$ interchanging $1$ and $2$. We also suppress mention of the dependence of the metrics on the choices of $\varepsilon_j$.  

\gap

Axiom [Expansion] and Proposition \ref{prop:elementary}(2)(a) imply the following result.  \\

\begin{lemma}\label{l1homeoimpqi} For $j=1,2$, integers $n_j$ exist such that
 \be
\item  for any $U_j\in S_j(n_j)$, there is $U_{j^*}\in S_{j^*}(1)$ which contains $U_j$. 
\item for any $U_{j^*}\in S_{j^*}(1)$, there is  $U_j\in S_j(n_j)$ contained in $U_{j^*}$.\eb \end{lemma}

\pf The finiteness of $S_1(1)\cup S_2(1)$  implies that there is some entourage $\Omega$ of $\XX_1$ such that
 any $U\in S_1(1)\cup S_2(1)$ 
contains $\Omega(x)$ for some $x\in X$, and, for any $x\in X$, 
$\Omega(x)$ is contained in some element of $ S_1(1)$ and of $ S_2(1)$.
\gap

We treat the case $j=1$. Since [Expansion] holds, there is some $n_1$ so that  any $W\in S_1(n_1)$
is contained in $\Omega (x)$ for any $x\in W\cap X$.

\be
\item If $U_1\in S_1(n_1)$, then consider $x\in U_1$ so that $U_1\subset \Omega(x)$. There is some $U_2\in S_2(1)$ 
such that $\Omega(x)\subset U_2$.
Therefore $U_1\subset U_2$.
\item If $U_2\in S_2(1)$, let $x\in U_2$ such that $\Omega(x)\subset U_2$. 
Let $U_1 \in S_1(n_1)$ contain $x$. Thus, $$U_1\subset \Omega(x)\subset U_2.$$\qed
\eb

\begin{lemma}\label{l2homeoimpqi} A constant $K\ge 1$ exists such that, for $j=1,2$, for any $x\in X$ and any $n\ge n_j$, 
there are some $U\in S_j(n)$, $W'\in S_{j^*}(n+n_{j^*}-1)$ and $W\in S_{j^*}(n-n_{j}+1)$ such that
 \be
\item  $x\in W'\subset U\subset W$. 
\item $\roundness_{_{j}}(\phi_j(U),\phi_j(x))\le K$,
$\roundness_{_{j^*}}(\phi_{j^*}(W'),\phi_{j^*}(x))\le K$,  and \\ 
$\roundness_{_{j^*}}(\phi_{j^*}(W),\phi_{j^*}(x))\le K$. \eb 
It follows that
$$\diam_{_{j^*}} \phi_{j^*}(W') \asymp\diam_{_{j^*}} \phi_{j^*}(U)\asymp  
\diam_{_{j^*}} \phi_{j^*}(W)\,.$$
\end{lemma}

\pf We let $j=1$.  Proceeding as usual, let us rename $x=\tx$.  

Let $\tx\in X$. Using the fact that $S_1(n_1)$ is finite, there is some $U\in S_1(n_1)$ such that 
$\roundness_{_{1}}  (\phi_1(U), \phi_1(f^{n-n_1}(\tx)))\le K_1'$
for some constant $K_1'\ge 1$.

Let $\wtU$ be the component of $f^{-(n-n_1)}(U)$ which contains $\tx$. Then $\wtU\in S_1(n)$ and Propositions \ref{ba2} and \ref{ba3}(1) 
imply $\roundness_{_{1}}  (\phi_1(\wtU), \phi_1(\tx))\le K_1$
for some constant $K_1\ge 1$.

\gap

By Lemma \ref{l1homeoimpqi}, there is some $W\in S_2(1)$ which contains $U$. It follows from compactness that there exists 
a constant $K_2'$ independent of $\tx$ such that $$\roundness_{_{2}}  
(\phi_2(W), \phi_2(f^{n-n_1}(\tx)))\le K_2'; $$
see the proof of Proposition \ref{prop:uniform_roundness}(1).

Let $\wtW$ be the component of $f^{-(n-n_1)}(W)$ which contains $\tx$. Then $\wtW \in S_1(n-n_1+1)$ and Proposition \ref{ba2}
implies $\roundness_{_{2}}  (\phi_2(\wtW), \phi_2(\tx))\le K_2$
for some constant $K_2\ge 1$.

\gap

By Lemma \ref{l1homeoimpqi}, the point $f^{n-1}(\tx)$ belongs to some $W'\in S_2(n_2)$ contained in $f^{n-1}(\wtU)$. 
Since $S_2(n_2)$ is finite, one can assume that 
$$\roundness_{_{2}}  (\phi_2(W'), \phi_2(f^{n-1}(\tx)))\le K_3'$$
for some constant $K_3'\ge 1$. 

 Let $\wtW'$ be the component of $f^{-(n-1)}(W')$ which contains $\tx$. Then $\wtW'\in S_1(n+n_2-1)$ and Proposition \ref{ba2}
implies $\roundness_{_{2}}  (\phi_2(\wtW'), \phi_2(\tx))\le K_3$
for some constant $K_3\ge 1$.

\gap

Let $K=\max\{K_1,K_2,K_3\}$. The lemma follows from Proposition \ref{ba3} once we have noticed that
$|W-W'|=n_1+n_2$. \qed

\gap

We now give the proof of Proposition \ref{homeoimpqi}.

\gap

\pf (Proposition \ref{homeoimpqi})  

Let $\mathbf{U}_j$, $j=1,2$ denote the corresponding collections of open sets defined by two different coverings at level zero.   For $j=1,2$ let $|\cdot|_j$ denote the metric on the repellor $X$ of $f: \XX_1 \to \XX_0$ obtained by pulling back the metric $d_{\varepsilon_j}$ on $\bdry_{\varepsilon_j}\Gamma$ via the homeomorphism $\phi_j$.   Roundness and diameters in these metrics will be denoted with subscripts.  We will show that the identity map is quasisymmetric: we want to find a homeomorphism $\eta:\R_+\to\R_+$ such that, given any  $x,y,z\in X$,
$$\frac{|x-z|_2}{|x-y|_2}\le\eta\left(\frac{|x-z|_1}{|x-y|_1}\right).$$  

\gap

Let $\Omega$ be an entourage such that, for any $x\in X$ any  $j=1,2$, there is some $U_j\in S_j(1)$, such that
$\Omega(x)\subset U_j$.

By the uniform continuity of $\phi_1$, $\phi_2$ and their inverses, it is enough to consider 
$x,y,z\in X$ such that $y,z\in \Omega(x)$. 

The strategy is the following. Let us assume that $z$ is closer to $x$ than $y$. Then, we may find neighborhoods $U_y, U_x \in \mathbf{U}$ of $x$ and $y$ respectively such that the ``ring'' $U_y\setminus U_x$ separates the set $\{x, z\}$ from $y$.   The 3-point condition will follow from a straightforward argument using what is known about the sizes of the neighborhoods in each of the two metrics.

By Proposition \ref{ba3}, there exists a neighborhood $U_y \in \mathbf{U}_1$ of $x$ not containing $y$ such that $|x-y|_1 \asymp \diam_1(U_y)$, and $\roundness_1(U_y, x) \leq K$, where $K$ is a uniform constant.  

Again by Proposition \ref{ba3}, there exists a neighborhood $U_z \in \mathbf{U}_1$ of $x$  containing $z$ such that $|x-z|_1 \asymp \diam_1(U_z)$ and $\roundness_1(U_z, z) \leq K$.  

\gap

Therefore, Lemma  \ref{l2homeoimpqi} implies the existence of $W_y'\in S_2(|U_y|_1+n_2-1)$ which
contains $x$ but is contained in $U_y$ such that $\roundness_{2} (W_y',x)\le K$ and
$\diam_{2} (W_y') \asymp\diam_{2} (U_y)$.

Similarly, a vertex $W_z\in S_2(|U_z|_1+n_1+1)$ which
contains $U_z$ exists such that $\roundness_{2} (W_z,x)\le K$ and
$\diam_{2} (W_z) \asymp\diam_{2}( U_z)$.

\gap

Since $\roundness_{2} (W_y',x)\le K$, it follows that
$$\frac{|x-z|_{2}}{|x-y|_{2}}
\lesssim \frac{\diam_{2} (W_z)}{\diam_{2} (W_y')}
\asymp e^{- {\varepsilon_2}(|W_z|_2-|W_y'|_2)}$$

But since $x\in W_z\cap W_y'$, $$|W_z|_2-|W_y'|_2= (|U_z|_1-|U_y|_1) +(n_2-n_1)$$
one obtains 
$$\frac{|x-z|_{2}}{|x-y|_{2}}
\lesssim \left(\frac{\diam_{1} (U_z)}{\diam_{1} (U_y)}\right)^{\varepsilon_2/\varepsilon_1}
\lesssim \left(\frac{|x-z|_{1}}{|x-y|_{1}}\right)^{\varepsilon_2/\varepsilon_1}
$$
and so the identity map is a quasisymmetry.   
 \qed

This concludes the proof of the second conclusion of Theorem \ref{homeoimphyp}.  

\gap

The last conclusion of Theorem \ref{homeoimphyp} is easily proved along the following lines.
There is a canonical inclusion $\iota$ from the vertices of $\Gamma(f^n,\UUU)$ to those of $\Gamma(f,\UUU)$ which sends a vertex of $\Gamma(f^n, \UUU)$, say $V \in \UUU_{nk}$ with $|V|=k$, to the vertex in $\Gamma(f,\UUU)$ called again $V \in \UUU_{nk}$ with now $|V|=nk$.   The image of $\iota$ is clearly $n$-cobounded, an isometry on horizontal paths, and multiplies the lengths of vertical paths by a factor of $n$.  Hence $\iota$ is $n$-Lipschitz.  Using these facts one proves easily that $\iota$ is in fact a quasi-isometry, and the proof of  Theorem \ref{homeoimphyp} is complete.  \qed

\section{Measure theory} 
\label{secn:measure_theory} 

In this section, we assume that $f: \XX_1 \to \XX_0$ is a degree $d$ fbc with repellor $X$ as in \S\,\ref{secn:spaces_associated}.  As in the previous section, we assume that we are  given a covering $\UUU$ of $X$ by connected open subsets of $\XX_1$ 
which satisfies [Expansion].   
Let $\Gamma=\Gamma(f,\UUU)$ be the Gromov hyperbolic graph associated to $f$ and $\UUU$ as in the previous section.  
Fix $\varepsilon>0$ small enough so that $\phi_f: X \to \bdry_\varepsilon\Gamma$ is a homeomorphism.  

We now assume that axiom [Irreducibility] holds as well.

The main result of this section is the following theorem.

\begin{thm}\label{maxent} Assume that Axioms [Expansion] and [Irreducibility] hold.  
Then there is a unique invariant \qc measure $\mu_f$; its dimension is $(1/\varepsilon)\log d$. 
This measure is also mixing and ergodic, and it  describes the distribution of preimages of  points and of periodic points.
Furthermore, the metric entropy and topological entropy satisfy the following bounds
$$0<\log d -\int \log d_F\; d\mu_f\le h_\mu(F)\le h_{top}(F) \le v\le \log d$$
and  $$\frac{h_\mu(F)}{\varepsilon}\le  \dim \mu_f\le \dim\partial_{\varepsilon}\Gamma \le \frac{v}{\varepsilon}\le \frac{\log d}{\varepsilon},$$
where $$v=\lim\frac{1}{n}\log|S(n)|.$$
 \end{thm}\index{Symbols}{$v$}

Precise statements and definitions are given in the next few subsections.

In the top chain of inequalities, 
the first one is a consequence of Rohlin's inequality, which always applies in our setting.  
The second follows from the Variational Principle, the third from generalities since $F$ is Lipschitz, and the last since $F$ is a degree $d$ fbc.

When $f$ is topologically cxc, then we will prove that both chains of inequalities in Theorem
\ref{maxent} are equalities, that $\mu_f$ is the unique measure of maximal entropy $\log d$, 
and that $(\partial_{\varepsilon}\Gamma,d_{\varepsilon},\mu_f)$ is an Ahlfors regular metric
measure space of dimension $(1/\varepsilon)\log d$, see Theorem \ref{thm:cor_cxc}.

In the remainder of this section, we dispense completely with the topological spaces $X, \XX_0, \XX_1$ and deal exclusively with $F: \cl{\Gamma}_\varepsilon \to \cl{\Gamma}_\varepsilon$.  

\subsection{Quasiconformal measures}
\label{subsecn:qc_measures}

Recall that $f$ induces a continuous surjective Lipschitz map $F: \cl{\Gamma}_\varepsilon \to \cl{\Gamma}_\varepsilon$ which maps vertices to vertices and edges (outside $\cl{B_\varepsilon(o,2)}$) homeomorphically onto edges.  \\
\gap 

{\noindent\bf Multiplicity function for $F$.} 
\index{Index}{multiplicity function}\index{Symbols}{$d_F$}%
Let $d_f(x)$ denote the local degree of $f$ at a point $x\in X$.
\bi
\item[--] If $\xi\in \partial_\varepsilon \Gamma$, let $d_F(\xi)=d_f(\phi_f^{-1}(\xi))$.
\item[--] If $W\in V$, $|W|\ge 2$, let $d_F(W)= \deg(f|_W)$.
\item[--] For each (open) edge $e=(W,W')$ with $|W|,|W'|\ge 1$, choose a point $x_e\in W\cap W'\cap X$.

If $e\subset \cl{\Gamma}_\varepsilon\setminus B(o,2)$,
set, for all $\xi\in e$, $$d_F(\xi)= \sum_{y\in (W\cap W'\cap  f^{-1}(x_{F(e)}))}d_f(y)\,.$$\ib

\noindent{\bf Remarks:}
\be
\item The definition depends on the choices of points $x_e$, but this is irrelevant for our purposes.

\item The function $d_F$ may vanish on certain edges.  For example, let $X=\XX_1=\XX_0 = \R/\Z$, let $f(x)=2x$ modulo $1$, and let $\UUU_0 = \{U,V\}$ where $U=X-\{1/4+\Z\}$ and $v=X-\{3/4+\Z\}$.   Note that $0+\Z \subset U \intersect V$ but that $U \intersect V$ is not connected.  The set $S(1)$ consists of the two vertices $U,V$ joined by a single edge $e$.  Choose $x_e=0+\Z$.  The four elements of $S(2)$ are the two preimages of $U$ given by the intervals (mod $\Z$) $(-3/8, 1/8)$ and $(1/8, 5/8)$ and the two preimages of $V$ are $(-1/8, 3/8)$ and $(3/8, 7/8)$.  According to the definition, the edge joining $(1/8, 5/8)$ and $(-1/8, 3/8)$ is given weight zero by $d_F$ since the intersection of these two intervals contains neither $0$ nor $1/2$, the preimages of the origin.   

\item If $d_F(\xi)\geq 2$ at a point $\xi$ in the interior of an edge $e$ (such as when the chosen point $x_e \in \bdry_\varepsilon\Gamma$ is a branch point of $F$ on the boundary), then $F$ is never a branched covering with degree function $d_F$, since $d_F$ is constant on interiors of edges, and an honest f.b.c. is unramified on a dense open set.  Conversely, if $d_F\equiv 1$ on $\bdry\Gamma_{\varepsilon}$, then  $F$ is an fbc in a neighborhood of $\partial_\varepsilon \Gamma$.

\eb

The following properties hold.

\begin{lemma} The multiplicity function behaves as a local degree function. More precisely,
\bi
\item[\rm (i)] for any $\xi\in\overline{\Gamma}_\varepsilon\setminus B(o,1)$, $$\sum_{F(\zeta)=\xi}d_F(\zeta)=d\,;$$
\item[\rm (ii)] for  any $\xi\in\overline{\Gamma}_\varepsilon\setminus B(o,2)$, there is a neighborhood $N$ such
that, for any $\zeta\in N$, $$d_F(\xi)=\sum_{\zeta'\in F^{-1}(\{F(\zeta)\})\cap N} d_F(\zeta').$$\ib
\end{lemma}

\pf \bi
\item[(i)] The statement is clear for vertices and points from the boundary. Let $e=(W,W')$ be an edge, and let us denote by
$\wtW_1,\ldots, \wtW_k$ the components of $f^{-1}(W)$, and by $\wtW_1',\ldots, \wtW_{k'}'$ the components of $f^{-1}(W')$.

If $f(y)=x_e$, then  there exists a unique edge  $\tilde{e}=(\wtW_y, \wtW_{y'})$ such that $y\in \wtW_y \intersect \wtW_y'$. Therefore
$$\sum_{F(\tilde{e})=e} d_F(\tilde{e})=\sum_{F(\tilde{e})=e} \sum_{y\in (\wtW_y \cap\wtW_y')\cap f^{-1}(x_e)}d_f(y)
=\sum_{f(y)=x_e} d_f(y)=d\,.$$
\item[(ii)] The statement is clear on $\Gamma\setminus B(o,2)$. Let $\xi\in\partial\Gamma$. There is some vertex $W_0$ such that 
$\phi_f(W_0)\ni \xi$, and $d_F(W_0)=d_F(\xi)$. Let $W_1\subset W_0$ small enough so that $\mho_{\infty}(F(W_1))\subset \phi_f(F(W_0))$.
Thus, for any $U\in\interior(\mho(F(W_1)))$, $U\subset F(W_0)$, so that 
$$\sum_{F(\wtU)=U,\wtU\subset W_0} d_F(\wtU)= d_f(W_0)=d_F(\xi)\,.$$
\ib\qed

Note that if we set  $d_{F^n}(\xi)=d_F(\xi)\ldots d_{F}(F^{n-1}(\xi))$, then the lemma remains true for $d_{F^n}$ as well. 

\gap

{\noindent\bf Action of $F$ on measures.} If $\varphi$ is a continuous test function defined on $\overline{\Gamma}_\varepsilon \setminus B_\epsilon(o,1)$, 
then  its {\em pullback}
\index{Index}{pullback, continuous functions}\index{Symbols}{$F^*\varphi$}\index{Symbols}{$F_*\nu$}%
  under $F$, given by the formula $F^*\varphi(\xi)=\varphi\circ F(\xi)$, defines a continuous function 
on $\overline{\Gamma}_\varepsilon\setminus B(o,2)$.
By duality, one may define for Borel probability measures $\nu$ with support in $\overline{\Gamma}_\varepsilon\setminus B(o,2)$ its {\em pushforward}
\index{Index}{pushforward, measure}%
 by
$\la F_*\nu,\varphi\ra= \la \nu,F^*\varphi\ra$.  Thus in particular, $(F_*\nu)(E)=\nu(F^{-1}(E))$ for all Borel sets $E$.  

The point of the construction of the multiplicity function $d_F$ is the following.  
If $\varphi$ is a continuous test\index{Index}{pushforward, continuous functions}\index{Symbols}{$F_*\varphi$}\index{Symbols}{$F^*\nu$}%
  function on $\cl{\Gamma}_\varepsilon \setminus B_\varepsilon(o,1)$, its {\em pushforward} under $F$ 
$$F_*\varphi(\xi)=\sum_{F(\zeta)=\xi} d_F(\zeta)\varphi(\zeta)$$ 
is again a continuous function on 
$\overline{\Gamma}\setminus B(o,1)$. 
By duality, we define the {\em pullback}\index{Index}{pullback, measure} of a Borel measure $\nu$ by the formula
$\la F^*\nu,\varphi\ra= \la \nu,F_*\varphi\ra$ (cf. \cite[\S\,2]{dinh:sibony:dynamique}).

\gap

{\noindent\bf Quasiconformal measures.} If $\mu, \nu$ are measures we write $\nu \ll \mu$ if $\nu$ is absolutely continuous with respect to $\mu$.  Let $\mu$ be a regular Borel probability measure on $\partial_\varepsilon \Gamma$. Inspired by the group setting \cite{coornaert:patterson-sullivan},
we say $\mu$ is a {\em  \qc measure of dimension $\alpha$}
\index{Index}{measure, quasiconformal}%
 if, for all $n\ge 1$, $(F^n)^*\mu \ll \mu$ and the Radon-Nikodym derivative
satisfies $$\frac{d (F^n)^*\mu}{d\mu}\asymp (e^{n\varepsilon})^\alpha \quad \mu-a.e.\,.$$

\gap

The quantity $e^{n\varepsilon}$ stands for the derivative of $F^n$ (cf. Proposition \ref{boulesF}). 

\gap

Let  $\mu$ be a quasiconformal measure on $\bdry_\varepsilon\Gamma$.  Fix $n \in \N$.  
Suppose $E \subset \bdry_\varepsilon\Gamma$ is a Borel subset of positive measure, $F^n|_E$ is injective, and  
the local degree of $F^n$ is constant on $E$, i.e.  for all $\xi \in E$, $d_{F^n}(\xi)=d_E$.  Then, it follows from the regularity of
the measure that 
\begin{equation}\label{eqn:F^*}
\la (F^n)^*\mu,\chi_E\ra =\la \mu, (F^n)_*\chi_E\ra =\int\sum_{F^n(\zeta)=\xi}d_F(\zeta)\chi_E(\zeta)d\mu(\xi)=d_E\mu (F^n(E))\,,
\end{equation}
and the quasiconformality of the measure implies $\la (F^n)^*\mu,\chi_E\ra \asymp e^{n\alpha\varepsilon}\mu(E)$.
Hence
\begin{equation}
\label{eqn:how_F_changes_measures}
\mu(F^n(E))\asymp \frac{e^{n\alpha\varepsilon}}{d_E}\mu(E).
\end{equation}

Axiom [Irreducibility] implies that the support of a quasiconformal measure is the whole set $\partial_\varepsilon \Gamma$.
Therefore, there is some $m>0$ such that, for all $x\in S(1)$, $\mu(\phi_f(W(x)))\ge m$. 

\gap

We let $d(W)$ be the degree of $f^{n-1}|_{W}$ for $W\in S(n)$. Since $\mu$ is a \qc measure, it follows
that $$0< m\le \mu(\phi_f(f^{n-1}(W))) =  \mu(F^{n-1} \phi_f(W)) \asymp \frac{e^{n\alpha \varepsilon} }{d(W)}\mu(\phi_f(W))\,.$$ 
This proves 

\begin{lemma}\label{shagen}{\em\bf of the shadow.} For any $W\in V$, $$\mu(\phi_f(W))\asymp d(W) e^{-\alpha\varepsilon|W|}\,.$$
\end{lemma}
\index{Index}{shadow, lemma of}

We use this lemma for the classification of quasiconformal measures.

\begin{thm}\label{trivmeas} Let $\mu$ be a \qc measure of dimension $\alpha$. The following are equivalent.
\begin{itemize}
\item[(i)] $\mu$ is atomic.
\item[(ii)] $\alpha =0$.
\item[(iii)] $\partial_\varepsilon\Gamma$ is a point.\end{itemize}
If $\alpha >0$ then $\alpha=\frac{1}{\varepsilon}\log d$, and any two such \qc measures are equivalent.\end{thm}

\pf  Consider first the constant function $\varphi=1$ on $\partial_\varepsilon \Gamma$. 
Then $F_*\varphi= d\varphi$ so that
$$\la (F^n)^*\mu,\varphi \ra =\la \mu, F^n_*\varphi\ra = \la \mu,d^n\varphi\ra = d^n \asymp e^{n\alpha\varepsilon}.$$Thus, $$\alpha =\frac{1}{\varepsilon}\log d.$$

It follows that $\alpha=0$ if and only if $d=1$, so that $\partial_\varepsilon \Gamma$ is a point since $f$
satisfies  [Expansion]. Hence (ii) implies (iii).

\gap

If $\mu$ is atomic, then there is some $\xi\in\partial_\varepsilon\Gamma$ such that $\mu\{\xi\}>0$. By definition, for all $n\ge 0$ and any $x \in \bdry_\varepsilon\Gamma$, 
$$((F^n)^*\mu)(\{x\})= d_{F^n}(x)\mu\{F^n(x)\}\,.$$ Since $\mu$ is \qc then by Equation (\ref{eqn:how_F_changes_measures})
$$d_{F^n}(\xi)\mu\{F^n(\xi)\}\asymp e^{n\alpha\varepsilon}\mu\{\xi\}=d^n\mu\{\xi\}\,.$$

\gap

But $d_{F^n}(\xi)\le d^n$ so $\mu\{F^n(\xi)\}\gtrsim \mu\{\xi\}$. Since the total mass of $\mu$ is finite, the orbit of $\xi$ has to be finite.
Let $\zeta=F^{\ell}(\xi)$ be periodic and let  $k$ be its period. Then $\mu\{\zeta\}>0$ and
$$(d_{F^k}(\zeta))^n\mu\{\zeta\}=((F^{kn})^*\mu)(\{\zeta\})\asymp  d^{nk}\mu\{\zeta\}$$ from which we deduce that
$d_{F^k}(\zeta)=d^k$. This means that the local degree at every point in its orbit is maximal, so that
its grand orbit is finite. Since $f$
satisfies [Irreducibility], the grand orbit of any point is dense in $X$ 
(Proposition \ref{prop:elementary}3(a)) and so $\partial_\varepsilon \Gamma$ is a point, $d=1$ and $\alpha=0$. So (i) implies (ii) and (iii).

\gap

The last implication (iii) implies (i) is obvious.

\gap

The Lemma of the Shadow (Lemma \ref{shagen}) and the assumption that quasiconformal measures are regular imply that two measures
of the same dimension are equivalent.\qed

\gap

We will now construct a \qc measure using the Patterson-Sullivan procedure \cite{coornaert:patterson-sullivan}. It turns out that this measure will
be invariant.

\gap

{\noindent\bf Poincar\'e series.} 
\index{Index}{Poincar\'e series}%
Let $$P(s)=|S(1)|\sum_{n\ge 1} d^{n-1}e^{-ns}=|S(1)|\frac{1}{e^s-d}\,.$$
It follows that $P(s)<\infty$ if and only if $s>\log d$. Let, for $s>\log d$, 
$$\mu_s=\frac{1}{P(s)}\sum_{n\ge 1} \sum_{\xi\in S(n)}e^{-ns}d(\xi)\delta_{\xi}\,.$$
For every $n\ge 1$, $F^n:S(n+1)\to S(1)$ has degree $d^n$.  Recall that for $\xi \in S(n)$, we denoted by $d(\xi)= d_{F^{n-1}}(\xi)$.    
So $$|S(n+1)|=d^n|S(1)|-\sum_{\xi\in S(n+1)} [d(\xi)-1]$$
and $\sum_{\xi\in S(n+1)} d(\xi)=d^n |S(1)|$. Therefore
$$\mu_s(\overline{\Gamma}_\varepsilon)=\frac{1}{P(s)}\sum_{n\ge 1} e^{-ns}\sum_{\xi\in S(n)}d(\xi)=\frac{1}{P(s)}\sum_{n\ge 1} e^{-ns}d^{n-1}|S(1)|=1\,.$$

\gap

Hence $\{\mu_s\}_{s >\log d}$ is a family of  probability measures on $\cl{\Gamma}_\varepsilon$. Let\index{Symbols}{$\mu_f$} $\mu_f$ be any weak limit of this family as $s$ decreases to $\log d$.
Since the Poincar\'e series diverges at $\log d$, it follows that the support of $\mu_f$ is contained in $\bdry \Gamma$.

\gap

If $\varphi$ is a continuous function with support close to $\partial \Gamma$, then
\[
\renewcommand{\arraystretch}{2}
\begin{array}{ll}
\la F^*\mu_s, \varphi\ra  & =\dis\frac{1}{P(s)}\dis\sum_{n\ge 1} e^{-ns}\dis\sum_{\xi\in S(n)}d(\xi)(F_*\varphi)(\xi)\\
& =\dis\frac{1}{P(s)}\dis\sum_{n\ge 1} e^{-ns}\dis\sum_{\xi\in S(n)}d(\xi)\sum_{F(\zeta)=\xi}d_F(\zeta) \varphi(\zeta)\\
& = \dis\frac{1}{P(s)}\dis\sum_{n\ge 1} e^{-ns}\dis\sum_{\zeta\in S(n+1)}d(\zeta)\varphi(\zeta)\\
&= e^s \la \mu_s, \varphi\ra + O(1/P(s)) \end{array}
\] where we have used that $d(\zeta)= d(F(\zeta))d_F(\zeta)$.

It follows that, as $s$ decreases to $\log d$, 
$$\la F^*\mu_f, \varphi\ra  =\la d\mu_f, \varphi\ra $$ and so  $F^*\mu_f= d\mu_f$. 
In other words, $\mu_f$ is a \qc measure of dimension $(1/\varepsilon)\log d$.

\gap

Let us look at $F_*\mu_f$:

$$\la F_*\mu_f, \varphi\ra = \la \mu_f, F^*\varphi\ra = (1/d)\la  F^*\mu_f, F^*\varphi\ra=(1/d)\la  \mu_f, F_*(F^*\varphi)\ra.$$
But $$F_*(F^*\varphi)(\xi)=  \left(\sum_{F(\zeta)=\xi} d_F(\zeta)(F^*\varphi)(\zeta)\right)= d\varphi(\xi)\,.$$ 

Therefore, $F_*\mu_f=\mu_f$, so $\mu_f$ is an invariant measure.

\gap

\noindent{\bf Ergodicity.} Let $(Z,\nu)$ be a probability space, and  $T:Z\to Z$ a transformation preserving the measure $\nu$,
i.e., $\nu(T^{-1}(A))=\nu(A)$ for every measurable subset $A$ of $Z$. The measure $\nu$ is {\it ergodic} 
\index{Index}{ergodic}%
if any invariant
measurable set $A$ has zero or full measure.  The measure $\nu$ is {\em mixing}
\index{Index}{mixing}%
 if for any two measurable subsets $A, B$, one has $\nu(T^{-n}(A)\intersect B) \to \nu(A)\nu(B)$ as $n \to \infty$; mixing implies ergodicity.  

A fundamental theorem of ergodic theory is the following {\em Birkhoff ergodic
theorem}
\index{Index}{Birkhoff ergodic theorem}%
 (see e.g. \cite[Thm 4.1.2 and Cor. 4.1.9]{hasselblatt:katok:dynamics}).

\begin{thm}[Birkhoff ergodic theorem]\label{thm:birkhoff} If $T:Z\to Z$ is a $\nu$-preserving ergodic transformation, then,
for any function $\varphi\in L^1(Z,\nu)$ and for $\nu$-almost every $z\in Z$, 
$$\frac{1}{n}\sum_{k=0}^{n-1} \varphi(T^k(z))=\int_Z \varphi d\nu\,.$$\end{thm}

Let us prove that $\mu_f$ is ergodic. Let $E$ be an invariant subset of $\partial_\varepsilon \Gamma$ with positive measure. 
Let $\nu=\mu_f|_E/\mu_f(E)$.
It follows that $\nu$ is also an invariant \qc measure. The Lemma of the shadow (Lemma \ref{shagen}) implies that $\mu_f(W)\asymp\nu(W)$
for all $W\in V$.  This implies that $\mu_f$ and $\nu$ are equivalent.  Hence $\mu_f(E)=1$.  
\gap

Since $\mu_f$ is an ergodic invariant measure, it follows that $\mu_s$ converges to $\mu_f$ in the weak-$*$ topology when $s$ decreases to
$\log d$.

\gap

\noindent{\bf Remark.}  On $\bdry_\varepsilon\Gamma$, the local degree function $d_F$ is multiplicative:  $d_{F^n}(\xi) = \prod_{i=0}^{n-1}d_F(f^i(\xi))$.  
From the Birkhoff ergodic theorem (Thm \ref{thm:birkhoff}) and the ergodicity of $F$ with respect to $\mu_f$, it follows that for $\mu_f$-almost any $\xi\in\partial_\varepsilon\Gamma$, 
$$\lim_{n\to\infty} \frac{1}{n} \log d_{F^n}(\xi)=\int\log d_F \;d\mu_f\,.$$
Thus, either the critical set has measure $0$ and the Jacobian of 
$F$ with respect to $\mu_f$ is constant and equal to $d$ almost everywhere, or almost every point visits the branch set so frequently that the  local degrees increase exponentially fast. 
Unfortunately, given the assumptions under which we are currently working, we have neither a proof that this latter possibility cannot occur, nor an example showing that it can occur.

\subsection{Entropy} \label{subsecn:entropy} 

We refer to \cite[Chap.\,3, \S\,3.1, Chap.\,4, \S\,3]{hasselblatt:katok:dynamics}, 
\cite{mane:hdim} and \cite[Chap.\,1]{przytycki:urbanski:book},  for background on entropy.

\gap

\noindent{\bf Topological entropy.}  Let $T: Z \to Z$ be a continuous map of a compact metric space $(Z,d)$ to itself.  The {\em dynamical distance} and the corresponding {\em dynamical balls at level $n$} are defined as\index{Symbols}{$S(\xi,n,r)$} 
 $$d_n(\xi,\zeta)=\max_{0\le j\le
n}\{d(T^j(\xi),T^j(\zeta))\} \quad\hbox{and}\quad S(\xi,n,r)=\{\zeta\in Z,\ d_n(\xi,\zeta)\le r\}\,.$$
Let $c_n(r)$ be the minimal number of dynamical balls $S(\cdot,n,r)$ at level $n$ needed to cover $Z$ and $s_n(r)$ the 
maximal number of disjoint  dynamical balls $S(\cdot,n,r)$.
The {\em topological entropy of $T$}
\index{Index}{entropy, topological}%
 may be defined  as  
$$h_{top}(T)= \lim_{r\to 0}\limsup_{n\to\infty}\frac{1}{n}\log c_n(r)=\lim_{r\to 0}\liminf_{n\to\infty}\frac{1}{n}\log s_n(r)\,.$$
\index{Symbols}{$h_{top}(T)$}

\gap

We now estimate $h_{top}(F)$, where $F$ denotes the restriction of $F: \cl{\Gamma}_\varepsilon \to \cl{\Gamma}_\varepsilon$ to the boundary $\partial_\epsilon \Gamma$. 
\gap 

Since $F$ is $e^\varepsilon$-Lipschitz, we have $d_n(\xi,\zeta)\le e^{n\varepsilon}d_\varepsilon(\xi,\zeta)$ and  hence 
$S(\xi,n,r)\supset B_\varepsilon(\xi,re^{-\varepsilon n})$.
For any $n\ge 1$, $\{\mho_\infty(\xi)\}_{\xi\in S(n)}$ is a covering of $\partial  
\Gamma$ by at
most $|S(n)|$ sets.  For any
$\xi\in S(n)$, $\diam\mho_\infty(\xi)\le Ce^{-\varepsilon |\xi|}$.  So,   
$\mho_{\infty}(\xi)\subset
S(\xi',p,Ce^{-\varepsilon(n-p)})$, for any $p\in\N$ and for any $\xi'\in\mho_{\infty}(\xi)$.
Hence $c_p(Ce^{-\varepsilon(n-p)}) \leq |S(n)|$.  

\gap

Recall that by definition $v=\lim\frac{1}{n}\log|S(n)|$; the limit exists since $|S(n+1)|\leq d|S(n)|$.  
Let $\eta >0$ be small.  For any $p\ge 1$, there is some $n\in\N$ such  
that
$\eta\asymp e^{-\varepsilon(n-p)}$, meaning that
$n$ is equal to $p +(1/\varepsilon) \log 1/\eta$ up to
a universal additive constant.
The discussion above implies that $\partial_\varepsilon \Gamma\subset \cup_{\xi\in
S(n)}S(\xi,p,Ce^{-\varepsilon(n-p)})$.  Since for any $\eta'>0$ and any $n$ large enough, 
$\log|S(n)|\le n(v+\eta')$ holds, we have $\log c_p(\eta)\le n(v+\eta')$ and
$$h_{top}(F)\le  \lim_{\eta\to 0}\limsup_{p\to\infty} \frac{p+(1/\varepsilon) \log  
1/\eta}{p}(v+\eta')$$ from which $h_{top}(F)\le v$ follows. Since $|S(n)|\le d^n$, one has $v\le \log d$.

\gap

\noindent{\bf Measure-theoretic entropy.}  
We recall first the definition of measure-theoretic entropy, more commonly referred to as {\em metric entropy}.
\index{Index}{metric entropy}%
\index{Index}{entropy, metric}%
Suppose  $(Z,\nu)$ is a probability space, and $T: Z \to Z$ preserves $\nu$.  If $\PPP$ is a partition of $Z$ 
into a countable collection of measurable sets, 
define its {\em entropy with respect to $\nu$} by \index{Symbols}{$H_{\nu}(\PPP)$}
$$H_{\nu}(\PPP)= \sum_{A \in \PPP}\nu(A)\log(1/\nu(A))\,.$$

If $\AAA$ and $\BBB$ are two measurable partitions, we define\index{Symbols}{$\AAA\vee\BBB$} $\AAA\vee\BBB$ as the partition given
by $\{A\cap B,\ A\in\AAA,\ B\in\BBB\}$. Furthermore, we say that $\AAA$ is {\it finer} than $\BBB$ if,
for any $A\in\AAA$, there is some $B\in\BBB$ such that $A\subset B$. 

For $n \in \N$ set 
\[ \PPP_n = P \vee T^{-1}(\PPP) \ldots \vee T^{-n}(\PPP).\]
Then \index{Symbols}{$h_{\nu}(T,\PPP)$}
$$h_{\nu}(T,\PPP)=\lim_{n \to \infty} \frac{1}{n}H_{\nu}(\PPP_n)$$
exists.  The supremum of $h_{\nu}(T,\PPP)$ over all partitions with finite
entropy defines
the metric entropy\index{Symbols}{$h_{\nu}(T)$} $h_{\nu}(T)$. If $(\PPP^n)$ is an increasing sequence of measurable partitions tending towards the
partition into points, then $\lim_{n\to\infty} h_{\nu}(T,\PPP^n)=h_{\nu}(T)$. 

A partition  $\PPP$ with finite entropy
is called a {\em generator}
\index{Index}{generator}%
 if it separates points i.e., for any distinct $z,z'\in Z$, there exist  some $n\ge 0$, and disjoint sets $A,A'\in\PPP_n$ such that $z\in A$ and $z'\in A'$; equivalently, there exist some $n\ge 0$, disjoint pieces $A,A'\in \PPP$ such that $T^n(z)\in A$ and $T^n(z')\in A'$.
For a generating partition $\PPP$, the entropy $h_{\nu}(T,\PPP)$ achieves the maximum of $h_{\nu}(T,\PPP')$ over all measurable partitions with finite entropy, 
so that $h_{\nu}(T,\PPP)=h_\nu(T)$ holds.  

The {\em  variational principle}
\index{Index}{variational principle}%
 (see \cite{walters:gtm}, Thm. 8.6) asserts that, when $T$ is continuous, then  \begin{equation}
\label{eqn:variational_principle}
h_{top}(T)=\sup_\mu h_\mu(T)
\end{equation}
 where $\mu$ varies over all invariant ergodic Borel measures.

\gap 
{\noindent\bf Jacobian.} Let $T:Z\to Z$ be  a continuous, countable-to-one  map and
$\nu$ an invariant regular Borel probability measure on $Z$. A {\em special
set}
\index{Index}{special set}%
 $A$ is a measurable subset of $Z$ such that $T|_A$ is injective.
 A {\em weak
Jacobian}
\index{Symbols}{$J_\nu$}\index{Index}{Jacobian, weak}%
 is a measurable function $J_\nu:Z\to\R_+$ such that there is some
set $Y$ such that $\nu(Y)=0$ and, for any special set $A$ disjoint from
$Y$, the equation
\begin{equation}
\label{eqn:jacobian_def}
\nu(T(A))=\int_A J_\nu d\nu
\end{equation}
 holds. The function $J_\nu$ is a {\em (strong)
Jacobian}
\index{Index}{Jacobian, strong}%
 if one can choose $Y=\emptyset$. 

Let us note that if $J_\nu$ is a weak Jacobian and $Y$ is the forbidden set as above, then the set  $Z'=Z\setminus \cup_{n\ge 0}T^{-n}(Y)$
has full measure, is forward invariant ($T(Z')\subset Z'$) and now, $J_\nu$ is a strong Jacobian for $T|_{Z'}$; see \cite[Lemma 1.9.3, Proposition
1.9.5]{przytycki:urbanski:book} for details. 

Weak Jacobians always exist
for finite branched coverings between compact spaces, and they are well-defined mod 0 sets. Let us sketch their construction in this case.

We start with a deep result of Rohlin for countable-to-one maps \cite{rohlin:measure_theory}, which turns out to be much easier in our setting.

\begin{prop}\label{prop:decomposition} Let $f:\XXX_1\to\XXX_0$ be an fbc of degree $d\ge 2$, with repellor $X$, which
satisfies [Expansion]. There exists a finite measurable partition $\PPP$ of $X$ into special sets.
\end{prop}

\pf Let $d_\varepsilon$ be the visual metric given by Theorem \ref{thm:construction} transported to $X$ via the conjugacy $\phi_f$.  In what follows, balls will be with respect to this metric.

Let $\DDD=\{d_1 < d_2<\ldots < d_N\}$ denote the set of integers which appear as local degrees of $f$ at points in $X$. For $p\in\DDD$, let  $X_p=\{x\in X,\ \deg(f;x)=p\}$.  Then $X$ is partitioned into finitely many sets $X_p$, $p \in \DDD$.  The semicontinuity properties of the local degree function imply that each set $X_p$ is measurable, so it suffices to show that for each $p \in \DDD$, the restriction $f|_{X_p}: X_p \to f(X_p)$ admits a partition into special sets.  

We will exploit the fact that $f$ maps balls onto balls as follows.  Let $x \in X$, and suppose $f^{-1}(x) = \{\tx_1, \ldots, \tx_l\}$ with $\deg(f, \tx_j) = k_j$.  By Lemma \ref{lemma:edmonds}, for each $x \in X$, there is an open connected neighborhood $U$ of $x$ in $\XX_1$ such that $f^{-1}(U)$ is a disjoint union of open connected sets $\wtU_j$, $j=1, \ldots, l$, such that $\deg(f: \wtU_j \to U) = k_j$.  
Let $r_x>0$ be so small that $B(x, r_x) \subset U$.  Then Proposition  \ref{boulesF} shows that  $f^{-1}(B_\varepsilon(x,r_x))$ is the disjoint union of the balls $B_\varepsilon(\tx_j,r_xe ^{-\varepsilon})$,
$j=1, \ldots, l$.    Furthermore, for any $y\in  B_\varepsilon(x,r_x )$, and any fixed $j=1, \ldots, l$, 
$$\sum_{y'\in f^{-1}(\{y\})\cap B_\varepsilon(\tx_j, r_xe ^{-\varepsilon})} \deg(f;y')=\deg(f;\tx_j)\,.$$
It follows that $B_{\varepsilon}(\tx_j,r_x e^{-\varepsilon})\cap X_p$ is  a special set, since if $\tx_j$ and $y'$ belong to this intersection, there can be at most one term in the sum.  
 
We next partition the image $f(X_p)$ into a countable collection $Q_1, Q_2, Q_3, \ldots$ of measurable pieces as follows.  
For each $x \in X_p$, let $B(x, r_x)$ be the ball constructed in the previous paragraph.  
By the $5r$-covering theorem \cite[Chapter 1]{heinonen:analysis}, there exists a set $\{x_i, i \in I\}$ of points in $f(X_p)$ for which the union of the balls $B(x_i, r_{x_i}), i \in I$, covers $f(X_p)$, and for which the balls $B(x_i, r_{x_i}/5)$ are pairwise disjoint.  Since the metric space $(X, d_\varepsilon)$ is separable, the index set $I$ can be taken to be countable.    
We are now ready to construct the elements of our partition inductively.  
Pick arbitrarily an element $i \in I$.  Set $x_1 = x_i$, set $r_1 = r_{x_1}$, and let $Q_1 = B(x_1, r_1)$ be the first element of our partition.  Suppose inductively that $Q_1, \ldots, Q_n$ have already been defined.  If $f(X_p) \subset Q_1 \union \ldots\union Q_n$, we stop; in this case our partition is finite.  
Otherwise, there is some $x_i, i \in I$, with $x_i \not\in Q_1 \union \ldots\union Q_n$; call this element $x_{n+1}$.  Let $r_{n+1}=r_{x_{n+1}}$ and  let $Q_{n+1}=B(x_{n+1}, r_{n+1})\setminus (Q_1 \union \ldots \union Q_n)$ be the next piece in our partition.   
 
In this paragraph, we show that for each piece of the partition constructed in the previous paragraph, its inverse image is a disjoint union of special sets.    Fix such a piece $Q=Q_n$,  and let $B(x, r)=B(x_n, r_n)$ be the corresponding ball.   Then  $f^{-1}(\{x\})\cap X_p= \{\tx_1,\ldots, \tx_l\}$ as above, where $l \le d/p$  depends on $x$.  Define 
$\widetilde{Q}_j=f^{-1}(Q)\cap B(\tx_j,r e ^{-\varepsilon})\intersect X_p$ for $1\le j\le l$ and for $j=l+1, \ldots d/p$, set $\widetilde{Q}_j=\emptyset$.  Then each set $\widetilde{Q}_j$, $j=1, \ldots, d/p$, is a special set, and the union of the $\widetilde{Q}_j$'s is the entire inverse image of $Q$ under $f|_{X_p}$.  

Finally, for $i=1, \ldots, d/p$, let $P_i = \union_n \widetilde{Q_{n,i}}$, where $\widetilde{Q_{n,i}}$ is as constructed in the previous paragraph with $Q=Q_n$.  Since a countable union of special sets with disjoint images is again special, each set $P_i, i=1, \ldots, d/p$, is special, and the proof is complete.  
\qed

\gap 

We go back to the construction of a weak Jacobian. Since we assume that $T$ is an fbc, we may consider a finite measurable partition  $\AAA$ of $Z$ into special sets according to
Prop.\,\ref{prop:decomposition}. 
Fix $A\in\AAA$. Since $T|_A$ is injective, the formula
$\nu_A(B)=\nu(T(B))$ applied to any  measurable subset $B\subset A$ defines
a measure $\nu_A$ on $A$. By the invariance of $\nu$, it follows that $\nu|_A$
is absolutely continuous with respect to $\nu_A$.  The Radon-Nikodym
theorem implies the existence of a measurable non-negative function $h_A$ defined
on $A$ such that $d\nu=h_Ad\nu_A$.

Let $Y_A=\{h_A=0\}$, and define $J_\nu$ on $A$ by
$$
J_\nu=
\left\{
\begin{array}{ll} 
0 &\hbox{on } Y_A,\\
1/h_A & \hbox{on } A\setminus Y_A\,.
\end{array}
\right.
$$
It follows that $\nu(Y_A)=0$ and that for any special set $B\subset A$ disjoint from $Y_A$, 
$$\nu(T(B))=\int_B d\nu_A=\int_B J_\nu( h_Ad\nu_A)=\int_B J_\nu d\nu\,.$$

Let us define $Y=\cup_{\AAA} Y_A$; it follows that $\nu(Y)=0$ and
for any special set $B$ disjoint from $Y$, 
$$\nu(T(B))=\sum_{\AAA}\nu(T(B\cap A))=\sum_{\AAA}\int_{B\cap A} J_\nu d\nu
=\int_B J_\nu d\nu\,.$$

{\noindent\bf Rohlin formula.}\index{Index}{Rohlin formula}
If $\nu$ is ergodic, and if it admits a countable generator of finite
entropy, then the so-called following {\em Rohlin formula} holds \cite[Thm 1.9.7]{przytycki:urbanski:book}: 
\begin{equation}
\label{eqn:rohlin}
h_\nu(T)=\int \log J_\nu d\nu
\end{equation}
 where $J_\nu$ is the Jacobian.

\gap

In general, for non-invertible transformations, the existence of a countable generator of finite entropy can be difficult to establish.
In our setting, we are able to obtain one inequality:

\begin{thm}[Rohlin's inequality]\label{thm:rohlin} Let $f:\XXX_1\to\XXX_0$ be an fbc of degree $d\ge 2$ with repellor $X$
which satisfies [Expansion]. Then, for any ergodic invariant probability measure $\nu$ on $X$, Rohlin's inequality always holds:
$$h_\nu(f)\ge  \int \log J_\nu d\nu$$
\end{thm} 

{\noindent\bf Remark.} The proof will show that if Rohlin's formula holds, 
then $h_\nu(f)=h_\nu(f,\PPP)$ for any finite partition by special sets.
\gap

Let $J_\nu$ be the Jacobian of $\nu$. We let $Y\subset X$ be of $\nu$-measure $0$ 
such that $\nu(f(E))=\int_E J_\nu d\nu$ for all special sets
$E$ with $E\cap Y=\emptyset$. 
Let us restrict $f$ to $X'=X\setminus \cup_{n\ge 0}f^{-n}(Y)$ so that $J_\nu$ becomes a strong Jacobian
for $f:(X',\nu)\to(X',\nu)$.

\gap

We start with a proposition essentially due to Rohlin, cf \cite[Thm 1]{parry:rohlin}.

\begin{prop}\label{prop:symb} Under the assumptions of Theorem \ref{thm:rohlin},
let $\PPP$ be a finite and measurable partition of $X'$ such that,
for any $P\in \PPP$, $f|_P$ is injective.
There exists a measurable map $\psi:(X',\BBB(X'),\nu)\to (\PPP^{\N},\FFF,\mu)$, where $\BBB(X')$ is 
the Borel $\sigma$-algebra of $X'$ and $\FFF$ is the $\sigma$-algebra generated by the cylinders of $Z$,  which satisfies the following
properties:
\bi
\item the map $\psi$ semiconjugates $f$ to the shift map $\sigma$ on $Z=\psi(X')$;
\item the space $Z$ is isomorphic to $X'/(\AAA)$ where $\AAA=\vee_{n\ge 0} f^{-n}(\PPP)$;
\item the probability measure $\mu$ is invariant under $\sigma$;
\item the following holds: $h_\nu(f,\PPP)= h_\mu(\sigma,\psi(\PPP))$.
\ib\end{prop}

\pf For any $x\in X'$, we set $\psi(x)=(P_n)$ where $f^n(x)\in P_n$. That is, $\psi$ sends the point $x$ to its itinerary with respect to $\PPP$ under forward iteration.
It follows that $x\in \cap_{n\ge 0} f^{-n}(P_n)$, so that $Z$ is naturally identified with
$X'/(\AAA)$. Clearly, the map $\psi$ is measurable, $Z$ is invariant under $\sigma$, $\psi$ semiconjugates $f$ to $\sigma$, and $\mu=\psi_*\nu$ is invariant under $\sigma$. 
Let us note that $\PPP_n$ is mapped into the partition of $Z$ by its cylinders of length $(n+1)$. 
Therefore, $\psi(\PPP_n)=(\psi(\PPP))_n=\vee_{0\le k\le n} \sigma^{-k}(\psi(\PPP))$.

For the entropy, one obtains
\begin{eqnarray*}
h(f,\PPP)& = & \lim \frac{-1}{n} \sum_{P\in\PPP_n}\nu(P)\log\nu(P)\\
&= & \lim \frac{-1}{n} \sum_{\psi(P)\in\psi(\PPP_n)}\mu(\psi(P))\log\mu(\psi(P))\\
&= & \lim \frac{-1}{n} \sum_{Q\in(\psi(\PPP))_n}\mu(Q)\log\mu(Q)\\
&=& h(\sigma,\psi(\PPP))\,.\end{eqnarray*}\qed

\pf (Theorem \ref{thm:rohlin}) Let $\III$ be a finite measurable partition of $X'$ obtained by Proposition \ref{prop:decomposition}.
We refine this partition into a finite partition $\PPP$ so that if $P\in\PPP$ and $P\cap W\ne\emptyset$ for some $W\in S(1)$, then $P\subset W$. 
Proposition \ref{prop:symb} implies the existence of $(Z,\sigma,\mu)$ and a measurable
map $\psi:X'\to (\PPP)^{\N}$. Define also $\AAA=\vee_{k\ge 0} T^{-k}(\PPP)$.

\gap

Let us assume that $\mu$ admits an atom $a\in(\PPP)^{\N}$. By invariance of $\mu$, it follows that the point
$a$ is periodic under $\sigma$ of some period $k\ge 1$. Let $E=\psi^{-1}(\{a\})$, and $\hat{E}=\cup_{0\le j\le k-1}f^{-j}(E)$.
It follows that $f^{-1}(\hat{E})=\hat{E}$ and that $\nu(\hat{E})>0$ so that the ergodicity of $\nu$ implies that
$\nu(\hat{E})=1$, and that $f$ shifts cyclically $f^{-j}(E)$, $j\ge 1$. But, by construction of $\psi$, $f|_{\hat{E}}$
is injective, since it corresponds to elements of the partition $\AAA$. Thus, $J_\nu=1$ almost everywhere so that $\int\log J_\nu d\nu =0$. 
Moreover, since $f|_{\hat{E}}$ is injective, there is some integer $N$ such that, for all $n\ge 1$, $\hat{E}$ is covered by $N$ sets from $S(n)$.
This implies with the axiom [Expansion] that $h_{top}(f|_{\hat{E}})=0$. 
Thus by the Variational Principle (\ref{eqn:variational_principle}), $h_\nu(f)=0=\int \log J_\nu d\nu$.

\gap

We may now assume that $\mu$ is non-atomic. By construction, the partition $\psi(\PPP)$ is a finite generator. So Rohlin's formula
holds for $\sigma$ \cite[Thm 2]{parry:rohlin}: $$h_{\mu}(\sigma)=h_{\mu}(\psi(\PPP),\sigma)=\int \log J_{\mu} d\mu\,.$$

Since $(X',\BBB(X'),\nu)$ is  a Lebesgue space \cite{rohlin:measure_theory},
there are conditional measures $\nu_A$ for almost every  $A\in\AAA$ (in the sense that the union of atoms where conditional measures exist is measurable and of 
full $\nu$-measure) and a measure $\nu_{\AAA}$ on $X'/\AAA$ such that, for any Borel set $E$,
$$\nu(E)=\int_{X'/\AAA} \nu_A(A\cap E)  d\nu_{\AAA}\,.$$

We note that $(X'/\AAA,\nu_{\AAA})$ is isomorphic to $(Z,\mu)$ by construction. Let us prove that $J_{\sigma}\ge \EE(J_\nu|\AAA)\circ\psi^{-1}$, 
where  
\begin{equation}
\label{eqn:EJ}
\EE(J_\nu|\AAA)(A)=\int_A J_\nu d\nu_A
\end{equation}
 by definition and where we consider $\psi$ defined on $X'/\AAA$. 

Let $E\subset \psi(P)$, for some $P\in\psi(\PPP)$, so that $\sigma|_E$ is injective, and $\mu(\sigma(E))=\int_E J_{\mu}d\mu$. 
Assume that $\psi^{-1}(E)$ is contained in the set of atoms
for which the conditional measures are well-defined (which has full measure). 
We note that $f|_{\psi^{-1}(E)}$ is injective as well,
since $\AAA$ is finer than $\III$.
By construction, 
$$\sigma(E)=\sigma\circ\psi\circ\psi^{-1}(E)= \psi\circ f\circ \psi^{-1}(E)\,.$$
Since $f(\psi^{-1}(E))\subset \psi^{-1}\circ\psi (f(\psi^{-1}(E)))$, 
it follows that
\begin{eqnarray*}
\mu(\sigma(E))& \ge & \nu(f(\psi^{-1}(E))\\
& = & \int_{\psi^{-1}(E)} J_\nu d\nu\\
& = & \int_{X/\AAA} d\nu_{\AAA}\int_A \chi_{\psi^{-1}(E)} J_\nu d\nu_A\\
&= & \int_E \EE(J_\nu|\AAA)\circ\psi^{-1} d\mu\,.\end{eqnarray*}

According to Proposition \ref{prop:symb} and Rohlin's formula, one obtains using (\ref{eqn:EJ})
$$h_\nu(f,\PPP)=h_{\mu}(\sigma,\psi(\PPP))\ge \int \log\EE(J_\nu|\AAA)d\nu_{\AAA}\,.$$
But, Jensen's formula implies that, for almost every $A\in\AAA$,
$$\log\EE(J_\nu|\AAA)\ge\int_A \log J_\nu d\nu_A$$ so that
$$ \int \log\EE(J_\nu|\AAA)d\nu_{\AAA}\ge \int\log J_\nu d\nu\,.$$
Therefore, 
$$h_\nu(f)\ge h_\nu(f,\PPP)\ge \int \log\EE(J_\nu|\AAA)d\nu_{\AAA}\ge \int\log J_\nu d\nu\,.$$
\qed

We now estimate $h_{\mu_f}(F)$ where 
$F: \partial_\varepsilon \Gamma \to \partial_\varepsilon\Gamma$ and $\mu_f$ 
is the quasiconformal measure constructed in the previous section.  

For $\nu=\mu_f$, one has $J_{\mu_f}(F)=d/ d_F$, cf.(\ref{eqn:F^*}), so that
$$
h_{\mu_f}(F) \ge\int_{\bdry_\varepsilon\Gamma} \log J_{\mu_f}(F)d\mu_f= \log d -\dis\int_{\bdry_\varepsilon\Gamma}\log d_F d\mu_f\,.
$$

\gap

The variational principle (\ref{eqn:variational_principle}) applied to $F$ then implies
$$\log d -\int_{\partial_\varepsilon\Gamma} \log d_F d\mu_f \le h_{\mu_f}(F) \le h_{top}(F)\le v\le  \log d\,.$$

\gap

As a corollary, we obtain the positivity of the topological entropy.

\begin{cor}\label{cor:entropy_positive}
If $d\ge 2$, then the metric entropy of $\mu_f$ is positive.\end{cor}

\pf If this was not the case, it would follow that $\int\log d_F d\mu_f=\log d$.
But $1\le d_F\le d$, so that $d_F=d$ $\mu_f$-a.e.. But, being of maximal degree,
the set $\{d=d_F\}$ is closed, and since the measure $\mu_f$ is supported by all the set $\partial\Gamma$,
this implies that $d_F=d$ everywhere, and by the definition of an fbc, $d=1$.\qed

Furthermore,

\begin{prop}
\label{prop:zero_measure_branch}
If the branch set $B_F$ has measure zero, then $\mu_f$ has maximal entropy $\log d$.  
\end{prop}

\subsection{Equidistribution}
In this subsection, we prove that iterated preimages of points and periodic points are equidistributed according to $\mu_f$.

Let us note that since $F_*F^*\varphi=d\varphi$, the operator $\nu\mapsto (1/d)F^*\nu$ has norm equal to $1$.

\begin{thm}\label{equidb}{\em\bf (Equidistribution of preimages)} For any probability measure $\nu$ whose support is disjoint from 
$o\in\overline{\Gamma}_\varepsilon$,
the sequence $(1/d^n)(F^n)^*\nu$ converges to $\mu_f$ in the weak-$*$ topology. In particular, for any 
$\xi\in\overline{\Gamma_\varepsilon}\setminus\{o\}$ and $n\ge 1$, the sequence of measures 
$$\mu_n^{\xi}=(1/d^n)\sum_{F^n(\zeta)=\xi} d_{F^n}(\zeta)\delta_{\zeta}=(1/d^n)  (F^n)^*\delta_{\xi}$$ 
converges to  $\mu_f$ in the weak-$*$ topology.\end{thm}
\index{Symbols}{$\mu_n^{\xi}$}

We may then deduce the following.

\begin{thm}\label{equidc}{\em\bf (Equidistribution of periodic points)} The sequence of measures supported on $\partial \Gamma_{\varepsilon}$
$$\hat{\mu}_n= \frac{1}{d^n}\sum_{F^n(\xi)=\xi} d_{F^n}(\xi)\delta_{\xi}$$
 converges to  $\mu_f$ in the weak-$*$ topology.\end{thm}
\index{Symbols}{$\hat{\mu}_n$}

 \noindent{\bf Remark:}  Since the number of cycles of period $n$ is not known, the measures $\hat{\mu}_n$ need not be probability measures.\\
 \gap

We start with a lemma (compare with the theory of primitive
almost periodic operators, e.g.  Theorem 3.9 in \cite{eremenko:lyubich}).

\begin{lemma} For any continuous function $\varphi:\overline{\Gamma}\setminus
B(o,1)\to\R$, the sequence of functions $(1/d^n)(F^n)_*\varphi$ is uniformly
convergent towards the constant function $$\int \varphi d\mu_f\,.$$
\end{lemma}

\pf Let us define $A(\varphi)=(1/d)F_*\varphi$.
Let us consider two points $\xi$ and $\zeta$ close enough so
that there exists a curve $\gamma$ joining them and avoiding $o$.  
It follows that the points of $F^{-n}(\{\xi\})$ and $F^{-n}(\{\zeta\})$ are joined together by subcurves of $F^{-n}(\gamma)$
of length bounded by $\ell_{\varepsilon}(\gamma)\cdot e^{-\varepsilon n}$. 

If $\varphi$ is a continuous function on $\overline{\Gamma}\setminus B(o,1)$ with modulus continuity $\omega_{\varphi}$, 
it follows that 
$$|A^n\varphi(\xi)- A^n\varphi(\zeta)|_{\varepsilon} \le 
\frac{1}{d^n}\sum_{F^n(\xi')=\xi}d_{F^n}(\xi')\omega_\varphi(\ell_\varepsilon (\gamma) e^{-\varepsilon n})
\le\omega_\varphi(\ell_\varepsilon (\gamma) e^{-\varepsilon n}).$$
This shows that the sequence $(A^n\varphi)$ is uniformly equicontinuous and
that any limit is locally constant. Thus, if $\overline{\Gamma}\setminus B(o,1)$ is connected, then any limit is constant. Furthermore, since 
$F^*\mu_f=d\mu_f$, it follows that, for any $n$,
$$\int A^n\varphi d\mu_f=\int\varphi d\mu_f$$ so that any constant limit has
to be $\int\varphi d\mu_f$.

If  $\overline{\Gamma}\setminus B(o,1)$ is not connected, one can
argue as follows. Adding a constant if necessary, we can 
assume that $\varphi\ge 0$. Then $(A^n\varphi)$ is a sequence
of non-negative functions, and 
$$\|A(\varphi)\|_{\infty}\le \|\varphi\|_\infty,$$ so that
the norms of $(A^n(\varphi))$ form a decreasing convergent sequence. Let $\varphi_\infty$
be any limit. One knows that it is locally constant; let
us assume that it is not constant. We let $k$ be any
iterate large enough so that, for any maximal open set $E$
such that $\varphi_\infty$ is constant, $F^k(E\cap\partial\Gamma)=\partial\Gamma$.
Then, for any $\xi\in\bdry\Gamma$,
$$|(F^k)_*\varphi_\infty(\xi)|=\sum_{F^k(\zeta)=\xi} \frac{d_{F^k}}{d^k} \varphi_\infty(\zeta)
< \|\varphi_\infty\|_\infty$$ since $\varphi_\infty$ is not locally constant,
but non-negative.
This contradicts the fact that 
$$\|\varphi_\infty\|_\infty=\inf_n \|A^n(\varphi)\|_\infty\,.$$
Thus $\varphi_\infty$ is constant.\qed

\gap

\begin{cor}\label{cor:mixing} The measure $\mu_f$ is mixing.\end{cor}

\pf For any continuous function $\varphi$, and almost every $\xi \in \partial_\varepsilon\Gamma$, the sequence $(1/d^n)(F^n)_* \varphi(\xi)$ tends  
to the value $\mu_f(\varphi)$ by the above lemma.   The operator $A$ has norm one,  so for all $\xi$,  $|(1/d^n)(F^n)_*\varphi(\xi)|\le ||\varphi||_\infty$. 
Hence $$\left|\frac{1}{d^n} (F^n)_*\varphi-\mu_f(\varphi)\right|^2\le 4 ||\varphi||_{\infty}^2$$
and the dominated convergence theorem implies that  $F^n_*\varphi \to \mu_f(\varphi)$ in $L^2(\partial_\varepsilon\Gamma, \mu_f)$.   
It follows from  Proposition 2.2.2 in \cite{dinh:sibony:dynamique} 
that $\mu_f$ is mixing.\qed

\gap

\pf (Theorem \ref{equidb}) Let $\nu$ be a measure supported off
the origin in $\Gamma$. For any continuous function $\varphi$,
one has $$\la (1/d^n)(F^n)^*\nu,\varphi\ra=\la\nu,(1/d^n)(F^n)_*\varphi\ra=
\int A^n\varphi d\nu\,.$$
It follows from dominated convergence that this sequence tends to
$$\int \left(\int\varphi d\mu_f\right)d\nu=\int\varphi d\mu_f=\la \mu_f,\varphi\ra$$
so that $(1/d^n)(F^n)^*\nu$ tends to $\mu_f$.\qed

We precede the proof of Theorem \ref{equidc} by a couple of intermediate results concerning periodic points, beginning with an index-type result.

\begin{prop}\label{eqp} Let $U$ be a Hausdorff connected, locally connected and locally compact open set,
$U'$ a relatively compact connected subset of $U$, and $g:U'\to U$  a finite branched covering of degree $d\ge 1$
which satisfies Axiom [Expansion] with respect to the covering $\UUU_0=\{U\}$. Then
$$d=\sum_{g(x)=x} d_g(x)\,.$$\end{prop}

\pf If, for every $n$, $g^{-n}(U)$ is connected, then Axiom  [Expansion] implies that $\cap f^{-n}(U)$ is a single point $x$,
which is fixed: thus $d=d_g(x)$.

Otherwise, let $k_0$ be the maximal integer such that $g^{-n}(U)$ is connected. Then $g^{-(k_0+1)}(U)$ is a finite
union of connected open sets  $U_1^0,\ldots, U_{m_0}^0$ where $m_0>1$. Each restriction $g_j:U_j^0 \to g^{-k_0}(U)$ is a finite
branched covering of degree $d_j<d$, and $d=\sum d_j$.

For each $g_j$, one may repeat this procedure until it stops.  The proposition follows easily.\qed


For an open set $U\in\mathbf{U}$ corresponding to a point $\xi \in \Gamma$, denote by 
$\mu_n^U=d^{-n}\sum_{F^n(\wtU)=U} d_{F^n}(\wtU)\delta_{\wtU}=\frac{1}{d^n}(F^n)^*\delta_\xi$ the measure appearing in the statement of Theorem \ref{equidb}.  
Recall that the measure $\hat{\mu}_n$ describing the distribution of periodic points is given by $\hat{\mu}_n=d^{-n}\sum_{F^n(\xi)=\xi} d_{F^n}(\xi)\delta_\xi$.

\begin{lemma}\label{eql2} Let $U$ be a vertex, and let us consider open subsets
$W_1$ and $W_2$ of $\overline{\Gamma}_\varepsilon\setminus\{o\}$ which intersect $\partial_\varepsilon\Gamma$ such that  
$\overline{W_1}\subset W_2$
and such that $\overline{W_2}\cap\partial_\varepsilon\Gamma\subset \phi_f(U)$.
For $n$ large enough, $$\mu_n^U(W_1)\le \hat{\mu}_n(W_2)\quad \mbox{and}\quad\hat{\mu}_n(W_1)\le\mu_n^U(W_2).$$
\end{lemma}

\pf  See Figure 
3.4.3.  

\begin{figure}[htb]
\label{fig:perptfig}
\psfrag{1}{\mbox{$U\cap X$}}
\psfrag{2}{\mbox{$\phi_f^{-1}(\cl{W}_2\cap \partial_\varepsilon \Gamma)$}}
\psfrag{3}{\mbox{$\widetilde{U}\cap X$}}
\psfrag{4}{\mbox{$\phi_f$}}
\psfrag{5}{\mbox{$\widetilde{U}$}}
\psfrag{6}{\mbox{$W_2$}}
\psfrag{7}{\mbox{$\overline{W}_1$}}
\psfrag{8}{\mbox{$X$}}
\psfrag{9}{\mbox{$\partial_\varepsilon \Gamma$}}
\psfrag{x}{\mbox{$\overline{\Gamma}_\varepsilon$}}
\includegraphics[width=5in]{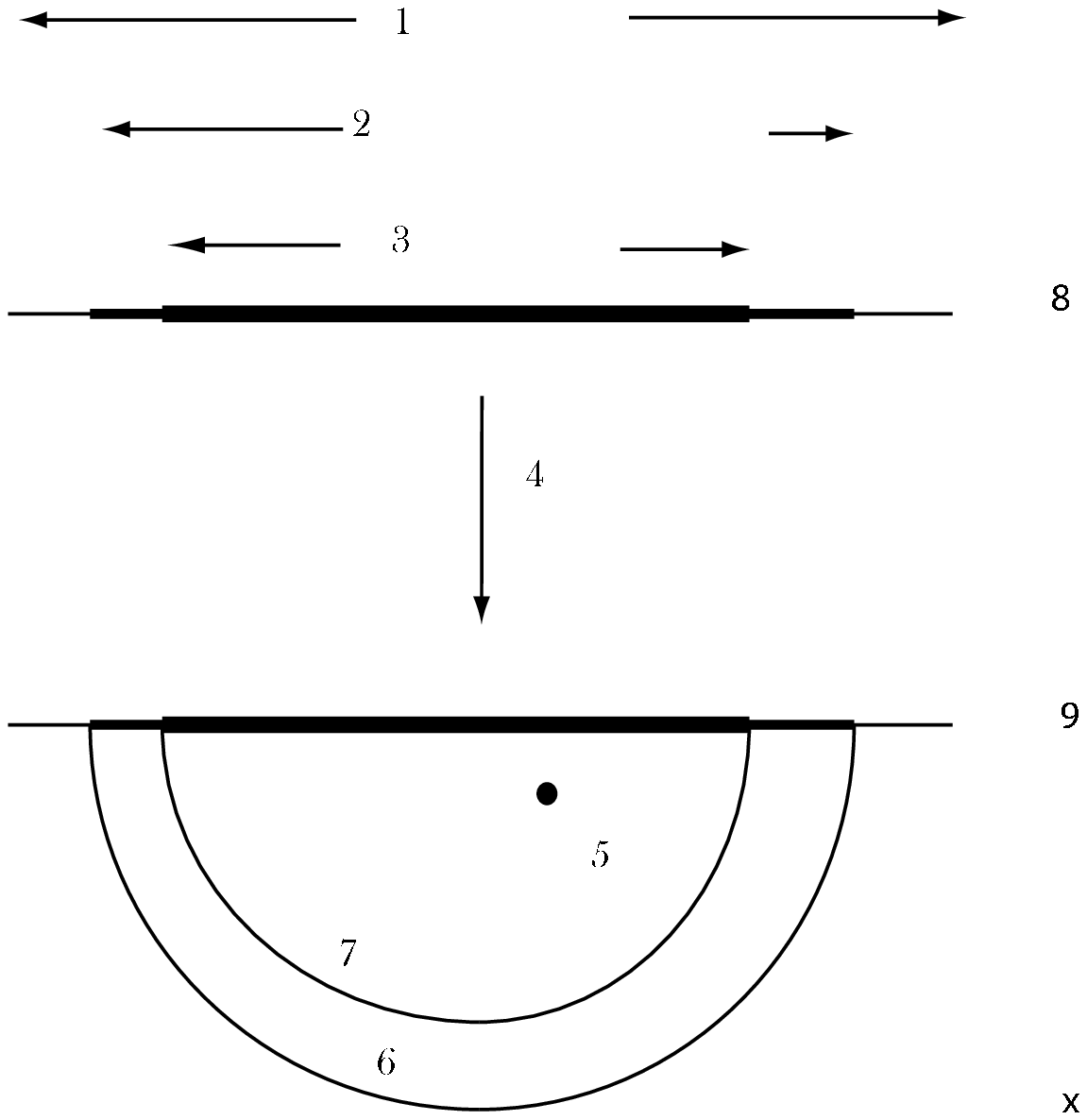}
\begin{center}{\sf Figure \ref{fig:perptfig}}\end{center}
\end{figure}

Suppose  $\wtU\in S(n+|U|)$ and, as a vertex point in $\cl{\Gamma}_\varepsilon$, belongs to $W_1$. Then $\phi_f(\wtU)\subset W_2$ if  $n$ is large enough,
since $d_{\varepsilon}(\wtU,\phi_f(\wtU))= (1/\varepsilon)e^{-|\wtU|} \to 0$. 
Hence $\wtU \subset \phi_f^{-1}(\cl{W}_2 \intersect \bdry_\varepsilon\Gamma) \subset U$.  
So, if moreover $f^n(\wtU)=U$ then $f^n|_{\wtU}: \wtU \to U$ satisfies the hypotheses of Proposition \ref{eqp} and so $$d_{f^n}(\wtU)=\sum_{f^n(x)=x, x\in\wtU} d_{f^n}(x).$$  For the periodic points $x$ appearing in the sum, $\phi_f(x)\in W_2$.
\gap

Therefore,
\[ 
\begin{array}{lcl}
\mu_n^U(W_1)& =& \displaystyle{\frac{1}{d^n}\sum_{f^n(\wtU)=U,\wtU\in W_1}d_{f^n}(\wtU)} \\[.5in]
\; & = & \displaystyle{\frac{1}{d^n}\sum_{f^n(\wtU)=U,\wtU\in W_1}\sum_{x\in \wtU, f^n(x)=x} d_{f^n}(x)}.
\end{array}
\]
Since for each such $\wtU$ appearing in the sum, (i) $\phi_f(\wtU) \subset W_2$, and (ii) for fixed $n$, the $\wtU$'s and their images under $\phi_f$ are  pairwise disjoint, we have 
\[ \mu_n^U(W_1) \leq \frac{1}{d^n}\sum_{x \in W_2, f^n(x)=x}d_{f^n}(x) = \hat{\mu}_n(W_2).\]

Similarly,  if $f^n(x)=x$ and $x\in W_1$, then there is a unique $\wtU\in f^{-n}(U)$ such that $x\in\wtU$.
Therefore, for $n$ large enough, $\phi_f(\wtU)$ has compact closure in $W_2$ and $\wtU\in W_2$. Thus
$$\hat{\mu}_n(W_1)\le\frac{1}{d^n}\sum_{f^n(\wtU)=U,\wtU\in W_2}\sum_{f^n(x)=x,\phi_f x\in\wtU}d_{f^n}(x)\leq \frac{1}{d^n}\sum_{f^n(\wtU)=U,\wtU\in W_1}d_{f^n}(\wtU)=\mu_n(W_2).$$\qed

\pf (Theorem \ref{equidc}) Note that the number of cycles is unknown. Nevertheless,
it follows from Lemma \ref{eql2} that $\{\hat{\mu}_n\}_n$ is relatively compact in the weak topology.
Let $\hat{\mu}$ be an accumulation point. We will prove that $\hat{\mu}=\mu_f$ using their Borel regularity.

\gap

Let $U$ be a vertex, and let us consider a compact subset $K$ of $\partial\Gamma$, open subsets
$W_1,W_2$ and $W_3$ of $\overline{\Gamma}\setminus\{o\}$ such that  
$$K\subset W_1\cap\partial\Gamma \subset\overline{W_1}\subset W_2\subset\overline{W_2}\subset W_3$$
and such that $\overline{W_3}\cap\partial\Gamma\subset \phi_f(U)$.

\gap

Let $\varphi_1$ and  $\varphi_2$ be two continuous functions such that 
$$\chi_K\le\varphi_1\le\chi_{W_1}\le\chi_{W_2}\le\varphi_2\le\chi_{W_3}.$$

Let us fix $\eta>0$; if $n$ is large enough then
$$\left\{\begin{array}{l} |\hat{\mu}(\varphi_j)-\hat{\mu}_n(\varphi_j)|\le \eta\\
|\mu_f(\varphi_j)-\mu_n^U(\varphi_j)|\le \eta\end{array}\right.$$ for $j=1,2$.

\gap

Therefore, by the preceding Lemma \ref{eql2} and the regularity of the measures,
$$\begin{array}{lcccccc}
\hat{\mu}(K) & \le &\hat{\mu}(\varphi_1) & \le&\hat{\mu}_n(\varphi_1)+\eta & \le & \hat{\mu}_n(W_1)+\eta\\
& \le  & \mu_n^U(W_2)+\eta & \le &  \mu_n^U(\varphi_2)+\eta &  \le  & \mu_f(\varphi_2)+2\eta\\
& \le  & \mu_f(U)+2\eta\,. & & &  &   \end{array}$$
Since this is true for any compact subset of $U$, the regularity of the measures imply $\hat{\mu}(U)\le \mu_f(U).$

\gap

Similarly, $$\begin{array}{lcccccc}
\mu_f(K) & \le  & \mu_f(\varphi_1) & \le & \mu_n^U(\varphi_1)+\eta & \le & \mu_n^U(W_1)+\eta\\
& \le  & \hat{\mu}_n(W_2)+\eta & \le &  \hat{\mu}_n(\varphi_2)+\eta & \le  & \hat{\mu}(\varphi_2)+2\eta\\
& \le & \mu_f(\phi_f(U))+2\eta & & &  & \end{array}$$
from which we deduce $\hat{\mu}(\phi_f(U))\ge \mu_f(\phi_f(U))$, so that $\hat{\mu}=\mu_f$.\qed

\subsection{Hausdorff dimension}   \label{subsecn:dim} 
Let $Z$ be a metric space. 
 Given $\delta>0$, a {\em $\delta$-cover} of $Z$ is a covering of $Z$ by sets of diameter at most $\delta$.
For $s \geq 0$, set 
\[ \HHH_\delta^s(Z)  = \inf  \sum_i (\diam U_i)^s  \]
where the infimum is over all $\delta$-coverings of $Z$ by sets $U_i$.  As $\delta$ decreases, $\HHH_\delta^s$ increases and so the {\em $s$-dimensional Hausdorff measure}
\index{Index}{Hausdorff measure}%
of $Z$ 
\[ \HHH^s(Z) = \lim_{\delta \to 0^+} \HHH^s_\delta(Z) \in [0,\infty]\]
exists.   The {\em Hausdorff dimension}
\index{Index}{Hausdorff dimension}%
 of $Z$  is given by 
\[ \dim_H(Z) = \inf\{ s: \HHH^s(Z)=0\} = \sup\{s: \HHH^s(Z) = \infty\}.\]
Using balls instead of arbitrary sets in the definition leaves the dimension unchanged.  See for instance \cite[\S\,5.1]{mattila:gmt}.

We now compute the Hausdorff dimension of the boundary $\partial_\varepsilon\Gamma$.   
Fix $s>0$.  By Lemma \ref{omdiam}, for any vertex $\xi \in S(n)$, $\diam_\varepsilon\mho_\infty(\xi) \leq Ce^{-\varepsilon n}$.  Therefore $\bdry_\varepsilon\Gamma$ is covered by at most $|S(n)|$ sets of diameter $\delta_n = Ce^{-\varepsilon n}$ and so 
\[ \HHH^s_{\delta_n} (\bdry_\varepsilon\Gamma) \leq |S(n)|e^{-\varepsilon n s}.\]
Suppose now that $s>\frac{v}{\varepsilon}$.  Recall that by definition, $v=\lim \frac{1}{n}\log |S(n)|$.
There exists $\eta>0$ with $v+\eta-\varepsilon s < 0$.  It follows that for all $n$ sufficiently large, 
\[ \HHH^s_{\delta_n}(\bdry_\varepsilon\Gamma) \leq |S(n)|e^{-\varepsilon n s} < e^{n(v+\eta-\varepsilon s)} < \infty.\]
Hence 
\[ \HHH^s(\bdry_\varepsilon\Gamma) < \infty \;\; \mbox{ for all } \;\; s>\frac{v}{\varepsilon}\]
and therefore $\dim_H(\bdry_\varepsilon\Gamma) \leq \frac{v}{\varepsilon}$.  
\gap

We now investigate lower bounds by  appealing to  the following  result, 
which is
similar to R.\,Ma\~n\'e's dimension formula \cite{mane:hdim}\,:

\begin{thm}\label{entmeath} If $\mu$ is an ergodic invariant measure  
with
positive entropy, then
$$\liminf \frac{\log \mu(B_\varepsilon(\xi,r))}{\log r}\ge  
h_\mu(F)/\varepsilon\,,$$ 
for $\mu$-almost every $\xi$.
\end{thm}

\pf Since $F$ is  
$e^\varepsilon$-Lipschitz,
it follows that
$B_\varepsilon(\xi,r)\subset S(\xi,n,re^{n\varepsilon})$. Since $\mu$ is invariant and ergodic,  
it follows from a formula of Brin and Katok \cite{brin:katok:local} that an equivalent definition of metric entropy is 
$$h_\mu(F)=\lim_{r\to 0}\limsup_{n\to\infty} -\frac{1}{n}\log\mu 
(S(\xi,n,r))\,,$$ for $\mu$-a.e. $\xi$.

Choose a generic point $\xi$ for $\mu$ and let $\eta>0$; we will write $B_\varepsilon(r)=B_\varepsilon(\xi,r)$. 
There are some $r_0>0$ and $n_0\in\N$ such that, if $r\le  
r_0$ and
$n\ge n_0$ then
$$\left|-\frac{1}{n}\log\mu S(n,r)- h_\mu(F)\right|\le 2\eta\,.$$ We  
choose
$r_n=r_0e^{-\varepsilon n}$ and we obtain 
$$\frac{\log \mu(B_\varepsilon(r_n))}{\log r_n}\ge -\frac{\log  
\mu(S(n,r_0))}{n(\varepsilon
-\log (r_0)/n)}\ge \frac{h_\mu(F)-2\eta}{\varepsilon -\log(r_0)/n}\,,$$
so $$\liminf \frac{\log \mu(B_\varepsilon(r_n))}{\log r_n}\ge
\frac{h_\mu(F)-2\eta}{\varepsilon}\,.$$
Given $r>0$, fix $n$ so that $B_\varepsilon(r_{n+1})\subset B_\varepsilon(r)\subset  
B_\varepsilon(r_n)$
and
$$\frac{\log \mu(B_\varepsilon(r))}{\log r}\ge\frac{\log  
\mu(B_\varepsilon(r_n))}{\left(\log
\frac{r_{n+1}}{r_n}\right) + \log r_n}\,.$$
Thus,
$$\liminf \frac{\log \mu(B_\varepsilon(r))}{\log r}\ge h_\mu(F)/\varepsilon\,.$$\qed

\gap

It follows that for the measure $\mu_f$ we have constructed and for any  
$\eta>0$ and $r$
small enough,
$$\mu_f(B_\varepsilon(r))\le r^{(1/\varepsilon)h_{\mu_f}(F) -\eta}\,.$$ 
 This implies that the local upper pointwise dimension of $\mu_f$ satisfies $\dim
\mu_f\ge(1/\varepsilon)h_{\mu_f}(F)$. Therefore 
$$\frac{h_{\mu_f}(F)}{\varepsilon}\le \dim\mu_f\le\dim\partial_\varepsilon \Gamma\le\frac{v}{\varepsilon}\le\frac{1}{\varepsilon}\log  d\,.$$

\gap

\pf (Theorem \ref{maxent}) By ergodicity and uniqueness of the class of \qc measures of given dimension, 
it follows that $\mu_f$ is unique (cf. Theorem \ref{trivmeas}). Theorem \ref{equidb} and Theorem \ref{equidc} prove
the equidistribution of preimages and periodic points according to $\mu_f$. The mixing property has also been
proved (Corollary \ref{cor:mixing}).  The claimed entropy and dimension estimates were proven in \S\,\ref{subsecn:entropy} and \S\,\ref{subsecn:dim}, respectively.  

\qed

\section{Properties for cxc maps following hyperbolicity}
\label{secn:properties_following_hyp}

In this section we assume that $f: \XX_1 \to \XX_0$ is a finite branched covering with repellor $X$ satisfying the conditions at the beginning of section \S\,\ref{secn:def_top_cxc}, 
and is topologically cxc with respect to some open cover $\UUU_0$.  Thus, the topological Axioms [Expansion], [Irreducibility], and [Degree] hold.

\subsection{Canonical gauge}
Let $\Gamma=\Gamma(f, \UUU_0)$ be the associated Gromov hyperbolic graph as in Section \ref{secn:spaces_associated}. 
 Recall that by Theorem \ref{thm:construction}, for $\varepsilon>0$ small enough, there is a homeomorphism $\phi_f: X \to \partial_\varepsilon\Gamma$ conjugating $f$ on $X$ to the Lipschitz map $F: \partial_\varepsilon\Gamma \to \partial_\varepsilon\Gamma$.  

\begin{thm}\label{cxcimphyp} If $f$ is topological cxc, then $\Gamma$ is hyperbolic for any covering $\UUU$ with sufficiently small mesh, and 
$F:\partial \Gamma\to\partial \Gamma$   satisfies axioms [Roundness distortion] and [Relative diameter distortion] with respect to the covering $\VVV_0 = \{\phi_f(\UUU_0\intersect X)\}_{U \in \UUU_0}$. If $f$ is furthermore metric cxc with respect to a metric $d$ on $\XX_0$, then the map $\phi_f$ is a quasisymmetry between the metric spaces $(X,d)$ and 
$\bdry_\varepsilon \Gamma$.\end{thm}

So, if $f$ is topologically cxc, and if $d_\varepsilon$ denotes the metric given by Theorem \ref{thm:construction}, then the dynamics on $X$, when equipped with the pulled-back metric $\phi_f^*(d_\varepsilon)$, is essentially metrically cxc.  The qualifier ``essential'' is necessary:  without further assumptions, we do not know that the dynamics of $f$ on the repellor itself is even topologically cxc with respect to the covering $\{U \intersect X : U \in \UUU_0\}$.  In particular, we do not know that $F: \bdry_\varepsilon \Gamma \to \bdry_\varepsilon\Gamma$ is metrically cxc, and so we are unable to apply Theorem \ref{thm:invariance_of_cxc} to prove the last conclusion.  

\gap

\pf  Since $f$ is topologically cxc, a metric exists so that the mesh of $S(n)$ has exponential decay (cf. Theorem \ref{thm:exp_decay}) and
$\phi_f:X\to\partial_\varepsilon \Gamma$ is a homeomorphism as soon as  $\varepsilon$ is small enough. 
Therefore Theorem \ref{homeoimphyp} implies that $\Gamma$ is hyperbolic and that its quasi-isometry class is well-defined. 

We let $\mathbf{V}=\phi_f(\mathbf{U})$. Axioms [Irreducibility], [Expansion] and [Degree] hold through the conjugation.
The [Diameter distortion] and [Roundness distortion] Axioms follow from Proposition \ref{prop:F_is_nearly_cxc}.

\gap

Let us assume from now on that $f$ is cxc.
 Our strategy is as follows.  We will first establish that $\phi_f$ is weakly quasisymmetric by the blowing up/down argument given in the proof of Theorem \ref{thm:invariance_of_cxc}.  The proof concludes by arguments similar to those given in the proof of Proposition  \ref{homeoimpqi}.  

Let $\delta$ be the Lebesgue number of $\UUU_0=S(1)$. Let $x\in X$ and let $r\in(0,\delta/L)$ where $L$ is  given by Proposition \ref{prop:BLC}. 

\gap

By Proposition \ref{prop:BLC}  we may find vertices $W',W$ such that 
$$B(x,r/L)\subset W'\subset B(x,r)\subset W\subset B(x,Lr).$$

Since $\diam W'\asymp \diam W$, we have  $|W-W'|=||W|-|W'||\le N$ for some constant $N$.
Let $n=\min\{|W|,|W'|\}-1$. It follows that $f^n( W')\subset f^n(  B(x,r))\subset f^n(  W)$ and that
the roundness of $f^n(W')$, $f^n(B(x,r))$ and $f^n(W)$ at $f^n(x)$ is bounded by $\rho_+(L)$.

By the uniform continuity of the conjugacy $\phi_f$ and its inverse, and the fact that all these sets have a definite size, 
there exists a constant $K$ such that
the roundness of $\phi_f(f^n(W'))$, $\phi_f(f^n(B(x,r)))$ and $\phi_f(f^n(W))$ at $\phi_f(f^n(x))$ is bounded by $K$.

\gap

Therefore, radii $r$ and $r'$ exist such that
$$B(F^n\phi_f(x), r' e^{\varepsilon p}) \subset \phi_f(f^n(W'))\subset B(F^n\phi_f(x), Kr' e^{\varepsilon n}) $$
and $$B(F^n\phi_f(x), r e^{\varepsilon p}) \subset \phi_f(f^n(W))\subset B(F^n\phi_f(x), Kr e^{\varepsilon n})\,. $$

Proposition \ref{ba2} implies that there is some finite constant $H$ such that
$$\roundness(\phi_f(B(x,r)), \phi_f(x))\lesssim \frac{\diam_{\varepsilon} \phi_f(W)}{\diam_{\varepsilon} \phi_f(W')}
\asymp e^{-\varepsilon (|W-W'|)}\le H.$$
Therefore $\phi_f$ is weakly quasisymmetric.

\gap

Using the uniform continuity of $\phi_f$ and its inverse again, it is enough to consider 
$x,y,z\in X$ such that $|x-y|_X,|x-z|_X\le \delta/L$. We argue as for Proposition \ref{homeoimpqi}.

\gap

Hence, we may find $W_y'$ and $W_z$ in $\Gamma$ such that
\be 
\item $y\notin W_y'$, $\diam W_y'\asymp |x-y|_X$ and $\roundness (W_y',x)\le K$, 
\item $z\in W_z$, $\diam W_z\asymp |x-z|$ and $\roundness (W_z,x)\le K$, 
\eb
for some universal $K$. 

It follows that $$\frac{|\phi_f(x)-\phi_f(y)|_\varepsilon}{|\phi_f(x)-\phi_f(z)|_\varepsilon} 
\asymp \frac{\diam_\varepsilon W_y'}{\diam_\varepsilon W_z}\,.$$

\gap

If $|x-y|_X$ and $|x-z|_X$ are equivalent, then Proposition \ref{prop:local_comparability} implies the bounds.

\gap

If $|x-y|_X$ is small compared to  $|x-z|_X$, then $W_y'\subset W_z$.
Therefore $$\frac{|x-y|_X}{|x-z|_X}\asymp \frac{\diam W_y'}{\diam W_z}\ge \delta_+^{-1}\left(\frac{c_{|W_y'|-|W_z|+1}}{d_1}\right),$$
where we recall that $c_n$ denotes the smaller diameter of sets in $S(n)$. 

This implies that $$\frac{|\phi_f(x)-\phi_f(y)|_\varepsilon}{|\phi_f(x)-\phi_f(z)|_\varepsilon}$$ is bounded
by a function of  $$\frac{|x-y|_X}{|x-z|_X}$$ which goes to zero as the ratio tends to zero.

\gap

If $|x-y|_X$ is large compared to  $|x-z|_X$, then $W_y'\supset W_z$.
Therefore $$\frac{|x-z|_X}{|x-y|_X}\asymp \frac{\diam W_z}{\diam W_y'}\ge \delta_-\left(\frac{d_{|W_z|-|W_y'|+1}}{c_1}\right).$$
We may conclude as above.

This proves that $\phi_f$ is quasisymmetric. 
 \qed

\gap

{\noindent\bf Remark.} Theorem \ref{thm:invariance_of_cxc} can be recovered with Theorem \ref{homeoimphyp}
and Theorem \ref{cxcimphyp}. 

\gap

As an application, we obtain the following result.  

\begin{defn}Let $f: \XX_1 \to \XX_0$ have  repellor $X$ and be topologically cxc with respect to some open covering $\UUU_0$.   The associated {\em conformal gauge} 
\index{Index}{conformal gauge}%
\index{Index}{gauge, conformal}%
$\GGG$ is the set of all metrics on $X$ which are quasisymmetrically equivalent to a metric of the form $\phi_f^*(d_\varepsilon)$, where $d_\varepsilon$ is the metric on $\bdry_\varepsilon\Gamma$ and $\phi_f: X \to \bdry_\varepsilon\Gamma$ is as above. 
\end{defn} 

\begin{thm}[Canonical gauge]
\label{thm:canonical_gauge}
\be
\item 
Let $f: \XX_1 \to \XX_0$ have  repellor $X$ and be topologically cxc with respect to some open covering $\UUU_0$.   Then the conformal gauge $\GGG$  is nonempty and depends only on its topological conjugacy class.  
\item If $U\intersect X$ is connected for every $U \in \UUU_0$, then the conformal gauges of $f: \XX_1 \to \XX_0$ and $f|_X: X \to X$ agree.
\item If the system is in addition metrically cxc with respect to some metric $d$ on $\XX_0$, the conformal gauge $\GGG$ of $f$ agrees with the conformal gauge of the metric space $(X, d|_X)$.   
\eb
\end{thm}

\pf  (1) follows from the uniqueness result Theorem  \ref{homeoimpqi} and (3) from the preceding Theorem \ref{cxcimphyp}.   The graph constructed using $f$ and $\UUU_0$ is naturally identified with that constructed using $f|_X$ and $\VVV_0=\{U\intersect X : U \in \UUU_0\}$ and the induced conjugacies respect this identification.  Therefore the metrics on $X$ obtained by pulling back the metrics on the boundaries of the two graphs coincide and (2) follows.  
\qed

\noindent{\bf Remark.}  The preceding theorem implies that the gauge depends only on the dynamics {\em near} the repellor.   One may surmise that it should really depend only on the dynamics {\em on} the repellor itself.   Conclusion (2) implies that this is true once $X$ is locally connected.  In the non-connected case, however, a proof remains elusive.

To illustrate the subtleties, fix $d \geq 2$,  let $X=\{1, 2, \ldots, d\}^\N$ be equipped with the metric $|x-y|=2^{-(x|y)}$ where $(x|y) = \min_i\{x_i \neq y_i\}$, and suppose $h: X \to X$ is a topological conjugacy, i.e. an automorphism of the one-sided shift on $d$ symbols.  If the gauge of $f$ depends only on the dynamics on $X$, then every such $h$ should be quasisymmetric.   This is indeed the case, and a proof may be given along the following lines.

\be
\item Start with a round closed disk $D \subset \C$.  
 For each $i=1, \ldots, d$, choose a similarity $g_i: \C \to \C$ such that $g_i(D) \intersect g_j(D) = \emptyset$ whenever $i \neq j$.  This defines a conformal iterated function system (IFS).  There is a unique nonempty compact set $K \subset D$ for which $K=\union_{i=1}^d g_i(K)$.   Using a blowing up/down argument, one shows that the attractor of this IFS is quasisymmetrically equivalent to $X$.

\item Using quasiconformal surgery, one builds a uniformly quasiregular map $G: \rs \to \rs$ such that $G|_{g_i(D)} = g_i^{-1}$ for each $i$, and such that $G=z^d + O(z^{d-1})$ as $z \to \infty$.   By Sullivan's Theorem \ref{thm:uqrsullivan}, $G$ is quasiconformally conjugate to a degree $d$ polynomial $p(z)$, and $K$ is quasiconformally (hence quasisymmetically) equivalent to the Julia set $J$ of $p$.

\item By results of Blanchard, Devaney, and Keen \cite{blanchard:devaney:keen:automorphisms}, every automorphism of the shift on $d$ symbols is realized as a monodromy in the shift locus of degree $d$ polynomials.  (The proof depends on the existence of a nice set of generators.)

\item As a polynomial varies in the shift locus, its Julia sets varies holomorphically \cite{ctm:ds:qciii}.  Hence the induced monodromy is quasiconformal, hence quasisymmetric.  In conclusion, we see that every automorphism is realized by a quasisymmetric map.   
\eb 

There seems to be a combinatorial obstacle to promoting topological conjugacies to quasisymmetric conjugacies for noninvertible expanding conformal dynamical systems with disconnected repellors.  Even for hyperbolic rational maps $f, g$ with disconnected Julia sets $J(f), J(g)$, it is not known if every topological conjugacy $h: J(f) \to J(g)$ is quasisymmetric. 

This is known in the following special cases.  First, if  $h$ extends to a conjugacy on a neighborhood of $J(f), J(g)$ then Theorem \ref{thm:invariance_of_cxc} applies and $h$ is quasisymmetric.  However, even for maps with connected Julia set, such an extension need not exist.   Second, if $f$ and $g$ are merely combinatorially equivalent in the sense of McMullen \cite{ctm:siegel} on a neighborhood of their Julia sets, then there is a quasiconformal conjugacy between $f$ and $g$ near their Julia sets.  In both cases, conditions on the dynamics near, not just on, the Julia sets are assumed.

In contrast, we have the following result in the setting of hyperbolic groups; see the appendix \ref{appendix:hyperbolic_groups} for definitions in what follows.  Suppose  $G_1, G_2$ are two hyperbolic groups, and suppose  $h: \bdry G_1 \to \bdry G_2$ conjugates the action of $G_1$ to the action of $G_2$.    By definition, this implies that there is some isomorphism $\Phi: G_1 \to G_2$ for which $h(g(x)) = \Phi(g)(h(x))$ for all $x \in \bdry G_1$ and all $g\in G_1$.  One has necessarily that $h$ arises as the boundary values of $\Phi$.  To see this, note that it is enough to verify that $h=\bdry\Phi$ on the dense set of fixed points of hyperbolic elements.  Suppose $g_1 \in G_1$ is hyperbolic with attracting fixed point $\omega_1$ and $g_2=\Phi(g_1)$ has attracting fixed point $\omega_2$.    Since $h$ is a continuous conjugacy we have $h(\omega_1)=\omega_2$.  But $\omega_i = \lim_n g_i^n$ and this forces $\bdry\Phi(\omega_1)=\omega_2=h(\omega_1)$.  Thus, every topological conjugacy on the boundary is induced from a combinatorial equivalence, i.e. from an isomorphism of the groups.  

This suggests that perhaps there is yet another essential difference between the setting of noninvertible cxc maps and of hyperbolic groups.  
\gap

\begin{cor}\label{topcxc_cxc1} If $f:X\to X$ is a topological cxc map, where $\XX_1=\XX_0=X$, then
$F:\partial\Gamma\to\bdry\Gamma$ is metrically cxc.  Therefore $X$ admits a metric, unique up to quasisymmetry, for which the dynamics is metrically cxc.\end{cor}

\pf The assumptions imply that $X$ is locally connected, and that $(S(n))_n$ is a basis of the topology
by connected open sets. Proposition \ref{ba3} implies the cxc property.\qed


\begin{cor}\label{topcxc_cxc2} If $f:(\XX_1,X)\to (\XX_0,X)$ is a topological cxc map with $f$ a non ramified covering, 
then there is some $R>0$, such that, if we set $\YY_0= \overline{\Gamma_{\varepsilon}}\setminus B_{\varepsilon}(o,R)$
and $\YY_1= F^{-1}(\YY_0)$, then 
$F:(\YY_1,\partial\Gamma)\to (\YY_0,\Gamma)$ is cxc.\end{cor}

\pf Since $f$ is a cover, there is some level $n_0$ such that, for any $n\ge n_0$, any $U\in S(n)$, the restriction
of $f$ to $U$ is injective. This implies that the local degree function for $F$ is $1$ at any point close enough
to $\partial\Gamma$.

Furthermore, if $n_0$ is large  enough, then $F^{-1}(\mho(U))$ will be a disjoint union of $d$ shadows based
at $F^{-1}(\{U\})$.  

Therefore, if we set $\YY_0=\overline{\Gamma_{\varepsilon}}\setminus B_{\varepsilon}(o,e^{-\varepsilon n_0})$ and
$\YY_1=F^{-1}(\YY_0)$, then
$F:\YY_1\to \YY_0$ is a degree $d$ covering.

For any $\xi\in\partial\Gamma$, let $V(\xi)$ be the connected component of the interior of $\mho(W)$ for some $W\in S(n_0)$ containing
$\phi_f^{-1}(\xi)$. Note that the interior of $\mho(W)$ is not empty since it contains $\phi_f(W)$. Since $\overline{\Gamma}$ is locally connected, $V(\xi)$ is open, and we may extract a finite subcover $\VVV$. Proposition \ref{ba3} implies
that $F$ is cxc. 
\qed

\subsection{Existence and uniqueness of the measure of maximal entropy}

Let $(X,d)$ be a metric space and $Q>0$.  A Radon measure $\mu$  is {\em Ahlfors-regular of dimension $Q$} 
\index{Index}{Ahlfors regular}%
if, for any $r\le \diam X$ and any ball $B(r)$, the measure of a ball of radius $r$ satisfies $\mu(B(r))\asymp r^Q$.
In this case, the measure $\mu$ is equivalent to the Hausdorff measure of dimension $Q$ on $X$. We may
then also speak of an Ahlfors regular metric space, keeping the measure implicit.
The measure is {\em doubling} 
\index{Index}{doubling}%
if there is some constant $C>1$ such that $\mu_f(2B)\le C\mu_f(B)$.  Ahlfors regularity implies doubling, but not conversely.  

\gap

\begin{thm}\label{thm:cor_cxc} Let $f:\XX_1 \to \XX_0$ be a topological cxc map of degree $d$ having repellor $X$. 
Let $d_\varepsilon$, $\phi_f$, and $\mu_f$ be the metric, conjugacy, and measure on $\bdry_\varepsilon\Gamma $ given by Theorems \ref{thm:construction} and \ref{maxent}, respectively.  
Then $\mu_f$ is the unique measure of maximal entropy $\log d$. 
It is Ahlfors regular of dimension $\frac{1}{\varepsilon}\log d$. 
If $f$ is furthermore metric cxc with respect to a metric $d$ on $\XX_0$, then on the metric space $(X,d)$, the measure $\phi_f^{*}(\mu_f)$ is doubling.  \end{thm}

\pf 
Ahlfors regularity follows from Axiom [Degree] and the Lemma of the shadow. Let us fix a ball $B_{\varepsilon}(\xi,r)\subset\partial \Gamma$. First, Proposition \ref{ba3} implies that we may find two vertices $W_1, W_2$ such that $\phi_f(W_1)\subset B_{\varepsilon}(\xi,r)
\subset\phi_f(W_2)$ and $$r\asymp e^{-\varepsilon |W_1|}\asymp  e^{-\varepsilon |W_2|}\,.$$

From the Lemma of the shadow (Lemma \ref{shagen}) and Axiom [Degree] follows 
$$ r^{\alpha}\asymp e^{-\varepsilon\alpha |W_1|}\lesssim \mu_f(B_{\varepsilon}(\xi,r))\lesssim e^{-\varepsilon\alpha |W_2|}\asymp r^{\alpha}\,.$$

The fact that the entropy is $\log d$ follows from Theorem \ref{maxent} and Axiom [Degree]:
since the degree is bounded along any pull-back, it follows that, for any $\xi\in\partial \Gamma$, we have
$$\lim\frac{1}{n}\log d_{F^n}(\xi)=0.$$
Hence Birkhoff's ergodic theorem (Thm \ref{thm:birkhoff}) implies that $$\int \log d_F d\mu_f=0$$ so that
$d_F(\xi)=1$ for $\mu_f$-a.e. every $\xi$.

If $f$ is metrically cxc, then the conjugacy $\phi_f$ is quasisymmetric, by Theorem \ref{cxcimphyp}, and quasisymmetric maps preserve the property of being doubling \cite[Chap.\,15]{heinonen:analysis}; this proves the last assertion. 

The proof of the uniqueness of $\mu_f$ occupies the rest of this section.
\qed

\gap

\begin{prop}\label{porous}
The critical set $C(F)$  and the set of critical values $V(F)=F(C(F))$ are porous i.e., 
there is some constant $c$, such that, any ball $B_\varepsilon(r)\subset \partial \Gamma$, $r\le \diam_\varepsilon \partial \Gamma$,
contains a ball of radius at least $c\cdot r$ disjoint from $C(F)$, or $V(F)$.\end{prop}

{\noindent\bf Proof of Proposition \ref{porous}:}
 Let us first prove that $C(F)$ is porous. We will use the fact that the critical set 
is nowhere dense.
If not, there would be a sequence of balls $(B(\xi_n,r_n))_n$ such that any ball of
radius $r_n/n$ contains critical points. 

\bigskip

Since $\Gamma$ is doubling, the sequences of pointed metric spaces 
$$(\overline{B(\xi_n,r_n)},\xi_n, (1/r_n)d_\varepsilon)) \quad \mbox{and} \quad 
(\overline{B(F(\xi_n),e^\varepsilon r_n)},\xi_n, (1/r_n)d_\varepsilon))$$ 
is compact in the Hausdorff-Gromov topology \cite{gromov:expanding}.
Hence, one can extract  convergent subsequences. Since the restrictions of $F$ are uniformly Lipschitz, of bounded multiplicity
and onto (cf. Lemma \ref{boulesF}), it is the case of any limit so the lemma above implies that any limit has a nowhere dense critical set which
yields a contradiction. Therefore, the critical set is porous. 

\bigskip

To see that the set of critical values is also porous, we pick a ball $B$ and write $F^{-1}(B)=B_1\cup B_2\ldots\cup B_k$, with
$k\le d$ (cf. Lemma \ref{boulesF}).  
We  first consider a ball $B_1'\subset B_1$ disjoint from $C(F)$ of definite size. Therefore $F(B_1')$ has definite size in $B$.
Let us pull it back in $B_2$ and define $B_2'$ inside this new ball with definite size.   By induction, since $k$ is bounded by $d$, one constructs in this way a ball $B'\subset B$ of size comparable to that of $B$ which is disjoint from the set of critical values.\qed

\bigskip

{\noindent\bf Topological entropy revisited.} 
We have proved that $\partial \Gamma$ is Ahlfors regular of dimension $\alpha=(1/\varepsilon)\log d$. 
We adapt the estimate of M.\,Gromov on the topological entropy to our setting
\cite{gromov:entropy} to obtain an upper bound for the relative entropy. Let $n\ge 1$ and endow $(\partial \Gamma)^n$ with 
the  metric $|(\xi_j)-(\zeta_j)|_\varepsilon=\max\{|\xi_j-\zeta_j|_\varepsilon\}$. Let $\pi_j:(\partial \Gamma)^n\to \partial \Gamma$
be the canonical projection to the $j$th factor and let us define 
$I_n:\partial \Gamma\to (\partial \Gamma)^n$ by 
$I_n(\xi)=(\xi,F(\xi),\ldots,F^{n-1}(\xi))$.
Set $\Gamma_n=I_n(\partial \Gamma)$, and let $H$ be the Hausdorff measure of dimension 
$\alpha$.
We let $H^r(Y)=\inf \left\{\sum (\diam U_j)^\alpha\right\}$, where the infimum is taken over all the coverings of a given set 
$Y$ by sets
$(U_j)$ of
diameter at most $r$.

\bigskip

Given a subset $Y\subset \partial \Gamma$ et $\eta\ge 0$, we let $I_n(Y)_\eta$ be the $\eta$-neighborhood of $I_n(Y)$ in $\Gamma_n$,
and we  define $$\begin{array}{l} 
\mbox{lov}(F|Y,\eta) =\dis\limsup_{n\to\infty} \dis\frac{1}{n}\log H(I_n(Y)_\eta)\\ \\
\mbox{lodn}(F|Y) =\dis\liminf_{r\to 0}\dis\liminf_{n\to\infty} \dis\frac{1}{n}\log \dis\inf_{\xi\in Y} H((B(I_n(\xi),r)\cap\Gamma_n)\end{array}$$
We set $\mbox{lov}(F)=\mbox{lov}(F|\partial \Gamma)$ and $\mbox{lodn}(F)=\mbox{lodn}(F|\partial \Gamma)$.

\bigskip

The main observation is that, for any $\xi\in \partial\Gamma$, $B(I_n(\xi),r)\cap \Gamma_n= I_n(S(\xi,n,r))$.
Let us fix $Y\subset \partial \Gamma$ and a maximal family $(S_j)_{1\le j\le s_n(Y,r)}$ of disjoint dynamical balls $S(\cdot,n,r)$
centered in $Y$. 
We fix $\eta>r$. 
It follows that $$H(I_n(Y)_\eta)\ge \sum_{1\le j\le s_n(Y,r)} H(I_n(S_j))\ge s_n(Y,r)\cdot\inf_{\xi\in Y} H((B(I_n(\xi),r)\cap\Gamma_n)\,.$$
Therefore, we obtain the formula $$h_{top}(F|Y)\le \mbox{lov}(F|Y,\eta)-\mbox{lodn}(F|Y)\,.$$

\bigskip

In our setting, for any point $\xi\in \partial \Gamma$, we have $H(B(I_n(\xi),r))\ge H(\pi_n B(I_n(\xi),r))$ 
since the projection decreases distances. 
But $\pi_n B(I_n(\xi),r)\supset B(F^{n-1}(\xi),r)$
so that $H(B(I_n(\xi),r))\ge H( B(F^{n-1}(\xi),r)\gtrsim r^\alpha$ since $\partial \Gamma$ is Ahlfors regular. Thus $\mbox{lodn}(F) =0$.

\bigskip

On the other hand, Proposition \ref{boulesF} implies that $I_n(B(\xi,r))\subset B(I_n(\xi),e^{\varepsilon (n-1)}r)$. Let $\delta>0$.
Let us cover $Y_\eta$ by sets $\{U_j\}_{j\in J}$ of diameter at most $r$ so that $\sum (\diam U_j)^\alpha\le H^r(Y_\eta)+\delta$.
It follows that we may pick points $\xi_j\in U_j$ so that $I_n(Y)_\eta$ is covered by the balls $B(I_n(\xi_j),e^{\varepsilon (n-1)}\diam U_j )$. Hence
$$H^{re^{\varepsilon (n-1)}}(I_n(Y)_\eta)\le 2^\alpha e^{\varepsilon\alpha (n-1)} (H^r(Y_\eta)+\delta)\le 2^\alpha d^{(n-1)} (H^r(Y_\eta)+\delta)\,.$$
Letting $\delta,r\to 0$, it follows that, for any set $Y\subset \partial \Gamma$ and any $n\ge 1$, $$H(I_n(Y)_\eta)\lesssim d^n H(Y_\eta)\,.$$
This shows that $\mbox{lov}(F|Y,\eta)\le \log d$ and the estimate on the relative entropy becomes $$h_{top}(F|Y)\le \log d\,.$$

\bigskip

{\noindent\bf Uniqueness of the measure of maximal entropy.} It follows from the ergodicity of $\mu_f$ and Theorem \ref{maxent}
that $\mu_f$ is the unique invariant measure of constant Jacobian. To prove that $\mu_f$ is the unique measure of maximal entropy, it
is enough to prove that if $\nu$ is an ergodic invariant measure with non-constant Jacobian, then $h_\nu(F)<\log d$.

\bigskip

We will  adapt the argument of M.\,Lyubich \cite{lyubich:entropy} following the ideas of J.Y.\,Briend and J.\,Duval \cite{briend:duval:deux}.
Let us first recall that if $Y$ has positive $\nu$-measure, then $h_\nu(F)\le h_{top}(F|Y)$ (cf. Lemma 7.1 in \cite{lyubich:entropy}).

\bigskip

If the measure $\nu$ charges the set of critical values $V(F)$, then, since this set is porous in an Ahlfors regular set
(Proposition \ref{porous}), its
box dimension is strictly smaller than $\alpha$ so that $h_{\nu}(F)\le e^{\varepsilon} \dim V(F) <\log d$.
 Therefore $\nu$ does not have maximal entropy. 

\bigskip

We assume from now on that the measure does not charge the set of critical values $V(F)$. For any point $\zeta \not\in V(F)$,
there is some radius $r_\zeta >0$ such that $F^{-1}(B(\zeta ,r_\zeta ))$ 
is the union of $d$ balls of radius $r_\zeta  e^{-\varepsilon}$.
We may assume that the measure  of the boundary is null. 
Let us extract a locally finite subcovering of $\partial \Gamma\setminus V(F)$.
This yields a measurable partition $\PPP=\{P_j\}$ of $\partial \Gamma$ such that, for any piece $P$, there are $d$ preimages 
$Q_1(P),\ldots, Q_d(P)$
such that $\diam Q_j(P)= e^{-\varepsilon}\diam P$ and $F|_{Q_j}:Q_j\to P$ is a homeomorphism. We label these preimages so that
$\nu(Q_j)\ge \nu(Q_{j+1})$, and we define $U_j=\cup_{P\in\PPP} Q_j(P)$. It follows that $\nu(U_1)>1/d$ since $\nu$ has non constant
Jacobian. 

\bigskip

Furthermore, $\sum\nu(U_j)=1$, so $\nu$-a.e. point has an itinerary defined by its visits to the sets $U_j$\,: for $\nu$ almost
every point $\xi$, for every iterate $n$, there is some index $j(n)$ such that $f^n(x)\in U_{j(n)}$. We let
$U_j^n=\{\xi\in U_j, B(\xi,1/n)\subset U_j\}$. We fix $N$ large enough so that $\nu(U_1^N)>1/d$, and we write $O=U_1^N$.

\bigskip

Let us define by $r_n(\xi)$ the number of iterates $0\le k\le n-1$ such that $f^k(\xi)\in O$. It follows from Birkhoff's ergodic 
theorem that $r_n(\xi)/n$ tends to $\nu(O)$ for almost every $\xi$. Let us fix $\sigma\in (1/d,\nu(O))$. By Egoroff's theorem, there is some
$m\ge 1$ such that the set $$Y=\{\xi\in \partial \Gamma,\ r_n(\xi)\ge \sigma n,\ \forall\,n\ge m\}$$ has positive $\nu$-measure. 
We will show that
$\mbox{lov}(f|Y,{(1/2N)})<\log d$\,: this will then imply by the remark above that $h_\nu(f)\le \mbox{lov}(f|Y,{(1/2N)})<\log d$.

\bigskip

Given $J=(j_1,\ldots,j_n)\in\{1,\ldots,d\}^n$, we let $U_J=\prod U_{j_i}$ and $\Gamma_n(J)=\Gamma_n\cap U_J$.
Let $\Sigma_n$ be the set of itineraries such that the number of occurrences of $1$'s is at least $\sigma n$. M.\,Lyubich has shown
that the cardinality of $\Sigma_n$ is bounded by $d^{\rho n}$ where $\rho<1$ (cf.\,Lemma 7.2 in \cite{lyubich:entropy}).

\bigskip

By definition, $I_n(Y)_{(1/2N)}\subset \cup_{J\in\Sigma_n}\Gamma_n(J)$, so that 
$$H(I_n(Y)_{(1/2N)})\le \sum_{J\in\Sigma_n}H(\Gamma_n(J)\cap I_n(Y)_{(1/2N)})\,.$$

\bigskip

Let us cover $\overline{Y_{(1/2N)}}$ by finitely many balls $B(\zeta ,r_\zeta  )$. 
It follows that $r_\zeta \ge r_0>0$ for these particular balls
and for some $r_0>0$.
Thus, if $U\subset Y_{(1/2N)}$ has diameter at most $r_0/2$, then we may define $d$ inverse branches $F_{j}$ so
that $\diam F_j(U)= e^{-\varepsilon}\diam U$.

\bigskip

Let $r\in (0,r_0/2)$, and let us cover 
$\pi_n(I_n(Y)_{(1/2N)}\cap\Gamma_n(J))$ by sets $E_j$ of diameter at most $r$ such that 
$$\sum (\diam E_j)^\alpha\le H^r(\pi_n(I_n(Y)_{(1/2N)}\cap\Gamma_n(J)))+\eta\,.$$
By construction, $F^{n-1}|_{\pi_1(\Gamma_n(J))}$ is injective. It follows that, for any $\ell=1,\ldots,n-1$, 
$$\diam f^{-\ell}(E_j)\cap\pi_{n-\ell}(\Gamma_n(J))=e^{-\varepsilon\ell}\diam E_j$$
so that we may pick points $\xi_j\in Y_{(1/2N)}\cap f^{-(n-1)}(E_j)$ such that $\Gamma_n(J)\cap I_n (Y)_{(1/2N)}\subset \cup_j B(I_n(\xi_j),r_j)$
with $r_j=\diam E_j$. Therefore, 
$$H^r(\Gamma_n(J)\cap I_n(Y)_{(1/2N)})\le 2^\alpha ( H^r(\pi_n(\Gamma_n(J)\cap I_n(Y)_{(1/2N)} ))+\eta )\lesssim 1$$
so that $H(I_n(Y)_{(1/2N)})\lesssim d^{\rho n}$ and $$h_\nu(f)\le \mbox{lov}(f|Y,{(1/2N)})\le \rho\log d <\log d\,.$$
This establishes the uniqueness of $\mu_f$ as a measure of maximal entropy.

\subsection{BPI-spaces} Following David and Semmes \cite{david:semmes:dreams}, a bounded space $(X,d,\mu)$ is called BPI (``Big pieces
of itself'') 
\index{Index}{BPI space}%
if $X$ is Ahlfors regular of dimension $\alpha$, and if the following homogeneity condition holds.  There are  constants $\theta<1$  and $C>1$
such that, given any balls $B(x_1,r_1)$ and $B(x_2,r_2)$ with $r_1,r_2\le\diam X$,  there exists a closed set $A\subset B(x_1,r_1)$
with $\mu(A)\ge \theta r_1^\alpha$ and an embedding $h:A\to B(x_2,r_2)$ such that $h$ is a $(C,r_2/r_1)$-quasisimilarity, i.e. $$C^{-1} \le \frac{|h(a)-h(b)|}{(r_2/r_1)|a-b|}\le C$$ for all $a, b \in A$.  

\begin{thm}\label{thm:BPI} Under the hypotheses of Theorem \ref{cxcimphyp}, the metric space $\bdry_\varepsilon\Gamma$ is BPI.\end{thm}

\pf We start with a preliminary step.

Suppose $\phi_f: X \to \bdry_\varepsilon\Gamma$ is the conjugacy given by Theorem \ref{thm:construction} and $d_\varepsilon$ is the metric on $\bdry_\varepsilon \Gamma$. 
 Let $d_{\varepsilon, X} = \phi_f^*(d_\varepsilon)$.   For convenience of notation, we will show $(X,  d_{\varepsilon, X})$ is BPI.

Recall that since $f$ is topologically cxc, there is a uniform (in $n$) upper bound $p$ on the degree $d(U)$ by which an element of $U\in\UUU_n$ maps under $f^n$.   
Choose $W \in \UUU_{n_0}$ arbitrarily so that the multiplicity $d(W)$ is maximal, so that any further preimages $\wtW$ of $W$ map onto $W$ by degree one i.e., are homeomorphisms. 
 It follows from Proposition \ref{ba3} that its image under $\phi_f$ contains some ball $B_{\varepsilon}(\xi,4r)$, such that, for any iterate $n$, any $\txi\in F^{-n}(\xi)$,
$F^n:B_{\varepsilon}(\txi,4re^{-\varepsilon n})\to B_{\varepsilon}(\xi,4r)$ is a homeomorphism. Therefore, Proposition \ref{prop:tentlike} shows that
$F^n:B_{\varepsilon}(\txi,re^{-\varepsilon n})\to B_{\varepsilon}(\xi,r)$ is a $(1,e^{\varepsilon n})$-quasisimilarity.

\gap

By Proposition \ref{prop:repellors_are_fractal}, for each $U_0 \in \UUU_0$, there exists $k \in \N$ and some $\wtW \in \UUU_{n_0+k}$ such that 
\bi
\item $\cl{\wtW} \subset U_0$,
\item $\wtW$ is a preimage of $W$ under $f^k$, and 
\item $\deg(f^k: \wtW \to W)=1$.  
\ib
 Since $\UUU_0$ is finite,  the $\wtW$'s considered above have a level bounded by some $n_0+k_0$.

Furthermore, for any $n$ and any $ U\in\UUU_n$, one has $f^n(U)\in\UUU_0$, so one may find a preimage  $W_U$ of $W$ so that 
$\cl{W_U}\subset U$, and $|W_U|=n +O(1)$. Thus, one can find a ball $B(\xi',re^{-\varepsilon (n+k)})\subset W_U$
so that $f^{n+k}:B(\xi',re^{-\varepsilon (n+k)})\to B(\xi,r)$ is a $(1,e^{-\varepsilon (n+k)})$-quasisimilarity.
Let us note that $re^{-\varepsilon (n+k)}\asymp \diam_\varepsilon U$. 

 Now suppose we are given $d_{\varepsilon,X}$-balls $B_i = B(\xi_i, r_i) \subset X$, $i=1,2$.  By Proposition \ref{ba3}, 
there exist $U_i', U_i \in \mathbf{U}$ with 
 \[ U_i' \intersect X  \subset B_i \subset U_i \intersect X\]
 such that $n_i = |U_i'| = \frac{1}{\varepsilon}\log \frac{1}{r_i} + O(1)$.  For each $i=1,2$, let $W_i$ be  a preimage  of $W$ so that 
$\cl{W_i}\subset U_i'$, and $|W_i|=n_i +k_i$ as in the previous paragraph.
Moreover, we consider balls $B(\xi'_i,re^{-\varepsilon (n_i+k_i)})\subset W_i$ as above. It follows from Ahlfors-regularity that $\mu_f(B(\xi'_i,re^{-\varepsilon (n_i+k_i)}))\asymp \mu_f(B_i)$. 
Let 
$h_i$ be the restriction of $f^{n_i+k_i}$ to the ball $B(\xi'_i,re^{-\varepsilon (n_i+k_i)})$,
 for $i=1,2$; the map
$h=h_2^{-1}\circ h_1$ is a quasisimilarity between big pieces of $B_1$ and $B_2$.
\qed

\chapter{Expanding non-invertible dynamics}

In this chapter, we give different classes of expanding, non-invertible, topological dynamical systems to which we may apply the theory developed in earlier chapters. We first compare our notion of cxc with classical conformal dynamical systems
on compact Riemannian manifolds\,: coverings on the circle  (\S\,\ref{sec:S1}), rational maps on the Riemann sphere (\S\,\ref{secn:semihyperbolic}) and
uniformly quasiregular mappings in higher dimension (\S\,\ref{secn:lattes}). We also provide other examples of cxc maps for which the conformal
structure is not given a priori. These classes  come from finite subdivision rules (\S\,\ref{secn:fsr}) and from expanding
maps on manifolds (\S\,\ref{secn:expanding_maps}). In \S\,\ref{secn:periodic_branch}, 
we provide two examples of maps which satisfy the axioms
[Expansion] and [Irreducibility], but not [Degree].   
We conclude the chapter in \S\,4.7 by comparing and contrasting our constructions with formally similar ones arising in $p$-adic dynamics.

\section{No exotic cxc systems on $\mathbf{\IS^1}$}   
\label{sec:S1}

Metric cxc systems on the Euclidean circle $\IS^1$ include the covering maps $z \mapsto z^d$, $|d| \geq 2$, and essentially nothing else.  

\begin{thm}[Cxc on $S^1$ implies qs conjugate to $z^d$]
Suppose $f: X \to X$ is a metric cxc dynamical system where $X$ is homeomorphic to 
$\IS^1$.   Then there exists a quasisymmetric homeomorphism $h: X \to \IS^1$ conjugating $f$  to the map $z \mapsto z^{\deg f}$.

\end{thm}

\pf An open connected subset of $\IS^1$ is an interval.  
Since $f$ is open, it 
sends small open intervals onto small open intervals.  
Moreover, if these intervals
are small enough, $f$ must be injective on such intervals, else 
there is a turning
point in the graph and openness fails.  Hence $f$ is a local homeomorphism.  
A local
homeomorphism on a compact space is a covering map, see
\cite[Thm. 2.1.1]{aoki:hiraide:topological}.  
In particular $f$ is strictly monotone.

Such a map admits a monotone factor map $\pi$ onto $g(z)=z^d$
where $d=\deg f$ \cite[Prop. 2.4.9]{hasselblatt:katok:dynamics}.  If $\pi$ is not a
homeomorphism, then there is an interval 
$I \subset \pi^{-1}(x)$ for some $x \in
\IS^1$.  Axiom [Irreducibility] implies  $f^N(I) = \IS^1$
for some $N$.  
Then 
\[ g^N(x) = g^N(\pi(I)) = \pi(f^N(I)) = \pi(\IS^1)= \IS^1\]
which is impossible.  
Thus $\pi$ is a homeomorphism and $f$ is topologically conjugate
to $g$.  Since $g$ is cxc with respect to the Euclidean metric, 
$\pi$ is
quasisymmetric, by Theorem \ref{thm:invariance_of_cxc}.  
\qed

\section{Semi-hyperbolic rational maps}
\label{secn:semihyperbolic}

We endow the Riemann sphere $\cbar$ with the spherical metric, and we will talk of disks $D(x,r)$\index{Symbols}{$D(x,r)$} rather than balls 
$B(x,r)$ in this context.

If $g$ is a rational map, the {\em Fatou set} $F(g)$
\index{Symbols}{$F(g)$}%
\index{Index}{Fatou set}%
 of $g$ is the set of points  $z\in\cbar$ having  a neighborhood $N(z)$ such 
that the set of restrictions of iterates $\{(f^n)|_{N(z)}\}_n$ forms a normal family. The {\em Julia set} $J(g)$ 
\index{Symbols}{$J(g)$}%
\index{Index}{Julia set}%
of $g$ is the complement of $F(g)$.
We shall say that $g$ is {\em chaotic}
\index{Index}{chaotic}%
  when $J(g)=\cbar$. 

The class of semi-hyperbolic rational maps
\index{Index}{semi-hyperbolic}%
\index{Index}{rational map, semi-hyperbolic}%
 has been introduced by L.\,Carleson, P.\, Jones and J.-C.\,Yoccoz in \cite{carleson:jones:yoccoz}.
In their paper, they provide several different characterizations, some of which we now recall.

\begin{thm}[definition of semi-hyperbolic rational maps]
\label{thm:def_semihyperbolic} 
Let $g$ be a rational map. The following conditions are equivalent and
define the class of semi-hyperbolic rational maps.
\be 
\item A radius $r>0$ and a maximal degree $p<\infty$ exist, such that, for any $z\in J(g)$, for any iterate $n\ge 1$
and any connected component $W_n$ of $g^{-n}(D(x,r))$, the degree of $g^n|_{W_n}$ is at most $p$.
\item A radius $r >0$,  a maximal degree $p <\infty$, and constants $c>0$ and $\theta<1$ exist such that, 
for any $x\in J(g)$, any
iterate $n\ge
1$, and
any component $W_n$ of $g^{-n}(D(x,r))$, the degree of the
restriction 
of
$g^n$ to
$W_n$ is at most $p$, and the diameter of $W_n$ is at most $c\theta^n$.
\item  The map $g$ has neither recurrent
critical point in the Julia set nor parabolic cycles.
\item A maximal degree $p_0$ exists such that, for any $r>0$  and any $x\in J(g)$, if we let $n$ be the least
iterate such $g^n(D(x,r)\cap J(g))=J(g))$ then $g^n|_{D(x,2r)}$ has degree at most $p_0$.
\eb
\end{thm}

We refer to \cite[Theorem 2.1]{carleson:jones:yoccoz} for the proofs of the equivalence above.

\begin{cor}[Topological cxc rational maps are semi-hyperbolic]\label{cor:topcxc_shyp}
A rational map is topological cxc if and only if it is semi-hyperbolic.\end{cor}

\pf Assume that $f$ is a topological cxc rational map. Then it satisfies conclusion (1) of Theorem \ref{thm:def_semihyperbolic}
with radius the Lebesgue number of the cover $\UUU$.

\gap 

Conversely, the classification of stable domains \cite[Chapter 16]{milnor:dynamics} implies that the complement of 
the Julia set $J(f)$ consists of points which converge to attracting cycles under iteration.  
Thus if $\XX_0$ is the complement of a suitable neighborhood of attacting cycles and their preimages, 
then $\XX_1 = f^{-1}(\XX_0)$ has closure in $\XX_0$ and the branch points of $f: \XX_0 \to \XX_1$ lie in $J(f)$.  

The axiom [Irreducibility] holds for any rational map 
in a neighborhood of its Julia set \cite[Thm. 4.10]{milnor:dynamics}.  

If $f$ is semihyperbolic, Theorem \ref{thm:def_semihyperbolic}  asserts that 
there
exist an $r >0$, $p <\infty$, $c>0$, and $\theta<1$  such that, 
for any $x\in J(f)$, any
iterate $n\ge
1$, and
any component $W_n$ of $f^{-n}(D(x,3r))$, the degree of the
restriction 
of
$f^n$ to
$W_n$ is at most $p$, and the diameter of $W_n$ is at most $c\theta^n$. 
We let $\UUU_0$ be a finite subcovering of $J(f)$
of
$\{D(x,r),\ x\in J(f)\}$,
and $\UUU_n$ 
be the set of components of $f^{-n}(U)$ when $U$ ranges
over
$\UUU_0$.
Therefore, Axioms [Degree] and [Expansion] hold, so that $f$ is a topological cxc map.\qed

\gap

Condition (4) of Theorem \ref{thm:def_semihyperbolic} says that, just as for convex cocompact groups), any point in the Julia set of a semihyperbolic rational map is
{\em conical}.
\index{Index}{conical point}%
That is, one may use the dynamics to go from
small scales to large scales and vice versa with bounded distortion (cf. Lemma \ref{lemma:pvalent} below).  Indeed, Lyubich and Minsky \cite{lyubich:minsky:lamination} call semihyperbolic maps {\em convex cocompact} and show that such maps are characterized by the following property:  the quotient (by the induced invertible dynamics of $f$) of the convex hull of the ``Julia set'' (the hull taken in their affine hyperbolic three-dimensional lamination associated to $f$) is compact.   

The aim of this section is first to prove that these maps are metrically cxc (Theorem \ref{thm:semihyperbolic_are_cxc}) 
and also to strenghten their
relationship to convex cocompact Kleinian groups within the dictionary by establishing new characterisations
of this class. Theorems \ref{thm:semihyperbolic_are_cxc}, \ref{thm:nothing_new}, and \ref{thm:car} below imply the following.

\begin{thm}[Characterizations of semi-hyperbolic rational maps]
\label{thm:newchar_semihyperbolic} 
Let $g$ be a rational map. The following propositions are equivalent.
\be 
\item $g$ is semi-hyperbolic.
\item $g$ is metric cxc on its Julia set, with respect to the spherical metric.
\item  There is a covering $\UUU$ of $J(g)$ such that the associated graph $\Gamma$ is quasi-isometric to
the convex hull of $J(g)$ in $\IH^3$ by a quasi-isometry which extends to $\phi_g:J(g)\to \partial\Gamma$.
\eb
\end{thm}

The map $\phi_g$ in the statement above is the one defined by Theorem \ref{thm:construction}.
The symbol $\IH^3$
\index{Symbol}{$\IH^3$}%
 denotes hyperbolic three-space; see \S 4.2.2.

\gap

The last subsections deal with the topological characterization of semi-hyperbolic maps in the
spirit of Cannon's conjecture for hyperbolic groups, which claims that a hyperbolic group $G$ with 
a topological $2$-sphere as boundary admits a faithful cocompact Kleinian action.

\subsection{Characterization of cxc mappings on the standard $2$-sphere}

\begin{thm}[semi-hyperbolic rational maps are cxc]
\label{thm:semihyperbolic_are_cxc} 
Let $f$ be a semihyperbolic rational map with Julia set $J(f)$.  
Then there are closed neighborhoods $\XX_0, \XX_1$ of $J(f)$  such that, in the spherical metric,  $f: \XX_1 \to \XX_0$ is metrically cxc with repellor $J(f)$ 
with respect to a finite collection $\UUU_0$ of open spherical disks.   
\end{thm}

By the Riemann mapping theorem, a simply-connected domain $V \subset \rs$ which is neither the whole sphere, nor the whole sphere with one point removed, is conformally isomorphic to the unit disk and so carries a unique hyperbolic metric $\rho_V$ of curvature $-1$.  We call such a domain a {\em simply-connected hyperbolic domain}.  

\gap
\noindent{\bf Notation.}  Let $\sigma = |dz|/(1+|z|^2)$ denote the spherical Riemannian metric on $\rs$.  For a simply-connected hyperbolic domain $V$ in $\rs$, let $\rho_V$ denote the hyperbolic metric on $V$ and $V^c$ its complement in $\rs$.   Given a metric $g$, we denote by $B(a, r; g)$ the ball of radius $r$ about $a$ and by $\diam(A; g)$ the diameter of a set $A$.  
\gap

\noindent{\bf Univalent functions.}   
\index{Index}{univalent functions} %
See \cite[Chap.\,5]{ahlfors:book:invariants}.  Let $\D$ denote the unit disk
\index{Symbols}{$\D$}%
 centered at the origin in $\C$, and let 
\[ \SSS = \{ f: (\D, 0) \to (\C, 0)\} : f \ \mbox{ is 1-1 and analytic}\}\]
denote the class of so-called {\em Schlicht functions}.  In the topology of local uniform convergence, $\SSS$ is compact.  This implies that restrictions to smaller balls are uniformly bi-Lipschitz, a fact known rather loosely as the {\em Koebe distortion principle}.
\index{Index}{Koebe distortion principle}%
  More precisely:  for all $0 < r < 1$, there is a constant $C(r)>1$ such that for all $|z|, |w| \leq r$ and all $f \in \SSS$,
\begin{equation}
\label{eqn:schlicht}
\frac{1}{C(r)} \leq \frac{|f(z)-f(w)|}{|z-w|} \leq C(r).
\end{equation}
The Koebe principle also implies that for any simply-connected hyperbolic domain $V \subset \C$,  
\begin{equation}
\label{eqn:like_oneoverd}
\rho_V (w)\asymp \frac{|dw|}{\dist(w, V^c)}.
\end{equation}
Finally, we will need the {\em Schwarz-Pick lemma}:  if $f: \wtV \to V$ is a proper, holomorphic map between hyperbolic domains, 
then with respect to the metrics $\rho_{\wtV}$ and $\rho_V$, for all tangent vectors $v$, 
\begin{equation}
\label{eqn:schwarz}
||df(v)||/||v|| \leq 1\ \mbox{with equality if and only if $f$ is an isometry}.
\end{equation}
While essentially classical, the specific versions given above may be found in  
\cite[Chapter 2]{ctm:renorm}.

\begin{lemma}[Comparing metrics]
\label{lemma:comparing_metrics}
There exists a universal constant $C$ such that the following holds.  Let $W \subset \rs$ be a simply-connected hyperbolic domain of spherical diameter $<\pi/4$, let $x \in W$, and set $D=B(x, 2; \rho_W)$.  Then, restricted to the domain $D$, the metrics $\rho_W$ and $\sigma/\diam(D; \sigma)$ are $C$-bi-Lipschitz equivalent.  
\end{lemma}

\pf By applying a rigid spherical rotation we may assume $W$ is contained in $\D$.  For such domains, by compactness, $\sigma$ and the Euclidean metric $|dw|$ are bi-Lipschitz equivalent.  By (\ref{eqn:like_oneoverd})
we have 
\[ \rho_W = \frac{1}{\dist(w, W^c)} |dw|\]
where $\dist(w, W^c)$ denotes the Euclidean distance from $w$ to the complement $W^c$ of $W$ in $\C$.  Suppose $\phi: (\D, 0) \to (W, x)$ is a holomorphic isomorphism.  By the Koebe principle (\ref{eqn:schlicht})   and the fact that $D$ is a hyperbolic ball of radius $2$, hence precompact, 
\begin{equation}
\label{eqn:diamDsigma}
\diam(D; \sigma) \asymp \diam(D; |dw|) \asymp |\phi'(0)|.
\end{equation}

Let $w \in D$.  The Schwarz-Pick lemma (\ref{eqn:schwarz}) implies $\dist(w, W^c) \leq \mbox{const}\cdot |\phi'(0)|$, and the Koebe principle (\ref{eqn:schlicht})  implies $\dist(w, W^c) \geq \mbox{const}\cdot |\phi'(0)|$, so that 
\begin{equation}
\label{eqn:distwwc}
\dist(w, W^c) \asymp |\phi'(0)|.
\end{equation}
Dividing (\ref{eqn:diamDsigma}) by (\ref{eqn:distwwc}) yields 
\[ \frac{\diam(D;\sigma)}{\dist(w, W^c)} \asymp 1\]
and so 
\[ \frac{\rho_W(w)}{\sigma(w)/\diam(D;\sigma)} \asymp \frac{\diam(D;\sigma)}{\dist(w, W^c)} \asymp 1.\]
\qed

While the Koebe distortion principle applies to univalent maps, there are variants for proper, noninjective maps as well.  See Figure 4.1
\begin{figure}
\label{fig:Koebe4P}
\psfragscanon
\psfrag{1}{$\tilde{w}$}
\psfrag{2}{$C_h(p,r)$}
\psfrag{4}{$\widetilde{W}$}
\psfrag{5}{$r$}
\psfrag{6}{$W$}
\psfrag{7}{$f$}
\psfrag{8}{$w$}
\psfrag{9}{$f(B(\tilde{w},r))$}
\psfrag{10}{$C_h^{-1}(p,r)$}
\begin{center}
\includegraphics[width=5in]{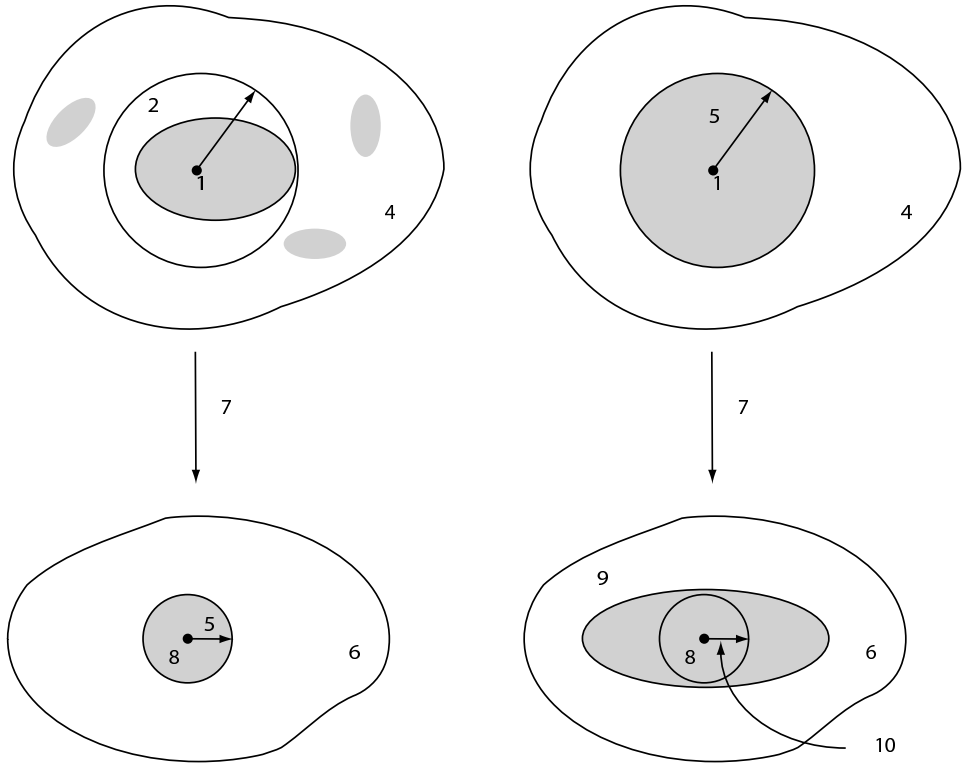}
\end{center}
\caption{{\sf At left:  the hyperbolic diameter of the component of $f^{-1}(B(w,r))$ containing $\tilde{w}$ cannot be too large.  At right:  the image of the hyperbolic ball of radius $r$ cannot be too skinny. }}
\end{figure}

\begin{lemma}[Distortion of $p$-valent maps]
\label{lemma:pvalent}
For $p \in \N$ and $\tr, r >0$, there exist real-valued functions $C_h(p,r)$ and $C^{-1}_h(p,\tr)$, tending to zero as $r, \tr$ tend to zero, with the following property.  Suppose $\wtW, W \subset \rs$ are hyperbolic simply-connected domains, $f: \wtW \to W$ is a proper, holomorphic map such that $\#f^{-1}(w) \leq p$ for all $w \in W$, and $f(\tw) = w$.  
\be
\item Let $B=B(w, r; \rho_W) \subset W$ and let $\wtB$ be the component of $f^{-1}(B)$ containing $\tw$.  Then 
\be
\item $B(\tw, r; \rho_{\wtW}) \subset \wtB \subset B(\tw, C_h(p,r); \rho_{\wtW})$.
\item If $B$ is replaced by an open connected set, then 
\[ \diam(B; \rho_W) \leq \diam(\wtB; \rho_{\wtW}) \leq C_h(p, \diam B).\]
\eb
\item Given $\tr>0$, 
\[ B(w, C_h^{-1}(p, \tr); \rho_W) \subset f(B(\tw, \tr; \rho_{\wtW} )) \subset B(w, \tr; \rho_W).\]
\eb
\end{lemma}
Basically, the above lemma says that for connected sets of a fixed size, preimages cannot be too large or too skinny, and images cannot be too large or too skinny.  

\pf 1(a) is the content of  \cite[Lemma 2.2]{carleson:jones:yoccoz} and implies the lower containment in (2) and the upper bound in 1(b); see also \cite{shishikura:tan:mane}.  The lower bound in 1(b) and the upper containment in (2) follow from the Schwarz-Pick lemma (\ref{eqn:schwarz}).  
\qed

\pf  (Theorem \ref{thm:semihyperbolic_are_cxc})  Suppose $f: \rs \to \rs$ is semihyperbolic.  Let  $r, p$ be the constants as in Theorem \ref{thm:def_semihyperbolic}.  Let $\XX_0$ be the complement of a forward-invariant neighborhood of the attractors as in the proof of  Corollary \ref{cor:topcxc_shyp} and $\XX_1$ its preimage, so that $f: \XX_1 \to \XX_0$ is an fbc satisfying the conditions at the start of \S 2.2.   To define the level zero good open sets $\UUU_0$ we proceed as follows.
\gap

\noindent{\bf Definition of $\UUU_0$.}  For $x \in J(f)$,
let $W(x)$ be the spherical ball whose radius is $r/2$.  By Lemma \ref{lemma:pvalent} there exists $r_0$ so small that  $C_h(p, r_0) < 1/2$ and let $U(x) = B(x, r_0; \rho_{W(x)})$.
Let $\UUU_0$ be a finite open cover of $J(f)$ by pointed sets of the form $(U(x), x)$.  Then we have a finite set of triples $(W(x), U(x), x)$.  By taking preimages, we obtain for each $n \in \N$ a covering $\UUU_n$ of $J(f)$ by Jordan domains $\wtU$ such that each has a preferred basepoint $\tx$ and is  compactly contained in a larger domain $\wtW$.  Moreover, 
\[ f^k: (\wtW, \wtU, \tx) \to (W, U, x)\]
whenever $U \in \UUU_n$, $\wtU \in \UUU_{n+k}$, and $f^k(\tx)=x$.  Note that by construction 
and Lemma \ref{lemma:pvalent}, 
for all $n$ and all $U \in \mathbf{U}=\union_n \UUU_n$ with basepoint $x$, 
\[ B(x, r_0; \rho_{W}) \subset U \subset B(x, 1/2; \rho_{W}).\]
In particular, 
\begin{equation}
\label{eqn:r0}
2r_0 \leq \diam(U; \rho_{W}) \leq 1.
\end{equation}  
\gap

\noindent{\bf Diameter distortion.}  Suppose $ f^k: (\wtU, \wtU') \to (U, U')$, and 
let $\wtW$ and $W$ be the larger sets given with $\wtU$ and $U$.  
We have by Lemma \ref{lemma:comparing_metrics}
\begin{equation}
\label{eqn:dd1} 
\frac{\diam(U';\sigma)}{\diam(U;\sigma)} = \diam(U'; \sigma/\diam(U; \sigma)) \asymp \diam(U'; \rho_{W})
\end{equation}
and similarly
\begin{equation}
\label{eqn:dd2}
\frac{\diam(\wtU'; \sigma)}{\diam(\wtU ;\sigma)} \asymp \diam(\wtU'; \rho_{\wtW}).
\end{equation}
By Lemma \ref{lemma:pvalent} 1(b), we have 
\begin{equation}
\label{eqn:dd3}
\diam(U'; \rho_{W}) \leq \diam(\wtU'; \rho_{\wtW}) \leq C_h(p, \diam(U'; \rho_{W})).
\end{equation}
Together, (\ref{eqn:dd1}), (\ref{eqn:dd2}), and (\ref{eqn:dd3}) imply 
\[ \frac{\diam(U'; \sigma)}{\diam(U ;\sigma)} \leq \mbox{const} \cdot \frac{\diam(\wtU'; \sigma)}{\diam(\wtU ;\sigma)}\]
and 
\[\frac{\diam(\wtU'; \sigma)}{\diam(\wtU ;\sigma)} \leq \mbox{const}\cdot C_h\left(p, \mbox{const}\frac{\diam(U'; \sigma)}{\diam(U ;\sigma)} \right) \]
as required.
\gap

\noindent{\bf Roundness distortion.}  We first estimate the distortion of roundness with respect to hyperbolic metrics, and then relate the hyperbolic to the spherical metric.  

Suppose $U \in \mathbf{U}=\union_n\UUU_n$, $a \in U$, and $\roundness(U,a) = K$ in the hyperbolic metric of $W$.
By the definition of roundness, there exists $s>0$ such that with respect to the hyperbolic metric on $W$, 
\[ B(a,s) \subset U \subset B(a, Ks)\]
and no smaller $K$ will do.  Thus 
\begin{equation}
\label{eqn:Ks1}
 \frac{1}{2} \diam(U) \leq Ks \leq \diam(U).
 \end{equation}
 Combining (\ref{eqn:r0}) and (\ref{eqn:Ks1}) we obtain 
 \begin{equation}
 \label{eqn:Ks2}
 r_0 \leq Ks \leq 1,
 \end{equation}
i.e. that  $K \asymp 1/s$  where the implicit constant is independent of $U \in \mathbf{U}$.  

Now suppose $f^k:(\wtU, \ta) \to (U,a)$.  
\gap

\noindent{\bf Backward roundness distortion.}  By Lemma \ref{lemma:pvalent} 1(a), with respect to the hyperbolic metric on $\wtW$, 
\[ B(\ta, s) \subset \wtU \subset B(\ta, C_h(p, Ks)) \subset B(\ta, C_h(p,1))\]
and so 
\[ \roundness(\wtU, \ta) \leq \frac{C_h(p,1)}{s} \leq \mbox{const}\cdot K\]
since $K \asymp 1/s$.  
Hence we obtain a linear backwards roundness distortion function.
\gap

\noindent{\bf Forward roundness distortion.}  Suppose now $\roundness(\wtU, \ta) = \widetilde{K}$.  Then with respect to the hyperbolic metric on $\wtW$, there exists $\tilde{s}>0$ such that 
\[ B(\ta, \tilde{s}) \subset \wtU \subset B(\ta, \widetilde{K}\tilde{s}) \subset B(\ta, 1)\]
so that $\widetilde{K} \asymp 1/\tilde{s}$.  Hence by Lemma \ref{lemma:pvalent}(2), with respect to the hyperbolic metric on $W$, 
\[ B(a, C_h^{-1}(p, \tilde{s})) \subset U \subset B(a, 1).\]
Therefore 
\[ \roundness(U,a) \leq 1/C_h^{-1}(p, \tilde{s})  \leq 1/C_h^{-1}(p, r_0/\tilde{K})\]
since $\tilde{s} \geq r_0/\tilde{K}$ by (\ref{eqn:Ks2}).  
\gap

It  remains only to transfer the roundness estimates from the hyperbolic to the spherical metric.  
Suppose $U \in \mathbf{U}$ has basepoint $x$, $a \in U$, and $\roundness(U,a)=K$ with respect to the hyperbolic metric on $W$.   
By construction, 
\[ U =B(x, r_0; \rho_W) \subset  B(x,1; \rho_W)\]
and we have already shown 
\[ B(a, Ks; \rho_W) \subset B(a, 1;\rho_W).\]
Therefore, the set $U$, its largest inscribed ball $B(a,s)$ about $a$, and its smallest circumscribing ball $B(a, Ks)$ about $a$ are all contained in the hyperbolic ball $D=B(x,2)$.  On the set $D$, Lemma \ref{lemma:comparing_metrics} implies that the hyperbolic metric $\rho_{W}$ is bi-Lipschitz equivalent to the metric $\sigma/\diam(D; \sigma)$.  Since roundness is   invariant under constant scalings of the metric, the factor $1/\diam(D; \sigma)$ is irrelevant.  It follows easily that the roundness computed with respect to the hyperbolic metric on $W$ is comparable to that computed with respect to the spherical metric $\sigma$.  

This completes the proof of Theorem \ref{thm:semihyperbolic_are_cxc}.
\qed

We now provide a converse statement.  

\begin{thm}[cxc on the Euclidean $\IS^2$ implies uniformly quasiregular]
\label{thm:nothing_new}
Suppose $f: \IS^2 \to \IS^2$ is an orientation-preserving metric cxc map with respect to the standard spherical
metric.  
Then $f$ is quasisymmetrically, hence
quasiconformally conjugate to a chaotic semi-hyperbolic rational map.
\end{thm}

We view this theorem as an analog of a theorem of Sullivan and Tukia 
which says that a group of uniformly $K$-quasiconformal homeomorphisms of $\IS^n$, $n\ge 2$, which acts as a uniform convergence group action
is quasiconformally conjugate to a cocompact Kleinian group; see \S\,4.4, Appendix \ref{appendix:hyperbolic_groups}, and \cite{DS3,tukia:quasiconformal_groups}.  
The proof of the theorem uses facts from quasiconformal analysis, see Appendix \ref{appendix:qrqc}.

\gap

\pf Any iterate
$f^n$ is also a finite  branched covering from the sphere to itself, hence has 
finitely many critical points.  
We will first prove that a constant $K$ exists such that $f^n$ is 
$K$-quasiregular for all $n\ge 1$.
Fix $n$ and consider a small disk $2D$ disjoint from the critical set of $f^n$. 
It follows
from Lemma \ref{prop:uwqr}
that there is a constant 
$H<\infty$ such that, 
for all $x\in 2D$, 
$$\limsup_{r\to 0} \frac{\max\{ |f^n(x)-f^n(y)|\; :\; |x-y|=r\}}{\min\{ 
|f^n(x)-f^n(y)|\;:\; |x-y|=r\}}\le H\,.$$
By Theorem \ref{thm:qc_via_H}, on the domain $D$, the iterates $f^n$ are uniformly $K=K(H)$-quasiconformal.  

This implies that $f^n$ 
is $K$-quasiregular off
the critical set of $f^n$. 
But finitely many points are removable for 
quasiregularity, hence
$f^n$ is $K$-quasiregular.
It follows from Theorem \ref{thm:uqrsullivan} that $f$ is quasiconformally conjugate to a rational map.
The semi-hyperbolicity follows from the property of bounded degree along 
pull-backs (cf. Theorem \ref{thm:def_semihyperbolic}). 
\qed

\subsection{Convex Hull of Julia sets}

We recall that the Euclidean sphere $\IS^2$ may be regarded as the sphere at infinity of hyperbolic three-space $\IH^3$ in the Poincar\'e ball model.  
Given a set $E \subset \IS^2$, the {\em convex hull} \index{Index}{convex hull}%
of $E$ is the smallest convex subset of $\IH^3$ which contains every hyperbolic geodesic joining pairs of distinct points in $E$.  

In this paragraph, we prove the following theorem.  Below, $\Gamma$ and $\phi_f$ are as in \S \ref{secn:spaces_associated}.

\begin{thm}\label{thm:car} Let $f$ be a rational map. Then $f$ is semihyperbolic if and only if there is a finite cover $\UUU$ of $J(f)$ such that
the space $\Gamma=\Gamma(f,\UUU)$ is quasi-isometric to the convex 
hull of the Julia set in $\HH^3$ by a quasi-isometry which extends as $\phi_f:J(f)\to\partial\Gamma$.
\end{thm}


We begin with the sufficiency.

\begin{prop}\label{hypimpsehyp} Let $f$ be a rational map, and $\Gamma$ the graph constructed from a finite covering $\UUU$.
If  $\Gamma$ is hyperbolic and $\phi_f:J(f)\to \partial_\varepsilon \Gamma$ is quasisymmetric for some
$\varepsilon>0$ then $f$ is semi-hyperbolic.
Furthermore, the measure of maximal entropy is doubling on $J(f)$.
\end{prop}


\gap

We start with a lemma.  In the statement, the notation $\diam$ refers to the Euclidean diameter.  

\begin{lemma}\label{zalchyp} Let $K$ be a compact subset of $J(f)-\{\infty\}$.  For any $\alpha =(w_0, W)$, where $W\in {\bf U}$ and $w_0\in W\cap K$, 
 let $A_{\alpha}(z) = w_0+\diam(W\cap J(f))z$.  
Under the assumptions of Proposition \ref{hypimpsehyp}, the family of maps $\{f^{|W|} \circ A_\alpha\}_{\alpha}$ 
is normal on $\C$. Furthermore,
any limiting map of this family is an open map.\end{lemma}

Any sequence $\{W_n\} \subset \mathbf{U}$ of neighborhoods of $w_0$ is 
is contained in some fixed Euclidean disk about $w_0$ on which the Euclidean 
and spherical metrics are comparable.   
Changing between comparable metrics does not affect normality, so we use whichever is most convenient.  

Below, the notation $\diam W$ will denote the diameter with respect to the spherical metric in $\cbar$, otherwise,
we will write $\diam_{\varepsilon}$ for the diameter with respect to the metric $d_\varepsilon$.

\gap

\pf Let us consider a sequence $(\alpha_n)$ of pointed sets; if the sequence of levels is bounded, then  the lemma is clearly true.
So we might as well assume that $\alpha_n=(w_n,W_n)$ with $|W_n|=n$. We shall then write $A_{\alpha_n}=A_n$.

 We note that $F^n$ is $e^{\varepsilon n}$-Lipschitz in $\partial_\varepsilon \Gamma$. 
But Proposition \ref{ba3} (which applies since $\phi_f$ is assumed to be a homeomorphism)
asserts that
$\diam_\varepsilon \phi_f(W\intersect J(f))\asymp e^{-\varepsilon|W|}$ for all $W\in \mathbf{U}$. Therefore, if $W\in S(n)$, 
then the  Lipschitz constant of
$$F^n \,: \left(\partial\Gamma, \frac{d_\varepsilon}{\diam_\varepsilon(\phi_f(W\intersect J(f))}\right)\to (\partial \Gamma,d_\varepsilon)$$ 
depends neither on the chosen element $W$ nor on $n$.

\gap

Suppose now that $\phi_f$ is $\eta$-quasisymmetric. We observe that the family 
$$ \phi_f\circ A_n \,:  A_n^{-1}(J(f))\to \left(\partial \Gamma, \frac{d_\varepsilon}{\diam_{\varepsilon}\phi_f(W\intersect J(f))}\right)$$ is equicontinuous. 
This follows from
Proposition 10.26 in \cite{heinonen:analysis}\,:
all these maps are $\eta$-quasisymmetric, and normalized: for any $z,w\in W$,
$$|A_n^{-1}(z)-A_n^{-1}(w)|\le 1$$ and
$$\frac{|\phi_f(z)-\phi_f(w)|_\varepsilon}{\diam_\varepsilon\phi_f(W\intersect J(f))}\le\eta\left(2\frac{|z-w|}{\diam(W\intersect J(f))}\right)\le
\eta (2|A_n^{-1}(z)-A_n^{-1}(w)|)\le\eta(2)\,.$$

\gap

This implies that, for all $R>0$,  all the maps $(f^n\circ A_n)|_{A_n^{-1}(J(f))\cap D(0,R)}$ share a common modulus of continuity $\omega_R$ since 
$$f^n\circ A_n = \phi_f^{-1}\circ F^n\circ (\phi_f\circ A_n)$$ and 
$$F^n \,: \left(\partial \Gamma, \frac{d_\varepsilon}{\diam_\varepsilon \phi_f(W\intersect J(f))}\right)\to (\partial \Gamma,d_\varepsilon)$$ 
is uniformly Lipschitz.

\gap

Let us now prove the lemma. If $(f^n\circ A_n)_{ n\ge 0}$ was not a normal family at a point
$z_\infty\in\C$, then Zalcman's lemma \cite{zalcman:bulletin} would imply
the existence of a convergent sequence of points $(z_k)$ with $z_\infty$ as a limit, a sequence $(\rho_k)$ of positive numbers
  decreasing to zero and a subsequence $(n_k)$ such that $f^{n_k}\circ A_{n_k}\circ B_k$ tends to an open map 
$g:\C\to\cbar$ where $B_k(z)= z_k+\rho_k z$.

\gap

Let $R=2|z_{\infty}|$,
and
let us choose $R'>2 d(0,g^{-1}(J(f))$. Then, for $k$ large enough, it follows that $B_k(D(0,R'))\subset D(0,R)$ and
$$\diam (f^{n_k}\circ A_{n_k})\circ B_k (D(0,R')\cap (A_{n_k}\circ B_k)^{-1}(J(f)))\le  \omega_R(2\rho_k R')$$
which tends to $0$. This contradicts the fact that $g$ is open since $\diam(g^{-1}(J(f)\cap D(0,R'))>0$. 

Therefore, $(f^n\circ A_n)_{n\ge 0}$ is a normal family on $\C$.

\gap

By construction, for all $k$ the domains $(f^{n_k} \circ A_{n_k})^{-1}(W_{0}\cap J(f))$ have diameter one, contain the 
origin, and map onto  $W_{0}\cap J(f)$ at level zero. 
Therefore any limiting map is nonconstant, hence open.
\qed

\gap

\pf (Proposition \ref{hypimpsehyp})
We will prove that the condition (1) in Theorem \ref{thm:def_semihyperbolic} follows from Lemma \ref{zalchyp}.  
Let $r>0$ be such that  any disk of radius $r$ 
centered at a point of $J(f)$ is contained in some open set defining the cover $\UUU$. 

If the condition were not satisfied, 
we would find a sequence of points $z_k\in J(f)$ and connected components $W_k$ of $f^{-n_k}(D(z_k,r))$ such that the degree of 
$f^{n_k}|_{W_k}$ would tend to infinity. We may assume that $(z_k)$ tends to some $z_\infty$. Let $w_k\in W_k$ be such that $f^{n_k}(w_k)=z_k$. 
It follows from Lemma \ref{zalchyp} that the sequence of maps $q_k(z)= f^{n_k}(w_k + \diam (W_k\cap J(f) z)$ is a normal sequence on $\C$ with open limits.  
Hence after passing to a subsequence we may assume $q_k \to q$ uniformly on the closed unit disk $\cl{\D}$.  
Since $q_k(\cl{\D}) \supset (D(z_k,r_k)\cap J(f))$ by construction, we have using Hurwitz' theorem that for all sufficiently large $k$, 
\[ \deg(f^{n_k}|W_k) \leq  \# \{q_k^{-1}(z_\infty)\cap\cl{\D}\} \leq  \# \{q^{-1}(z_\infty)\cap \cl{\D}\},\]
where $\#$ counts with multiplicity.
So the degree has to be eventually bounded, contradicting our assumption. Therefore, $f$ is semi-hyperbolic.

\gap

The statement on the measure follows from the following argument. 
Since $f$ is semi-hyperbolic, $f$ is also cxc (Theorem \ref{thm:semihyperbolic_are_cxc}), 
so Theorem \ref{thm:cor_cxc} implies that
$\mu_f$ is the unique measure of maximal entropy and that $\mu_f$ is also Ahlfors-regular of 
dimension $(1/\varepsilon)\log d$. 

 In particular $\mu_f$ is doubling\,:
there is a constant $C>0$ such that, for any ball $B(x,r)$, with $r\le\diam_{\varepsilon}\partial\Gamma$,
$\mu_f(B(x,2r))\le C\mu_f(B(x,r))$. Since this
condition is preserved under the application of quasisymmetric mappings, 
the same is true for $\phi_f^*\mu_f$ (cf. Cor.\,4.15 from \cite{heinonen:analysis}). Furthermore,
metric entropy is invariant under Borelian isomorphisms (Proposition 4.3.16 in \cite{hasselblatt:katok:dynamics}),
and in particular under homeomorphisms.  So, we recover the fact that $f$ admits a unique measure of maximal entropy\,:
the pull-back under $\phi_f$ of $\mu_f$, and this measure is doubling.

\qed

\gap

We may now prove Theorem \ref{thm:car}\,:

\gap

\pf (Theorem \ref{thm:car})
Suppose $\Gamma$ is quasi-isometric to the convex hull of $J(f)$ via a map which extends 
as $\phi_f$. First, since quasi-isometries between proper geodesic spaces preserve hyperbolicity, it follows
at once that $\Gamma$ is hyperbolic. Alternatively, the comment following the
statement of Proposition \ref{hypimpsehyp} shows that one may apply Theorem \ref{homeoimphyp} 
to conclude the hyperbolicity of $\Gamma$. Second, since quasi-isometries extend as quasisymmetric maps, 
$\phi_f$ is quasisymmetric. 

Therefore, Proposition \ref{hypimpsehyp} applies and  shows that $f$
is semi-hyperbolic. 
\gap

Conversely, if $f$ is semi-hyperbolic, then Theorem \ref{cxcimphyp} 
shows that $\phi_f$ is quasisymmetric and that $\Gamma$ is hyperbolic. 
Since both  $\Gamma$ and 
the convex hull of $J(f)$
are quasi-starlike Gromov spaces, the  quasisymmetry $\phi_f$ extends as
a quasi-isometry between $\Gamma$ and the convex hull of $J(f)$ (Theorem \ref{thm:bsch_funct}).
\qed

\subsection{Topological characterizations of chaotic semihyperbolic rational maps}

In this subsection, we prove the following theorem.  Below, the notation $\IS^2$ denotes the Riemann sphere $\rs$ equipped with the spherical metric.  

\begin{thm}[Characterization of chaotic semi-hyperbolic rational maps]
\label{thm:characterization}
Let $f: S^2 \to S^2$ be an orientation-preserving
finite branched covering map defined on a topological
$2$-sphere which satisfies [Expansion] with respect to some covering $\UUU$, and suppose $\varepsilon$ is small enough so that Theorem \ref{thm:construction} applies.   Then the following are equivalent:
\be
\item $f$ is topologically conjugate to a semi-hyperbolic rational function $R: \IS^2 \to
\IS^2$ with $J(R)=\IS^2$.
\item $\bdry_\varepsilon\Gamma$ is quasisymmetrically equivalent to $\IS^2$.
\item $\Gamma$ is quasi-isometric to hyperbolic three-space $\IH^3$.
\item The conformal gauge of $\bdry_\varepsilon \Gamma$ contains a $2$-Ahlfors regular metric.
\item The map $f$ is topological cxc and the sequence $\{\UUU_n\}_n$ of coverings of $S^2$ is conformal 
in the sense of Cannon.
\eb
\end{thm}

Recall that the conformal gauge of a metric space $(X,d)$ is the set of metrics $\hat{d}$ on $X$
such that the identity map $(X,d)\to (X,\hat{d})$ is quasisymmetric (see Chapter 15 in \cite{heinonen:analysis}).

\gap

The equivalence with (5) will be proved after developing some needed background. 
This
corresponds to a theorem  of Cannon and Swenson
for hyperbolic groups whose boundary 
is homeomorphic to the two-sphere \cite{cannon:swenson:characterization}.

\gap

\pf Suppose $\UUU=\UUU_0$ satisfies the axiom [Expansion].    Then there exists
some $N$ such that for all $n \geq N$, elements of $\UUU_n$ contain
at most one branch value of $f$.  Let $\VVV$ be a finite covering of the sphere by
Jordan domains $V$ which is finer than $\UUU_N$ and such that $\bdry V$ avoids 
the countable set of forward orbits of critical points.  The elements of $\mathbf{V}=
\union_n\VVV_n$ are then homeomorphic to Jordan domains, since they are
coverings of disks ramified over at most one point
and their boundaries are unramified covers of the Jordan curve
boundaries of elements of $\VVV$.
Since the quasi-isometry class of $\Gamma$ does not depend on the
cover (Theorem \ref{homeoimphyp}), we may assume at the outset that $\UUU=\VVV$ and hence
that elements $U$ of $\mathbf{U}$ and their complements in the sphere
are connected.  
\gap

$\mathbf{(1) \implies (2).}$  Suppose $h_1: S^2 \to \IS^2$ conjugates $f$
to a  semi-hyperbolic  rational 
function $R$ and $h_2: S^2 \to \bdry_\varepsilon \Gamma$
conjugates $f$ to the dynamics $F$ on the boundary of 
$\Gamma$. Since $R$ is semi-hyperbolic,
$f$ is topologically cxc  (Corollary \ref{cor:topcxc_shyp}) and so $F$ is topologically cxc as well. 
By Corollary \ref{topcxc_cxc1}, $F$ is metric cxc. 
The rigidity theorem, Theorem  \ref{thm:invariance_of_cxc}, implies that  $h_2\circ h_1^{-1}:\IS^2\to\partial_\varepsilon\Gamma$ is quasisymmetric.

\gap

$\mathbf{(2) \implies (1).}$  
Suppose $h: \bdry_\varepsilon \Gamma \to \IS^2$ is
a quasisymmetric map.  By Propositions \ref{prop:F_is_nearly_cxc} and \ref{prop:uwqr}, $F$ is uniformly weakly quasiregular (in the sense that it satisfies the conclusion of Proposition \ref{prop:uwqr}).   Since this condition is preserved under quasisymmetric conjugacies, so is $G=hFh^{-1}$.     By 
Theorem \ref{thm:qc_via_H}, the iterates of the map $G$ are uniformly quasiregular.  Sullivan's Theorem \ref{thm:uqrsullivan} implies that $G$ is quasiconformally conjugate to a rational map $R$.
Since $h$ is a quasisymmetry, Theorems \ref{thm:car} and Theorem \ref{thm:bsch_funct} together 
imply that $R$ has to semi-hyperbolic.

\gap

$\mathbf{(3) \iff (2).}$ This follows from Theorem \ref{thm:bsch_funct}:   
boundary
values of quasi-isometries are quasisymmetries and, conversely,
quasisymmetric maps of 
boundaries extend to quasi-isometries.  

$\mathbf{(4) \iff (2).}$ The fact that (2) implies (4) follows from the fact that $\IS^2$ is naturally a $2$-Ahlfors
regular metric space.

For the converse, since all elements of $\mathbf{U}$ are Jordan domains, Proposition \ref{prop:BLC} 
shows that $\partial_\varepsilon \Gamma$ is linearly locally connected.
 Since linear local connectivity is a 
quasisymmetry invariant \cite{heinonen:analysis}, there exists a metric in the gauge of  
$\partial_\varepsilon \Gamma$ which is both 
linearly locally connected and, by hypothesis, Ahlfors 2-regular.  By 
M.\,Bonk and B.\,Kleiner's characterization of the standard two-sphere 
\cite{bonk:kleiner:qsparam}, this implies that $\bdry_\varepsilon\Gamma$
quasisymmetrically equivalent to the standard Euclidean two-sphere.
\qed

\gap

These statements mimic similar theorems for Gromov hyperbolic groups in the context of
Cannon's conjecture. Statement (2) is concerned with Sullivan/Tukia's straightening theorem
of \qc groups \cite{DS3, tukia:quasiconformal_groups}; statement (3) is due to J.\,Cannon
and D.\,Cooper \cite{cannon:cooper} in the context of groups; statement (5) is due to M.\,Bonk and
B.\,Kleiner, and can be deduced
either  from \cite{bonk:kleiner:rigidity}, or from \cite{bonk:kleiner:qsparam} and \cite{bonk:kleiner:conf_dim}.

In \cite{bonk:kleiner:conf_dim}, M.\,Bonk and B.\,Kleiner also prove that a Gromov hyperbolic group 
admits a cocompact Kleinian action on $\rs$ if the {\em Ahlfors-regular conformal dimension}
\index{Index}{conformal dimension}%
 of
the gauge of its boundary is attained. The Ahlfors-regular conformal dimension of $(X,d)$ is the infimum of the Hausdorff 
dimensions over all Ahlfors-regular metrics in its gauge. In our setting of non-invertible dynamical systems, however, the analogous statement does not hold:

\begin{prop} There is a metric $d$ on the $2$-sphere $\IS^2$ and a metric cxc map $f:(\IS^2,d)\to (\IS^2,d)$  such that
the Ahlfors-regular conformal dimension is attained by $d$, but $f$ is not topologically conjugate to a rational map.\end{prop}

\pf Let us consider $F:\C\to\C$ be defined by $F(x+iy)= 2x + 3iy$. Let us consider the metric
$\hat{d}(x+iy,x'+iy')= |x-x'|+|y-y'|^\alpha$ where $\alpha=\log 2/\log 3$. One may check that $(\C,\hat{d})$ is
Ahlfors regular of dimension $1+1/\alpha$ and that this dimension is also its Ahlfors-regular conformal dimension since
the $(1+1/\alpha)$-modulus of the family of horizontal curves is clearly positive (cf. Theorem 15.10  in \cite{heinonen:analysis}).

Since $F(\Z[i])\subset \Z[i]$ and $F(-z)=-z$, this map descends to a map $f:\rs\to\rs$, onto which one can push down
the metric $\hat{d}$ to a metric $d$. It follows that this metric satisfies the same properties as $\hat{d}$.
Furthermore, since  $\hat{d}(F(z),F(z'))= 2 \hat{d}(z,z')$ clearly holds for any $z,z'\in \C$, it follows that
$f$ is cxc with the metric $d$. 

But since the conformal dimension of $(\rs,d)$ is strictly larger than $2$, Theorem \ref{thm:characterization}
shows that $f$ is not equivalent to a rational map.\qed

\noindent{\bf Remark.}  The preceding proposition implies that the metric space $(S^2, d)$ need not be a so-called {\em Loewner space} even if the Ahlfors-regular conformal dimension is attained (see \cite{heinonen:analysis}).
Also, if $\partial\Gamma$ admits an Ahlfors regular Loewner metric in its gauge, then Theorem \ref{thm:characterization} above together with a result of Bonk and Kleiner \cite{bonk:kleiner:qsparam} imply $f$ is conjugate to a rational map.

\subsection{Cannon's combinatorial Riemann mapping theorem}
\label{scn:cannoncomb}

Before we prove the equivalence with (5), 
we first review the notions that 
are needed to understand the statement and
the proof.   The basic idea is to estimate the classical modulus $\mod(A)$ of an annulus  using combinatorial data coming from a sequence of finer and finer coverings (see Appendix \ref{appendix:qrqc} for the definition of classical moduli).  
\gap

{\noindent\bf Combinatorial moduli.}
\index{Index}{modulus, combinatorial}%
Let $\SSS$ 
\index{Symbols}{$\SSS$}%
 be a covering of a  topological 
surface $X$.  
Denote by $\MMM(\SSS)$ the set
of maps  $\rho:\SSS\to \R_+$ such that $0<\sum_{s \in \SSS} \rho(s)^2<\infty$ 
which we call  {\sl admissible metrics}.  
Let $K\subset X$\,; the {\sl $\rho$-length} of $K$ is by definition 
$$\ell_\rho 
(K)=\sum_{s\in\SSS,\,s\cap K\ne\emptyset} \rho(s)$$
and its {\em $\rho$-area} is  
$$A_{\rho}(K)=\sum_{s\in\SSS,\,s\cap K\ne\emptyset} \rho(s)^2\,.$$
If $\Gamma$ is a family of curves in $X$ and if $\rho\in\MMM(\SSS)$, we define
$$L_\rho(\Gamma,\SSS)=\inf_{\gamma\in\Gamma} \ell_\rho(\gamma),$$
and its {\em combinatorial modulus} by
$$\mod(\Gamma,\SSS) = \inf_{\rho\in\MMM(\SSS)} 
\frac{A_{\rho}(X)}{L_\rho(\Gamma,\SSS)^2}= \inf_{\rho\in\MMM(\SSS)} 
\mod(\Gamma,\rho,\SSS).$$

Let $A$ be an annulus in $X$.  Let  $\Gamma_t$ be the set of curves in $A$ which join the boundary components  of $A$, and  $\Gamma_s$ those which  separate the boundary components of $A$. 
Define
$$\mod_{\sup} (A,\SSS)= \frac{1}{\mod\,(\Gamma_t,\SSS)}\quad and \quad 
\mod_{\inf} (A,\SSS)= 
\mod\,(\Gamma_s,\SSS)\,.$$
The classical moduli of $\Gamma_s, \Gamma_t$ are mutually reciprocal.  
In the combinatorial setting, this is no longer quite true.  However 
J.\,Cannon, W.\,Floyd and W.\,Parry have proved that always  
$\mod_{\inf}(A,\SSS)\le\mod_{\sup}(A,\SSS)$  \cite{cannon:floyd:parry:sq_rect}.

A covering $\SSS$ has {\em $N$-bounded overlap} if, for all $x\in X$, $$\sum_{s\in\SSS} \chi_s(x)\le N$$ where
$\chi_s$ denotes the characteristic function of $s$. Two coverings are said to be {\em $N$-equivalent}, or to have {\em $N$-bounded overlap}, if each piece of one intersects at most $N$ pieces of the other, and vice-versa.

\bigskip

{\noindent\bf Sequence of coverings.} 
In order to state  J.W.\,Cannon's combinatorial Riemann mapping 
theorem, we  
introduce a couple of new notions.

\bigskip

\noindent{\bf Definition.} A {\em shingle}
\index{Index}{shingle}%
 is a connected compact subset of $X$, and a {\em shingling} is 
a covering of $X$ by shingles.  
\bigskip

\noindent{\bf Definition.} A sequence of coverings $(\SSS_n)$ of $X$
is {\em $K$-conformal} 
($K\ge 1$) if 
\be
\item the mesh of $(\SSS_n)$ tends to zero;
\item for any annulus $A$ in $X$, there exist an integer $n_0$ and a positive constant $m=m(A)>0$ 
such that, for all
$n\ge n_0$, $$\mod_{\sup}(A,\SSS_n),\mod_{\inf}(A,\SSS_n)\in [m/K,Km];$$
\item for any $x\in X$, any $m>0$ and any neighborhood $V$, there is an annulus 
$A\subset V$ which separates $x$ from
$X\setminus V$ such that $\mod_{*} (A,\SSS_n)\ge m$ for all large $n$, where
$*\in\{\inf,\sup\}$.\eb
In the preceding definition, we assume neither that the elements of $\SSS_n$ are connected, nor that they are compact.

The quantity $m(A)$ will be referred to as the {\em combinatorial modulus}
\index{Index}{combinatorial modulus}%
 of $A$ with respect to the sequence $(\SSS_n)$. 
If $\SSS'=(\SSS_n')$ is another sequence of coverings whose elements $\SSS_n'$ are $N$-equivalent with $\SSS_n$, 
where $N$ is independent of $n$, then the combinatorial moduli computed with respect to $\SSS$ and $\SSS'$ are known 
to be comparable (Theorem 4.3.1 in \cite{cannon:swenson:characterization}).  Hence $\SSS$ is conformal iff $\SSS'$ is conformal.  

\begin{thm}[combinatorial Riemann mapping theorem \cite{cannon:rmt}] 
\label{thm:comb_rmt}
If $(\SSS_n)$ is a conformal 
sequence of shinglings,
on a topological surface $X$, then $X$ admits a complex structure such 
that the 
analytic moduli of annuli are comparable
with their combinatorial moduli.\end{thm}

There is also a converse\,:

\begin{thm}\label{thm:shi_is_conf} A sequence $(\SSS_n)$ of shinglings on 
the Riemann sphere is conformal if all of the following conditions are satisfied:
\be
\item the maximum diameter of an element of $\SSS_n$ tends to zero as $n \to \infty$,
\item each covering $\SSS_n$ has overlap bounded by some universal
constant $N$, and 
\item there exists a constant $K >1$ such that 
for any $n$ and any $s\in \SSS_n$, there exist two concentric disks $D_s$ and $\Delta_s$ such that $D_s\subset 
s\subset \Delta_s$, and such that  
 $\diam \Delta_s\le K\diam D_s$. 
\eb
\end{thm}

This is slightly different from Theorem 7.1 in \cite{cannon:rmt}.  There, the smaller disks $D_s$ are required to be pairwise disjoint.  There is no such assumption here,  so we provide a proof.
We will use the following lemma of J.\,Str\"omberg and A.\,Torchinsky \cite{st}.  Below, disks are spherical, and integrals are over the whole sphere.

\begin{lemma}\label{lmclan} Let $\BBB$ be a family of disks $B$, each equipped with a weight $a_B>0$.
For any $p>1$ and any $\lambda\in (0,1)$, there exists a constant $C=C(p,\lambda)>0$, independent of 
the family and of the weights, such that
$$\int \left(\sum a_B\chi_B\right)^p \le C \int \left(\sum a_B\chi_{\lambda B}\right)^p \,.$$
\end{lemma}

\pf (Theorem \ref{thm:shi_is_conf}) It suffices to prove that there is some constant $C>0$ such that, for any annulus $A$, 
there is some $n(A)$ such that, if $n\ge n(A)$, then $\mod_{\inf}(A,\SSS_n)\ge (1/C)\mod A$ and $\mod_{\sup}(A,\SSS_n)\le C\mod A$.

Fix an annulus $A$. Since the mesh of $\SSS_n$ tends to $0$, we may find some $n(A)$ and $\kappa>0$ ($\kappa$ independent
from $A$ and $n$) such that, 
for any $n\ge n(A)$, any piece $s\in\SSS_n$ which intersects $A$ and
any curve $\gamma\in\Gamma_t\cup\Gamma_s$ which intersects $s$, the length of $\gamma\cap 2\Delta_s$ is at least $\kappa\diam s$.

Let $\Gamma$ denote $\Gamma_s$ or $\Gamma_t$ and $\SSS=\SSS_n$ for some $n\ge n(A)$. 
If $\gamma\in\Gamma$, the family of pieces $s\in\SSS$ which intersects
$\gamma$ is denoted by $\SSS(\gamma)$.

If $\rho:\SSS\to \R_+$ is an admissible metric for $\Gamma$, we define a classical test metric 
$$\hat{\rho}=\sum_{s\in\SSS}\frac{\rho(s)}{\diam s} \chi_{2\Delta_s}\,,$$
where $\chi_{2\Delta s}$ denotes the characteristic  function of $2\Delta_s$.
Therefore, if $\gamma\in\Gamma$, then the definitions of $\hat{\rho}$ and $\kappa$ imply 
$$\begin{array}{ll}\ell_{\hat{\rho}}(\gamma) &
 \ge \dis\sum_{s\in\SSS(\gamma)}\frac{\rho(s)}{\diam s} \ell(\gamma \cap 2\Delta s)\\ & \\
& \ge \kappa\dis\sum_{s\in\SSS(\gamma)} \rho(s)\\ &\\
& \ge \kappa L_\rho(\Gamma,\SSS)\,\end{array}$$
and so 
\begin{equation}
\label{eqn:lengthest}
L_{\hat{\rho}}(\Gamma) \geq \kappa L_\rho(\Gamma, \SSS).
\end{equation}  
\medskip

On the other hand,
$$\begin{array}{ll} \mbox{Area}(\cbar,\hat{\rho}) & = \dis\int_{\cbar} \left(\dis\sum_{s\in\SSS} \frac{\rho(s)}{\diam s}\chi_{2\Delta_s}\right)^2\\&\\
& \le C\dis\int_{\cbar} \left(\dis\sum_{s\in\SSS} \frac{\rho(s)}{\diam s}\chi_{D_s}\right)^2\end{array}$$
by Lemma \ref{lmclan}. Since
 $\SSS$ has bounded overlap, 

$$\begin{array}{ll}\left(\dis\sum_{s\in\SSS} \dis\frac{\rho(s)}{\diam s}\chi_{D_s}\right)^2  &
\le N^2 \left(\max \left\{\dis\frac{\rho(s)}{\diam s}\chi_{D_s}\right\}\right)^2 \\&\\
& \le N^2\dis\sum_{s\in\SSS}  \left( \dis\frac{\rho(s)}{\diam s}\chi_{D_s}\right)^2\,.\end{array}$$
Therefore
$$\begin{array}{ll} \mbox{Area}(\cbar,\hat{\rho}) 
& \le CN^2\dis\sum_{s\in\SSS}\dis\int_{s} \left( \dis\frac{\rho(s)}{\diam s}\right)^2\\&\\
& \le CN^2K^2\pi\dis\sum_{s\in\SSS}\rho(s)^2\,.\end{array}$$ Hence 
\begin{equation} 
\label{eqn:areaest}
\mbox{Area}(\cbar,\hat{\rho}) \le CN^2K^2 \pi  A_\rho(\cbar)\,.
\end{equation}

Combining (\ref{eqn:lengthest}) and (\ref{eqn:areaest}) yields 
\[ \mod(\Gamma) \leq c \; \mod(\Gamma, \SSS)\]
where $c=CN^2K^2\pi/\kappa^2$ is independent of the level $n$ of $\SSS = \SSS_n$.   Taking $\Gamma = \Gamma_s$, we obtain
\[ \mod_{\inf}(A, \SSS) = \mod(\Gamma_s, \SSS) \geq \frac{1}{c} \mod(\Gamma_s) = \frac{2\pi}{c} \mod(A).\]
Taking $\Gamma = \Gamma_t$, we obtain 
\[ \mod_{\sup}(A, \SSS) = \frac{1}{\mod(\Gamma_t, \SSS)} \leq \frac{c}{\mod(\Gamma_t)} = \frac{c}{2\pi}\mod(A).\]\qed

\subsection{Proof of rational iff Cannon-conformal}

We now conclude the proof of Theorem \ref{thm:characterization}.  

\pf $[\mathbf{ (5)\implies (1)}]$
Assume that $f$ is topological cxc with respect to a covering $\UUU_0$ and that the sequence 
$\{\UUU_n\}_n$  is conformal. 
Let $\VVV_0$ be a finite covering of $S^2$ by Jordan domains so small that for each $V \in \VVV_0$, the closure of $V$ is contained in an element of $\UUU_0$, and let $\{\VVV_n\}_n$ be the corresponding sequence of coverings obtained by pulling back under iterates of $f$.  For $n=0, 1, 2, \ldots$ let $\SSS_n$ be the shingling of $S^2$ whose elements are the closures of the elements of $\VVV_n$.  Axiom [Degree] implies that the coverings $\UUU_n$ and $\SSS_n$ have bounded overlap.  
Since the sequence $\{\UUU_n\}$ is conformal, so is the sequence $\{\SSS_n\}_n$.   
By the combinatorial 
Riemann mapping theorem, Theorem \ref{thm:comb_rmt}, 
the sphere $S^2$ has 
a complex structure compatible with its 
combinatorial structure. In other words,
there is a homeomorphism $h:S^2\to \IS^2$ such that, for any 
annulus $A$ and for $n$ large enough,
$$\mod_*(A,\SSS_n)\asymp \mod (h(A))$$
where as before $* \in \{\inf, \sup\}$.  
The map $G=h\circ f \circ h^{-1}: \IS^2 \to \IS^2$ is a finite branched covering. 
We will prove that iterates of $G$ are uniformly $K$-quasiregular (see \S 4.4 and Appendix \ref{appendix:qrqc}).  This will establish {\bf (1)} 
by Sullivan's straightening theorem (Theorem \ref{thm:uqrsullivan}), and the fact that $G$ is topological cxc (Corollary \ref{cor:topcxc_shyp}). 
Fix $k$  and 
$z\in\IS^2$ off the (finite) branching set $B(G^k)$ 
of $G^k$. Let $V$ be a neighborhood of $z$ disjoint from  $B(G^k)$. Therefore, 
$G^k|V$ is injective so, if $A\subset V$ is an annulus,
then, for all $n$ large enough
$$\mod G^k(A)\asymp \mod_{\sup} (f ^k(h^{-1}(A)),S(n))
=\mod_{\sup} 
(h^{-1}(A),S(n+k))\asymp \mod A$$
hence $G^k|_V$ is $K$-quasiconformal for some universal $K$. Therefore, $G^k$ is 
$K$-quasiregular since $B(G^k)$ is finite,
hence removable; see Appendix\ref{appendix:qrqc}. 
\gap

 $[\mathbf{ (2) \implies (5)}]$ Since we have already proved that (2) $\implies$ (1), assumption (2) implies that $f$ is topologically conjugate to a semi-hyperbolic rational map.  
So, we may assume $f$ is topologically cxc by Corollary \ref{cor:topcxc_shyp}.  
In particular, Axiom [Degree] holds.  
By assumption, there exists a quasisymmetric homeomorphism $h:\partial_\varepsilon\Gamma\to \IS^2$.   
Let $\phi_f:S^2\to\partial_\varepsilon\Gamma$ be the conjugacy given by Theorem \ref{thm:construction}. Since $h$ is quasisymmetric, 
Proposition \ref{ba3} implies that the roundness of $h\circ\phi_f(W)$ is uniformly bounded for any $W\in\mathbf{U}$.   
Axiom [Degree] implies that  the sequence of coverings $\{\UUU_n\}$  has uniformly bounded overlap, and that the pairs $\UUU_n, \cl{\UUU}_n$ 
of coverings also have uniformly bounded overlap.  
Theorem \ref{thm:shi_is_conf} implies that the sequence of shinglings $\{\cl{\UUU}_n\}_n$ is conformal, and we conclude by bounded overlap that the sequence $\{\UUU_n\}_n$ is conformal.
\qed

\section{Finite subdivision rules} 
\label{secn:fsr}

Finite subdivision rules
\index{Index}{finite subdivision rule}%
 have been intensively studied since they give natural concrete examples with which to study Cannon's problem of determining when a sequence of finer and finer combinatorial structures yields a compatible conformal structure; see \cite{cfp:fsr}, \cite{cfkp:rationalmaps}, \cite{cfp:xcI} and the discussion in the preceding section.  

A {\it finite subdivision rule} (abbreviated fsr) $\RRR$
\index{Symbols}{$\RRR$}%
  consists of a finite 2-dimensional CW complex  $S_\RRR$, 
  \index{Symbols}{$S_\RRR$}%
  a subdivision $\RRR(S_\RRR)$ of $S_\RRR$,  and a continuous cellular map $\phi_\RRR: \RRR(S_\RRR) \to S_\RRR$ 
  \index{Symbols}{$\phi_\RRR$}%
   whose restriction to each open cell is a homeomorphism.  We assume throughout this section that the underlying space of $S_\RRR$ is homeomorphic to the two-sphere $S^2$ and $\phi_\RRR$ is orientation-preserving.  In this case, $\phi_\RRR$ is a postcritically finite branched covering of the sphere with the property that pulling back the tiles effects a recursive subdivision of the sphere.   That is, for each $n \in \N$, there is a subdivision $\RRR^n(S_\RRR)$ 
   \index{Symbols}{$\RRR^n(S_\RRR)$}%
   of the sphere such that $f$ is a cellular map from the $n$th to the $(n-1)$st subdivisions.   Thus, we may speak of {\em tiles} (which are closed 2-cells), {\em faces} (which are the interiors of tiles), {\em edges}, {\em vertices}, etc. at {\em level } $n$.  It is important to note that formally, a finite subdivision rule is {\em not} a combinatorial object, since the map $\phi_\RRR$, which is part of the data, is assumed given.   In other words:  as a dynamical system on the sphere, the topological conjugacy class of $\phi$ is well-defined.  

Let $\RRR$ be a finite subdivision rule on the sphere such that $\phi_\RRR$ is orientation-preserving. The fsr  $\RRR$ has {\em mesh going to zero} if for every open cover of $S_\RRR$, there is some integer $n$ for which each tile at level $n$ is contained in an element of the cover.  A {\em tile type} is a tile at level zero equipped with the cell structure induced by the first subdivision.  The fsr $\RRR$ is {\em irreducible} if, given any pair of tile types, an iterated subdivision of the first contains an isomorphic copy of the second.  If $\RRR$ has mesh going to zero, then it is easy to see that $\RRR$ is irreducible:  any two tile types are joined by a path of edges of some bounded length.  $\RRR$ is of {\em bounded valence} if there is a uniform upper bound on the valence of any vertex at any level.  

\begin{thm}
\label{thm:fsr_cxc}
Suppose $\RRR$ is a finite subdivision rule for which $S_\RRR$ is the two-sphere and the subdivision map $\phi_\RRR$ is orientation-preserving.  

If $\RRR$ has mesh going to zero, then there exists an open covering $\UUU_0$ such that $\phi_\RRR$ satisfies Axiom [Expansion] and [Irreducibility].

If in addition $\RRR$ has bounded valence, then $\phi_\RRR$ satisfies Axiom [Degree], and so $\phi_\RRR$ is topologically cxc.  
\end{thm}

\pf To define the covering $\UUU_0$, we recall a few notions from \cite{cfp:xcI}.   Given a subcomplex $Y$ of a CW complex $X$ the {\em star} of $Y$ in $X$, denoted $\cwstar(Y,X)$,  is the union of all closed tiles intersecting $Y$. Let $X$ denote the CW structure on the sphere at level zero, and set $X_n=\RRR^n(X)$.   

\begin{lemma}  Suppose $\RRR$ has mesh going to zero.  Then there exist $n_0, n_1 \in \N$ with the following property.  For each closed 2-cell  $t \in X_{n_0}$, the set $D_t=\cwstar(t, X_{n_0+n_1})$ is a closed disk which, if it meets the postcritical set $P$ of $\phi_\RRR$, does  so in at most one point, and this point lies in the interior $U_t$ of $D$. 
\end{lemma}

\pf {\bf (of Lemma)}  Mesh going to zero implies that for some $n_0$, each 2-cell $t$ of $X_{n_0}$ meets $P$ in at most one point.  It also implies that for some $n_1$, for any 2-cell $s$ of $X_0$, and any two 0-cells $x,y$ of $s$, no 2-cell of $X_{n_1}$ contains both $x$ and $y$.  Together, these two observations imply that for any 2-cell $t$ of $X_{n_0}$, the set $D_t=\cwstar(t, X_{n_0+n_1})$ is a cell complex which contains $t$ in its interior $U_t$, and which, if it intersects $P$, does so in its interior.  Since $D_t$ is the closure of $U_t$ and its boundary is a simple closed curve, $D_t$ is a disk.  
\qed

Let $\UUU_0$ be the finite open covering of the two-sphere underlying $X$ given by the Jordan domains $U_t$ constructed in the Lemma above, and consider the topological dynamical system $f=\phi_\RRR: X \to X$ together with $\UUU_0$.  Since $\RRR$ is irreducible, Proposition \ref{prop:elementary} 3(a) implies that axiom [Irreducibility] in the definition of topologically cxc holds.  For any $k \in \N$, the restriction of $f^k$  to an element $\wtU$ of $\UUU_k$ is a branched covering onto its image $U$ which is ramified at at most one $0$-cell $c$ which maps onto some $0$-cell $v$ .  Let $\tilde{w} \in \cl{U}$ be a 0-cell and put $w=f^k(\tilde{w})$.  Then $w$ is joined by an edge-path (i.e. a union of 1-cells) to  $v$ whose interior avoids $v$, and the length of this edge path (i.e. the number of 1-cells comprising it) is at most some constant $q$.  Since $f^k: \wtU \to U$ is ramified only at $c$, this edge-path lifts to an edge path of length at most $q$ joining $\tilde{w}$ to $c$.  It follows that the combinatorial diameter of the zero-skeleton of $\cl{U}$ is uniformly bounded.  Since $\RRR$ has mesh going to zero, it follows that axiom [Expansion] holds.  

Moreover, if in addition $\RRR$ has bounded valence,  then the ramification of $f^k$ at $c$ is uniformly bounded.  This implies that $\cl{\wtU}$ comprises a uniformly bounded number of cells and hence that the degree of $f^k: \wtU \to U$ is uniformly bounded, so that axiom [Degree] holds.  Axiom [Irreducibility] follows immediately from the irreducibility of $\RRR$.  

Hence, $f: X \to X$ together with $\UUU_0$ yields a topologically cxc system on the sphere.  
\qed

Under the hypotheses of Theorem \ref{thm:fsr_cxc}, if $\RRR$ has mesh going to zero and bounded valence, then the covering $\SSS_n$ by closed tiles 
at level $n$ and the covering $\UUU_n$ have bounded overlap independent of $n$.  It follows that the sequence $\{\SSS_n\}_n$ is conformal if and only if the sequence 
$\{\UUU_n\}_n$ is conformal.  Combining this observation with Theorem \ref{thm:characterization}, we conclude that for the subdivision maps $\sigma_\RRR$ 
of such subdivision rules, yet another, equivalent characterization of rational maps is the conformality 
of the sequence $\{\SSS_n\}_n$.  
Compare \cite{meyer:origami} and \cite[Thm. 3.1]{cfkp:rationalmaps}.

\section{Uniformly quasiregular dynamics}
\label{secn:lattes}

Let $M$ be a compact $C^\infty$ Riemannian manifold of dimension $n \geq 2$, 
and suppose $f: M \to M$ is a nonconstant quasiregular map (Appendix \ref{appendix:qrqc}).  
This condition implies in particular that $f$ is a finite branched covering. 
We say that $f$ is
{\em uniformly quasiregular}
\index{Index}{uniformly quasiregular}%
 if all its iterates are $K$-quasiregular for a fixed $K$.

When $n= 2$, D.\,Sullivan proved the following theorem \cite{DS2} in parallel with
a similar statement for quasiconformal groups on the $2$-sphere \cite{DS3}:

\begin{thm}[D.\,Sullivan]\label{thm:uqrsullivan} A uniformly quasiregular map 
of the standard Euclidean two-sphere to itself is quasiconformally conjugate to a rational map.
\end{thm}

The iteration of uniformly quasiregular maps on the standard two-sphere therefore reduces to the iteration of rational maps. 

In higher dimension $n\ge 3$, uniformly quasiregular maps generalize one-dimensional 
holomorphic dynamics, and have been introduced in this setting
 by T.\,Iwaniec and G.\,Martin in \cite{iwaniec:martin:uqr}. 
Uniformly quasiregular maps on space-forms have been classified
in \cite{martin:mayer:peltonen}. They can be seen as analogs of
quasiconformal groups.

For such maps,  Fatou sets are defined as the set of normality, and
Julia sets as the set of non-normality.

\gap

In \cite{mayer:lattes}, V.\,Mayer proposes a generalization of classical Latt\`es 
examples (\S \ref{subsecn:rational_maps}) to higher dimensions. 
\index{Index}{Latt\`{e}s map}%
They are uniformly quasiregular maps of finite degree $f:M\to M$, 
where $M$ is a compact Riemannian manifold,
which are defined as follows. 

There are a crystallographic group $\Gamma$ and an onto $\Gamma$-automorphic quasiregular map 
$h:\R^n\to M$ such that
$h(x)=h(y)$ if and only if there is some element $\gamma\in\Gamma$ so that 
$y=\gamma(x)$, and there are
a matrix $U\in {\rm SO}_n(\R)$ and a constant $\lambda >1$ such that, if we set 
$A=\lambda U$,
then $A\Gamma A^{-1}\subset \Gamma$ and such that the following diagram 
commutes
$$\begin{array}{rcl} \R^n & \stackrel{A}{\longrightarrow} & \R^n\\
h\downarrow & &\downarrow h\\
M & \stackrel{f}{\longrightarrow} & M\end{array}$$

For more precise statements, we refer to V.\,Mayer's article 
\cite{mayer:lattes}.

\gap

Let us recall the following compactness result (cf. Theorem 2.4 
\cite{martio:srebro:vaisala:normal}).  Here, $|\cdot|$ denotes the Euclidean metric on $\R^n$, and $B^n(R)$ the Euclidean ball of radius $R$ about the origin.

\begin{thm}[normality of qr mappings]\label{thm:norm-qr}
Suppose that $0<r<R\le\infty$, $0<r'<\infty$, $1\le K<\infty$, $N\ge 1$, and 
that $\FFF$
is a family of $K$-quasiregular mappings $f:B^n(R)\to\R^n$ such that every point 
has 
at most $N$ preimages, $f(0)=0$, and such that for each $f\in\FFF$ there is a 
continuum
$A(f)$ with the properties
$$0\in A(f),\quad \max\{|x|,\ x\in A(f)\}=r,\quad \max\{|f(x)|,\ x\in 
A(f)\}=r'.$$
Then $\FFF$ is a normal family and any limit map is $K$-quasiregular, and any 
point in the
range has at most $N$ preimages.\end{thm}

This implies that, under the assumptions of Theorem \ref{thm:norm-qr}, assuming $R=1$, there are functions
$d_+$ and $d_-$ such that $d_{\pm}(t)\to 0$ with $t$ and such that, for any $f\in\FFF$, and any
set $U\subset B^n(1)$, $\diam f(U)\le d_+(\diam U)$ and, for any compact connected subset $V$ of the image
of $f$ which contains the origin, $\diam W\le d_-(\diam V) $ where $W$ denotes the component of $f^{-1}(V)$
which contains the origin. 

\begin{thm} 
\label{thm:Lattes_are_cxc}
Latt\`es maps are cxc.\end{thm}

\pf Axiom [Irreducibility] clearly  holds.

Fix $r_0>0$; for any $x\in \R^n$, we denote by
$W(x)$ the connected component of $h^{-1}(B(h(x),r_0))$ which contains $x$.
It follows from the quasiregularity and the fact that $h$ is automorphic with 
respect to a cocompact
group of Euclidean motions that we may choose $r_0>0$ such that a constant
$N<\infty$ exists
so that, for all $x\in\R^n$,
the degree of $h|_{W(x)}$ is  bounded by $N$ 
(Lemma III.4.1 in \cite{rickman:qrmappings}).

We fix some size $r_1>0$ small enough so that, for any $x,y\in\R^n$, if $x$ 
belongs to 
the component $V(y)$ of $h^{-1}(B(h(y),r_1))$ containing $y$, then $B(x,2\diam 
V(y))\subset W(x)$. 

We define $\UUU_0$ as a finite subcover of $\{B(x,r_1),\ x\in M\}$. 
Then $\UUU_0$ satisfies the axioms [Degree] and [Expansion].
It remains to prove  [Roundness distortion] and [Diameter distortion]. 
Since $f$ is semi-conjugate to a conformal map,
we need only verify (i) $h$ distorts the roundness of (small) sets by a controlled amount, and 
(ii) $h$ distorts ratios of diameters of nested sets by a controlled amount.

We note that since $M$ is compact, one may find uniformly quasiconformal charts which map
balls of radius $3r_1$ in $M$ onto the unit ball of $\R^n$. Therefore, we may assume that $h$
takes its image into $\R^n$ in the sequel. 

\gap

For each $x_0\in\R^n$ and each connected open set $V$ contained in some $V(y)$ which contains $x_0$,  we consider 
the map
$$h_{x_0,V} :x\in B^n(0,1)\mapsto \frac{1}{\diam h(V)}\cdot
h(x_0+x\cdot \diam V ).$$
All these maps define a compact family $\FFF$ of degree at most $N$ 
according to 
Theorem \ref{thm:norm-qr} since $h_{x_0,V}(B(0,1))$ contains at least one point at distance $1/2$ from the origin.

\gap

If $W\subset V\subset \R^n$, then it follows that $W\subset B(x_0,\diam V)$ and 
$$\frac{\diam h(W)}{\diam h(V)}\le d_+\left(\frac{\diam W}{\diam V}\right)\,.$$

Similarly, if $V'\subset M$ is small enough, if $W'\subset V'$ and if $V$ and $W$ denote
connected components of $h^{-1}(V)$ and $h^{-1}(W)$ such that $W\subset V$, then 
$$\frac{\diam W}{\diam V}\le d_-\left(\frac{\diam W'}{\diam V'}\right)\,.$$
This establishes (ii).

\gap

Let $V\subset \R^n$ contained in some $V(y)$, and let $x_0\in V$. Denote by $K=\roundness(V,x_0)$ its roundness. Then 
$B(x_0,\diam V/(2K))\subset V$ so
$$B(h(x_0),d_+(1/(2K))\diam h(V))\subset h(V).$$ This proves that $\roundness(h(V),x_0)\le 1/d_+(1/(2K))$.

 Let us denote by $\FFF(K')$ the
subset of $\FFF$ obtained from pairs $(V,x_0)$ such that $\diam h(V)\le r_1$ and $\roundness(h(V),h(x_0))\le K'$.
This family is also compact, so the roundness of $V$ at $x_0$ depends only on $K'$. Hence (i) holds.

This ends the proof that a Latt\`es example is cxc.\qed

Conversely, one has\,:

\begin{thm} 
\label{thm:cxc_implies_Lattes}
Let $f:M\to M$ be an orientation preserving metric cxc mapping, where 
$M$ is a compact Riemannian manifold of dimension at least $3$.
Then $f$ is a Latt\`es map.\end{thm}

\pf It follows from Proposition \ref{prop:uwqr} and Theorem \ref{thm:qc_via_H} that
$f$ is uniformly quasiregular. Furthermore, it follows from compactness 
properties of qr mappings
that every point is conical: for any $x_0\in M$, a sequence of sizes $r_n\to 
0$ and a sequence
of iterates $k_n$ exist such that $x\in B(0,1)\mapsto f^{k_n}(x_0 +r_n x)$ 
defines a convergent
sequence to a non constant map. Therefore, Theorem 1.3 in 
\cite{martin:mayer:rigidity} implies that
$f$ is a Latt\`es map.\qed

\gap

Let us note that V.\,Mayer has also generalized the notion of power maps
in \cite{mayer:lattes}\,: these are uniformly quasiregular self-maps of the Euclidean sphere 
$f:\IS^n\to\IS^n$, $n\ge 3$, such that the Fatou set consists of two
totally invariant attracting basins, the Julia set is a sphere $\IS^{n-1}$,
and the dynamics on the Julia set is of Latt\`es type. These maps
are also clearly cxc, if one restricts the dynamics to suitable neighborhoods
of the Julia sets.

\gap

For all other known examples of uniformly quasiregular maps, the Julia set
is a Cantor set, and the Fatou set is the basin of an attracting or of
a parabolic fixed point \cite{iwaniec:martin:uqr, martin:branch, peltonen, hinkkanen:martin:mayer,
martin:extending,martin:mayer:peltonen}. In the former case, when
$f$ has degree $d$, then there are $d+1$ embedded balls $B_0,\ldots, B_d$,
such that $B_1,\ldots,B_d$ have pairwise disjoint closures, all of them
contained in $B_0$, and the restriction to each $B_j$, $j=1,\ldots,d$,
is a homeomorphism onto $B_0$: the Julia set is contained in these balls,
and the restriction of $f$ to these balls is clearly cxc.

\section{Expanding maps on manifolds}
\label{secn:expanding_maps}
\renewcommand{\arraystretch}{1}
If $X$ is metric space, and $f:X\to X$ is continuous, we say that $f$ is {\em expanding} if, for any $x\in X$, there is a neighborhood $U$ such that, for any distinct $y,z\in U$, one has $|f(y)-f(z)|>|y-z|$; cf. \cite[\S\,1]{gromov:expanding}.

\gap

\noindent{\bf A baby example.}  Let $X=\T^2=\R^2/\Z^2$ be the two-torus and 
$f:X \to
X$ the degree twelve covering map induced by $v \mapsto \phi v$ where
$\Phi: \R^2 \to \R^2$ is the linear map given by $\Phi(x,y)=(3x,4y)$.  
Equip
$\R^2$ with the norm $|\cdot|$ given by
\[ |(x,y)| = \max\{|x|,|y|^\lambda\} \]
where $\lambda = \log 3/\log 4$.  Then for all $v \in \R^2$,
\[ |\Phi(v)|=3|v|.  \]
It follows that for all $a, b \in X$ sufficiently close,
\[ |f(a)-f(b)| = 3|a-b|\]
and it follows easily that $(X,f)$ is cxc.

We now (greatly) generalize this example.

\begin{thm}[From expanding to homothety]
\label{thm:from_expanding_to_homothety}
Let $f: M \to M$ be an expanding map of a compact connected
manifold to itself.
Then there exists a distance function on $d$ on $M$ and constants
$\delta>0$ and $\rho>1$ such that for all $x,y \in M$, 
\[ d(x,y)<\delta \implies d(f(x),f(y)) = \rho \cdot d(x,y)\]
and such that balls of radius $\leq \delta$ are connected and
contractible.    
\end{thm}

\begin{cor}[Expanding implies cxc]
The dynamical system $((M,d),f)$ is  cxc.  Hence the metric $d$ is unique, up to quasisymmetry.
\end{cor}

\pf  (of Corollary) We remark that $f: M \to M$ is necessarily a
covering map of degree $D=\deg f$. 

Since $f$ is expanding on a compact manifold, Axiom [Irreducibility] holds.  

Let $\UUU_0$ be a finite open cover of
$M$ by open balls of radius $\delta$.   If $U \in
\UUU$ then since $U$ is contractible we have 
\[ f^{-1}(U) = \bigcup_{1}^D \wtU_i\]
where the union is disjoint and where each $f|\wtU_i : \wtU_i \to U$ is a
homeomorphism which multiplies distances by exactly the factor $\rho$. 
Thus for each $i$ there is an inverse branch $g_i: U \to \wtU_i$ which 
is a  homeomorphism and which contracts distances exactly by the factor
$\rho^{-1}$.   By induction, for each $n$ and each $U \in \UUU_0$ there
are $D^n$ inverse branches of $f^n$ over $U$ which are homeomorphisms and
which contract distances by $\rho^{-n}$.  Verification of the axioms is
now straightforward.   The last claim follows from Theorem \ref{thm:invariance_of_cxc}.  
\qed

The proof of Theorem \ref{thm:from_expanding_to_homothety} occupies the
remainder of this section.
\gap

\noindent{\bf Sketch of proof.}  One way to prove the theorem is to apply
the geometric constructions of the previous chapter. We prefer however to give a self-contained proof
using the algebra hidden behind expanding covers of Riemannian manifolds.

\be

\item[I.]  By a celebrated result of Gromov
\cite{gromov:expanding}, $f$ is topologically conjugate to the action of
an expanding endomorphism on an infra-nilmanifold.  Thus we may assume 
$M$ is an  infra-nilmanifold modelled on a simply connected nilpotent Lie
group $G$ and $f$ is such an endomorphism.  

\item [II.]  Let $\tilde{f}$ denote the
lift of $f$ to the universal  cover 
$G$.   We shall show that there exists an associated $\tilde{f}$-{\em
homogeneous norm} $|\cdot |: G \to [0,\infty)$ satisfying the following
properties for all $x \in G$:
\be
\item[1] $|x|=0 \iff x=1_G$,
\item[2] $|x^{-1}|=|x|$,
\item[3] $\exists \rho>1$ such that $|\tilde{f}(x)| = \rho|x|$ 
\item[4] $|\cdot |$ is proper and continuous.  
\eb 

\item[III.] 
For some $0<\epsilon \leq 1$, the function 
\[ x,y \mapsto |x^{-1}y|^\epsilon\]
is bi-Lipschitz equivalent to a left-invariant metric $d=d_\epsilon$ on
$G$.  In the metric $d$, the map $\tilde{f}$ expands distances by the
constant factor $\rho^\epsilon$, and thus $d$ descends to a distance on
$M$ with the desired properties.
\eb
\gap

We now begin the proof of Theorem \ref{thm:from_expanding_to_homothety}.
\gap

\noindent{\bf Infra-nilmanifolds.}  For background, see
\cite{dekimpe:svln}.  Let $G$ be a real, simply connected, finite
dimensional, nilpotent Lie group. 
\index{Index}{Lie group}%
Then
$G \rtimes \Aut(G)$ acts on $G$ on the left via 
\[ ^{(g,\Phi)}x = g\cdot \Phi(x).\]
An {\em almost-Bieberbach group}
\index{Index}{almost-Bieberbach group}
 is a torsion-free subgroup $E < G \rtimes
\Aut(G)$ of the form $L \rtimes F$ where $L<G$ is discrete and cocompact
and $F<\Aut(G)$ is finite.  Recalling that $E$ then acts freely on $G$,
the quotient $E\backslash G$ (which is not a coset space) is called an
{\em infra-nilmanifold modelled on $G$}.
\index{Index}{infra-nilmanifold}%

\gap 

\noindent{\bf Expanding endomorphisms.}  
Suppose $E$ is an  almost-Bieberbach group, $M=E\backslash G$, and
$(g,\Phi) \in G \rtimes \Aut(G)$ satisfies $(g,\Phi)E(g,\Phi)^{-1}
\subset E$.  Define 
\[ \tilde{f}: G \to G\]
by 
\[ \tilde{f}(x) = \;^{(g,\Phi)}x.\]
Then $\tilde{f}$ descends to a map 
\[ f: M \to M\]
which is called an {\em
endomorphism}
\index{Index}{endomorphism}%
 of the infra-nilmanifold $M$.  It is called {\em expanding}
\index{Index}{expanding}%
if all eigenvalues of the differential $d\Phi:
\mathfrak{g} \to \mathfrak{g}$ lie outside the closed unit disk, where
$\mathfrak{g}$ is the Lie algebra of $G$. 

We remark that 
\begin{equation}
\label{eqn:cancels} \tilde{f}(x)^{-1}\cdot \tilde{f}(y) =
\left(^{(g,\Phi)}x\right)^{-1}\cdot
\left(^{(g,\Phi)}y\right) =
\Phi(x^{-1}y).
\end{equation}
\gap

\noindent{\bf Homogeneous norms.}   If $\Psi \in G \rtimes \Aut(G)$, a
function $|\cdot |: G \to [0,\infty)$ will  be called a $\Psi$-{\em
homogeneous norm} if it satisfies properties (1)-(4) in (II) with
$\tilde{f}$ replaced by $\Psi$ in (3).   Equation (\ref{eqn:cancels})
implies that if $\tilde{f}$ is given by the action of $(g,\Phi)$,
then $|\cdot |$ is a $\tilde{f}$-homogeneous norm if and only if it is a
$\Phi$-homogeneous norm.  
\index{Index}{homogeneous norm}%

Since $G$ is simply connected, the exponential map $\exp: \mathfrak{g}
\to G$ is a diffeomorphism.  Hence we may identify $\mathfrak{g}$ and
$G$.  In this identification, $\Phi$ becomes $d\Phi$, which we again
denote by $\Phi$.  Thus, we may assume that $\Phi: \mathfrak{g} \to
\mathfrak{g}$ is a linear map and search for $\Phi$-homogeneous norms on
$\mathfrak{g}$.  
\gap

The case when $\Phi$ is semisimple is treated in detail in
\cite{folland:stein:homogeneous}. In general, we need the following
development.  
\gap

\noindent{\bf Linear algebra.}  

\begin{lemma}
\label{lemma:linear_algebra}
Let $\VVV$ be a finite-dimensional real vector space and $\Phi \in
\Aut(\VVV)$ have all eigenvalues strictly outside the closed unit disk.  Then
there exists a function 
\[ |\cdot|: \VVV \to [0,\infty),\]
and a real 1-parameter family $\Phi_t \subset \Aut(\VVV)$ with $\Phi =
\Phi_{1}$ such that for all $v \in \VVV$ and all $t \in \R$
\be
\item $|v|=0 \iff v=0$
\item $|-v|=|v|$
\item $|\Phi_t(v)| = e^t|v|$
\item $|\cdot |$ is proper and continuous.  
\eb
\end{lemma}

Assuming the lemma, we proceed as in \cite{folland:stein:homogeneous}. 
Let $|\cdot|$ be the homogeneous norm on $\mathfrak{g}$ given by Lemma
\ref{lemma:linear_algebra} and transfer this via the exponential map to a
homogeneous norm on $G$ satisfying conditions (1)-(4) in the sketch of
the proof.  Now we are done with Lie algebras and work only on $G$. 
Condition (4) implies that 
\[ \{(x,y) \in G \times G : |x|+|y|=1\}\]
is compact.  Therefore 
\[ Q = \sup\{ |xy| : |x|+|y|=1 \}\]
exists.  For any $x,y \in G$, let $t$ be so that $e^t = |x|+|y|$.  Then 
\[
\begin{array}{ccl}
 |xy| &  = & e^{t}e^{-t}|xy|\\
\; & = & e^t|\Phi_{-t}(xy)| \\
\; & = & e^t|\Phi_{-t}(x)
\Phi_{-t}(y)| \\
\; & \leq & e^t Q \\
\; & = & Q(|x|+|y|)
\end{array}\]
since $|\Phi_{-t}(x)| +|\Phi_{-t}(y)| = 1$ by construction.  In summary,
the norm $|\cdot|$ satisfies the additional property 
\be
\setcounter{enumi}{4}
\item $|xy| \leq 2Q\max\{|x|, |y|\}$.    
\eb
for some constant $Q>0$.  
\gap

\noindent{\bf Quasi-ultrametrics.} The function 
\[ \varrho(x,y) = |x^{-1}y|\]
satisfies the symmetry and nondegeneracy conditions of a distance
function by properties (1) and (2) of the norm
$|\cdot|$.  
By (5) above, we have that for $C=2Q$
\[ \varrho(x,z) \leq C\max\{ \varrho(x,y), \varrho(y,z)\}.\]
This turns $\varrho$ into a so-called {\em quasi-ultrametric}.  
\index{Index}{quasi-ultrametric}%

Given any quasi-ultrametric $\varrho$, there are
constants $C', \alpha>0$ such that $C'\varrho^\alpha$ defines a metric.  
We outline the construction below and refer to 
e.g. \cite{ghys:delaharpe:groupes}, \S 7.3 for details.  
Define 
\[ \varrho_\epsilon(x,y) =
|x^{-1}y|^\epsilon\]
which is now a quasimetric with constant $Q^\epsilon$.  
Moreover, it satisfies the homogeneity property  
\begin{equation}
\label{eqn:homog}
 \varrho_\epsilon(\Phi(x),\Phi(y)) = \rho^\epsilon
\varrho_\epsilon(x,y).
\end{equation}

 Given $x,y \in G$
a {\em chain} $C$ from $x$ to $y$ is a sequence 
\[
x=x_0, x_1, \ldots, x_n=y\]
of elements of $G$; its {\em length} is given by 
\[ l_\epsilon(C) = \sum_{i=1}^{n} \varrho_\epsilon(x_{i-1}, x_i).\]
The set of
chains from
$x$ to
$y$ is denoted
$\CCC_{xy}$.  Define a new function on pairs of points by 
\[ d_\epsilon(x,y) = \inf\{l_\epsilon(C) : C \in \CCC_{xy}\}.\]
The function $d_\epsilon$ is symmetric and trivially
satisfies the triangle inequality.  Since
$\varrho_\epsilon$ satisfies Equation (\ref{eqn:homog}), so does
$d_\epsilon$.   Moreover, if $Q^\epsilon < \sqrt{2}$ then for all
$x,y \in G$, one has (\ibid, Prop. 10)
\[ (3-2Q^\epsilon)\varrho_\epsilon(x,y) \leq d_\epsilon(x,y) \leq
\varrho_\epsilon(x,y)\]
so that the nondegeneracy condition holds and the functions 
$d_\epsilon,
\varrho_\epsilon$ are bi-Lipschitz equivalent.  

This completes the proof, modulo the proof of Lemma
\ref{lemma:linear_algebra}.
\gap

\pf (Lemma \ref{lemma:linear_algebra})  Assume first
that $\Phi$ lies on a 1-parameter subgroup 
\[ \Phi_t = \exp(\phi t)\]
for some $\phi \in \mbox{\rm End}(\VVV)$.  Then the real parts of the
eigenvalues of $\phi$ have strictly positive real parts.  
\gap

\noindent{\bf Claim.}  {\em There exists a basis for $\VVV$ such that if
$|| \cdot ||$ is the corresponding Euclidean norm, then for all
$0 \neq v \in \VVV$ the function 
\[ t \mapsto || \Phi_t x|| \]
is strictly increasing. } 
\gap

The claim implies that for nonzero $v$, there is exactly one $t(v)$ such
that 
\[ ||\Phi_{t(v)}(v)|| = 1.\]
Define $|0|=0$ and for $v\neq 0$ define 
\[ |v| = e^{-t(v)}.\]
Conclusions (1) and (2) are clearly satisfied.  To prove (3), note that
the conclusion is obvious if $v=0$ and if $v\neq 0$ we have 
\[ 1=||\Phi_{t(v)}v||=||\Phi_{t(v)-s}\Phi_s(v)||\]
hence 
\[ t(\Phi_s(v))=e^{t(v)-s} \implies |\Phi_s(v)|=e^s|v|.\]
Clearly $|\cdot |$ is continuous.  To prove properness,  note that the
Claim implies that for all $t\leq 0$, and for all $v$ with $||v||=1$,
$||\Phi_t(v)|| \leq 1$.  Thus 
\[ B=\{v : \; |v| \leq 1\}\]
is compact.  Therefore, given any $r=e^t$ we have by (3) that the set 
\[ \{v : \; |v| \leq r\} = \Phi_t(B)\]
is also compact.  It follows easily that $|\cdot |$ is proper.  
\gap

\noindent{\bf Proof of Claim.}
To prove the claim, let $\VVV = \oplus_i \VVV_i$ be the real Jordan
decomposition of $\VVV$ given by $\phi$ (not $\Phi$), and choose a basis
of $\VVV$ such that each Jordan block is either of the form 
\[
\left(
\begin{array}{cccccc}
\lambda_i & 1 & \\
\; & \lambda_i & \ldots \\
\; & \; & \ldots & 1 & \\
\; & \; & \; & \lambda_i & 1 \\
\; & \; & \; & \; & \lambda_i 
\end{array}
\right) 
\mbox{ or } 
\left(
\begin{array}{cccccc}
\rho_i R_{\theta_i}& I & \\
\; & \rho_i R_{\theta_i} & \ldots \\
\; & \; & \ldots & I & \\
\; & \; & \; & \rho_iR_{\theta_i} & I \\
\; & \; & \; & \; & \rho_i R_{\theta_i}
\end{array}
\right) 
\] 
where $\lambda_i,\rho_i>0$, $I$ is the 2-by-2 identity matrix, and
$R_{\theta} = \mtwo{\cos\theta}{-\sin\theta}{\sin\theta}{\cos\theta}$.
If $i$ corresponds to a block of the second kind we set $\lambda_i =
\rho_i \cos \theta_i$; this is positive since this is the real part of
the corresponding complex eigenvalue.   By making a coordinate change of
the form 
\[
\left(
\begin{array}{cccccc}
1 & \; & \\
\; & \delta^{-1}& \; \\
\; & \; & \ldots & \; & \\
\; & \; & \; & \delta^{-(m-1)}
\end{array}
\right) 
\]
for an $m$-by-$m$ block we may assume that the off-diagonal elements are
$\delta$ in the first case and $\delta I$ in the second, where
\[ 0 < \delta < \lambda_i.\]
Thus if $\phi_i = \phi|_{\VVV_i}$ then 
\[ \phi_i = \lambda_iI + \delta N_i + K_i\]
where $N_i$ is the nilpotent matrix with ones just above the diagonal,
$K_i$ is skew-symmetric, and the three terms commute pairwise.  

So setting 
\[ \Phi^i_t = \exp(\phi_i t) \]
we have 
\[ \Phi^i_t = \exp( (\lambda_i I + \delta N_i) t) \cdot \exp(K_it)\]
where the second factor is orthogonal.  

Let $\langle \cdot , \cdot \rangle_i$ denote the inner product on
$\VVV_i$ corresponding to the above basis on $\VVV_i$ and extend to
$\VVV$ in the obvious way so that the $\VVV_i$ are orthogonal.  
The claim is proved once we show that for each $i$ 
\[ t \mapsto ||\Phi^i_t(v)||_i \]
is strictly increasing.  

We have for all $t_0 \in \R$ and all $v \neq 0$ 
\renewcommand{\arraystretch}{2}
\[ \begin{array}{ccl}
\frac{d}{dt}\left|_{t=t_0}\right.\left< \Phi_t(v),
\Phi_t(v)\right> & = &
\frac{d}{dt}\left|_{t=t_0}\right.\left<e^{(\lambda_iI_i+\delta N_i)t}v
\;,\;
e^{(\lambda_iI_i+\delta N_i)t}v\right> \\
\; & = & 2 \left<e^{(\lambda_iI+\delta N_i)t_0},
\frac{d}{dt}\left|_{t=t_0}\right.e^{(\lambda_iI+\delta N_i)t}v\right>\\
\; & = & 2 \left< e^{(\lambda_i I + \delta N_i)t_0}v, (\lambda_i I +
\delta N_i)e^{(\lambda_i I + \delta N_i)t_0}v\right> \\
\; & = & 2\lambda_i \left< y , y \right> + 2\delta \left< y , N_i y
\right> 
\end{array}
\]
where $y = e^{(\lambda_i I + \delta N_i)t_0}v$.  The Cauchy-Schwarz
inequality shows that $|\left< y , N_i y \right>| < \langle y, y\rangle$ and so since
$\delta < \lambda_i $ we have that the derivative at $t_0$ is strictly
positive and the claim is proved.
\gap

If $\Phi$ does not lie on a 1-parameter subgroup we proceed as
follows.  It is well-known that $\Phi$ lies on a
1-parameter subgroup if and only if the Jordan blocks with negative
real eigenvalues occur in identical pairs.  If this is not the
case, we first change notation so $\Phi = \Phi'$.  Next, let $M_i: \VVV_i
\to \VVV_i$ be given by $-\id$ if the $i$th Jordan block of $\Phi$  is
real with negative eigenvalue and by $\id$ otherwise, and set
$M=\oplus_i M_i: \VVV \to \VVV$.  Then $\Phi'$ commutes with $M$ and we
set $\Phi = M \Phi'$.    Then $\Phi$ lies on a 1-parameter
subgroup $\Phi_t = \exp(\phi t)$ and we set $\Phi'_t = M
\Phi_t$.  Since 
\[ || \Phi'_t (v)|| = ||M\Phi_t(v)|| = 1 \iff ||\Phi_t(v)|| = 1\]
we have 
\[ |\Phi'_t(v)| = e^t|v|\]
for all nonzero $v$ and the proof is complete.
\qed 

\noindent{\bf Remarks:}  In many cases, raising to a power in step (III) of the construction of $d$ is unnecessary and a representative metric $d$ can either be written down explicitly or is a well-studied object.  

For example, suppose $\mathfrak{g}$ is abelian (i.e. all brackets are trivial) and $\Phi$ is diagonalizable over $\R$.  This is a generalization of the baby example and one can write the metric $d$ explicitly.  The resulting gauges on the universal cover $\R^n$ are studied by Tyson \cite{tyson:mgqc}, \S 15.  If not equivalent to the Euclidean gauge, these gauges are highly anisotropic:  there exist a flag $V_0 \subset V_1 \subset \ldots \subset V_m = \R^n$ such that any qs automorphism $h$ satisfies $h(V_k)=V_k$, $k=1, \ldots, m$.  

Another well-studied situation arises in the {\em Carnot-Carath\'eodory case},
\index{Index}{Carnot-Carath\'eodory space}%
 i.e. when $\Phi|\HHH = \lambda \id _\HHH$  on a subalgebra $\HHH$ which generates $\mathfrak{g}$ as a Lie algebra.   In this case any two points are joined by a smooth curve with tangent in the distribution defined by $\HHH$.  The resulting length space is a so-called {\em smooth Carnot-Carath\'eodory} metric space; cf. Pansu \cite{pansu:cc}.    The prototypical example is the map $(x,y,z) \mapsto (2x, 2y, 4z)$ on the {\em Heisenberg manifold} 
 \index{Index}{Heisenberg manifold}%
 $M=H/\Gamma$ where $H$ is the three-dimensional Heisenberg group of upper triangular matrices with ones on the diagonal and $\Gamma$ is the lattice consisting of such matrices with integer entries.  

In both cases, the conformal dimension (i.e. the infimum of the Hausdorff dimension over all quasisymmetrically equivalent spaces) is given by 
\[ \frac{1}{\lambda_1} (\lambda_1 + \lambda_2 + \ldots + \lambda_n) \]
where $\lambda_1 \leq \lambda_2 \leq \ldots \leq \lambda_n$ are the eigenvalues of $\Phi$.   

The classification of Lie algebras admitting expanding endomorphisms is still in progress; see \cite{dekimpe:lee:expandingaffine}.

\section{Expanding maps with periodic branch points}
\label{secn:periodic_branch}

\subsection{Barycentric subdivision} Given a Euclidean triangle $T$, its {\em barycentric subdivision} 
\index{Index}{barycentric subdivision}%
is the collection of six smaller triangles formed by the three medians.  Barycentric subdivision is natural with respect to Euclidean affine maps:  if $A: \R^2 \to \R^2$ is affine, then the small triangles comprising the barycentric subdivision of $T$ are sent by $A$ to those comprising $A(T)$.  If $T$ is equilateral, the six smaller triangles are congruent.  Suppose $T$ has side length one,  and let $B$ be an orientation-preserving affine map sending $T$ to one of the six smaller triangles in its barycentric subdivision.  Then 
\[ B = L \circ S \circ K\]
where $K$ is an linear isometry, $L$ is a translation, and 
\renewcommand{\arraystretch}{1}
\[ S=\mtwo{1/2}{\sqrt{3}/6}{0}{1/3}.\]  
Using the naturality of barycentric subdivision and the fact that the Euclidean operator norm of $S$ is $(\sqrt{7}+1)/6 \approx .608<1$, it follows easily that for any triangle, under iterated barycentric subdivision, the diameters of the smaller triangles after $n$ subdivisions tend to zero exponentially in $n$.  

\begin{figure}[htb]
\label{fig:bary}
\psfrag{1}{\mbox{$T_1$}}
\psfrag{T}{\mbox{$T$}}
\psfrag{3}{\mbox{$f$}}
\psfrag{5}{\mbox{$c$}}
\psfrag{6}{\mbox{$a$}}
\psfrag{7}{\mbox{$b$}}
\includegraphics[width=5in]{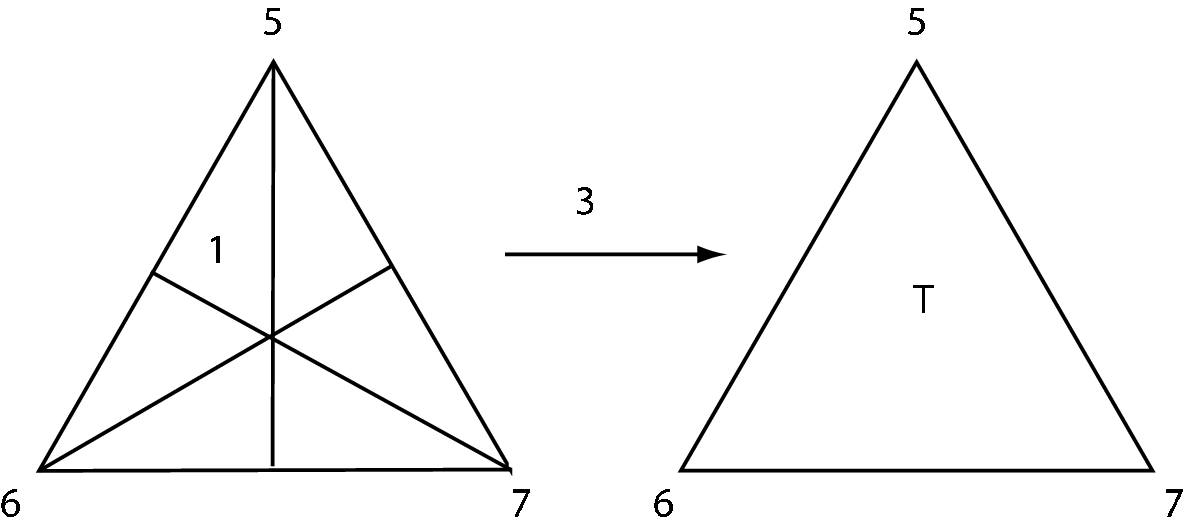}
\begin{center}{\sf Figure \ref{fig:bary}}\end{center}
\end{figure}

Let $T$ as above be the Euclidean equilateral unit triangle.  Equip $T$ with an orientation and label the vertices of $T$ as $a, b, c$ as shown.  Let $T_1$ be one of the two smaller triangles in the first subdivision meeting at the vertex $c$, and let $\phi: T_1 \to T$ be the restriction of the unique orientation-preserving Euclidan affine map fixing $c$ and sending $T_1$ onto $T$.  Regard now the two-sphere $S^2$ as the double $T \union \cl{T}$ of the triangle $T$ across its boundary.  Equip $S^2$ with the complete length structure inherited from the Euclidean metric on $T$ and its mirror image, so that the sphere becomes a CW complex $X$ equipped with a path metric.  By composing with reflections, there is a unique affinely natural extension of $\phi$ to an orientation-preserving degree six branched covering map $f=\phi_\RRR: S^2 \to S^2$ sending each of the twelve smaller triangles at level one onto $T$ or $\cl{T}$; see Figure \ref{fig:bary}.  

The twelve smaller triangles give a CW structure $\RRR(X)$ on $X$ subdividing the original one, and we obtain a finite subdivision rule (in the sense of \S \ref{secn:fsr}) with mesh going to zero.  Notice, however, that this f.s.r. does not have bounded valence, since the branch point $c$ of $\phi_\RRR$ is a fixed 0-cell.  

Let $\UUU_0$ be the finite open cover of the sphere whose elements are given by the construction in Section \ref{secn:fsr}.  The discussion there implies that together, $f:  S^2 \to S^2$ and $\UUU_0$ satisfy axioms [Irreducibility] and [Expansion], but not [Degree] in the definition of topologically cxc, and that the diameters of the elements of $\UUU_n$ tend to zero exponentially in $n$.  

Let $\Gamma_f=\Gamma(f,\UUU_0)$ be the associated graph constructed in \S 
\ref{secn:spaces_associated}.  By theorem \ref{thm:construction}, for some $\varepsilon>0$, there is a homeomorphism 
$ \phi_f: S^2 \to \bdry_\varepsilon \Gamma_f$
conjugating $f$ to the induced map $F$ on the boundary.  Since $P_f$ consists of a finite set of points, Proposition \ref{prop:degree_fails} applies and hence $\bdry_\varepsilon \Gamma_f$ fails to be doubling.  

The map $f$ is not the only dynamical system naturally associated to the barycentric subdivision rule.  Let $\IH \subset \C$ denote the upper half-plane and let $\rho: T \to \IH$ be the unique Riemann map sending $a \mapsto 0$, $b \mapsto 1$, $c \mapsto \infty$.  By Schwarz reflection, this defines a conformal isomorphism $\rho: S^2 \to \rs$, where now $S^2$ is the sphere endowed with the conformal structure of the path metric defined above. Let $\psi: T_1 \to T$ be given by the unique Riemann map fixing $c$ and sending vertices to vertices.   As before, this determines an f.s.r. $\SSS$ with an associated map $\psi_\SSS: S^2 \to S^2$.  By the symmetry of the construction, the map $g: \rs \to \rs$ given by $\rho \circ \phi_\SSS \circ \rho^{-1}$ is a rational map; it is given by 
\[ g(z) = \frac{4}{27}\frac{(z^2-z+1)^3}{z^2(z-1)^2}.\]
See \cite{cfkp:rationalmaps}.  

The composition $h'=(\phi_\RRR|_{T_1})^{-1}\circ (\phi_\SSS|_{T_1}): T_1 \to T_1$ extends by reflection to a homeomorphism $h_1': (S^2, a, b, c) \to (S^2, a, b, c)$.  Letting $h_0'=\id$, then
\[ h_0' \circ \phi_\SSS = \phi_\RRR \circ h_1'\]
and $h_0'$ is isotopic to $h_1'$ relative to the set $\{a, b, c\}$.  That is, as postcritically finite branched coverings of $S^2$, $\phi_\RRR$ and $\phi_\SSS$ are combinatorially equivalent.  Letting $h_1 = h_1' \circ \rho^{-1}$ gives $h_0 \circ g = f \circ h_1$ with $h_0, h_1$ isotopic relative to the set $\{0, 1, \infty\}$.  

By lifting under the dynamics, we obtain for each $n \in \N$ a homeomorphism $h_n: \rs \to S^2$ such that $h_n\circ g = f \circ h_{n+1}$ with $h_n \sim h_{n+1}$ relative to $\{0,1,\infty\}$.  Since $f$ is uniformly expanding with respect to the length metric on $S^2$, the sequence of maps $\{h_n\}$ converges uniformly to a map $h: \rs \to S^2$ for which $hg=gh$.  Since $g$ is locally contracting near infinity, the diameters of the  preimages of the two half planes $\IH^{\pm}$ under $g^{-k}$ which meet the point at infinity remain bounded from below as $k \to \infty$.  Therefore $h$ is not injective.  Indeed, it is easy to see that $h$ collapses the closure of each Fatou component to a point.   

Let $\VVV_0 = \{ h^{-1}(U):  U \in \UUU_0\}$ be the open covering of $\rs$ given by pulling back the elements of $\UUU_0$ under $h^{-1}$ and let $\Gamma_g=\Gamma(g, \VVV_0)$.   Then $h$ induces an isometry $h_\Gamma: \Gamma_g \to \Gamma_f$.  The natural map $\phi_g: J_g \to \bdry\Gamma_g$   satisfies $\phi_f \circ h|J_g = \bdry h_\Gamma \circ \phi_g|J_g$ and collapses the closure of every Fatou component to a point.

\subsection{Expanding polymodials}
For $z = re^{i\theta}$ let $f: \C \to \C$ be given by
$f(z) = 1-are^{i2\theta}$ where $a = (1+\sqrt{5})/2$ is the golden
ratio.   This map is an {\em expanding polymodial}
\index{Index}{polymodial}%
 in the sense of
\cite{bcm:polymodialsI} and is studied in their Example 5.2.   

The origin is a  critical point which is periodic of period
three, hence for each $n \in \N$, $f^{3n}$ is locally $2^n$-to-one on
neighborhoods of the origin.   
The point $-\beta \approx -1.7$ is a repelling fixed
point with preimage $\beta$.  Let $I=[-\beta, \beta]$, $I^-=[-\beta, 0], 
I^+=[0,\beta]$.  Then each map $f|_{I^\pm}: I^\pm \to I$ is a homeomorphism onto its image which 
expands Euclidean lengths by the factor $a$.  

Let $T_0$ be the metric tree with underlying space $I$ and length metric given by the
Euclidean length metric $\sigma_0$.  It is easy to see that for all $n \in \N$, the set
$f^{-1}(T_n)$ is a tree $T_{n+1}$ which is the union of $T_n$ together with a finite collection of
smooth closed arcs $J_i$. Each such $J_i$ is attached to $T_n$ at a single
endpoint which lies in
$f^{-n+1}(\{0, f(0), f(f(0))\})$, and $f|_{J_i}$ is a homeomorphism onto
its image.  Inductively, define a length metric $\sigma_n$ on  $T_n$ by setting 
\[ \sigma_{n+1}|_{J_i} = a^{-1}(f|_{J_i})^*(\sigma_n).\]
Then $f:T_{n+1} \to T_n$ multiplies the lengths of curves by the factor $a$.

Let $\pi_{n+1}: T_{n+1} \to T_n$ be the map which collapses each such ``new'' interval
$J_i$ to the point on $T_n$ to which it is attached.  Clearly, $\pi_n$ is
distance-decreasing for all $n$.  Let 
\[ X = T_0 \stackrel{\pi_1}{\leftarrow} T_1 \stackrel{\pi_2}{\leftarrow} T_2 \ldots  \]
denote the inverse limit.   Metrize $X$ as follows.  The
diameters of the $T_n$ are bounded by the partial sums of a convergent geometric series
and thus are uniformly bounded.  Hence for all $x=(x_n), y=(y_n) \in T$, 
\[ \sup_n \sigma_n(x_n, y_n)\]
is bounded and increasing, hence convergent.  It follows easily that $T$ inherits a
length metric $\sigma$ such that the map $f$ of $T$ induced by $f|_{T_{n+1}}:
T_{n+1} \to T_n$ multiplies the lengths of curves by the factor $a$.    

It is easy to see that for each $k \geq 1$, near the origin,
$X$ contains an isometrically embedded copy of the one-point union of
$2^k$ copies of a Euclidean interval of length $a^{-k}$ where
the common vertex is the origin $(0,0,0,\ldots)$.  This implies that $X$
is not doubling, since (i) doubling is hereditary under passing to
subspaces, and (ii) at least $2^k$ balls of radius $a^{-k}/2$ will be
needed to cover the ball of radius $a^{-k}$ centered at the origin.
On the other hand, it is also easy to see that $(X,f)$ satisfies the other 
axioms for a cxc system.

\section{Some comparisons with $p$-adic dynamics}

The construction of the graph $\Gamma$ is reminiscent of certain constructions in $p$-adic dynamics.
Below, we  give a quick and partial account of $p$-adic
dynamics 
\index{Index}{p-adic dynamics}%
in order to point out some formal similarities and major differences between our setting and the $p$-adic setting. References include \cite{baker:hsia, benedetto:reduction,
rivera:thesis}.

\gap

The main object of $p$-adic dynamics is to understand the iterates of rational
maps with $p$-adic coefficients acting on $\PP^1(\C_p)$, where $\C_p$ is the metric completion 
of the algebraic closure
of $\Q_p$ endowed with the $p$-adic norm.  $\C_p$ is an algebraically closed, 
non-Archimedean valued field, and a complete {\it
non-locally compact} totally disconnected ultra-metric space.

Let us note that the first difference with our setting is that $\C_p$ is neither locally compact
nor connected ! 

\gap

Since the metric on $\C_p$ is an ultrametric, two balls are either disjoint, or one is contained in the other.  
In turn this induces a tree structure on the family of balls: 
the vertices are the balls of rational radii,  and  the edges originating
from such a vertex are parametrised by the residual field $\cl{\mathbb{F}_p}$. 
If $B\subset B'$ are two balls, then the edge joining them is made of the intermediate balls,
and if $B\cap B'=\emptyset$, then the edge joining these balls is made of the two edges joining
these balls to
the smallest ball which contains both of them.
This is the {\it $p$-adic hyperbolic space} $\mathcal{H}$ \cite{rivera:padichyp}.
$\mathcal{H}$ 
can be metrized to become a complete $\R$-tree i.e., a $0$-hyperbolic  metric
space. The projective space $\PP^1(\C_p)$ is a part of 
the boundary of $\mathcal{H}$.  The tree $\mathcal{H}$ is isometric to the Bruhat-Tits building 
for ${\rm SL}(2,\C_p)$ and is closely related to the Berkovich line \cite{berkovich:book}.

We emphasize that the  boundary at infinity of $\mathcal{H}$  is larger 
than $\PP^1(\C_p)$, since some intersections of balls with radii 
not converging to $0$ may be empty, yielding points  of $\bdry \mathcal{H}$ not in $\PP^1(\C_p)$.

\gap

It turns out that rational maps send balls to balls, and that rational maps always  act on the tree $\mathcal{H}$. So in the $p$-adic setting, the natural
hyperbolic space $\mathcal{H}$ on which any rational map $f$ acts does not depend on the dynamics:
it is a universal object independent of $f$. Another difference is that dynamics can be
tame on the boundary, but never on the tree.  That is, 
in $\PP^1(\C_p)$ the chaotic locus may be empty, but in $\mathcal{H}$ it is always nonempty. 
In contrast, in our setting, the dynamics
is always chaotic on the boundary $\bdry \Gamma$, while the induced dynamics on the tree $\Gamma$ itself is transient.

Finally, one can define, for rational maps $R$ of degree $d$ on $\C_p$, 
an invariant measure $\mu$ such that $R^*\mu=d\mu$. While its metric entropy
is at most $\log d$, there are examples for which the inequality is
strict.  This happens when the Julia set is contained in the hyperbolic
space $\mathcal{H}$ and the topological degree of the map on the Julia set is 
 also stricty smaller than $d$; see for instance \cite{favre:rivera:brolin}.

\gap

Of course, there are examples of rational maps over the $p$-adics for which the dynamics on the Julia set is  conjugate to a full shift. In such cases, one obtains cxc maps.
But generally, the $p$-adic case is rather different from ours, and the
similarities are merely formal.


\begin{appendix}
\chapter{Quasiconformal analysis}
\label{appendix:qrqc}
 
In this chapter, we summarize concepts and basic facts about quasiconformal and quasiregular mappings.  Standard references are the texts by 
V\"ais\"al\"a \cite{vaisala:lectures_qc} for quasiconformal
maps, and Rickman \cite{rickman:qrmappings} for quasiregular maps.  
This chapter does not deal with dynamics.  

\gap
The following results are summarized from \cite[\S\S I.4 and II.6]{rickman:qrmappings}.  Let $\R^n$ denote Euclidean space with the usual metric $ds$ and Lebesgue measure $m$,  and let $U \subset \R^n$ be a domain.  Fix $n \geq 2$ and let $f: U \to \R^n$ be a continuous, not necessarily invertible, nonconstant map.   The map $f$ is called {\em quasiregular} provided $f$ belongs to the Sobolev space $W^{1,n}_{loc}$, and, for some $K<\infty$, satisfies 
\begin{equation}
\label{eqn:qr}
|Df(x)|^n \leq K \cdot J_f(x)\;\;\; a.e.
\end{equation}
where $|Df(x)|$ is the Euclidean operator norm of the derivative and $J_f(x)$ is the Jacobian derivative.  In this case we say also that $f$ is {\em $K$-quasiregular}.  
\index{Index}{quasiregular}%
If in addition $f: U \to f(U)$ is  a homeomorphism, $f$ is said to be {\em $K$-quasiconformal}.  
\index{Index}{quasiconformal}%
A quasiregular map is discrete and open.  The {\em branch set} $B_f$ is the set of points at which $f$ fails to be a local homeomorphism.  The branch set $B_f$ and its image $f(B_f)$ have measure zero.  Also, $f$ is differentiable almost everywhere.  The condition (\ref{eqn:qr}) implies that at almost every point $x$ in $U\setminus B_f$, the derivative sends round balls to ellipsoids of uniformly bounded  eccentricity.  
The composition of a $K$-quasiregular map with a conformal map is again $K$-quasiregular.  The inverse of a quasiconformal map is quasiconformal, and the composition of quasiregular maps is quasiregular.   

The following theorems give alternative useful characterizations.  To set up the statements, let 
\[ H(x,f) = \limsup_{r\to 0} \frac{\max\{ |f^n(x)-f^n(y)|\; :\; |x-y|=r\}}{\min\{ 
|f^n(x)-f^n(y)|\;:\; |x-y|=r\}}.\]
If $\Gamma$ is a set of rectifiable paths in $U$ and $p \geq 1$, the {\em $p$-modulus} 
\index{Index}{modulus}%
of $\Gamma$ is  
\[ \mod_p(\Gamma) = \sup_\rho \int_U\rho^p \; dm\]
where the supremum is over all Borel measurable functions which are {\em admissible} in the sense that 
\[ \int_\gamma \rho \ ds \geq 1 \;\;\;\mbox{for all } \gamma \in \Gamma.\]
Any annulus $A$ in the Riemann sphere is conformally isomorphic to an annulus of the form $\{z : 1 < |z| < R\}$.  The $2$-modulus of the path family $\Gamma_s$ separating the boundary components coincides with the classical modulus defined by $\mod(A) = (1/2\pi)\ln(R)$.  

One characterization is given in terms of infinitesimal pointwise distortion of the roundness of balls:\begin{thm}\cite[Thm. II.6.2]{rickman:qrmappings}  
\label{thm:qc_via_H}A nonconstant mapping $f: U \to \R^n$ is quasiregular if and only if it satisfies
\be
\item $f$ is orientation-preserving, discrete, and open;
\item $H(x,f)$ is finite at every point, except maybe on a countable
set;
\item there exists $a<\infty$ such that $H(x,f) \leq a$ for almost every $x \in U\setminus B_f$.
\eb
\end{thm}

Another is geometric and is given geometrically in terms of the distortion of moduli of path families.  For finite branched coverings, the quantity $N(f,A)$ in the statement is at most the degree of $f$.
\begin{thm}\cite[Thm. II.6.7]{rickman:qrmappings}
\label{thm:qc_via_mod}
A continuous, orientation-preserving, discrete, open map $f: U \to \R^n$ is $K$-quasiregular if and only if 
\[ \mod_n(\Gamma) \leq K\cdot N(f,A) \cdot \mod_n(\Gamma)\]
for all path families $\Gamma$ in $A$ and all Borel sets $A \subset U$, where 
\[ N(f,A) = \sup_{y \in \R^n}\#\{f^{-1}(y) \intersect A\}.\]
\end{thm}

In dimension two, there is yet another characterization:  a theorem of Sto{\"{\i}}low \cite{stoilow:factor} implies that a quasiregular map $f: \rs \to \rs$ is the composition of a rational map and a quasiconformal map.    Moreover, a homeomorphism $h: \rs \to \rs$ is quasiconformal if and only if it is quasisymmetric \cite[Thm. 11.14]{heinonen:analysis}.  

The definition of a quasiregular map extends readily to maps between oriented, $C^\infty$ Riemannian manifolds.  A nonconstant quasiregular map between compact Riemannian manifolds is a finite branched covering \cite[Prop. 4.4]{rickman:qrmappings}.

The definition of a quasiconformal map extends in the obvious way to the setting of  homeomorphisms between arbitrary metric spaces.

\chapter{Hyperbolic groups in a nutshell}
\label{appendix:hyperbolic_groups}

We briefly summarize concepts and results regarding hyperbolic and convergence groups.  Our treatment is intented to highlight similarities between cxc dynamics and group actions.

We refer to \S\,\ref{scn:compactification} for the definitions
of a hyperbolic space and of its boundary. More generally, references
on hyperbolic metric spaces and hyperbolic groups include 
\cite{coornaert:delzant:papadopoulos, ghys:delaharpe:groupes, bridson:haefliger:book}.
One may consult \cite{coornaert:patterson-sullivan} for quasiconformal measures on their boundaries,
\cite{gehring:martin:qcgroupsi, bowditch:characterization,bowditch:convergence_groups}
for convergence group actions.

\section{Definition of a hyperbolic group}

A metric space is {\em proper}
\index{Index}{proper metric space}%
 if its closed balls are compact.  
An action of a group $ G$ on a proper metric space $X$  is said to be {\it geometric} 
\index{Index}{geometric action}%
\index{Index}{action, geometric}%
if
\bi
\item[(GA1)] each element acts as an isometry;
\item[(GA2)] the action is properly discontinuous;
\item[(GA3)] the action is cocompact.
\ib

Recall that a group $ G$ 
of isometries acts properly discontinuously
\index{Index}{properly discontinuous}%
 on $X$ if,  for any compact sets $K$ and $L$, the number of
group elements $ g\in G$ such that $ g(K)\cap L\ne\emptyset$ is finite.

\gap

 For example, if $ G$ is a finitely generated group and $S$ is a finite  set of generators for which $s \in S$ implies $s^{-1} \in S$, one may consider the Cayley graph
 \index{Index}{Cayley graph}%
  $\Gamma$ associated with $S$: the set of vertices are the elements of
the group, and pairs $( g, g')\in G\times G$ define an edge if $ g^{-1} g'\in S$. Endowing $\Gamma$ with the length
metric which makes each edge isometric to the segment $[0,1]$ defines the {\it word metric associated with $S$}.
\index{Index}{word metric}%
It turns $\Gamma$ into a geodesic proper metric space on which $ G$ acts geometrically by left-translation. A different generating set yields
a quasi-isometric graph.

We recall \v{S}varc-Milnor's lemma which provides a sort of converse statement \cite[Prop.\,3.19]{ghys:delaharpe:groupes}:

\begin{lemma}\label{svmi} Let $X$ be a geodesic proper metric space, and $ G$ a group which acts geometrically
on $X$. Then $ G$ is finitely generated and $X$ is quasi-isometric to any locally finite Cayley graph
of $ G$.\end{lemma}

\bigskip

A group $ G$ is {\it hyperbolic}
\index{Index}{hyperbolic group}%
 if it acts geometrically on a geodesic proper
hyperbolic metric space $X$ (e.g. a locally finite Cayley graph). 
Then \v{S}varc-Milnor's lemma above implies that $ G$ is finitely generated.

By definition, we will say that a metric space
$(X,d)$ is {\it quasi-isometric}
\index{Index}{quasi-isometry}%
 to a group $ G$ if
it is quasi-isometric to a locally finite Cayley graph of $ G$.

A hyperbolic group is said to be {\it elementary}
\index{Index}{elementary group}%
 if it is finite or quasi-isometric to $\Z$.
We will only consider non-elementary hyperbolic groups.

\section{Action on the boundary}

Let $ G$ be a hyperbolic group acting geometrically
on $(X,d)$.
From \v{S}varc-Milnor's lemma, the homeomorphism type of the boundary $\partial X$
depends only on the group, and so
$\partial G$ is well-defined up to homeomorphisms.  

Moreover, the Cayley graph of $G$ in its word metric is quasi-starlike, and so 
the boundary $\partial G$ can be  equipped with a metric $d_\varepsilon$ 
compatible with its topology by means of the compactification procedure 
given in \S \ref{scn:compactification}.  
Different choices of generating set and sufficiently small parameter 
$\varepsilon$ yield metrics which differ by quasisymmetries.  
Therefore, the {\em conformal gauge} of $G$, defined
\index{Index}{conformal gauge, of a group}%
 as the set of all 
metrics on $\partial G$ which are quasisymmetric to such a metric 
$d_\varepsilon$,  depends only on the (quasi-isometry class of the) group $G$.

The action of $ G$ extends to the boundary by homeomorphisms.
It defines a {\it convergence group action}:
\index{Index}{convergence group action}%
 let $\Theta$
denote the set of distinct triples of points of $\partial X$; then
$ G$ acts properly discontinuously on $\Theta$.

This action on $\Theta$ is even cocompact and so defines
 a so-called {\it uniform convergence action}. 
 \index{Index}{uniform convergence group action}%
 \index{Index}{convergence group, uniform}%
 Conversely, 
 Bowditch proved that
a group $G$ admits a uniform convergence action on a metrizable
perfect compact space $Z$  if and only 
if $G$ is hyperbolic (Theorem \ref{thm:bowditch_characterization}). In this
case $Z$ is homeomorphic to $\partial G$. The action is thus also
canonical. 

It follows that the action is minimal on $\partial G$ (every
orbit is dense in $\partial G$). 

\gap

With respect to a visual metric $d_{\varepsilon}$,  the action of $G$
is {\it uniformly quasi-M\"obius}:  
\index{Index}{uniformly quasi-M\"obius}%
there is an increasing
homeomorphism $\eta:\R_+\to\R_+$ such that, for any distinct quadruples
$(a,b,c,d)$ in $\partial X$, and any $g\in G$,
$$\frac{d_{\varepsilon}(g(a),g(b))}{d_{\varepsilon}(g(a),g(c))}
\frac{d_{\varepsilon}(g(c),g(d))}{d_{\varepsilon}(g(b),g(d))}
\le \eta\left( \frac{d_{\varepsilon}(a,b)}{d_{\varepsilon}(a,c)}
\frac{d_{\varepsilon}(c,d)}{d_{\varepsilon}(b,d)}\right)\,.$$

This property is preserved under conjugation by a quasisymmetric map.

\gap

In  summary: given any uniform convergence  
action of a group $G$ on a perfect metrizable compact space, 
there is a canonically associated gauge of metrics in which 
the dynamics is uniformly quasi-M\"obius.

\section{Quasiconformal measures}

{\noindent\bf Busemann functions.} Let $(X,w)$ be a pointed hyperbolic geodesic proper space.
Let $a\in\partial X$, $x,y\in X$ and suppose $h:\R_+\to X$ is a  ray such that  $h(0)=y$ and $\lim_{t \to \infty}h(t)=a$.
Let $\beta_a(x, y; h) = \limsup_{t \to \infty} (|x-h(t)|-|y-h(t)|)$.  For fixed $a$, 
the {\em Busemann function}
\index{Index}{Busemann function}%
 $\beta_a(\cdot, \cdot)$ {\it at the point $a$} 
is the function of two variables $x, y \in X$ defined by 
$$\beta_a(x,y)=\sup\{\beta_a(x, y; h),\mbox{ with }h\mbox{ as above}\}\,.$$

Busemann functions, visual metrics and the action of $G$ are related by the following property:  
for any $a\in\partial X$ and any $ g\in G$, there
is some neighborhood $V$ of $a$ such that, for any $b,c\in V$,
$$d_\varepsilon(g(b),g(c))\asymp L_g(a) d_\varepsilon(b,c)\,$$  where $L_g(a) = e^{\varepsilon\beta_a(w, g^{-1}(w))}$.
Moreover, $ G$ also acts on measures on  $\partial X$
through the usual rule: $ (g^*\rho)(A) = \rho( g A)$.

\bigskip

The next theorem,  proved by  M.\,Coornaert \cite{coornaert:patterson-sullivan}, 
summarizes properties of quasiconformal measures
\index{Index}{quasiconformal measure}%
on the boundary of $X$.  
The symbol $G(w)$ denotes the orbit of $w$ under the action of $G$.  

\begin{thm}\label{reg} Let $ G$ be a non-elementary hyperbolic group acting 
geometrically on a geodesic proper hyperbolic metric space $(X,d)$.
For any small enough $\varepsilon>0$,  one has 
$0<\mbox{\em dim}_H\,(\partial X,d_\varepsilon)<\infty$, and
$$v = \limsup_{R \to \infty} \frac{1}{R}\log\left|\{G(w)\cap B(w,R)\}\right|=\varepsilon\cdot\mbox{\em dim}_H\,(\partial X,d_\varepsilon)\,.$$
Let $\rho$ be the  Hausdorff measure on $\partial X$ of dimension $\alpha= v/\varepsilon$.  Then 
\bi

\item[(i)] the measure $\rho$ is Ahlfors-regular of dimension $\alpha$: for any $a\in\partial X$,
for any $r\in (0,\diam\partial X)$, we have $\rho(B_{\varepsilon}(a,r))\asymp r^\alpha$.
In particular, $0<\rho(\partial X)<\infty$;
\item[(ii)] the measure $\rho$ is a $ G$-quasiconformal measure 
of dimension $\alpha$:
for any $g\in  G$, we have $\rho\ll g^*\rho\ll\rho$, and $\rho$-almost everywhere
$$\frac{ d g^*\rho}{d\rho}\asymp (L_g)^\alpha;$$\\
\item[(iii)] the action of $ G$ is ergodic for $\rho$: 
for any $ G$-invariant Borel set $A$ of $\partial X$,
$$
\rho (A) =0 \;\; \mathrm{or} \;\; \rho (\partial X \backslash A)=0 \, .
$$\ib

Moreover, if $\rho'$ is another  $ G$-quasiconformal measure, then 
its dimension is also $\al$,
 $\rho\ll\rho'\ll\rho$ and $\dis\frac{ d\rho}{d\rho'}\asymp 1$ a.e. and
$$\#\{G(w)\cap B(w,R)\}\asymp e^{vR}\,.$$ \end{thm}

The class of measures thus defined on $\partial X$ is called the
{\it Patterson-Sullivan} class.
It does not depend on the choice of the parameter $\varepsilon$, 
but it does depend on the metric $d$.

The study of quasiconformal measures yields the following key estimate for the measure of shadows. 
Below, $d(\cdot, \cdot)$ denotes the original metric on $X$.  Recall that 
The shadow $\mho(x,R)$ consists of those points $y \in \cl{X}_\varepsilon$ 
for which there exists a geodesic ray in $X$ emanating from the basepoint 
$w$ and passing through the closed $d$-ball $\cl{B}(x,R)$; 
compare \S \ref{scn:compactification}.

\begin{lemma} \label{ashadow} {\em\bf (Lemma of the shadow)} Under the assumptions of Theorem \ref{reg}, 
there  exists $R_0$, such that
if $R>R_0$, then, for any $x\in X$,  $$\rho(\mho(x,R))\asymp  e^{-v d(w,x)}$$
where the implicit constants do not depend on $x$. \end{lemma}

\section{Cannon's conjecture}

Let $G$ be a hyperbolic group. J.\,Cannon's conjecture
\index{Index}{Cannon's conjecture}%
 states
that if $\partial G$ is a topological $2$-sphere, then
it acts geometrically on hyperbolic three-space $\IH^3$, i.e. it is virtually a cocompact Kleinian group.
There have been essentially two approaches to prove this conjecture.

\gap

The first one is combinatorial and due to Cannon et al.
\cite{cannon:rmt,cannon:swenson:characterization}: let $\Gamma$ be a locally
finite Cayley graph of $G$, and let us consider the sequence
of covers $\{\UUU_n\}_n$ by shadows of balls centred at vertices
at distance $n$ from the identity. They prove the conjecture under
the assumption that the sequence
is conformal in the sense of Cannon (see \S \ref{scn:cannoncomb} for the definition).

\gap

The second approach is analytical and due to Bonk and Kleiner 
\cite{bonk:kleiner:qsparam, bonk:kleiner:rigidity, bonk:kleiner:conf_dim}.
They prove the conjecture under the assumption that the gauge of the group contains a $2$-Ahlfors
regular metric or a $Q$-regular $Q$-Loewner metric.  
\end{appendix}

\bibliographystyle{math}
\bibliography{cxcrefs-v2}


\twocolumn{\chapter*{Index}}
\markboth{INDEX}{INDEX}
  \input{Index.ind}

\twocolumn{\chapter*{Table of symbols}} 
\markboth{TABLE OF SYMBOLS}{TABLE OF SYMBOLS}
 \input{Symbols.ind} 

\end{document}